\pdfoutput=1
\RequirePackage{ifpdf}
\ifpdf 
\documentclass[pdftex]{sigma}
\else
\documentclass{sigma}
\fi

\numberwithin{equation}{section}

\newtheorem{Theorem}{Theorem}[section]

\newtheorem{Lemma}[Theorem]{Lemma}
\newtheorem{Proposition}[Theorem]{Proposition}
{ \theoremstyle{definition}
\newtheorem{Definition}[Theorem]{Definition}

\newtheorem{Example}[Theorem]{Example}
\newtheorem{Remark}[Theorem]{Remark}
\newtheorem{Notation}[Theorem]{Notation}
\newtheorem{Assumptions}[Theorem]{Assumptions}
}

\usepackage{mathtools}

\usepackage{tikz}
\usetikzlibrary{decorations.pathreplacing, calc, shapes.geometric, shapes.misc, cd, positioning}
\usepackage{tabularx}

\usepackage[only,llbracket,rrbracket]{stmaryrd}

\DeclareFontFamily{U}{MnSymbolC}{}
\DeclareSymbolFont{MnSyC}{U}{MnSymbolC}{m}{n}
\DeclareMathSymbol{\diamondplus}{\mathbin}{MnSyC}{"7C}
\DeclareFontShape{U}{MnSymbolC}{m}{n}{
    <-6>  MnSymbolC5
   <6-7>  MnSymbolC6
   <7-8>  MnSymbolC7
   <8-9>  MnSymbolC8
   <9-10> MnSymbolC9
  <10-12> MnSymbolC10
  <12->   MnSymbolC12}{}

\usepackage{enumitem}
\setlist[enumerate]{nosep}
\setlist[enumerate,1]{label=$(\alph*)$}
\setlist[enumerate,2]{label=(\roman*)}
\setlist[enumerate,3]{label=(\arabic*)}
\setlist[itemize]{nosep}
\setlist[itemize,1]{label=$\square$}

\begin{document}

\allowdisplaybreaks

\newcommand{\arXivNumber}{2305.12494}

\renewcommand{\PaperNumber}{082}

\FirstPageHeading

\ShortArticleName{First Cohomology with Trivial Coefficients of All Unitary Easy Quantum Group Duals}

\ArticleName{First Cohomology with Trivial Coefficients\\ of All Unitary Easy Quantum Group Duals}

\Author{Alexander MANG}

\AuthorNameForHeading{A. Mang}

\Address{Hamburg University, Bundesstra{\ss}e 55, 20146 Hamburg, Germany}
\Email{\href{mailto:alex@alexandermang.net}{alex@alexandermang.net}}
\URLaddress{\url{https://alexandermang.net}}

\ArticleDates{Received September 24, 2023, in final form August 22, 2024; Published online September 12, 2024}

\Abstract{The first quantum group cohomology with trivial coefficients of the discrete dual of any unitary easy quantum group is computed. That includes those potential quantum groups whose associated categories of two-colored partitions have not yet been found.}


\Classification{20G42; 05A18}

\newcommand\Tstrut{\rule{0pt}{2.6ex}}
\newcommand{\identitymatrix}{I}
\newcommand{\freealg}[2]{#1\langle #2\rangle}
\newcommand{\univalg}[3]{#1\langle #2 \classpredicate #3\rangle}
\newcommand{\HSalgebraconstructdeform}[4]{#2\diamondplus^{#1}_{#4}#3}
\newcommand{\HScocycles}[3]{Z^{#1}_{\mathrm{HS}}(#2,#3)}
\newcommand{\HScoboundaries}[3]{B^{#1}_{\mathrm{HS}}(#2,#3)}
\newcommand{\HScohomology}[3]{H^{#1}_{\mathrm{HS}}(#2,#3)}
\newcommand{\HSdifferential}[1]{\partial^{#1}}
\newcommand{\HSuacocycles}[4]{Z^{#1}_{#2,#3,#4}}
\newcommand{\HSuacoboundaries}[4]{B^{#1}_{#2,#3,#4}}
\newcommand{\HSuacohomology}[4]{H^{#1}_{#2,#3,#4}}
\newcommand{\HSuafunc}[4]{F^{#1}_{#2,#3,#4}}
\newcommand{\HSuafuncx}[5]{F^{#1,#5}_{#2,#3,#4}}
\newcommand{\HSuafuncxx}[6]{F^{#1,#5}_{#2,#3,#4}(#6)}
\newcommand{\HSuafuncxxx}[7]{F^{#1,#5}_{#2,#3,#4}(#6)(#7)}
\newcommand{\thefunc}{F}
\newcommand{\thefuncx}[1]{F_{#1}}
\newcommand{\thefuncxx}[2]{F_{#1}(#2)}
\newcommand{\Scomultiplication}{\Delta}
\newcommand{\Scounit}{\epsilon}
\newcommand{\Sstar}{\ast}
\newcommand{\tiff}{\text{if and only if}}
\newcommand{\tif}{\text{if }}
\newcommand{\tothw}{\text{otherwise}}
\newcommand{\tand}{\text{and}}
\newcommand{\tor}{\text{or}}
\newcommand{\twhere}{\text{where}}
\newcommand{\twith}{\text{with}}
\newcommand{\todd}{\text{odd}}
\newcommand{\teven}{\text{even}}
\newcommand{\Sand}{\,\wedge\,}
\newcommand{\Sor}{\,\vee\,}
\newcommand{\Sthen}{\,\Rightarrow\,}
\newcommand{\Simplies}{\,\implies\,}
\newcommand{\Siff}{\Leftrightarrow}
\newcommand{\Strue}{\top}
\newcommand{\Sfalse}{\bot}
\newcommand{\iffpd}{{\ \colon\!\!\!\iff}}
\newcommand{\funcdef}{\colon}
\newcommand{\classpredicate}{\mid}
\newcommand{\quantorpredicate}{\colon }
\newcommand{\Sidentity}{\mathrm{id}}
\newcommand{\SidentityX}[1]{\Sidentity_{#1}}
\newcommand{\Scomposition}{\circ}
\newcommand{\Smonoidalproduct}{\otimes}
\newcommand{\Sinternalhom}{[\argph,\argph]}
\newcommand{\SinternalhomX}[2]{[#1,#2]}

\newcommand{\nnint}{{\mathbb{N}_0}}
\newcommand{\pint}{\mathbb{N}}
\newcommand{\integers}{\mathbb{Z}}
\newcommand{\comps}{\mathbb{C}}
\newcommand{\reals}{\mathbb{R}}
\newcommand{\dwi}[1]{\llbracket{#1}\rrbracket}
\newcommand{\fromto}[2]{#1 \to #2}
\newcommand{\xfromto}[3]{#1 \colon  #2 \to #3}
\newcommand{\fromtomaps}[4]{#1 \to #2, \, #3 \mapsto #4}
\newcommand{\xfromtomaps}[5]{#1 \colon  #2 \to #3, \, #4 \mapsto #5}
\newcommand{\sran}{\mathrm{ran}}
\newcommand{\sranx}[1]{\mathrm{ran}(#1)}
\newcommand{\bfunc}{,\,}

\newcommand{\ASdualX}[1]{\overline{#1}}
\newcommand{\subproof}[1]{\textit{#1}}
\newcommand{\Sisomorphic}{\cong}
\newcommand{\Sdirectsum}{\oplus}
\newcommand{\Sdirectproduct}{\times}
\newcommand{\Ssetmonoidalproduct}{\times}
\newcommand{\Sleftmoduleaction}{\triangleright}
\newcommand{\Srightmoduleaction}{\triangleleft}
\newcommand{\Sadj}{^\ast}
\newcommand{\Stra}{^{\mathrm{t}}}
\newcommand{\argph}{\hspace{0.5pt}\cdot\hspace{0.5pt}}

\newcommand{\upp}[1]{\prescript{\resizebox{0.008\hsize}{!}{$\blacksquare$}}{}{#1}}
\newcommand{\lop}[1]{\prescript{}{\resizebox{0.008\hsize}{!}{$\blacksquare$}}{#1}}
\newcommand{\tsopx}[2]{\Pi^{#1}_{#2}}
\newcommand{\blpoint}{\mathnormal{\bullet}}
\newcommand{\whpoint}{\mathnormal{\circ}}
\newcommand{\blaw}{\{\whpoint,\blpoint\}}
\newcommand{\tcsx}[2]{\Sigma^{#1}_{#2}}
\newcommand{\csx}[3]{\sigma^{#1}_{#2}(#3)}
\newcommand{\csfunx}[2]{\sigma^{#1}_{#2}}
\newcommand{\simg}[1]{#1_{\rightarrow}}
\newcommand{\simgx}[2]{#1_{\rightarrow}(#2)}
\newcommand{\spimg}[1]{#1^{\leftarrow}}
\newcommand{\spimgx}[2]{#1^{\leftarrow}(#2)}
\newcommand{\finerthan}{\mathrel{\leq}}
\newcommand{\nfinerthan}{\mathrel{\not\leq}}
\newcommand{\finerthanrel}{\mathop{\leq}}
\newcommand{\Yc}[1]{\mathfrak{#1}}
\newcommand{\Yp}[1]{\mathtt{#1}}
\newcommand{\Yb}[1]{\mathsf{#1}}
\newcommand{\blo}[2]{\pi_{#1}(#2)}
\newcommand{\blofun}[1]{\pi_{#1}}
\newcommand{\kerp}[1]{\ker(#1)}
\newcommand{\djp}[2]{#1\begin{smallmatrix}\resizebox{0.008\hsize}{!}{$\blacksquare$}\hspace{4pt}\\\hspace{4pt}\resizebox{0.008\hsize}{!}{$\blacksquare$}\end{smallmatrix}#2}
\newcommand{\zetf}{\zeta}
\newcommand{\zetfx}[2]{\zetf(#1,#2)}
\newcommand{\squarematrices}[2]{M_{#1}(#2)}
\newcommand{\kron}[2]{\delta_{#1,#2}}
\newcommand{\quomap}[2]{\sfrac{#1}{#2}}
\newcommand{\tsquomap}[2]{#1 / #2}
\newcommand{\onehalf}{{ \frac{1}{2}}}
\newcommand{\rightwardproduct}[2]{\underset{#1}{\overset{#2}{\overrightarrow{\prod}}}}
\newcommand{\tsrightwardproduct}[2]{{\overrightarrow{ \prod}}\vphantom{ \prod}_{#1}^{#2}}
\newcommand{\leftwardproduct}[2]{\underset{#1}{\overset{#2}{\overleftarrow{\prod}}}}
\newcommand{\tsleftwardproduct}[2]{{\overleftarrow{ \prod}}\vphantom{ \prod}_{#1}^{#2}}

\newcommand{\DQGgroupalgebra}[1]{\comps [#1 ]}
\newcommand{\DQGcocycles}[3]{Z^{#1} (#2,#3 )}
\newcommand{\DQGcoboundaries}[3]{B^{#1} (#2,#3 )}
\newcommand{\DQGcohomology}[3]{H^{#1} (#2,#3 )}
\newcommand{\DQGcocyclesTC}[2]{Z^{#1} (#2 )}
\newcommand{\DQGcoboundariesTC}[2]{B^{#1} (#2 )}
\newcommand{\DQGcohomologyTC}[2]{H^{#1} \bigl(#2 \bigr)}
\newcommand{\CQGdualC}[1]{\widehat{#1}}

\newcommand{\ueqgU}[1]{U_{#1}}
\newcommand{\ueqgUplus}[1]{U^{+}_{#1}}
\newcommand{\ueqgUstar}[2]{U^{\ast}_{#1,#2}}
\newcommand{\ueqgUtimes}[2]{U^{\times}_{#1,#2}}
\newcommand{\ueqgUtimesplus}[2]{U^{{\times}{+}}_{#1,#2}}

\newcommand{\Slength}[1]{|#1|}
\newcommand{\therelpoly}[5]{\mathrm{r}^{#1}_{#2}(#3)_{#4,#5}}

\newcommand{\alterer}[3]{#1\downarrow_{#2}#3}

\newcommand{\thedim}{n}
\newcommand{\theqg}{G}
\newcommand{\thegens}{E}
\newcommand{\theuni}{u}
\newcommand{\theunim}[1]{\theuni^{#1}}
\newcommand{\theunix}[3]{\theuni^{#1}_{#2,#3}}
\newcommand{\therels}{R}
\newcommand{\thepartrels}[1]{\therels_{#1}}
\newcommand{\thecat}{\mathcal{C}}
\newcommand{\theideal}{J}
\newcommand{\thepartideal}[1]{\theideal_{#1}}
\newcommand{\thefield}{\mathbb{K}}
\newcommand{\themodule}{X}
\newcommand{\theresmodule}{Y}
\newcommand{\thecounit}{\Scounit}
  \newcommand{\thepredicate}{A}
  \newcommand{\thepredicatex}[2]{A(#1,#2)}

  \newcommand{\scenarioone}{R1}
  \newcommand{\scenariotwo}{R2}
  \newcommand{\scenariothree}{R3}
  \newcommand{\scenariofour}{R4}

  \newcommand{\conditionone}{P1}
  \newcommand{\conditiontwo}{P2}
  \newcommand{\conditionthree}{P3}

\newcommand{\oeqgS}[1]{S_{#1}}
\newcommand{\oeqgSprime}[1]{S^{\prime}_{#1}}
\newcommand{\oeqgSplus}[1]{S^{+}_{#1}}
\newcommand{\oeqgSprimeplus}[1]{S^{\prime +}_{#1}}

\newcommand{\oeqgO}[1]{O_{#1}}
\newcommand{\oeqgOstar}[1]{O^{\ast}_{#1}}
\newcommand{\oeqgOplus}[1]{O^{+}_{#1}}

\newcommand{\oeqgB}[1]{B_{#1}}
\newcommand{\oeqgBprime}[1]{B^{\prime}_{#1}}
\newcommand{\oeqgBhashstar}[1]{B^{\#\ast}_ {#1}}
\newcommand{\oeqgBplus}[1]{B^{+}_{#1}}
\newcommand{\oeqgBprimeplus}[1]{B^{\prime{+}}_{#1}}
\newcommand{\oeqgBhashplus}[1]{B^{\#{+}}_{#1}}

\newcommand{\oeqgH}[1]{H_{#1}}
\newcommand{\oeqgHstar}[1]{H^{\ast}_{#1}}
\newcommand{\oeqgHplus}[1]{H^{+}_{#1}}
\newcommand{\oeqgHcurly}[2]{H^{\{#1\}}_{#2}}
\newcommand{\oeqgHcurlyinfinity}[1]{H^{\{\infty\}}_{#1}}
\newcommand{\oeqgHsquare}[2]{H^{[#1]}_{#2}}
\newcommand{\oeqgHsquareinfinity}[1]{H^{[\infty]}_{#1}}
\newcommand{\oeqgHround}[2]{H^{(#1)}_{#2}}
\newcommand{\oeqgHangle}[2]{H^{\langle #1\rangle}_{#2}}

  \newcommand{\onlyneutralnonsingletonblocks}{\mathsf{NNSB}}
  \newcommand{\onlyneutralpartitions}{\mathsf{NP}}

\newcommand{\PartIdenB}{%
  \begin{tikzpicture}[scale=0.25,baseline=0.015cm]
    \def\xdist{1}
    \def\ydist{1}
    \node [scale=0.33, circle, draw=black, fill=black] (a1) at ({0*\xdist},{0*\ydist}) {};
    \node [scale=0.33, circle, draw=black, fill=black] (b1) at ({0*\xdist},{1*\ydist}) {};
    \draw (a1) to (b1);
  \end{tikzpicture}
}
\newcommand{\PartIdenW}{%
  \begin{tikzpicture}[scale=0.25,baseline=0.015cm]
    \def\xdist{1}
    \def\ydist{1}
    \node [scale=0.33, circle, draw=black, fill=white] (a1) at ({0*\xdist},{0*\ydist}) {};
    \node [scale=0.33, circle, draw=black, fill=white] (b1) at ({0*\xdist},{1*\ydist}) {};
    \draw (a1) to (b1);
  \end{tikzpicture}
}
\newcommand{\PartIdenUpBW}{%
  \begin{tikzpicture}[scale=0.25,baseline=-0.015cm]
    \def\xdist{1}
    \def\ydist{1}
    \def\ypar{1}
    \node [scale=0.33, circle, draw=black, fill=black] (a1) at ({0*\xdist},{1*\ydist}) {};
    \node [scale=0.33, circle, draw=black, fill=white] (a2) at ({1*\xdist},{1*\ydist}) {};
    \draw (a1) -- ++ (0,{-\ypar}) -| (a2);
  \end{tikzpicture}
}
\newcommand{\PartIdenUpWB}{%
  \begin{tikzpicture}[scale=0.25,baseline=-0.015cm]
    \def\xdist{1}
    \def\ydist{1}
    \def\ypar{1}
    \node [scale=0.33, circle, draw=black, fill=white] (a1) at ({0*\xdist},{1*\ydist}) {};
    \node [scale=0.33, circle, draw=black, fill=black] (a2) at ({1*\xdist},{1*\ydist}) {};
    \draw (a1) -- ++ (0,{-\ypar}) -| (a2);
  \end{tikzpicture}
  }
\newcommand{\PartIdenLoBW}{%
  \begin{tikzpicture}[scale=0.25,baseline=-0.015cm]
    \def\xdist{1}
    \def\ydist{1}
    \def\ypar{1}
    \node [scale=0.33, circle, draw=black, fill=black] (a1) at ({0*\xdist},{0*\ydist}) {};
    \node [scale=0.33, circle, draw=black, fill=white] (a2) at ({1*\xdist},{0*\ydist}) {};
    \draw (a1) -- ++ (0,{\ypar}) -| (a2);
  \end{tikzpicture}
}
\newcommand{\PartIdenLoWB}{%
  \begin{tikzpicture}[scale=0.25,baseline=-0.015cm]
    \def\xdist{1}
    \def\ydist{1}
    \def\ypar{1}
    \node [scale=0.33, circle, draw=black, fill=white] (a1) at ({0*\xdist},{0*\ydist}) {};
    \node [scale=0.33, circle, draw=black, fill=black] (a2) at ({1*\xdist},{0*\ydist}) {};
    \draw (a1) -- ++ (0,{\ypar}) -| (a2);
  \end{tikzpicture}
}
\newcommand{\PartIdenLoWW}{%
  \begin{tikzpicture}[scale=0.25,baseline=-0.015cm]
    \def\xdist{1}
    \def\ydist{1}
    \def\ypar{1}
    \node [scale=0.33, circle, draw=black, fill=white] (a1) at ({0*\xdist},{0*\ydist}) {};
    \node [scale=0.33, circle, draw=black, fill=white] (a2) at ({1*\xdist},{0*\ydist}) {};
    \draw (a1) -- ++ (0,{\ypar}) -| (a2);
  \end{tikzpicture}
}
\newcommand{\PartFourWBWB}{%
  \begin{tikzpicture}[scale=0.25,baseline=-0.015cm]
    \def\xdist{0.666}
    \def\ydist{1}
    \def\ypar{1}
    \node [scale=0.33, circle, draw=black, fill=white] (a1) at ({0*\xdist},{0*\ydist}) {};
    \node [scale=0.33, circle, draw=black, fill=black] (a2) at ({1*\xdist},{0*\ydist}) {};
    \node [scale=0.33, circle, draw=black, fill=white] (a3) at ({2*\xdist},{0*\ydist}) {};
    \node [scale=0.33, circle, draw=black, fill=black] (a4) at ({3*\xdist},{0*\ydist}) {};
    \draw (a1) -- ++ (0,{\ypar}) -| (a4);
    \draw (a2) -- ++ (0,{\ypar});
    \draw (a3) -- ++ (0,{\ypar});
  \end{tikzpicture}
}
\newcommand{\PartSinglesWB}{%
  \begin{tikzpicture}[scale=0.25,baseline=-0.015cm]
    \def\xdist{1}
    \def\ydist{1}
    \def\ypar{1}
    \node [scale=0.33, circle, draw=black, fill=white] (a1) at ({0*\xdist},{0*\ydist}) {};
    \node [scale=0.33, circle, draw=black, fill=black] (a2) at ({1*\xdist},{0*\ydist}) {};
    \draw[->] (a1) -- ++ (0,{\ypar});
    \draw[->] (a2) -- ++ (0,{\ypar});
  \end{tikzpicture}
}
\newcommand{\PartSingleW}{%
    \begin{tikzpicture}[scale=0.25,baseline=-0.015cm]
    \def\xdist{1}
    \def\ydist{1}
    \def\ypar{1}
    \node [scale=0.33, circle, draw=black, fill=white] (a1) at ({0*\xdist},{0*\ydist}) {};
    \draw[->] (a1) -- ++ (0,{\ypar});
  \end{tikzpicture}
}
\newcommand{\PartSingleB}{%
    \begin{tikzpicture}[scale=0.25,baseline=-0.015cm]
    \def\xdist{1}
    \def\ydist{1}
    \def\ypar{1}
    \node [scale=0.33, circle, draw=black, fill=black] (a1) at ({0*\xdist},{0*\ydist}) {};
    \draw[->] (a1) -- ++ (0,{\ypar});
  \end{tikzpicture}
}
\newcommand{\PartCoSingleWBWB}{%
    \begin{tikzpicture}[scale=0.25,baseline=-0.015cm]
    \def\xdist{0.666}
    \def\ydist{1.25}
    \def\ypar{1}
    \node [scale=0.33, circle, draw=black, fill=white] (a1) at ({0*\xdist},{0*\ydist}) {};
    \node [scale=0.33, circle, draw=black, fill=black] (a2) at ({1*\xdist},{0*\ydist}) {};
    \node [scale=0.33, circle, draw=black, fill=white] (a3) at ({2*\xdist},{0*\ydist}) {};
    \node [scale=0.33, circle, draw=black, fill=black] (a4) at ({3*\xdist},{0*\ydist}) {};
    \draw (a1) -- ++ (0,{1*\ydist}) -| (a3);
    \draw[->] (a2) -- ++ (0,{1*\ypar});
    \draw[->] (a4) -- ++ (0,{1*\ypar});
    \end{tikzpicture}
  }
\newcommand{\PartCrossWW}{%
  \begin{tikzpicture}[scale=0.25,baseline=0.015cm]
    \def\xdist{1}
    \def\ydist{1}
    \node [scale=0.33, circle, draw=black, fill=white] (a1) at ({0*\xdist},{0*\ydist}) {};
    \node [scale=0.33, circle, draw=black, fill=white] (a2) at ({1*\xdist},{0*\ydist}) {};
    \node [scale=0.33, circle, draw=black, fill=white] (b1) at ({0*\xdist},{1*\ydist}) {};
    \node [scale=0.33, circle, draw=black, fill=white] (b2) at ({1*\xdist},{1*\ydist}) {};
    \draw (a1) to (b2);
    \draw (a2) to (b1);
  \end{tikzpicture}
}
\newcommand{\UCPartIden}{%
  \begin{tikzpicture}[baseline=0.0cm]
    \def\xdist{0.2cm}
    \def\ydist{0.32cm}
    \def\ypar{0.25cm}
    \coordinate (a1) at ({0*\xdist},{0*\ydist});
    \draw (a1) -- ++ (0,{\ypar});
  \end{tikzpicture}
}
\newcommand{\UCPartIdenLo}{%
  \begin{tikzpicture}[scale=0.25,baseline=0.015cm]
    \def\xdist{1}
    \def\ydist{1}
    \def\ypar{1}
    \coordinate (a1) at ({0*\xdist},{0*\ydist}) {};
    \coordinate (a2) at ({1*\xdist},{0*\ydist}) {};
    \draw (a1) -- ++ (0,{\ypar}) -| (a2);
  \end{tikzpicture}
  }
\newcommand{\UCPartIdenUp}{%
  \begin{tikzpicture}[scale=0.25,baseline=0.015cm]
    \def\xdist{1}
    \def\ydist{1}
    \def\ypar{1}
    \coordinate (a1) at ({0*\xdist},{1*\ydist}) {};
    \coordinate (a2) at ({1*\xdist},{1*\ydist}) {};
    \draw (a1) -- ++ (0,{-\ypar}) -| (a2);
  \end{tikzpicture}
  }
\newcommand{\UCPartFour}{%
  \begin{tikzpicture}[baseline=0.0cm]
    \def\xdist{0.2cm}
    \def\ydist{0.32cm}
    \def\ypar{0.2cm}
    \coordinate (a1) at ({0*\xdist},{0*\ydist});
    \coordinate (a2) at ({1*\xdist},{0*\ydist});
    \coordinate (a3) at ({2*\xdist},{0*\ydist});
    \coordinate (a4) at ({3*\xdist},{0*\ydist});
    \draw (a1) -- ++ (0,{\ypar}) -| (a4);
    \draw (a2) -- ++ (0,{\ypar});
    \draw (a3) -- ++ (0,{\ypar});
  \end{tikzpicture}
}
\newcommand{\UCPartCoSingle}{%
  \begin{tikzpicture}[baseline=0.0cm]
    \def\xdist{0.2cm}
    \def\ydist{0.32cm}
    \def\ypar{0.2cm}
    \coordinate (a1) at ({0*\xdist},{0*\ydist});
    \coordinate (a2) at ({1*\xdist},{0*\ydist});
    \coordinate (a3) at ({2*\xdist},{0*\ydist});
    \coordinate (a4) at ({3*\xdist},{0*\ydist});
    \draw (a1) -- ++ (0,{\ydist}) -| (a3);
    \draw[->] (a2) -- ++ (0,{\ypar});
    \draw[->] (a4) -- ++ (0,{\ypar});
  \end{tikzpicture}
}
\newcommand{\UCPartSingle}{%
  \begin{tikzpicture}[baseline=0.0cm]
    \def\xdist{0.2cm}
    \def\ydist{0.32cm}
    \def\ypar{0.25cm}
    \coordinate (a1) at ({0*\xdist},{0*\ydist});
    \draw[->] (a1) -- ++ (0,{\ypar});
  \end{tikzpicture}
}
\newcommand{\UCPartSinglesTensor}{%
  \begin{tikzpicture}[baseline=0.0cm]
    \def\xdist{0.2cm}
    \def\ydist{0.32cm}
    \def\ypar{0.25cm}
    \coordinate (a1) at ({0*\xdist},{0*\ydist});
    \coordinate (a2) at ({1.75*\xdist},{0*\ydist});
    \draw[->] (a1) -- ++ (0,{\ypar});
    \draw[->] (a2) -- ++ (0,{\ypar});
    \node [scale=0.66] at ($(a1)+({0.5*1.75*\xdist},{0.325*\ydist})$) {$\otimes$};
  \end{tikzpicture}
}
\newcommand{\UCPartSingles}{%
  \begin{tikzpicture}[baseline=0.0cm]
    \def\xdist{0.2cm}
    \def\ydist{0.32cm}
    \def\ypar{0.2cm}
    \coordinate (a1) at ({0*\xdist},{0*\ydist});
    \coordinate (a2) at ({1*\xdist},{0*\ydist});
    \draw[->] (a1) -- ++ (0,{\ypar});
    \draw[->] (a2) -- ++ (0,{\ypar});
  \end{tikzpicture}
}

\section{Introduction}
\label{section:introduction}

\subsection{Background and context}
\label{section:introduction-background-and-context}
In \cite{BanicaSpeicher2009}, Banica and Speicher provided a way of constructing compact quantum groups (in the sense of \cite{Woronowicz1987b, Woronowicz1991, Woronowicz1998}) by solving infinite combinatorics puzzles:  They introduced three operations on the collection of all equivalence relations of disjoint unions of finite sets and showed that each subset which is closed under these operations gives rise to a compact quantum group. An uncountable number of such sets and of the resulting so-called ``easy''  quantum groups and, in fact, all there can be, have since been found in \cite{BanicaCurranSpeicher2010,BanicaSpeicher2009, RaumWeber2014, RaumWeber2016a, RaumWeber2016b, Weber2013}. In \cite{TarragoWeber2018}, Tarrago and Weber extended Banica and Speicher's operations to the collection of all ``two-co\-lo\-red'' partitions, thus providing even more quantum groups to find. The classification program they initiated to determine all sets closed under the operations is still ongoing (see \cite{ Gromada2018,Maassen2021, Mang2022, MangWeber2020, MangWeber2021a, MangWeber2021b, MangWeber2021c, TarragoWeber2018}). The construction has since been further extended to two-co\-lo\-red partitions with arbitrarily many ``colors'' by Freslon in \cite{Freslon2017}, to ``three-dimensional'' sets  by C\'ebron and Weber in \cite{CebronWeber2016} and to  equivalence relations on graphs instead of sets by Man\v{c}inska and Roberson in \cite{MancinskaRoberson2020}.

An issue that all these constructions share is that it is difficult to tell which of the resulting compact quantum groups are new and which are isomorphic to already known ones. In particular, each solution to the combinatorics puzzle does not only provide one quantum group but an entire countably infinite series, one for each dimension of its fundamental representation. And already  Banica and Speicher themselves observed in \cite[Proposition~2.4\,(4)]{BanicaSpeicher2009} that, at least in some cases, the quantum groups of one solution are isomorphic to those of another, just shifted by one dimension. That underlines the importance of studying quantum group invariants with the potential of telling easy quantum groups apart. Of course, these are often very difficult to compute like, e.g., the $L^2$-cohomology of \cite{Kyed2008a} of discrete quantum groups. But perhaps at least the cohomology with trivial coefficients is a reasonable goal to strive for.

The present article computes the first order of the quantum group cohomology with trivial coefficients of the discrete duals of all of  Tarrago and Weber's so-called unitary easy quantum groups. That includes even the potential  ones whose combinatorics puzzles have not been solved yet. Said cohomology can be realized as the first Hochschild cohomology of the trivial bimodule of an augmented algebra presented in terms of generators and relations. As with any augmented algebra the space of $1$-coboundaries is then trivial and the task thus boils down to solving the generally infinite system of linear equations in the finitely many generators determining the $1$-cocycles.

The results of the present article might be useful for the computation of the second order begun by Bichon, Das, Franz, Gerhold, Kula and Skalski in \cite{BichonFranzGerhold2017, DasFranzKulaSkalski2018} as well as Wendel in \cite{Wendel2020}. The former six investigated the cohomology of certain easy quantum groups out of a different motivation. In particular, they were interested in the Calabi--Yau property of \cite{Ginzburg2006}, a generalization of  Poincar\'e duality, and the classification of Sch\"urmann triples. Namely, a quantum group whose second cohomology vanishes has the AC property, defined in \cite{FranzGerholdThom2015}, which is important in the study of quantum L\'evy processes because it guarantees the existence of an associated Sch\"urmann triple.

In \cite{BichonFranzGerhold2017, DasFranzKulaSkalski2018}, Bichon, Das, Franz, Gerhold, Kula and Skalski had already laid out a potential strategy for computing the second cohomology of any easy quantum group (later refined in~\cite{DasFranzKulaSkalski2023} to address  universal unitary quantum groups). This strategy is based on two key insights and goes as follows. They interpreted quantum group cohomology as Hochschild cohomology and chose the Hochschild complex as their resolution. Thus, they were faced with having to compute the quotient of the $2$-cocycles by the $2$-coboundaries. By a very clever use of the universal property of the quantum groups in question,  they managed to solve the linear system of equations determining the space of $2$-coboundaries. This use of the universal property is the first key tool (see \cite[Lemma~5.4]{BichonFranzGerhold2017} and \cite[Lemma~4.1]{DasFranzKulaSkalski2018}).

Understanding the $2$-cocycles then allowed them to define a ``defect map'', an injective linear map  from  $2$-cohomology to a certain finite-dimensional vector space of matrices. Thus, at this point they only needed to determine the image of this defect map in order to compute the second cohomology. This is where their second key insight  comes into play. Namely, although being interested only in the second-order cohomology, they incidentally also computed the first. That is because they wanted to make  use of the multiplicative structure of the cohomology ring. They showed that, at least for the specific quantum groups they investigated, each $2$-cocycle was cohomologous to a linear combination of cup products of $1$-cocycles. Thus, rather than having to probe the potentially infinite-dimensional vector space of all $2$-cocycles as the domain of the defect map they could confine themselves to determining the image of the restriction to cup products, a finite-dimensional space.

In short, when trying to compute the second cohomology of any easy quantum group it might be helpful, perhaps even necessary, to know the first cohomology. Hence, the main result of the present article might also constitute an intermediate step in computing the higher cohomologies of all easy quantum groups.

\subsection{Main result}
\label{section:introduction-main-results}
{
  \newcommand{\themat}{v}
  \newcommand{\indi}{i}
  \newcommand{\indj}{j}
  \newcommand{\thecoc}{\eta}
  \newcommand{\thecocx}[1]{\eta(#1)}
  \newcommand{\thebetti}{\beta_1(\widehat{G})}
  \newcommand{\thesum}{\lambda}
  \newcommand{\neutralblocks}{3}
  \newcommand{\neutralpartitions}{4}
  \newcommand{\smallblocks}{1}
  \newcommand{\largeblocks}{2}

  Let  $\squarematrices{\thedim}{\comps}$ be the $\comps$-vector space of $(\thedim\times\thedim)$-matrices with complex entries and $\identitymatrix$ the identity $(\thedim\times\thedim)$-matrix. Moreover, call a matrix  ``small'' if each of its rows and each of its columns sums to $0$.

Then, the below theorem extends the results  of  \cite{BichonFranzGerhold2017,DasFranzKulaSkalski2018} as well as \cite{Wendel2020}.
  \begin{theorem*}
Let $\thedim\in\pint$, let $\theqg$ be any unitary easy compact $(\thedim\times\thedim)$-matrix quantum group, let  $\theuni$ be its fundamental representation and let $\thecat$ be the category of two-co\-lo\-red partitions associated with $\theqg$. Say that $\thecat$ has property
  \begin{itemize}\itemsep=0pt
  \item[$(1)$] if and only if  each block of each two-co\-lo\-red partition of $\thecat$ has at most two points,
    \item[$(2)$] if and only if  each block of each two-co\-lo\-red partition of $\thecat$ has at least two points,
  \item[$(3)$] if and only if  each block of each two-co\-lo\-red partition of $\thecat$ with at least two points contains as many  white lower and black upper points as it does black lower and white upper points,
  \item[$(4)$]    if and only if each two-co\-lo\-red partition of $\thecat$ has as many white lower and black upper points as it has black lower and white upper points.
  \end{itemize}
An isomorphism of complex vector spaces from the first quantum group cohomology of $\CQGdualC{\theqg}$ with trivial coefficients to the subspace
  \[
 \{\themat\in\squarematrices{\thedim}{\comps}\Sand \thepredicatex{\thecat}{\themat}\}
\cong \mathbb{C}^{\oplus \beta_1(\widehat{G})}
  \]
  of $\squarematrices{\thedim}{\comps}$ is defined by the rule which assigns to $($the one-elemental cohomology class of$)$ any $1$-cocycle $\thecoc$ the matrix $ (\thecocx{\theunix{}{\indj}{\indi}} )_{(\indj,\indi)\in \{1,\ldots,\thedim\}^{\Ssetmonoidalproduct 2}}$,
 where
  $\beta_1 \bigl(\widehat{G} \bigr)$ and for any $\themat\in \squarematrices{\thedim}{\comps}$ the predicate $\thepredicatex{\thecat}{\themat}$ are as follows:
  \begin{table}[h]\centering\renewcommand{\arraystretch}{1.2}
  \begin{tabular}{p{0.2\textwidth}|p{0.48\textwidth}|p{0.2\textwidth}}
    If $\thecat$ is \ldots, & then $\thepredicatex{\thecat}{\themat}$ is ``\ldots'' & and $\beta_1 \bigl(\widehat{G} \bigr)$ is \ldots \\ \hline
$\smallblocks \Sand \largeblocks\Sand\neutralblocks$                            & $\top$ & $\thedim^2$  \Tstrut \\
$\smallblocks \Sand \neg\largeblocks\Sand\neutralblocks\Sand\neutralpartitions$                           & $\exists\thesum\in\comps\quantorpredicate\themat-\thesum\identitymatrix$ is small & $(\thedim-1)^2+1$  \\
$\smallblocks \Sand \neg\largeblocks\Sand\neutralblocks\Sand\neg\neutralpartitions$                           &  $\themat$ is small & $(\thedim-1)^2$  \\
$\smallblocks \Sand \largeblocks\Sand\neg\neutralblocks\Sand \neutralpartitions$                           & $\exists\thesum\in\comps\quantorpredicate\themat-\thesum\identitymatrix$ is skew-sym\-met\-ric & $\onehalf\thedim(\thedim-1)+1$  \\
$\smallblocks \Sand \largeblocks\Sand\neg\neutralblocks\Sand \neg\neutralpartitions$                           &  $\themat$ is skew-sym\-met\-ric & $\onehalf\thedim(\thedim-1)$  \\
$\smallblocks \Sand \neg\largeblocks\Sand\neg\neutralblocks\Sand \neutralpartitions$                           & $\exists\thesum\in\comps\quantorpredicate\themat-\thesum\identitymatrix$ is skew-sym\-met\-ric and small  & $\onehalf(\thedim-1)(\thedim-2)+1$  \\
$\smallblocks \Sand \neg\largeblocks\Sand\neg\neutralblocks\Sand \neg\neutralpartitions$                           & $\themat$ is skew-sym\-met\-ric and small & $\onehalf(\thedim-1)(\thedim-2)$  \\
$\neg\smallblocks \Sand \largeblocks\Sand\neutralblocks\Sand \neutralpartitions$                           & $\themat$ is diagonal & $\thedim$  \\
$\neg\smallblocks \Sand\neg\neutralblocks\Sand \neutralpartitions$                        &  $\exists\thesum\in\comps\quantorpredicate\themat-\thesum\identitymatrix=0$ & $1$  \\
$\neg\smallblocks \Sand\neg\neutralblocks\Sand \neg\neutralpartitions$      & $\themat=0$ & $0$
  \end{tabular}
  \end{table}

  \noindent
  And these are all the cases that can occur.
\end{theorem*}
  }

\subsection{Structure of the article}
\label{section:introduction-structure-of-the-article}
Excluding the introduction, the article is divided into five sections.
\begin{enumerate}[label=$\square$]\itemsep=0pt
\item Section~\ref{section:quantum-groups} recalls the definitions of compact quantum groups and the quantum group cohomology with trivial coefficients of their discrete duals.
\item
Following that, Section~\ref{section:easy} provides particular examples of compact quantum groups by presenting the definitions of categories of two-co\-lo\-red partitions and unitary easy quantum groups.
\item
For the convenience of the reader, the definition of the first Hochschild cohomology and important results about it are recalled in Section~\ref{section:hochschild}.
\item
Section~\ref{section:dimensions} defines the vector spaces of matrices appearing in the main result and computes their dimensions.
\item
The proof of the main theorem is contained in Section~\ref{section:main-one}. Starting from a characterization of the first cohomology of a universal algebra recalled in Section~\ref{section:hochschild}  the first quantum group cohomology as defined in Section~\ref{section:quantum-groups} is computed of the discrete duals of the quantum groups defined in Section~\ref{section:easy}.
\end{enumerate}

\subsection{Notation}
\label{section:introduction-notation}

\makeatletter
\providecommand*{\cupdot}{%
  \mathbin{%
    \mathpalette\@cupdot{}%
  }%
}
\newcommand*{\@cupdot}[2]{%
  \ooalign{%
    $\m@th#1\cup$\cr
    \hidewidth$\m@th#1\cdot$\hidewidth
  }%
}
\makeatother

\makeatletter
\providecommand*{\bigcupdot}{%
  \mathbin{%
    \mathpalette\@bigcupdot{}%
  }%
}
\newcommand*{\@bigcupdot}[2]{%
  \ooalign{%
    $\m@th#1\bigcup$\cr
    \hidewidth$\m@th#1\cdot$\hidewidth
  }%
}
\makeatother

{
  \newcommand{\domspace}{V}
  \newcommand{\codspace}{W}
  \newcommand{\someset}{E}
  \newcommand{\somespace}{X}
  \newcommand{\somefield}{\mathbb{K}}
  \newcommand{\somerels}{R}
  \newcommand{\somemap}{f}
  \newcommand{\anyset}{S}
In the following, $0\notin \pint$. Rather, $\nnint=\pint\cupdot\{0\}$. Let $\dwi{k}:= \{i\in \pint\Sand i\leq k  \}$ for any $k\in\nnint$, in particular, $\dwi{0}=\varnothing$. The symbol $\Sdirectproduct$ will denote the Cartesian product of sets, with the convention $\anyset^{\Ssetmonoidalproduct 0}:= \{\varnothing\}$ for any set $\anyset$. Throughout, all algebras are meant to be associative and unital. The symbols $\Sleftmoduleaction$ and $\Srightmoduleaction$ are used to denote the left  respectively right actions of any algebra on any bimodule. Moreover, given any vector spaces $\domspace$ and $\codspace$ over any field the symbol $\SinternalhomX{\domspace}{\codspace}$ will stand for the vector space of linear maps from $\domspace$ to $\codspace$. Furthermore, for any vector space~$\somespace$ and any (possibly infinite) set  $\someset$ the notation $\somespace^{\Sdirectproduct\someset}$ will be used for the $\someset$-fold direct product  vector space of $\somespace$ (not to be confused with the direct sum $\somespace^{\Sdirectsum\someset}$). For any field $\somefield$ and any set~$\someset$, the free $\somefield$-algebra over $\someset$ will be denoted by $\freealg{\somefield}{ \someset}$. For any $\somerels\subseteq \freealg{\somefield}{ \someset}$, we will write $\univalg{\somefield}{\someset}{\somerels}$ for the universal $\somefield$-algebra with generators $\someset$ and relations $\somerels$.
}

\section{Quantum groups and their cohomology}
\label{section:quantum-groups}
The most general kind of ``quantum group'' in the sense considered here are the locally compact quantum groups introduced by Kustermans and Vaes in \cite{Kustermans2001, Kustermans2005, KustermansVaes2000, KustermansVaes2003}. Two subcategories of these are Woronowicz's compact quantum groups defined in \cite{Woronowicz1987b, Woronowicz1991, Woronowicz1998} and Van Daele's discrete quantum groups studied in \cite{VanDaele1996b, VanDaele1998}.

While of those two each is equivalent to the dual category of the other via Pontryagin duality, it is customary to ascribe the cohomology discussed in the present article to the discrete quantum group rather than its compact dual in order to preserve the analogy with the group case. At the same time, the particular quantum groups treated in this article are usually considered to be compact rather than discrete.

And it is in fact most convenient for the purpose of the present article to adopt the latter perspective and work with compact quantum groups. The fact that the quantum group cohomology is actually that of discrete quantum groups will be glossed over by only giving the definition of the composition of the cohomology functor with the Pontryagin transformation. However, the custom will be respected when it comes to notation.

\subsection{Compact quantum groups}
\label{section:quantum-groups-definition}
{
  \newcommand{\thealgebra}{A}
  \newcommand{\thefund}{\theunim{}}
  \newcommand{\thecomodule}{M}
  \newcommand{\thefundwh}{\theunim{\whpoint}}
  \newcommand{\thefundbl}{\theunim{\blpoint}}
  Quantum groups can be defined both on an analytic, namely von-Neumann- or $C^\ast$-algebraic level, and on a purely algebraic level. For the purposes of discussing quantum group cohomology, it fully suffices to consider the latter definition, given in \cite{DijkhuizenKoornwinder1994}. In that sense, an (\emph{algebraic}) \emph{compact quantum group} $\theqg$ is the formal dual of a Hopf $\ast$-algebra $\comps\bigl[\CQGdualC{\theqg}\bigr]$ which admits a faithful positive integral. A big class of examples is provided in Section~\ref{section:easy}. For any $\theqg$, in the present article $\comps \bigl[\CQGdualC{\theqg} \bigr]$ is generated as a $\ast$-algebra by the matrix coefficients of a single finite-dimensional unitary comodule $\thecomodule$. The coefficient matrix $\thefund$ of a choice of such an $\thecomodule$ is often called a \emph{fundamental representation}.
The axioms imply in particular that, if $\thefundbl$ is the matrix of the conjugate comodule of $\thecomodule$ and if $\thefundwh:= \thefund$, then $\comps \bigl[\CQGdualC{\theqg} \bigr]$ is generated as an algebra (as opposed to a $\ast$-algebra) by the union of the entries of $\thefundwh$ and $\thefundbl$. If $\thealgebra$ is the underlying algebra and $\thecounit$ the counit of the Hopf $\ast$-algebra  $\comps \bigl[\CQGdualC{\theqg} \bigr]$, it is entirely sufficient to think of $\theqg$ as the augmented algebra $(\thealgebra,\thecounit)$ and keep in mind that for the examples in this article a generating set of $\thealgebra$ can be given consisting of the entries of two matrices $\thefundwh$ and $\thefundbl$ of the same size.
}
\subsection{Quantum group cohomology}
\label{section:quantum-groups-cohomology}
{
  \newcommand{\thederivs}{\mathrm{Der}_\epsilon(A)}
  \newcommand{\anycoc}{\eta}
  \newcommand{\leftop}{a_1}
  \newcommand{\rightop}{a_2}
  \newcommand{\anygen}{e}
  \newcommand{\thetuple}{x}
  \newcommand{\thetuplex}[1]{\thetuple_{#1}}
  \newcommand{\genels}[1]{a_{#1}}
  \newcommand{\theorder}{p}
  \newcommand{\thealgebra}{A}
  \newcommand{\themod}{{}_\thecounit\comps_\thecounit}
  \newcommand{\somescal}{\lambda}
  \newcommand{\someop}{a}
One of many equivalent ways of introducing  quantum group cohomology is via Hochschild cohomology.
  For any compact quantum group $\theqg$ and any $\theorder\in\nnint$, if $\thealgebra$ is the underlying algebra and~$\thecounit$ the counit of  $\DQGgroupalgebra{\CQGdualC{\theqg}}$ and if $\themod$ denotes the $\thealgebra$-bimodule given by the $\comps$-vector space $\comps$ equipped with the left and right $\thealgebra$-actions defined by $\someop\Smonoidalproduct\somescal\mapsto \thecounit(\someop)\somescal$ respectively $\somescal\Smonoidalproduct\someop\mapsto \somescal\thecounit(\someop)$ for any $\someop\in\thealgebra$ and $\somescal\in\comps$,
  then the \emph{$\theorder$-th quantum group cohomology with trivial coefficients} of the discrete dual $\CQGdualC{\theqg}$ of $\theqg$ is defined as
  \[
    \DQGcohomologyTC{\theorder}{\CQGdualC{\theqg}}:=\HScohomology{\theorder}{\thealgebra}{\themod},
  \]
  the $\theorder$-th Hochschild cohomology of $\thealgebra$ with coefficients in $\themod$.
}

\section[Categories of two-colored partitions and unitary easy quantum groups]{Categories of two-colored partitions\\ and unitary easy quantum groups}
\label{section:easy}

The quantum groups whose quantum group cohomology is investigated in the present article are the discrete duals of so-called easy quantum groups. They can be defined via Tannaka--Krein duality (see \cite{Woronowicz1988}) using the combinatorics of so-called two-colored partitions (which will be explained in Definition~\ref{definition:two-colored-partition}). Throughout the article, it will be important to distinguish the notions of a two-co\-lo\-red partition and a set-theo\-re\-ti\-cal partition in the following sense.
\begin{Notation}
  \label{notation:set-theoretical-partitions}
  \newcommand{\theset}{X}
  Let $\theset$ be any set.
  \begin{enumerate}
  \item\label{notation:set-theoretical-partitions-1}
    {
      \newcommand{\thepartition}{p}
      A \emph{set-theo\-re\-ti\-cal partition} of $X$ is the quotient set of any equivalence relation on $X$ or, equivalently,  any set of non-empty pairwise disjoint subsets of $X$ whose union is $X$.
      }
  \item\label{notation:set-theoretical-partitions-2}
    {
      \newcommand{\smallpart}{p}
      \newcommand{\largepart}{q}
      \newcommand{\smallblock}{\Yb{B}}
      \newcommand{\largeblock}{\Yb{C}}
      \newcommand{\leftpartition}{p_1}
      \newcommand{\rightpartition}{p_2}
      \newcommand{\thejoin}{s}
      \newcommand{\thepartition}{p}
      \newcommand{\anyel}{\Yp{x}}
      \newcommand{\block}{\Yb{B}}
      Given any two set-theo\-re\-ti\-cal partitions $\smallpart$ and $\largepart$ of  $\theset$, write $\smallpart\finerthan\largepart$ if $\smallpart$ is \emph{finer} than $\largepart$, i.e., if for any $\smallblock\in\smallpart$ there exists $\largeblock\in\largepart$  with $\smallblock\subseteq \largeblock$.
      }
  \item\label{notation:set-theoretical-partitions-3}
    {
      \newcommand{\smallpart}{p}
      \newcommand{\largepart}{q}
      \newcommand{\smallblock}{\Yb{B}}
      \newcommand{\largeblock}{\Yb{C}}
      \newcommand{\leftpartition}{p_1}
      \newcommand{\rightpartition}{p_2}
      \newcommand{\thejoin}{s}
      \newcommand{\thepartition}{p}
      \newcommand{\anyel}{\Yp{x}}
      \newcommand{\block}{\Yb{B}}
      For any two set-theo\-re\-ti\-cal partitions $\smallpart$ and $\largepart$ of $\theset$, let $\zetfx{\smallpart}{\largepart}:= 1$ if $\smallpart\finerthan\largepart$ and  let $\zetfx{\smallpart}{\largepart}:= 0$ otherwise.
    }
  \item\label{notation:set-theoretical-partitions-4}
    {
      \newcommand{\smallpart}{p}
      \newcommand{\largepart}{q}
      \newcommand{\smallblock}{\Yb{B}}
      \newcommand{\largeblock}{\Yb{C}}
      \newcommand{\leftpartition}{p_1}
      \newcommand{\rightpartition}{p_2}
      \newcommand{\thejoin}{s}
      \newcommand{\thepartition}{p}
      \newcommand{\anyel}{\Yp{x}}
      \newcommand{\block}{\Yb{B}}
      Furthermore, for any set-theo\-re\-ti\-cal partitions $\leftpartition$ and $\rightpartition$ of $\theset$ the \emph{join} of $\leftpartition$ and $\rightpartition$ is the unique  set-theo\-re\-ti\-cal partition $\thejoin$ of $\theset$ which satisfies $\leftpartition\finerthan\thejoin$ and $\rightpartition\finerthan\thejoin$ and which is minimal with that property with respect to the partial order $\finerthanrel$.
    }
    \item\label{notation:set-theoretical-partitions-5}
    {
      \newcommand{\themap}{f}
      \newcommand{\domset}{X}
      \newcommand{\codset}{Y}
      \newcommand{\domsubset}{A}
      \newcommand{\codsubset}{B}
      \newcommand{\domel}{x}
      \newcommand{\codel}{y}
      For any mapping $\xfromto{\themap}{\domset}{\codset}$ from $\domset$ to any set $\codset$ and for any subset $\codsubset\subseteq \codset$, let $\spimgx{\themap}{\codsubset}:= \{\domel\in\domset\Sand \themap(\domel)\in\codsubset\}$ denote the \emph{pre-image of $\codsubset$ under $\themap$}. Moreover, let $\sranx{\themap}:= \{\themap(\domel)\classpredicate \domel\in\domset\}$ and $\kerp{\themap}:= \{\spimgx{\themap}{\{\codel\}}\classpredicate \codel\in\sranx{\themap}\}$ be the \emph{image} and \emph{kernel} of $\themap$, respectively.
    }
  \item\label{notation:set-theoretical-partitions-6}
    {
      \newcommand{\themap}{f}
      \newcommand{\domset}{X}
      \newcommand{\codset}{Y}
      \newcommand{\anyel}{\Yp{x}}
      \newcommand{\thepartition}{p}
      \newcommand{\block}{\Yb{B}}
      For any set-theo\-re\-ti\-cal partition  $\thepartition$ of $\domset$, write $\blofun{\thepartition}$ for the associated \emph{projection}, the mapping $\fromto{\domset}{\thepartition}$ which maps any  $\anyel\in\domset$ to the unique $\block\in\thepartition$ with $\anyel\in\block$. And for any second set~$\codset$ and any mapping $\xfromto{\themap}{\domset}{\codset}$ with $\thepartition\finerthan\kerp{\themap}$, let $\tsquomap{\themap}{\thepartition}$ denote the \emph{quotient mapping}, the unique mapping $\fromto{\thepartition}{\codset}$ with $(\tsquomap{\themap}{\thepartition})\Scomposition \blofun{\thepartition}=\themap$.
    }
  \end{enumerate}
\end{Notation}

\begin{Example}
      \newcommand{\themap}{f}
      \newcommand{\domset}{X}
      \newcommand{\codset}{Y}
      \newcommand{\thepartition}{p}
      \newcommand{\othermap}{g}
      \newcommand{\block}{\Yb{C}}
      If $\domset=\{1,2,3,4,5,6\}$, then $\thepartition=\{\{1\},\{2,4\},\{3,5,6\}\}$ is a set-theoretical partition of $\domset$ and the projection $\blofun{\thepartition}$ is the mapping $\domset\to \thepartition$, which sends $1$ to $\{1\}$, sends both~$2$ and~$4$ to $\{2,4\}$ and sends each of $3$, $5$ and $6$ to $\{3,5,6\}$.

      Moreover, if  $\codset=\{a,b,c,d\}$ and $|\codset|=4$ and if $\xfromto{\themap}{\domset}{\codset}$ maps each of $1$, $2$ and $4$ to $a$ and each of $3$, $5$ and $6$ to $c$, then the kernel of $\themap$ is $\kerp{\themap}=\{\{1,2,4\},\{3,5,6\}\}$. Since then $\thepartition\finerthan\kerp{\themap}$ the quotient mapping $\tsquomap{\themap}{\thepartition}$ exists and maps both  $\{1\}$ and $\{2,4\}$ to $a$ and~$\{3,5,6\}$ to~$c$.

      In contrast, if $\themap(4)$ was not given by $a$ but by $b$ instead, then $\kerp{\themap}$ would equal $\{\{1,2\},\{4\},\break\{3,5,6\}\}$, in which case $\thepartition$ would not be finer than  $\kerp{\themap}$ since there would be no $\block\in\kerp{\themap}$ with $\{2,4\}\subseteq \block$. There would be no $\tsquomap{\themap}{\thepartition}$ with $(\tsquomap{\themap}{\thepartition})\Scomposition\blofun{\thepartition}=\themap$.
\end{Example}

\subsection{Two-colored partitions and their categories}
\label{section:easy-partitions}
Two-co\-lo\-red partitions can be defined as follows. For further details see \cite{TarragoWeber2018}, where they  were first introduced, generalizing the (uncolored) ``partitions'' considered in \cite{BanicaSpeicher2009}.
\begin{Assumptions}\quad
  \begin{enumerate}\itemsep=0pt
  \item Let $\upp{(\argph)}$ and $\lop{(\argph)}$ be any two injections with common domain  $\pint$ and with disjoint ranges.
  \item  Let $\whpoint$ and $\blpoint $ be arbitrary with $\whpoint\neq \blpoint$.
  \end{enumerate}
\end{Assumptions}

\begin{Definition}\quad
  \label{definition:points-colors}
  \newcommand{\inlen}{k}
  \newcommand{\outlen}{\ell}
  \newcommand{\inind}{a}
  \newcommand{\outind}{b}
  \newcommand{\codset}{X}
  \newcommand{\upmap}{g}
  \newcommand{\lomap}{j}
  \newcommand{\incol}{\Yc{c}}
  \newcommand{\incolx}[1]{\incol_{#1}}
  \newcommand{\outcol}{\Yc{d}}
  \newcommand{\outcolx}[1]{\outcol_{#1}}
  \newcommand{\thepartition}{p}
  \begin{enumerate}
\item\label{definition:points-colors-1} For any $\{\inlen,\outlen\}\subseteq \nnint$, we call $\tsopx{\inlen}{\outlen}:= \{\upp{\inind},\lop{\outind}\classpredicate \inind\in\dwi{\inlen}\Sand \outind\in\dwi{\outlen}\}$ the set of $\inlen$ upper and $\outlen$ lower \emph{points}.
\item\label{definition:points-colors-2} Given any $\{\inlen,\outlen\}\subseteq \nnint$, any set $\codset$ and any mappings $\xfromto{\upmap}{\dwi{\inlen}}{\codset}$ and $\xfromto{\lomap}{\dwi{\outlen}}{\codset}$ denote by~$\djp{\upmap}{\lomap}$ the mapping $\fromto{\tsopx{\inlen}{\outlen}}{\codset}$ with $\upp{\inind}\mapsto \upmap(\inind)$ for any $\inind\in\dwi{\inlen}$ and $\lop{\outind}\mapsto \lomap(\outind)$ for any $\outind\in\dwi{\outlen}$.
\item\label{definition:points-colors-3} $\whpoint$ and $\blpoint$ are called the two \emph{colors} and are said to  be \emph{dual} to each other, in symbols, $\overline{\whpoint}:= \blpoint$ and $\overline{\blpoint}:= \whpoint$. They moreover have the \emph{color values} $\sigma(\whpoint):= 1$ and $\sigma(\blpoint):= -1$.
  \item\label{definition:points-colors-4} For any $\{\inlen,\outlen\}\subseteq \nnint$, any $\xfromto{\incol}{\dwi{\inlen}}{\blaw}$ and any $\xfromto{\outcol}{\dwi{\outlen}}{\blaw}$ the \emph{color sum} of $(\incol,\outcol)$ is the $\integers$-valued measure $\csfunx{\incol}{\outcol}$ on $\tsopx{\inlen}{\outlen}$ with density $-\sigma(\incolx{\inind})$ on $\upp{\inind}$ for any $\inind\in\dwi{\inlen}$ and density $\sigma(\outcolx{\outind})$ on $\lop{\outind}$ for any $\outind\in\dwi{\outlen}$. Moreover, $\tcsx{\incol}{\outcol}:= \csx{\incol}{\outcol}{\tsopx{\inlen}{\outlen}}$ is called the \emph{total color sum} of  $(\incol,\outcol)$.
  \item\label{definition:points-colors-5} For brevity, let $\Slength{\incol}:= \inlen$ for any $\inlen\in \nnint$ and any $\xfromto{\incol}{\dwi{\inlen}}{\blaw}$.
  \end{enumerate}
\end{Definition}

\begin{Example}
  \newcommand{\anypoints}{\Yb{S}}
  \newcommand{\inlen}{k}
  \newcommand{\outlen}{\ell}
  \newcommand{\incol}{\Yc{c}}
  \newcommand{\incolx}[1]{\incol_{#1}}
  \newcommand{\outcol}{\Yc{d}}
  \newcommand{\outcolx}[1]{\outcol_{#1}}
  \newcommand{\anypoint}{\Yp{x}}
  Consider $\xfromto{\incol}{\dwi{3}}{\blaw}$ and $\xfromto{\outcol}{\dwi{4}}{\blaw}$ with $\incolx{2}=\incolx{3}=\circ$ and $\incolx{1}=\bullet$ and $\outcolx{1}=\outcolx{2}=\outcolx{4}=\circ$  and $\outcolx{3}=\bullet$.
  $$
\begin{tikzpicture}
    \def\scp{0.666}
    \def\linksize{\scp*0.075cm}
    \def\pointsize{\scp*0.25cm}
    \def\dd{\scp*0.5cm}
    \def\dx{\scp*1cm}
    \def\cx{\scp*0.4cm}
    \def\txu{2*\dx}
    \def\txl{3*\dx}
    \def\dy{\scp*1cm}
    \def\cy{\scp*0.4cm}
    \def\ty{2*\dy}
    \tikzset{whp/.style={circle, inner sep=0pt, text width={\pointsize}, draw=black, fill=white}}
    \tikzset{blp/.style={circle, inner sep=0pt, text width={\pointsize}, draw=black, fill=black}}
    \tikzset{lk/.style={regular polygon, regular polygon sides=4, inner sep=0pt, text width={\linksize}, draw=black, fill=black}}
    \draw[dotted] ({0-\dd},{0}) -- ({\txl+\dd},{0});
    \draw[dotted] ({0-\dd},{\ty}) -- ({\txu+\dd},{\ty});
    \coordinate (l1) at ({0+0*\dx},{0+0*\ty}) {};
    \coordinate (l2) at ({0+1*\dx},{0+0*\ty}) {};
    \coordinate (l3) at ({0+2*\dx},{0+0*\ty}) {};
    \coordinate (l4) at ({0+3*\dx},{0+0*\ty}) {};
    \coordinate (u1) at ({0+0*\dx},{0+1*\ty}) {};
    \coordinate (u2) at ({0+1*\dx},{0+1*\ty}) {};
    \coordinate (u3) at ({0+2*\dx},{0+1*\ty}) {};
    \node[whp] at (l1) {};
    \node[whp] at (l2) {};
    \node[blp] at (l3) {};
    \node[whp] at (l4) {};
    \node[blp] at (u1) {};
    \node[whp] at (u2) {};
    \node[whp] at (u3) {};
    \node[above =7pt of u1] {$\incolx{1}$};
    \node[above =7pt of u2] {$\incolx{2}$};
    \node[above =7pt of u3] {$\incolx{3}$};
    \node[below =7pt of l1] {$\outcolx{1}$};
    \node[below =7pt of l2] {$\outcolx{2}$};
    \node[below =7pt of l3] {$\outcolx{3}$};
    \node[below =7pt of l4] {$\outcolx{4}$};
    \node[below =2pt of u1] {$\scriptstyle 1$};
    \node[below =2pt of u2] {$\scriptstyle -1$};
    \node[below =2pt of u3] {$\scriptstyle -1$};
    \node[above =2pt of l1] {$\scriptstyle 1$};
    \node[above =2pt of l2] {$\scriptstyle 1$};
    \node[above =2pt of l3] {$\scriptstyle -1$};
    \node[above =2pt of l4] {$\scriptstyle 1$};
    \draw[densely dashed] ($(u3)+({-1*\cx},{1*\cy})$) to ($(l3)+({-1*\cx},{-1*\cy})$) to ($(l4)+({1*\cx},{-1*\cy})$) to ($(l4)+({1*\cx},{1*\cy})$) to ($(u3)+({1*\cx},{-1*\cy})$) to ($(u3)+({1*\cx},{1*\cy})$) -- cycle;
    \node at ($(l4)+({.75*\dx},{1*\dy})$) {$\anypoints$};
    \node at ({-.5*\dx-\dd},{1*\dy}) {$\csfunx{\incol}{\outcol}$};
  \end{tikzpicture}
$$
  The color sum $\csfunx{\incol}{\outcol}$ has density $1$ at each of $\upp{1}$, $\lop{1}$, $\lop{2}$ and $\lop{4}$ and density $-1$ at each of  $\upp{2}$, $\upp{3}$ and~$\lop{3}$. Consequently, the subset $\anypoints=\{\upp{3},\lop{3},\lop{4}\}$ of $\tsopx{3}{4}$  has color sum $\csx{\incol}{\outcol}{\anypoints}=\csx{\incol}{\outcol}{\{\upp{3}\}}+\csx{\incol}{\outcol}{\{\lop{3}\}}+\csx{\incol}{\outcol}{\{\lop{4}\}}=-1-1+1=-1$. The total color sum is $\tcsx{\incol}{\outcol}=1$.
\end{Example}

\begin{Definition}\quad
  \label{definition:two-colored-partition}
  \begin{enumerate}
  \item {
  \newcommand{\inlen}{k}
  \newcommand{\outlen}{\ell}
  \newcommand{\inind}{i}
  \newcommand{\outind}{j}
  \newcommand{\codset}{X}
  \newcommand{\upmap}{g}
  \newcommand{\lomap}{j}
  \newcommand{\incol}{\Yc{c}}
  \newcommand{\incolx}[1]{\incol_{#1}}
  \newcommand{\outcol}{\Yc{d}}
  \newcommand{\outcolx}[1]{\outcol_{#1}}
  \newcommand{\thepartition}{p}
  A \emph{two-co\-lo\-red partition} is any triple $(\incol,\outcol,\thepartition)$ for which there exist $\{\inlen,\outlen\}\subseteq \nnint$ such that $\incol$ and $\outcol$ are mappings from $\dwi{\inlen}$ respectively $\dwi{\outlen}$ to $\blaw$, the upper and lower \emph{colorings}, and such that $\thepartition$, the collection of \emph{blocks}, is a  set-theo\-re\-ti\-cal partition of
  the set~$\tsopx{\inlen}{\outlen}$
  of points.
  
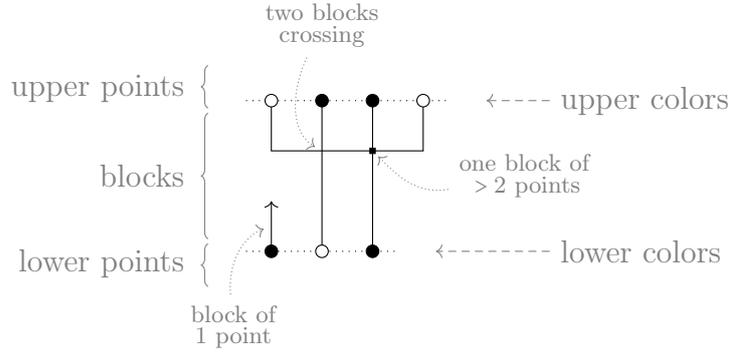
\begin{center}
          \hspace{-3em}
\begin{tikzpicture}[baseline=0.91cm]
    \def\scp{0.666}
    \def\linksize{\scp*0.075cm}
    \def\pointsize{\scp*0.25cm}
    \def\dd{\scp*0.5cm}
    \def\dx{\scp*1cm}
    \def\cx{\scp*0.3cm}
    \def\txu{3*\dx}
    \def\txl{2*\dx}
    \def\dy{\scp*1cm}
    \def\cy{\scp*0.3cm}
    \def\ty{3*\dy}
    \tikzset{whp/.style={circle, inner sep=0pt, text width={\pointsize}, draw=black, fill=white}}
    \tikzset{blp/.style={circle, inner sep=0pt, text width={\pointsize}, draw=black, fill=black}}
    \tikzset{lk/.style={regular polygon, regular polygon sides=4, inner sep=0pt, text width={\linksize}, draw=black, fill=black}}
    \draw[dotted] ({0-\dd},{0}) -- ({\txl+\dd},{0});
    \draw[dotted] ({0-\dd},{\ty}) -- ({\txu+\dd},{\ty});
    \coordinate (l1) at ({0+0*\dx},{0+0*\ty}) {};
    \coordinate (l2) at ({0+1*\dx},{0+0*\ty}) {};
    \coordinate (l3) at ({0+2*\dx},{0+0*\ty}) {};
    \coordinate (u1) at ({0+0*\dx},{0+1*\ty}) {};
    \coordinate (u2) at ({0+1*\dx},{0+1*\ty}) {};
    \coordinate (u3) at ({0+2*\dx},{0+1*\ty}) {};
    \coordinate (u4) at ({0+3*\dx},{0+1*\ty}) {};
    \node[lk] at ({2*\dx},{2*\dy}) {};
    \draw (u1) -- ++ ({0*\dx},{-1*\dy}) -| (u4);
    \draw[->] (l1) -- ++({0*\dx},{1*\dy});
    \draw (l2) -- (u2);
    \draw (l3) -- (u3);
    \node[blp] at (l1) {};
    \node[whp] at (l2) {};
    \node[blp] at (l3) {};
    \node[whp] at (u1) {};
    \node[blp] at (u2) {};
    \node[blp] at (u3) {};
    \node[whp] at (u4) {};
    \node[anchor=west, gray] (upper-colors) at ({5.5*\dx},{3*\dy}) {upper colors $\incol$};
    \node[anchor=east, gray] (upper-points) at ({-1.5*\dx},{3.25*\dy}) {$\inlen=4$ upper points};
    \node[anchor=west, gray] (lower-colors) at ({5.5*\dx},{0*\dy}) {lower colors $\outcol$};
    \node[anchor=east, gray] (lower-points) at ({-1.5*\dx},{-0.25*\dy}) {$\outlen=3$ lower points};
    \node[anchor=east, gray, align=left] (blocks) at ({-1.5*\dx},{1.5*\dy}) { $\thepartition=\{\{\lop{1}\},\{\upp{2},\lop{2}\},$ \\ \hphantom{$\thepartition=\{$}$\{\upp{1},\upp{3},\upp{4},\lop{3}\}\}$};
    \node[align=center, anchor=center, gray] (crossing-blocks) at ({1*\dx},{4.5*\dy}) {\scriptsize  blocks $\{\upp{1},\upp{3},\upp{4},\lop{3}\}$\\[-0.5em] \scriptsize and $\{\upp{2},\lop{2}\}$ crossing};
    \node[align=center, anchor=center, gray] (big-block) at ({6*\dx},{1.5*\dy}) {\scriptsize block $\{\upp{1},\upp{3},\upp{4},\lop{3}\}$};
    \node[align=center, anchor=center, gray] (singleton-block) at ({-0.75*\dx},{-1.5*\dy}) {\scriptsize block $\{\lop{1}\}$};
    \draw[->, densely dotted, out=95, in=200, shorten >= 3pt, gray] (singleton-block) to ({-0*\dx},{0.5*\dy});
    \draw[->, densely dotted, out=180, in=315, shorten >= 3pt, gray] (big-block) to ({2*\dx},{2*\dy});
    \draw[->, densely dotted, out=245, in=150, shorten >= 3pt, gray] (crossing-blocks) to ({1*\dx},{2*\dy});
    \draw[xshift={-0.75*\dx}, decorate, decoration={brace}, gray] ({-0.5*\dx},{0.25*\dy}) to ({-0.5*\dx},{2.75*\dy});
    \draw[xshift={-0.75*\dx}, decorate, decoration={brace}, gray] ({-0.5*\dx},{-0.7*\dy}) to ({-0.5*\dx},{0.15*\dy});
    \draw[xshift={-0.75*\dx}, decorate, decoration={brace}, gray] ({-0.5*\dx},{2.85*\dy}) to ({-0.5*\dx},{3.7*\dy});
    \draw[->, densely dashed, out=180, in=0, shorten >= 0pt, gray] (upper-colors) to ({3*\dx+2.5*\dd},{3*\dy});
    \draw[->, densely dashed, out=180, in=0, shorten >= 0pt, gray] (lower-colors) to ({2*\dx+2.5*\dd},{0*\dy});
  \end{tikzpicture}
  \end{center}
 }
\item
  {
    Any set $\thecat$ of two-co\-lo\-red partitions meeting the following conditions is called a \emph{category of two-co\-lo\-red partitions}:
  \begin{enumerate}\itemsep=0pt
  \item $\thecat$ contains $\PartIdenW$, $\PartIdenB$, $\PartIdenLoWB$, $\PartIdenUpBW$, $\PartIdenLoBW$ and $\PartIdenUpWB$.
  \item
       {
          \newcommand{\inlen}{k}
          \newcommand{\outlen}{\ell}
            \newcommand{\incol}{\Yc{c}}
            \newcommand{\outcol}{\Yc{d}}
          \newcommand{\thepartition}{p}
          \newcommand{\block}{\Yb{B}}
          \newcommand{\inind}{a}
          \newcommand{\outind}{b}
          $\thecat$ is closed under forming adjoints, that is, horizontal reflection. More precisely,
          $(\outcol,\incol,\thepartition\Sadj)\in\thecat$ for any  $(\incol,\outcol,\thepartition)\in \thecat$, where, if $\{\inlen,\outlen\}\subseteq \nnint$ are such that $\thepartition$ is a set-theo\-re\-ti\-cal partition of $\tsopx{\inlen}{\outlen}$, then
$\thepartition\Sadj:=\{ \{\upp{\outind}\classpredicate \outind\in\dwi{\outlen}\Sand \lop{\outind}\in\block\}\cupdot\{\lop{\inind}\classpredicate \inind\in\dwi{\inlen}\Sand \upp{\inind}\in\block\} \}_{\block\in\thepartition}$ is the \emph{adjoint}
of $\thepartition$.
$$
  \left(
              \begin{tikzpicture}[baseline=0.666cm]
    \def\scp{0.666}
    \def\linksize{\scp*0.075cm}
    \def\pointsize{\scp*0.25cm}
    \def\dd{\scp*0.5cm}
    \def\dx{\scp*1cm}
    \def\cx{\scp*0.3cm}
    \def\txu{2*\dx}
    \def\txl{0*\dx}
    \def\dy{\scp*1cm}
    \def\cy{\scp*0.3cm}
    \def\ty{2*\dy}
    \tikzset{whp/.style={circle, inner sep=0pt, text width={\pointsize}, draw=black, fill=white}}
    \tikzset{blp/.style={circle, inner sep=0pt, text width={\pointsize}, draw=black, fill=black}}
    \tikzset{lk/.style={regular polygon, regular polygon sides=4, inner sep=0pt, text width={\linksize}, draw=black, fill=black}}
    \draw[dotted] ({0-\dd},{0}) -- ({\txl+\dd},{0});
    \draw[dotted] ({0-\dd},{\ty}) -- ({\txu+\dd},{\ty});
    \coordinate (l1) at ({0+0*\dx},{0+0*\ty}) {};
    \coordinate (u1) at ({0+0*\dx},{0+1*\ty}) {};
    \coordinate (u2) at ({0+1*\dx},{0+1*\ty}) {};
    \coordinate (u3) at ({0+2*\dx},{0+1*\ty}) {};
    \draw (u2) -- ++ ({0*\dx},{-1*\dy}) -| (u3);
    \draw (l1) -- (u1);
    \node[whp] at (l1) {};
    \node[whp] at (u1) {};
    \node[whp] at (u2) {};
    \node[blp] at (u3) {};
  \end{tikzpicture}
  \right)^\ast=
                \begin{tikzpicture}[baseline=0.666cm]
    \def\scp{0.666}
    \def\linksize{\scp*0.075cm}
    \def\pointsize{\scp*0.25cm}
    \def\dd{\scp*0.5cm}
    \def\dx{\scp*1cm}
    \def\cx{\scp*0.3cm}
    \def\txu{0*\dx}
    \def\txl{2*\dx}
    \def\dy{\scp*1cm}
    \def\cy{\scp*0.3cm}
    \def\ty{2*\dy}
    \tikzset{whp/.style={circle, inner sep=0pt, text width={\pointsize}, draw=black, fill=white}}
    \tikzset{blp/.style={circle, inner sep=0pt, text width={\pointsize}, draw=black, fill=black}}
    \tikzset{lk/.style={regular polygon, regular polygon sides=4, inner sep=0pt, text width={\linksize}, draw=black, fill=black}}
    \draw[dotted] ({0-\dd},{0}) -- ({\txl+\dd},{0});
    \draw[dotted] ({0-\dd},{\ty}) -- ({\txu+\dd},{\ty});
    \coordinate (l1) at ({0+0*\dx},{0+0*\ty}) {};
    \coordinate (l2) at ({0+1*\dx},{0+0*\ty}) {};
    \coordinate (l3) at ({0+2*\dx},{0+0*\ty}) {};
    \coordinate (u1) at ({0+0*\dx},{0+1*\ty}) {};
    \draw (l2) -- ++ ({0*\dx},{1*\dy}) -| (l3);
    \draw (l1) -- (u1);
    \node[whp] at (l1) {};
    \node[whp] at (l2) {};
    \node[blp] at (l3) {};
    \node[whp] at (u1) {};
  \end{tikzpicture}
$$
        }
          \item {
              \newcommand{\inlenI}[1]{k_{#1}}
              \newcommand{\outlenI}[1]{\ell_{#1}}
              \newcommand{\incolI}[1]{\Yc{c}_{#1}}
            \newcommand{\outcolI}[1]{\Yc{d}_{#1}}
            \newcommand{\thepartitionI}[1]{p_{#1}}
            \newcommand{\inind}{a}
            \newcommand{\outind}{b}
            \newcommand{\twoind}{t}
          \newcommand{\block}{\Yb{B}}            $\thecat$ is closed under tensor products, i.e., horizontal concatenation. Formally,  $({\incolI{1}\Smonoidalproduct\incolI{2}},\allowbreak \outcolI{1}\Smonoidalproduct\outcolI{2},\thepartitionI{1}\Smonoidalproduct\thepartitionI{2})\in\thecat$  for any $(\incolI{1},\outcolI{1},\thepartitionI{1})\in\thecat$ and $(\incolI{2},\outcolI{2},\thepartitionI{2})\in\thecat$, where, if $\inlenI{\twoind}$ and $\outlenI{\twoind}$ are such that $\thepartitionI{\twoind}$ is a set-theo\-re\-ti\-cal partition of $\tsopx{\inlenI{\twoind}}{\outlenI{\twoind}}$ for each $\twoind\in\dwi{2}$, then
          $\incolI{1}\Smonoidalproduct\incolI{2}\in\blaw^{\Ssetmonoidalproduct(\inlenI{1}+\inlenI{2})}$ is defined by  $\inind\mapsto \incolI{1}(\inind)$ if $\inind\leq \inlenI{1}$ and $\inind\mapsto \incolI{2}(\inind-\inlenI{1})$ if $\inlenI{1}<\inind$ and, analogously, $\outcolI{1}\Smonoidalproduct\outcolI{2}\in \blaw^{\Ssetmonoidalproduct(\outlenI{1}+\outlenI{2})}$ is defined by $\outind\mapsto \outcolI{1}(\outind)$ if $\outind\leq \outlenI{1}$ and $\outind\mapsto \outcolI{2}(\outind-\outlenI{1})$ if $\outlenI{1}<\outind$, and where $\thepartitionI{1}\Smonoidalproduct\thepartitionI{2}:=\thepartitionI{1}\cupdot \{ \{\upp{(\inlenI{1}+\inind)}\classpredicate \inind\in\dwi{\inlenI{2}}\Sand \upp{\inind}\in \block\}\cupdot \big\{\lop{(\outlenI{1}+\outind)}\classpredicate \outind\in\dwi{\outlenI{2}}\Sand \lop{\outind}\in \block\}\}_{\block\in\thepartitionI{2}}$ is the \emph{tensor product} of $(\thepartitionI{1},\thepartitionI{2})$.
          $$
    
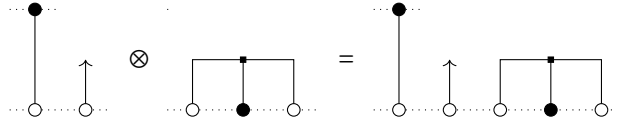
\begin{tikzpicture}[baseline=0.575cm]
    \def\scp{0.666}
    \def\linksize{\scp*0.075cm}
    \def\pointsize{\scp*0.25cm}
    \def\dd{\scp*0.5cm}
    \def\dx{\scp*1cm}
    \def\cx{\scp*0.3cm}
    \def\txu{0*\dx}
    \def\txl{1*\dx}
    \def\dy{\scp*1cm}
    \def\cy{\scp*0.3cm}
    \def\ty{2*\dy}
    \tikzset{whp/.style={circle, inner sep=0pt, text width={\pointsize}, draw=black, fill=white}}
    \tikzset{blp/.style={circle, inner sep=0pt, text width={\pointsize}, draw=black, fill=black}}
    \tikzset{lk/.style={regular polygon, regular polygon sides=4, inner sep=0pt, text width={\linksize}, draw=black, fill=black}}
    \draw[dotted] ({0-\dd},{0}) -- ({\txl+\dd},{0});
    \draw[dotted] ({0-\dd},{\ty}) -- ({\txu+\dd},{\ty});
    \coordinate (l1) at ({0+0*\dx},{0+0*\ty}) {};
    \coordinate (l2) at ({0+1*\dx},{0+0*\ty}) {};
    \coordinate (u1) at ({0+0*\dx},{0+1*\ty}) {};
    \draw[->] (l2) -- ++({0*\dx},{1*\dy});
    \draw (l1) -- (u1);
    \node[whp] at (l1) {};
    \node[whp] at (l2) {};
    \node[blp] at (u1) {};
  \end{tikzpicture}
\ {}  \Smonoidalproduct{}\
              \begin{tikzpicture}[baseline=0.575cm]
    \def\scp{0.666}
    \def\linksize{\scp*0.075cm}
    \def\pointsize{\scp*0.25cm}
    \def\dd{\scp*0.5cm}
    \def\dx{\scp*1cm}
    \def\cx{\scp*0.3cm}
    \def\txu{-1*\dx}
    \def\txl{2*\dx}
    \def\dy{\scp*1cm}
    \def\cy{\scp*0.3cm}
    \def\ty{2*\dy}
    \tikzset{whp/.style={circle, inner sep=0pt, text width={\pointsize}, draw=black, fill=white}}
    \tikzset{blp/.style={circle, inner sep=0pt, text width={\pointsize}, draw=black, fill=black}}
    \tikzset{lk/.style={regular polygon, regular polygon sides=4, inner sep=0pt, text width={\linksize}, draw=black, fill=black}}
    \draw[dotted] ({0-\dd},{0}) -- ({\txl+\dd},{0});
    \draw[dotted] ({0-\dd},{\ty}) -- ({\txu+\dd},{\ty});
    \coordinate (l1) at ({0+0*\dx},{0+0*\ty}) {};
    \coordinate (l2) at ({0+1*\dx},{0+0*\ty}) {};
    \coordinate (l3) at ({0+2*\dx},{0+0*\ty}) {};
    \node[lk] at ({1*\dx},{1*\dy}) {};
    \draw (l1) -- ++ ({0*\dx},{1*\dy}) -| (l2);
    \draw  ({1*\dx},{1*\dy}) -| (l3);
    \node[whp] at (l1) {};
    \node[blp] at (l2) {};
    \node[whp] at (l3) {};
  \end{tikzpicture}
  \ {}={} \
               \begin{tikzpicture}[baseline=0.575cm]
    \def\scp{0.666}
    \def\linksize{\scp*0.075cm}
    \def\pointsize{\scp*0.25cm}
    \def\dd{\scp*0.5cm}
    \def\dx{\scp*1cm}
    \def\cx{\scp*0.3cm}
    \def\txu{0*\dx}
    \def\txl{4*\dx}
    \def\dy{\scp*1cm}
    \def\cy{\scp*0.3cm}
    \def\ty{2*\dy}
    \tikzset{whp/.style={circle, inner sep=0pt, text width={\pointsize}, draw=black, fill=white}}
    \tikzset{blp/.style={circle, inner sep=0pt, text width={\pointsize}, draw=black, fill=black}}
    \tikzset{lk/.style={regular polygon, regular polygon sides=4, inner sep=0pt, text width={\linksize}, draw=black, fill=black}}
    \draw[dotted] ({0-\dd},{0}) -- ({\txl+\dd},{0});
    \draw[dotted] ({0-\dd},{\ty}) -- ({\txu+\dd},{\ty});
    \coordinate (l1) at ({0+0*\dx},{0+0*\ty}) {};
    \coordinate (l2) at ({0+1*\dx},{0+0*\ty}) {};
    \coordinate (l3) at ({0+2*\dx},{0+0*\ty}) {};
    \coordinate (l4) at ({0+3*\dx},{0+0*\ty}) {};
    \coordinate (l5) at ({0+4*\dx},{0+0*\ty}) {};
    \coordinate (u1) at ({0+0*\dx},{0+1*\ty}) {};
    \node[lk] at ({3*\dx},{1*\dy}) {};
    \draw[->] (l2) -- ++({0*\dx},{1*\dy});
    \draw (l1) -- (u1);
    \draw (l3) -- ++ ({0*\dx},{1*\dy}) -| (l4);
    \draw  ({3*\dx},{1*\dy}) -| (l5);
    \node[whp] at (l1) {};
    \node[whp] at (l2) {};
    \node[whp] at (l3) {};
    \node[blp] at (l4) {};
    \node[whp] at (l5) {};
    \node[blp] at (u1) {};
  \end{tikzpicture}
 $$
          }
      \item
        {
          \newcommand{\inlen}{k}
          \newcommand{\midlen}{\ell}
          \newcommand{\outlen}{m}
            \newcommand{\incol}{\Yc{c}}
            \newcommand{\midcol}{\Yc{d}}
            \newcommand{\outcol}{\Yc{e}}
          \newcommand{\firstpartition}{p}
          \newcommand{\secondpartition}{q}
          \newcommand{\midsup}{s}
          \newcommand{\inblock}{\Yb{A}}
          \newcommand{\midblock}{B}
          \newcommand{\outblock}{\Yb{C}}
          \newcommand{\inind}{i}
          \newcommand{\outind}{j}
          $\thecat$ is closed under composition, i.e., vertical concatenation in the following sense. If~for two set-theoretical partitions   the lower coloring of the first agrees with the upper coloring of the second, then the composition has the same upper coloring as the first and the same lower coloring as the second. Any blocks of the first which only include upper points are inherited by the composition, as are any blocks of the second which only include lower points. The remaining blocks of the composition are formed by the following procedure. The collection of all non-empty intersections of  blocks of the first two-colored partition  with the set of lower points is a set-theoretical partition of the latter. Likewise, a set-theoretical partition of the set of upper points of the second two-colored partition is given by the collection of all non-empty intersections of blocks of the second two-colored partition with it.  If the lower points of the first two-colored partition and the upper points of the second are identified according to the numbering, the two set-theoretical partitions just described have a join. For each element of the join, consider the union  of the following two sets. The first is the (possibly empty) set of upper points of the first two-colored partition which are contained in a block of the first two-colored partition which intersects the element of the join if the latter is interpreted as a set of lower points of the first two-colored partition. Similarly, the second is the (possibly empty) set of lower points of the second two-colored partition which are contained in a block of the second two-colored partition which intersects the element of the join if the latter is interpreted as a set of upper points of the second two-colored partition. Provided that the union of these two sets is not empty, it constitutes a block of the composition of the two-colored partitions. And all blocks of the composition arise in one of the three aforementioned ways.
          In formulas:
          $(\incol,\outcol,\secondpartition\firstpartition)\in \thecat$ for any $(\incol,\midcol,\firstpartition)\in\thecat$ and $(\midcol,\outcol,\secondpartition)\in\thecat$, where if $\{\inlen,\midlen,\outlen\}\subseteq \nnint$ are such that~$\firstpartition$ is a set-theo\-re\-ti\-cal partition of $\tsopx{\inlen}{\midlen}$ and $\secondpartition$ one of $\tsopx{\midlen}{\outlen}$, and if $\midsup$ is the join of the two set-theo\-re\-ti\-cal partitions $\{\{\outind\in\dwi{\midlen}\Sand \lop{\outind}\in \inblock\}\}_{\inblock\in\firstpartition}\backslash\{\varnothing\}$ and $\{\{\inind\in\dwi{\midlen}\Sand \upp{\inind}\in \outblock\}\}_{\outblock\in\secondpartition}\backslash\{\varnothing\}$ of $\dwi{\midlen}$, then $\secondpartition\firstpartition:= \{\inblock\in\firstpartition\Sand \inblock\subseteq \tsopx{\inlen}{0} \}\cupdot  \{\outblock\in\secondpartition\Sand \outblock\subseteq \tsopx{0}{\outlen} \}\cupdot \{  \bigcupdot\{\inblock\cap\tsopx{\inlen}{0}\classpredicate \inblock\in \firstpartition\Sand \exists \outind\in\midblock\quantorpredicate\lop{\outind}\in\inblock \}\cupdot  \bigcupdot\{\outblock\cap \tsopx{0}{\outlen}\classpredicate \outblock\in\secondpartition\Sand \exists \inind\in\midblock\quantorpredicate\upp{\inind}\in\outblock\} \}_{\midblock\in\midsup}\backslash \{\varnothing\} $ is the \emph{composition} of~$(\secondpartition,\firstpartition)$.
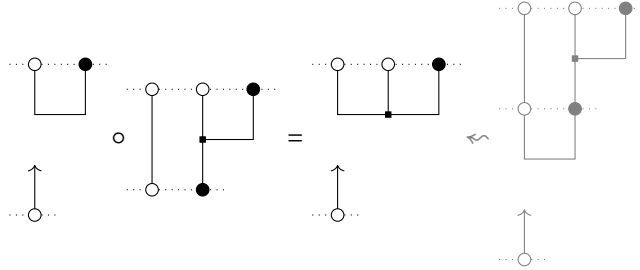
          $$
  \begin{tikzpicture}[baseline=0.91cm]
    \def\scp{0.666}
    \def\linksize{\scp*0.075cm}
    \def\pointsize{\scp*0.25cm}
    \def\dd{\scp*0.5cm}
    \def\dx{\scp*1cm}
    \def\cx{\scp*0.3cm}
    \def\txu{1*\dx}
    \def\txl{0*\dx}
    \def\dy{\scp*1cm}
    \def\cy{\scp*0.3cm}
    \def\ty{3*\dy}
    \tikzset{whp/.style={circle, inner sep=0pt, text width={\pointsize}, draw=black, fill=white}}
    \tikzset{blp/.style={circle, inner sep=0pt, text width={\pointsize}, draw=black, fill=black}}
    \tikzset{lk/.style={regular polygon, regular polygon sides=4, inner sep=0pt, text width={\linksize}, draw=black, fill=black}}
    \draw[dotted] ({0-\dd},{0}) -- ({\txl+\dd},{0});
    \draw[dotted] ({0-\dd},{\ty}) -- ({\txu+\dd},{\ty});
    \coordinate (l1) at ({0+0*\dx},{0+0*\ty}) {};
    \coordinate (u1) at ({0+0*\dx},{0+1*\ty}) {};
    \coordinate (u2) at ({0+1*\dx},{0+1*\ty}) {};
    \draw (u1) -- ++ ({0*\dx},{-1*\dy}) -| (u2);
    \draw[->] (l1) -- ++({0*\dx},{1*\dy});
    \node[whp] at (l1) {};
    \node[whp] at (u1) {};
    \node[blp] at (u2) {};
  \end{tikzpicture}
  \whpoint
  \begin{tikzpicture}[baseline=0.575cm]
    \def\scp{0.666}
    \def\linksize{\scp*0.075cm}
    \def\pointsize{\scp*0.25cm}
    \def\dd{\scp*0.5cm}
    \def\dx{\scp*1cm}
    \def\cx{\scp*0.3cm}
    \def\txu{2*\dx}
    \def\txl{1*\dx}
    \def\dy{\scp*1cm}
    \def\cy{\scp*0.3cm}
    \def\ty{2*\dy}
    \tikzset{whp/.style={circle, inner sep=0pt, text width={\pointsize}, draw=black, fill=white}}
    \tikzset{blp/.style={circle, inner sep=0pt, text width={\pointsize}, draw=black, fill=black}}
    \tikzset{lk/.style={regular polygon, regular polygon sides=4, inner sep=0pt, text width={\linksize}, draw=black, fill=black}}
    \draw[dotted] ({0-\dd},{0}) -- ({\txl+\dd},{0});
    \draw[dotted] ({0-\dd},{\ty}) -- ({\txu+\dd},{\ty});
    \coordinate (l1) at ({0+0*\dx},{0+0*\ty}) {};
    \coordinate (l2) at ({0+1*\dx},{0+0*\ty}) {};
    \coordinate (u1) at ({0+0*\dx},{0+1*\ty}) {};
    \coordinate (u2) at ({0+1*\dx},{0+1*\ty}) {};
    \coordinate (u3) at ({0+2*\dx},{0+1*\ty}) {};
    \node[lk] at ({1*\dx},{1*\dy}) {};
    \draw ({1*\dx},{1*\dy}) -| (u3);
    \draw (l1) -- (u1);
    \draw (l2) -- (u2);
    \node[whp] at (l1) {};
    \node[blp] at (l2) {};
    \node[whp] at (u1) {};
    \node[whp] at (u2) {};
    \node[blp] at (u3) {};
  \end{tikzpicture}
  =\begin{tikzpicture}[baseline=0.91cm]
    \def\scp{0.666}
    \def\linksize{\scp*0.075cm}
    \def\pointsize{\scp*0.25cm}
    \def\dd{\scp*0.5cm}
    \def\dx{\scp*1cm}
    \def\cx{\scp*0.3cm}
    \def\txu{2*\dx}
    \def\txl{0*\dx}
    \def\dy{\scp*1cm}
    \def\cy{\scp*0.3cm}
    \def\ty{3*\dy}
    \tikzset{whp/.style={circle, inner sep=0pt, text width={\pointsize}, draw=black, fill=white}}
    \tikzset{blp/.style={circle, inner sep=0pt, text width={\pointsize}, draw=black, fill=black}}
    \tikzset{lk/.style={regular polygon, regular polygon sides=4, inner sep=0pt, text width={\linksize}, draw=black, fill=black}}
    \draw[dotted] ({0-\dd},{0}) -- ({\txl+\dd},{0});
    \draw[dotted] ({0-\dd},{\ty}) -- ({\txu+\dd},{\ty});
    \coordinate (l1) at ({0+0*\dx},{0+0*\ty}) {};
    \coordinate (u1) at ({0+0*\dx},{0+1*\ty}) {};
    \coordinate (u2) at ({0+1*\dx},{0+1*\ty}) {};
    \coordinate (u3) at ({0+2*\dx},{0+1*\ty}) {};
    \node[lk] at ({1*\dx},{2*\dy}) {};
    \draw (u1) -- ++ ({0*\dx},{-1*\dy}) -| (u2);
    \draw ({1*\dx},{2*\dy}) -| (u3);
    \draw[->] (l1) -- ++({0*\dx},{1*\dy});
    \node[whp] at (l1) {};
    \node[whp] at (u1) {};
    \node[whp] at (u2) {};
    \node[blp] at (u3) {};
  \end{tikzpicture}
    {\color{gray} \rotatebox{180}{\text{$\rightsquigarrow$}} 
  \begin{tikzpicture}[baseline=1.5cm]
    \def\scp{0.666}
    \def\linksize{\scp*0.075cm}
    \def\pointsize{\scp*0.25cm}
    \def\dd{\scp*0.5cm}
    \def\dx{\scp*1cm}
    \def\cx{\scp*0.3cm}
    \def\txu{2*\dx}
    \def\txm{1*\dx}
    \def\txl{0*\dx}
    \def\dy{\scp*1cm}
    \def\cy{\scp*0.3cm}
    \def\tyl{3*\dy}
    \def\tyu{2*\dy}
    \tikzset{whp/.style={circle, inner sep=0pt, text width={\pointsize}, draw=gray, fill=white}}
    \tikzset{blp/.style={circle, inner sep=0pt, text width={\pointsize}, draw=gray, fill=gray}}
    \tikzset{lk/.style={regular polygon, regular polygon sides=4, inner sep=0pt, text width={\linksize}, draw=gray, fill=gray}}
    \draw[dotted] ({0-\dd},{0}) -- ({\txl+\dd},{0});
    \draw[dotted] ({0-\dd},{\tyl}) -- ({\txm+\dd},{\tyl});
    \draw[dotted] ({0-\dd},{\tyu+\tyl}) -- ({\txu+\dd},{\tyu+\tyl});
    \coordinate (l1) at ({0+0*\dx},{0+0*\tyl+0*\tyu}) {};
    \coordinate (m1) at ({0+0*\dx},{0+1*\tyl+0*\tyu}) {};
    \coordinate (m2) at ({0+1*\dx},{0+1*\tyl+0*\tyu}) {};
    \coordinate (u1) at ({0+0*\dx},{0+1*\tyl+1*\tyu}) {};
    \coordinate (u2) at ({0+1*\dx},{0+1*\tyl+1*\tyu}) {};
    \coordinate (u3) at ({0+2*\dx},{0+1*\tyl+1*\tyu}) {};
    \node[lk] at ({1*\dx},{1*\dy+1*\tyl+0*\tyu}) {};
    \draw ({1*\dx},{1*\dy+1*\tyl+0*\tyu}) -| (u3);
    \draw (m1) -- ({0*\dx},{-1*\dy+1*\tyl+0*\tyu}) -| (m2);
    \draw[->] (l1) -- ++({0*\dx},{1*\dy});
    \draw (m1) -- (u1);
    \draw (m2) -- (u2);
    \node[whp] at (l1) {};
    \node[whp] at (m1) {};
    \node[blp] at (m2) {};
    \node[whp] at (u1) {};
    \node[whp] at (u2) {};
    \node[blp] at (u3) {};
  \end{tikzpicture}          }
      $$
          }
  \end{enumerate}
  }
\item
  {
    \newcommand{\thegenset}{\mathcal{G}}
    For any set $\thegenset$ of two-co\-lo\-red partition, we write $\langle \thegenset\rangle$ for the intersection of all categories of two-co\-lo\-red partitions containing $\thegenset$ and we say that $\thegenset$ \emph{generates} $\langle \thegenset\rangle$.
    }
  \end{enumerate}
\end{Definition}

\begin{Example}\quad
            \newcommand{\inlen}{k}
          \newcommand{\outlen}{\ell}
            \newcommand{\incol}{\Yc{c}}
            \newcommand{\outcol}{\Yc{d}}
          \newcommand{\thepartition}{p}
          \newcommand{\block}{\Yb{B}}
          \newcommand{\inind}{a}
          \newcommand{\outind}{b}
          \begin{enumerate}
            \item Of course, the set of all two-co\-lo\-red partitions forms the maximal category of two-co\-lo\-red partitions. It follows from \cite[Theorem~8.3]{TarragoWeber2018} that it coincides with $\langle \PartCrossWW,\PartFourWBWB,\PartIdenLoWW,\PartSingleW\rangle$.
  \item Another category of two-co\-lo\-red partitions is given by  $\langle \PartCrossWW\rangle$, the category of two-co\-lo\-red pair partitions with neutral blocks,  i.e., all  $(\incol,\outcol,\thepartition)$ with $|\block|=2$ and $\csx{\incol}{\outcol}{\block}=0$    for any $\block\in \thepartition$, (see \cite[Proposition~5.3]{MangWeber2020}).
  \item By \cite[Theorem~7.2]{TarragoWeber2018}, the minimal category of two-co\-lo\-red partitions $\langle \varnothing\rangle$ is the subset of all elements of  $\langle \PartCrossWW\rangle$ which are \emph{non-crossing}. The precise definition of being non-crossing is unimportant here; informally, it means that blocks can be  ``drawn without intersections''.
  \end{enumerate}
\end{Example}
Not much familiarity with two-co\-lo\-red partitions and their categories is required in order to prove the main result. In particular, the full classification of all categories of two-co\-lo\-red partitions can remain open. However, we will need to divide the landscape of all possible categories as follows.
\begin{Definition}
  \label{definition:obsh_cases}
  \newcommand{\incol}{\Yc{c}}
  \newcommand{\outcol}{\Yc{d}}
  \newcommand{\thepartition}{p}
  \newcommand{\anyblock}{\Yb{B}}
  We say that any category  $\thecat$  of two-co\-lo\-red partitions is
  \begin{enumerate}
  \item \emph{case $\mathcal{O}$} if $\PartSinglesWB\notin\thecat$ and $\PartFourWBWB\notin\thecat$,
  \item \emph{case $\mathcal{B}$}  if $\PartSinglesWB\in\thecat$ and $\PartFourWBWB\notin\thecat$,
  \item \emph{case $\mathcal{H}$}  if $\PartSinglesWB\notin\thecat$ and $\PartFourWBWB\in\thecat$,
  \item \emph{case $\mathcal{S}$}     if $\PartSinglesWB\in\thecat$ and $\PartFourWBWB \in\thecat$,
  \item \emph{class $\onlyneutralnonsingletonblocks$}
    if $\csx{\incol}{\outcol}{\anyblock}=0$ for any $(\incol,\outcol,\thepartition)\in\thecat$ and any $\anyblock\in\thepartition$ with $2\leq |\anyblock|$,
  \item \emph{class $\onlyneutralpartitions$}
    if $\tcsx{\incol}{\outcol}=0$ for any $(\incol,\outcol,\thepartition)\in\thecat$.
  \end{enumerate}
\end{Definition}
The names $\onlyneutralnonsingletonblocks$ and $\onlyneutralpartitions$ reflect the defining conditions  of having only \underline{n}eutral \underline{n}on-\underline{s}ingleton \underline{b}locks respectively only \underline{n}eutral two-colored \underline{p}artitions, where ``neutral'' means vanishing color sum.
For the motivation behind the names $\mathcal{O}$, $\mathcal{B}$, $\mathcal{H}$ and $\mathcal{S}$ see Remark~\ref{remark:obsh_cases} below.
\begin{Remark}
  \label{remark:two-colored-categories}
  \newcommand{\anycat}{\mathcal{R}}
  For any of the known categories of two-colored partitions, it is easy to determine whether it has a given property in Definition~\ref{definition:obsh_cases} or not. Any known category which is not case~$\mathcal{H}$ is covered by \cite{MangWeber2021b} (see Section~7 there for the correspondence to the results of \cite{Gromada2018, MangWeber2020, MangWeber2021a, TarragoWeber2018}) and any known case-$\mathcal{H}$ category by \cite{Gromada2018,Maassen2021,TarragoWeber2018} or \cite[Chapter~1]{Mang2022}.

\emph{Cases~$\mathcal{O}$, $\mathcal{B}$, $\mathcal{S}$.} Any category $\mathcal{R}_{f,v,s,l,k,x}$ in the main theorem of \cite{MangWeber2021b} is case $\mathcal{O}$ if and only if $f=\{2\}$, case $\mathcal{B}$ if and only if $f=\{1,2\}$ and case $\mathcal{S}$ if and only if $f=\pint$ (and never case~$\mathcal{H}$). It is class $\onlyneutralnonsingletonblocks$ if and only if $v=\{0\}$ or $v=\pm\{0,1\}$ and class $\onlyneutralpartitions$ if and only if $s=\{0\}$.

\emph{Case~$\mathcal{H}$.} In \cite{TarragoWeber2018}, neither of the categories $\mathcal{H}_{\mathrm{glob}}(k)$ of Theorem~7.1 and $\mathcal{H}_{\mathrm{grp},\mathrm{glob}}(k)$ of Theorem~8.3 is class $\onlyneutralnonsingletonblocks$. And, each is class $\onlyneutralpartitions$ if and only if $k=0$. The category  $\mathcal{H}'{}_{\mathrm{loc}}$ in  Theorem~7.2 is both class $\onlyneutralnonsingletonblocks$ and class $\onlyneutralpartitions$. Each of $\mathcal{H}_{\mathrm{loc}}(k,d)$ from Theorem~7.2 and $\mathcal{H}_{\mathrm{grp},\mathrm{loc}}(k,d)$ from Theorem~8.3 is class $\onlyneutralnonsingletonblocks$  if and only if $k=d=0$ and is class $\onlyneutralpartitions$ if and only if $k=0$.

The  case-$\mathcal{H}$ categories in \cite[Table~1]{Gromada2018} which are not already covered by  \cite{TarragoWeber2018} are
$\mathcal{H}_{\mathrm{hl},\mathrm{glob}}(k,0)$,  $\mathcal{H}_{\mathrm{hl},\mathrm{glob}}(k,s)$,  $\mathcal{H}_{\pi}(k,s)$,  $\mathcal{H}_{\pi}(k,\infty)$ and $\mathcal{H}_{A}(k)$ (where the categories $\mathcal{H}_{\mathrm{hl},\mathrm{glob}}(k,0)$ and  $\mathcal{H}_{\mathrm{hl},\mathrm{glob}}(k,s)$ can each also be written as $\mathcal{H}_{A}(k)$ for certain $A$). None one of these are class $\onlyneutralnonsingletonblocks$. And any one is class $\onlyneutralpartitions$ if and only if $k= 0$.

In \cite{Maassen2021}, no group-theoretical category $\anycat$ of two-colored partitions in the sense of Definition~4.1.5 is class $\onlyneutralnonsingletonblocks$. And, any such category is class $\onlyneutralpartitions$ if and only if $F_\infty(\anycat)$ in the sense of Definition~4.3.21 contains no word with different numbers of generators and inverses of generators.

Lastly, the category $\mathcal{W}_{\mathcal{R}}$ in the sense of the main result of \cite[Chapter~1]{Mang2022} is both class $\onlyneutralnonsingletonblocks$ and class $\onlyneutralpartitions$ for any parameter $\mathcal{R}$.
\end{Remark}
Beyond those case distinctions, we will also need to know the following elementary facts about categories of two-co\-lo\-red partitions.
\begin{Definition}
  \newcommand{\inlen}{k}
  \newcommand{\outlen}{\ell}
  \newcommand{\thepartition}{p}
  \newcommand{\incol}{\Yc{c}}
  \newcommand{\outcol}{\Yc{d}}
  \newcommand{\inind}{i}
  \newcommand{\outind}{j}
  \newcommand{\incolx}[1]{\incol_{#1}}
  \newcommand{\outcolx}[1]{\outcol_{#1}}
  \newcommand{\block}{\Yb{B}}
Dual two-colored partitions are obtained by simultaneous horizontal reflection, vertical reflection and color inversion. More precisely, given any $\{\inlen,\outlen\}\subseteq \nnint$, any $\incol\in\blaw^{\Ssetmonoidalproduct\inlen}$, any $\outcol\in\blaw^{\Ssetmonoidalproduct\outlen}$ and any set-theo\-re\-ti\-cal partition $\thepartition$ of $\tsopx{\inlen}{\outlen}$ the \emph{dual} of $(\incol,\outcol,\thepartition)$ is the triple $(\ASdualX{\outcol},\ASdualX{\incol},\ASdualX{\thepartition})$, where $\ASdualX{\outcol}\in \blaw^{\Ssetmonoidalproduct\outlen}$ is defined by $\outind\mapsto \overline{\outcolx{\outlen-\outind+1}}$, where $\ASdualX{\incol}\in \blaw^{\Ssetmonoidalproduct\inlen}$ is defined by $\inind\mapsto \overline{\incolx{\inlen-\inind+1}}$, and where $\ASdualX{\thepartition}:=\{\{\upp{(\outlen-\outind+1)}\classpredicate \outind\in\dwi{\outlen}\Sand \lop{\outind}\in\block\}\cupdot\{\lop{(\inlen-\inind+1)}\classpredicate \inind\in\dwi{\inlen}\Sand \upp{\inind}\in\block\}\}_{\block\in\thepartition}$ is the \emph{dual} of~$\thepartition$.
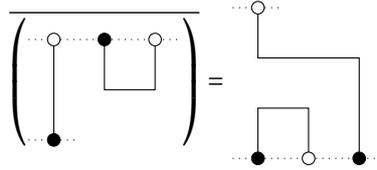
  $$
  \overline{\left(
              \begin{tikzpicture}[baseline=0.666cm]
    \def\scp{0.666}
    \def\linksize{\scp*0.075cm}
    \def\pointsize{\scp*0.25cm}
    \def\dd{\scp*0.5cm}
    \def\dx{\scp*1cm}
    \def\cx{\scp*0.3cm}
    \def\txu{2*\dx}
    \def\txl{0*\dx}
    \def\dy{\scp*1cm}
    \def\cy{\scp*0.3cm}
    \def\ty{2*\dy}
    \tikzset{whp/.style={circle, inner sep=0pt, text width={\pointsize}, draw=black, fill=white}}
    \tikzset{blp/.style={circle, inner sep=0pt, text width={\pointsize}, draw=black, fill=black}}
    \tikzset{lk/.style={regular polygon, regular polygon sides=4, inner sep=0pt, text width={\linksize}, draw=black, fill=black}}
    \draw[dotted] ({0-\dd},{0}) -- ({\txl+\dd},{0});
    \draw[dotted] ({0-\dd},{\ty}) -- ({\txu+\dd},{\ty});
    \coordinate (l1) at ({0+0*\dx},{0+0*\ty}) {};
    \coordinate (u1) at ({0+0*\dx},{0+1*\ty}) {};
    \coordinate (u2) at ({0+1*\dx},{0+1*\ty}) {};
    \coordinate (u3) at ({0+2*\dx},{0+1*\ty}) {};
    \draw (u2) -- ++ ({0*\dx},{-1*\dy}) -| (u3);
    \draw (l1) -- (u1);
    \node[blp] at (l1) {};
    \node[whp] at (u1) {};
    \node[blp] at (u2) {};
    \node[whp] at (u3) {};
  \end{tikzpicture}
\right)}=
                \begin{tikzpicture}[baseline=0.91cm]
    \def\scp{0.666}
    \def\linksize{\scp*0.075cm}
    \def\pointsize{\scp*0.25cm}
    \def\dd{\scp*0.5cm}
    \def\dx{\scp*1cm}
    \def\cx{\scp*0.3cm}
    \def\txu{0*\dx}
    \def\txl{2*\dx}
    \def\dy{\scp*1cm}
    \def\cy{\scp*0.3cm}
    \def\ty{3*\dy}
    \tikzset{whp/.style={circle, inner sep=0pt, text width={\pointsize}, draw=black, fill=white}}
    \tikzset{blp/.style={circle, inner sep=0pt, text width={\pointsize}, draw=black, fill=black}}
    \tikzset{lk/.style={regular polygon, regular polygon sides=4, inner sep=0pt, text width={\linksize}, draw=black, fill=black}}
    \draw[dotted] ({0-\dd},{0}) -- ({\txl+\dd},{0});
    \draw[dotted] ({0-\dd},{\ty}) -- ({\txu+\dd},{\ty});
    \coordinate (l1) at ({0+0*\dx},{0+0*\ty}) {};
    \coordinate (l2) at ({0+1*\dx},{0+0*\ty}) {};
    \coordinate (l3) at ({0+2*\dx},{0+0*\ty}) {};
    \coordinate (u1) at ({0+0*\dx},{0+1*\ty}) {};
    \draw (l1) -- ++ ({0*\dx},{1*\dy}) -| (l2);
    \draw (l3) -- ++ ({0*\dx},{2*\dy}) -| (u1);
    \node[blp] at (l1) {};
    \node[whp] at (l2) {};
    \node[blp] at (l3) {};
    \node[whp] at (u1) {};
  \end{tikzpicture}
      $$
\end{Definition}

{
  \newcommand{\incol}{\Yc{c}}
  \newcommand{\outcol}{\Yc{d}}
  \newcommand{\thepartition}{p}
  \newcommand{\param}{t}
  \newcommand{\someblock}{\Yb{B}}
  \begin{Lemma}
    \label{lemma:scalar-article-categories-helper}
    Let $\thecat$ be any category of two-co\-lo\-red partitions.
    \begin{enumerate}
    \item\label{lemma:scalar-article-categories-helper-0}  $(\ASdualX{\outcol},\ASdualX{\incol},\ASdualX{\thepartition})\in\thecat$ for any $(\incol,\outcol,\thepartition)\in\thecat$.
    \item\label{lemma:scalar-article-categories-helper-1} $\PartSinglesWB\in\thecat$ if and only if there exist $(\incol,\outcol,\thepartition)\in\thecat$ and $\someblock\in\thepartition$ such that $|\someblock|<2$.
    \item\label{lemma:scalar-article-categories-helper-2} $\PartFourWBWB\in\thecat$ if and only if there exist $(\incol,\outcol,\thepartition)\in\thecat$ and $\someblock\in\thepartition$ such that $|\someblock|>2$.
    \item\label{lemma:scalar-article-categories-helper-3} If $\PartSinglesWB\in\thecat$ and $\PartFourWBWB\in\thecat$, then $\PartCoSingleWBWB\in\thecat$.
    \item\label{lemma:scalar-article-categories-helper-4} If $\PartSinglesWB\in\thecat$, then $\{\PartSingleW^{\Smonoidalproduct |\tcsx{\incol}{\outcol}|},\PartSingleB^{\Smonoidalproduct |\tcsx{\incol}{\outcol}|}\}\subseteq\thecat$ for any $(\incol,\outcol,\thepartition)\in \thecat$ with $\tcsx{\incol}{\outcol}\neq 0$.
    \end{enumerate}
  \end{Lemma}
  \begin{proof}
    Part \ref{lemma:scalar-article-categories-helper-0} is implied by \cite[Lemma~1.1\,(a)]{TarragoWeber2018}.
Parts \ref{lemma:scalar-article-categories-helper-1} and \ref{lemma:scalar-article-categories-helper-2} are
 \cite[Lemmas~1.3\,(b) and~2.1\,(a)]{TarragoWeber2018}
and \cite[Lemmas~1.3\,(d) and 2.1\,(b)]{TarragoWeber2018}, respectively.
Part \ref{lemma:scalar-article-categories-helper-3} follows immediately from \cite[Lemma~1.3\,(b)]{TarragoWeber2018}.

In order to see part \ref{lemma:scalar-article-categories-helper-4}, use \cite[Lemma~1.1\,(a)]{TarragoWeber2018} to first ``rotate'' any potential upper point of $(\incol,\outcol,\thepartition)$ down (in an arbitrary direction). ``Disconnect'' then each and every point from its block with the help of  \cite[Lemma~1.3\,(b)]{TarragoWeber2018}. Following that, keep ``erasing'' neighboring points of different colors,  as  \cite[Lemma~1.1\,(b)]{TarragoWeber2018} permits, until no such points remain. None of these transformations have affected the total color sum. The resulting two-colored partition is either~$\PartSingleW^{\Smonoidalproduct |\tcsx{\incol}{\outcol}|}$ or $\PartSingleB^{\Smonoidalproduct |\tcsx{\incol}{\outcol}|}$. Passing to the adjoint of the dual  as allowed by \ref{lemma:scalar-article-categories-helper-0} hence shows the claim.
  \end{proof}
}

{
  \begin{Lemma}\quad
    \label{lemma:impossible-cases}
    \begin{enumerate}
    \item\label{lemma:impossible-cases-1} Any case-$\mathcal{O}$ or case-$\mathcal{H}$ category of two-colored partitions that is class $\onlyneutralnonsingletonblocks$ is class $\onlyneutralpartitions$.
    \item\label{lemma:impossible-cases-2} No case-$\mathcal{S}$ category of two-colored partitions  is class $\onlyneutralnonsingletonblocks$.
    \end{enumerate}
  \end{Lemma}
  \begin{proof}
      \newcommand{\anyblock}{\Yb{B}}
      \newcommand{\incol}{\Yc{c}}
      \newcommand{\outcol}{\Yc{d}}
      \newcommand{\thepartition}{p}
      (a) Let $\thecat$ be     case-$\mathcal{O}$ or case-$\mathcal{H}$ and class $\onlyneutralnonsingletonblocks$. Then, for any $(\incol,\outcol,\thepartition)\in \thecat$ and any $\anyblock\in \thepartition$ on the one hand  $2\leq |\anyblock|$ by the first assumption and thus on the other hand $\csx{\incol}{\outcol}{\anyblock}=0$ by the second assumption. Since that demands  $\tcsx{\incol}{\outcol}=\sum_{\anyblock\in\thepartition}\csx{\incol}{\outcol}{\anyblock}=0$ the category $\thecat$ is necessarily class~$\onlyneutralpartitions$.

(b) Since any case-$\mathcal{S}$ category contains both $\PartSinglesWB$ and $\PartFourWBWB$, it must also contain $\PartCoSingleWBWB$ by Lem\-ma~\hyperref[lemma:scalar-article-categories-helper-3]{\ref*{lemma:scalar-article-categories-helper}\,\ref*{lemma:scalar-article-categories-helper-3}}. The fact that $\{\lop{1},\lop{3}\}\in \UCPartCoSingle$ and $\csx{\varnothing}{\whpoint\blpoint\whpoint\blpoint}{\{\lop{1},\lop{3}\}}=2\neq 0$ hence shows that such a~category is not class $\onlyneutralnonsingletonblocks$.
  \end{proof}
}

\subsection{Unitary easy quantum groups}
\label{section:easy-relations}
``Easy'' quantum groups are now defined by transforming the elements of a given category of two-co\-lo\-red partitions into relations for the generators of a universal algebra that can be given the structure of a compact quantum group. To be more precise, an entire  series of compact quantum groups indexed by $\pint$ arises in this way.

\begin{Assumptions}
  \label{assumptions:easy}
  \newcommand{\indj}{j}
  \newcommand{\indi}{i}
In the following, fix any  $\thedim\!\in\!\pint$ and any $2\thedim^2$-elemental set
    $\thegens \!=\! \{\theunix{\whpoint}{\indj}{\indi},\theunix{\blpoint}{\indj}{\indi}\}_{\indi,\indj=1}^\thedim$ and define the two families $\theunim{\whpoint}:= (\theunix{\whpoint}{\indj}{\indi})_{(\indj,\indi)\in\dwi{\thedim}^{\Ssetmonoidalproduct 2}}$ and  $\theunim{\blpoint}:= (\theunix{\blpoint}{\indj}{\indi})_{(\indj,\indi)\in\dwi{\thedim}^{\Ssetmonoidalproduct 2}}$.
\end{Assumptions}
The transformation of two-co\-lo\-red partitions into relations is accomplished by the following formula, where  $\zetf$ was defined in Notation~\hyperref[notation:set-theoretical-partitions-3]{\ref*{notation:set-theoretical-partitions}\,\ref*{notation:set-theoretical-partitions-3}}.
\begin{Notation}
  \label{notation:easy}
  \newcommand{\indj}{j}
  \newcommand{\indi}{i}
  \newcommand{\incol}{\Yc{c}}
  \newcommand{\outcol}{\Yc{d}}
  \newcommand{\incolx}[1]{\Yc{c}_{#1}}
  \newcommand{\outcolx}[1]{\Yc{d}_{#1}}
  \newcommand{\inlen}{k}
  \newcommand{\outlen}{\ell}
  \newcommand{\infree}{g}
  \newcommand{\outfree}{j}
  \newcommand{\inbound}{h}
  \newcommand{\outbound}{i}
  \newcommand{\infreex}[1]{g_{#1}}
  \newcommand{\outfreex}[1]{j_{#1}}
  \newcommand{\inboundx}[1]{h_{#1}}
  \newcommand{\outboundx}[1]{i_{#1}}
  \newcommand{\incounter}{a}
  \newcommand{\outcounter}{b}
  \newcommand{\thepartition}{p}
    For any   $ \{\inlen,\outlen\}\subseteq \nnint$, any $\incol\in\blaw^{\Ssetmonoidalproduct\inlen}$ and $ \outcol\in\blaw^{\Ssetmonoidalproduct\outlen}$, any set-theo\-re\-ti\-cal partition $\thepartition$ of $\tsopx{\inlen}{\outlen}$ and any $\infree\in\dwi{\thedim}^{\Ssetmonoidalproduct\inlen}$ and $\outfree\in\dwi{\thedim}^{\Ssetmonoidalproduct\outlen}$, let in $\comps\langle \thegens\rangle$
    \[
      \therelpoly{\incol}{\outcol}{\thepartition}{\outfree}{\infree}:=      \sum_{\outbound\in\dwi{\thedim}^{\Ssetmonoidalproduct \outlen}}
      \zetfx{\thepartition}{\ker(\djp{\infree}{\outbound})}\;
\rightwardproduct{\outcounter=1}{\outlen}{\theunix{\outcolx{\outcounter}}{\outfreex{\outcounter}}{\outboundx{\outcounter}}}
      -\sum_{\inbound\in\dwi{\thedim}^{\Ssetmonoidalproduct \inlen}}  \zetfx{\thepartition}{\ker(\djp{\inbound}{\outfree})}\; \rightwardproduct{\incounter=1}{\inlen}{\theunix{\incolx{\incounter}}{\inboundx{\incounter}}{\infreex{\incounter}}}.
    \]
  \end{Notation}
For example, the two-colored partitions  $\PartIdenW$ and $\PartIdenB$ induce the trivial relation $0$. The relations induced by $\PartIdenLoWB$, $\PartIdenUpBW$, $\PartIdenLoBW$ and $\PartIdenUpWB$ will be of the utmost importance.

{
  \newcommand{\infree}{g}
  \newcommand{\infreex}[1]{\infree_{#1}}
  \newcommand{\inbound}{h}
  \newcommand{\inboundx}[1]{\inbound_{#1}}
  \newcommand{\outbound}{i}
  \newcommand{\outboundx}[1]{\outbound_{#1}}
  \newcommand{\outfree}{j}
  \newcommand{\outfreex}[1]{\outfree_{#1}}
  \newcommand{\inlen}{k}
  \newcommand{\outlen}{\ell}
  \newcommand{\incol}{\Yc{c}}
  \newcommand{\outcol}{\Yc{d}}
  \newcommand{\incolx}[1]{\Yc{c}_{#1}}
  \newcommand{\outcolx}[1]{\Yc{d}_{#1}}
  \newcommand{\thepartition}{p}
  \newcommand{\therel}{r}
  \newcommand{\incounter}{a}
  \newcommand{\outcounter}{b}
  \begin{Lemma}
    \label{lemma:unitarity-relations}
    For any $\infree\in\dwi{\thedim}^{\Ssetmonoidalproduct 2}$ and $\outfree\in\dwi{\thedim}^{\Ssetmonoidalproduct 2}$, the following hold:
    \begin{alignat*}{3}
&\therelpoly{\varnothing}{\whpoint\blpoint}{\UCPartIdenLo}{\outfree}{\varnothing}= \sum_{\outbound=1}^\thedim \theunix{\whpoint}{\outfreex{1}}{\outbound}\theunix{\blpoint}{\outfreex{2}}{\outbound}-\kron{\outfreex{1}}{\outfreex{2}}1,
     \qquad && \therelpoly{\blpoint\whpoint}{\varnothing}{\UCPartIdenUp}{\varnothing}{\infree}= \kron{\infreex{1}}{\infreex{2}}1-\sum_{\inbound=1}^\thedim \theunix{\blpoint}{\inbound}{\infreex{1}}\theunix{\whpoint}{\inbound}{\infreex{2}},&
      \\    &  \therelpoly{\varnothing}{\blpoint\whpoint}{\UCPartIdenLo}{\outfree}{\varnothing}= \sum_{\outbound=1}^\thedim       \theunix{\blpoint}{\outfreex{1}}{\outbound}\theunix{\whpoint}{\outfreex{2}}{\outbound}-\kron{\outfreex{1}}{\outfreex{2}}1,
\qquad &&       \therelpoly{\whpoint\blpoint}{\varnothing}{\UCPartIdenUp}{\varnothing}{\infree}= \kron{\infreex{1}}{\infreex{2}}1-\sum_{\inbound=1}^\thedim \theunix{\whpoint}{\inbound}{\infreex{1}}\theunix{\blpoint}{\inbound}{\infreex{2}}.&
    \end{alignat*}
  \end{Lemma}
  \begin{proof}
    Only the proof for $\therelpoly{\varnothing}{\whpoint\blpoint}{\UCPartIdenLo}{\outfree}{\varnothing}$ is given. With the names of Notation~\ref{notation:easy} then $\inlen=0$ and $\outlen=2$ and $\incol=\varnothing$ and $\outcol=\circ\bullet$ and $\thepartition=\{\{\lop{1},\lop{2}\}\}$ and $\infree=\varnothing$. On the one hand, for any $\outbound\in \dwi{\thedim}^{\Ssetmonoidalproduct 2}$ the set-theoretical partition $\kerp{\djp{\infree}{\outbound}}$ can only take two values, namely    $\{\{\lop{1},\lop{2}\}\}$ if $\outboundx{1}= \outboundx{2}$ and $\{\{\lop{1}\},\{\lop{2}\}\}$ if $\outboundx{1}\neq \outboundx{2}$.  Whereas $\kerp{\djp{\infree}{\outbound}}$ even agrees with $\thepartition$ in the former case, $\thepartition$ is not finer than $\kerp{\djp{\infree}{\outbound}}$ in the latter case. Hence, only if $\outboundx{1}= \outboundx{2}$ does  $\zetfx{\thepartition}{\ker(\djp{\infree}{\outbound})}$ evaluate to $1$. Consequently, the first of the two sums in the definition of $      \therelpoly{\incol}{\outcol}{\thepartition}{\outfree}{\infree}$ effectively runs only over the pairs $(\outbound,\outbound)$ for $\outbound\in\dwi{\thedim}$. That explains the term $ \sum_{\outbound=1}^\thedim \theunix{\whpoint}{\outfreex{1}}{\outbound}\theunix{\blpoint}{\outfreex{2}}{\outbound}$ in the claim. On the other hand, $\inlen=0$ by convention implies $\dwi{\thedim}^{\Ssetmonoidalproduct \inlen}=\{\varnothing\}$. Thus, $\varnothing$ is the only $\inbound$ over which the second sum in the definition of $\therelpoly{\incol}{\outcol}{\thepartition}{\outfree}{\infree}$ runs. Just like before, $\kerp{\djp{\inbound}{\outfree}}$ is then either $\{\{\lop{1},\lop{2}\}\}$ or $\{\{\lop{1}\},\{\lop{2}\}\}$, depending on whether  if $\outfreex{1}= \outfreex{2}$ or not. In other words, $\zetfx{\thepartition}{\ker(\djp{\inbound}{\outfree})}=\kron{\outfreex{1}}{\outfreex{2}}$. By $\inlen=0$ the set of indices $\incounter$ over which the product \smash{$\tsrightwardproduct{\incounter=1}{\inlen}{\theunix{\incolx{\incounter}}{\inboundx{\incounter}}{\infreex{\incounter}}}$} runs is the empty set $\varnothing$ (whereas the set of indices $\inbound$ before was not $\varnothing$ but $\{\varnothing\}$).  By common convention, a product with empty index set is~$1$. That is how the second term $-\kron{\outfreex{1}}{\outfreex{2}}1$ in the claim comes about.
  \end{proof}
  }

In general, the relations can become quite complicated.
\begin{Example}
  \label{example:relations}
  \newcommand{\indj}{j}
  \newcommand{\indi}{i}
  \newcommand{\incol}{\Yc{c}}
  \newcommand{\outcol}{\Yc{d}}
  \newcommand{\incolx}[1]{\Yc{c}_{#1}}
  \newcommand{\outcolx}[1]{\Yc{d}_{#1}}
  \newcommand{\inlen}{k}
  \newcommand{\outlen}{\ell}
  \newcommand{\infree}{g}
  \newcommand{\outfree}{j}
  \newcommand{\inbound}{h}
  \newcommand{\outbound}{i}
  \newcommand{\infreex}[1]{g_{#1}}
  \newcommand{\outfreex}[1]{j_{#1}}
  \newcommand{\inboundx}[1]{h_{#1}}
  \newcommand{\outboundx}[1]{i_{#1}}
  \newcommand{\incounter}{a}
  \newcommand{\outcounter}{b}
  \newcommand{\thepartition}{p}
    If $\inlen=4$ and $\outlen=5$ and $\incol=\mathnormal{\whpoint}\mathnormal{\blpoint}\mathnormal{\blpoint}\mathnormal{\whpoint}$ and $\outcol=\blpoint\whpoint\blpoint\blpoint\blpoint$ and $\thepartition=\{\{\lop{1}\},\{\upp{2},\lop{2}\},\allowbreak\{\lop{4},\lop{5}\},\allowbreak\{\upp{1},\upp{3},\upp{4},\lop{3}\}\}$,
$$
      (\incol,\outcol,\thepartition)
    \begin{tikzpicture}[baseline=0.91cm]
    \def\scp{0.666}
    \def\linksize{\scp*0.075cm}
    \def\pointsize{\scp*0.25cm}
    \def\dd{\scp*0.5cm}
    \def\dx{\scp*1cm}
    \def\cx{\scp*0.3cm}
    \def\txu{3*\dx}
    \def\txl{4*\dx}
    \def\dy{\scp*1cm}
    \def\cy{\scp*0.3cm}
    \def\ty{3*\dy}
    \tikzset{whp/.style={circle, inner sep=0pt, text width={\pointsize}, draw=black, fill=white}}
    \tikzset{blp/.style={circle, inner sep=0pt, text width={\pointsize}, draw=black, fill=black}}
    \tikzset{lk/.style={regular polygon, regular polygon sides=4, inner sep=0pt, text width={\linksize}, draw=black, fill=black}}
    \draw[dotted] ({0-\dd},{0}) -- ({\txl+\dd},{0});
    \draw[dotted] ({0-\dd},{\ty}) -- ({\txu+\dd},{\ty});
    \coordinate (l1) at ({0+0*\dx},{0+0*\ty}) {};
    \coordinate (l2) at ({0+1*\dx},{0+0*\ty}) {};
    \coordinate (l3) at ({0+2*\dx},{0+0*\ty}) {};
    \coordinate (l4) at ({0+3*\dx},{0+0*\ty}) {};
    \coordinate (l5) at ({0+4*\dx},{0+0*\ty}) {};
    \coordinate (u1) at ({0+0*\dx},{0+1*\ty}) {};
    \coordinate (u2) at ({0+1*\dx},{0+1*\ty}) {};
    \coordinate (u3) at ({0+2*\dx},{0+1*\ty}) {};
    \coordinate (u4) at ({0+3*\dx},{0+1*\ty}) {};
    \node[lk] at ({2*\dx},{2*\dy}) {};
    \draw (u1) -- ++ ({0*\dx},{-1*\dy}) -| (u4);
    \draw[->] (l1) -- ++({0*\dx},{1*\dy});
    \draw (l2) -- (u2);
    \draw (l3) -- (u3);
    \draw (l4) -- ++ ({0*\dx},{1*\dy}) -| (l5);
    \node[blp] at (l1) {};
    \node[whp] at (l2) {};
    \node[blp] at (l3) {};
    \node[blp] at (l4) {};
    \node[blp] at (l5) {};
    \node[whp] at (u1) {};
    \node[blp] at (u2) {};
    \node[blp] at (u3) {};
    \node[whp] at (u4) {};
  \end{tikzpicture}
$$
then for any $\infree\in\dwi{\thedim}^{\Ssetmonoidalproduct\inlen}$ and $\outfree\in\dwi{\thedim}^{\Ssetmonoidalproduct\outlen}$,
\[
      \therelpoly{\incol}{\outcol}{\thepartition}{\outfree}{\infree}  =\kron{\infreex{1}}{\infreex{3}}\kron{\infreex{1}}{\infreex{4}}\; \left(\sum_{\indi=1}^\thedim\theunix{\blpoint}{\outfreex{1}}{\indi}\right)\theunix{\whpoint}{\outfreex{2}}{\infreex{2}}\theunix{\blpoint}{\outfreex{3}}{\infreex{1}}\left(\sum_{\indi=1}^\thedim\theunix{\blpoint}{\outfreex{4}}{\indi}\theunix{\blpoint}{\outfreex{5}}{\indi}\right)-\kron{\outfreex{4}}{\outfreex{5}}\; \theunix{\whpoint}{\outfreex{3}}{\infreex{1}}\theunix{\blpoint}{\outfreex{2}}{\infreex{2}}\theunix{\blpoint}{\outfreex{3}}{\infreex{3}}\theunix{\whpoint}{\outfreex{3}}{\infreex{4}}.
\]
  \end{Example}

{
  \newcommand{\thealg}{A}
  \newcommand{\themultiplication}{\Smultiplication}
  \newcommand{\thecomultiplication}{\Scomultiplication}
  \newcommand{\theunit}{\Sunit}
  \newcommand{\thestar}{\Sstar}
  \newcommand{\thecatgen}{\mathcal{G}}
  \newcommand{\thecatgenextended}{\mathcal{R}}
  \newcommand{\indj}{j}
  \newcommand{\indi}{i}
  \newcommand{\indg}{g}
  \newcommand{\indh}{h}
  \newcommand{\incol}{\Yc{c}}
  \newcommand{\outcol}{\Yc{d}}
  \newcommand{\incolx}[1]{\Yc{c}_{#1}}
  \newcommand{\outcolx}[1]{\Yc{d}_{#1}}
  \newcommand{\inlen}{k}
  \newcommand{\outlen}{\ell}
  \newcommand{\infree}{g}
  \newcommand{\outfree}{j}
  \newcommand{\inbound}{h}
  \newcommand{\outbound}{i}
  \newcommand{\infreex}[1]{g_{#1}}
  \newcommand{\outfreex}[1]{j_{#1}}
  \newcommand{\inboundx}[1]{h_{#1}}
  \newcommand{\outboundx}[1]{i_{#1}}
  \newcommand{\incounter}{a}
  \newcommand{\outcounter}{b}
  \newcommand{\thepartition}{p}
  \newcommand{\sumind}{s}
  \newcommand{\anycol}{\Yc{c}}
  \newcommand{\thepartitionset}{\mathcal{P}}
  \newcommand{\theflip}{\SbraidingX{\thealg}{\thealg}}
  \newcommand{\thepolys}{\freealg{\comps}{\thegens}}
      The following definition of ``easy'' quantum groups is the algebraic version of \cite[Definition~5.1]{TarragoWeber2016}. Recall that $\thedim\in\pint$ is fixed per Assumptions~\ref{assumptions:easy}.
      \begin{Notation}
        \label{notation:relations_of_partition_set}
 For any set $\thepartitionset$ of two-colored partitions, let
    \[
      \thepartrels{\thepartitionset} :=\big\{ \therelpoly{\incol}{\outcol}{\thepartition}{\outfree}{\infree}    \vert \  (\incol,\outcol,\thepartition)\in \thecat\Sand \infree\in\dwi{\thedim}^{\Ssetmonoidalproduct\Slength{\incol}}\Sand\outfree\in\dwi{\thedim}^{\Ssetmonoidalproduct\Slength{\outcol}}\big\}
    \]
    and let $\thepartideal{\thepartitionset}$ be the two-sided ideal of $\thepolys$ generated by $\thepartrels{\thepartitionset}$.
      \end{Notation}
      \begin{Definition}
        \label{definition:easy-quantum-group}
      For any category $\thecat$  of two-co\-lo\-red partitions, the \emph{unitary easy compact quantum group} of $(\thecat,\thedim)$ is given by $(\univalg{\comps}{\thegens}{\thepartrels{\thecat}},\thestar,\thecomultiplication)$,
    where $\thestar$ and $\thecomultiplication$ are respectively the unique anti-multiplicative anti-linear self-map of $\univalg{\comps}{\thegens}{\thepartrels{\thecat}}$ and the unique multiplicative linear map from $\univalg{\comps}{\thegens}{\thepartrels{\thecat}}$ to the tensor product algebra of $\univalg{\comps}{\thegens}{\thepartrels{\thecat}}$ with itself which satisfy respectively%
    \begin{gather*}
      (\theunix{\anycol}{\indj}{\indi}{}+\thepartideal{\thecat})^\thestar= \theunix{\overline{\anycol}}{\indj}{\indi}+\thepartideal{\thecat} \qquad\tand \qquad \thecomultiplication(\theunix{\anycol}{\indj}{\indi}+\thepartideal{\thecat})=\sum_{\sumind=1}^\thedim(\theunix{\anycol}{\indj}{\sumind}+\thepartideal{\thecat})\Smonoidalproduct(\theunix{\anycol}{\sumind}{\indi}+\thepartideal{\thecat})
    \end{gather*}
    for any $\{\indi,\indj\}\subseteq \dwi{\thedim}$ and  $\anycol\in\blaw$.
    \end{Definition}

    \begin{Remark}
    \newcommand{\somelen}{m}
    \newcommand{\someind}{e}
    \newcommand{\someindx}[1]{e_{#1}}
    \newcommand{\dualincolx}[1]{\ASdualX{\incol}_{#1}}
    \newcommand{\dualoutcolx}[1]{\ASdualX{\outcol}_{#1}}
    \newcommand{\dualinfreex}[1]{\ASdualX{\infree}_{#1}}
    \newcommand{\dualinboundx}[1]{\ASdualX{\inbound}_{#1}}
    \newcommand{\dualoutfreex}[1]{\ASdualX{\outfree}_{#1}}
    \newcommand{\dualoutboundx}[1]{\ASdualX{\outbound}_{#1}}
      The definition of unitary easy quantum groups is usually given in terms of universal $\ast$-algebras, not universal algebras, cf.\ \cite[Definition~5.1]{TarragoWeber2016}. The variant given above is equivalent, as explained hereafter. For any $\somelen\in\nnint$ and any     $\someind\in\dwi{\thedim}^{\Ssetmonoidalproduct\somelen}$ let $\ASdualX{\someind}\in\dwi{\thedim}^{\Ssetmonoidalproduct\somelen}$ be defined by $\indi\mapsto \someindx{\somelen-\indi+1}$ for any $\indi\in\dwi{\somelen}$. Let   $ \{\inlen,\outlen\}\subseteq \nnint$, let $\incol\in\blaw^{\Ssetmonoidalproduct\inlen}$, let $ \outcol\in\blaw^{\Ssetmonoidalproduct\outlen}$, let $(\incol,\outcol,\thepartition)\in \thecat$, let $\infree\in\dwi{\thedim}^{\Ssetmonoidalproduct\inlen}$ and let $\outfree\in\dwi{\thedim}^{\Ssetmonoidalproduct\outlen}$. Then,  with respect to the $\ast$-map in De\-fi\-ni\-tion~\ref{definition:easy-quantum-group},
    \begin{align*}
      (\therelpoly{\incol}{\outcol}{\thepartition}{\outfree}{\infree})\Sadj&{}=      \sum_{\outbound\in\dwi{\thedim}^{\Ssetmonoidalproduct \outlen}}
      \zetfx{\thepartition}{\ker(\djp{\infree}{\outbound})}
\leftwardproduct{\outcounter=1}{\outlen}{(\theunix{\outcolx{\outcounter}}{\outfreex{\outcounter}}{\outboundx{\outcounter}})\Sadj}
-\sum_{\inbound\in\dwi{\thedim}^{\Ssetmonoidalproduct \inlen}}  \zetfx{\thepartition}{\ker(\djp{\inbound}{\outfree})} \leftwardproduct{\incounter=1}{\inlen}{(\theunix{\incolx{\incounter}}{\inboundx{\incounter}}{\infreex{\incounter}})\Sadj}\\
&{}=      \sum_{\outbound\in\dwi{\thedim}^{\Ssetmonoidalproduct \outlen}}
      \zetfx{(\ASdualX{\thepartition})\Sadj}{\ker(\djp{\ASdualX{\infree}}{\ASdualX{\outbound}})}
\rightwardproduct{\outcounter=1}{\outlen}{\theunix{\dualoutcolx{\outcounter}}{\dualoutfreex{\outcounter}}{\dualoutboundx{\outcounter}}}
-\sum_{\inbound\in\dwi{\thedim}^{\Ssetmonoidalproduct \inlen}}  \zetfx{(\ASdualX{\thepartition})\Sadj}{\ker(\djp{\ASdualX{\inbound}}{\ASdualX{\outfree}})} \rightwardproduct{\incounter=1}{\inlen}{\theunix{\dualincolx{\incounter}}{\dualinboundx{\incounter}}{\dualinfreex{\incounter}}}\\
&{}=\therelpoly{\ASdualX{\incol}}{\ASdualX{\outcol}}{(\ASdualX{\thepartition})\Sadj}{\ASdualX{\outfree}}{\ASdualX{\infree}}.
\end{align*}
Since $\thecat$ also contains the two-colored partition $(\ASdualX{\incol},\ASdualX{\outcol},(\ASdualX{\thepartition})\Sadj)$ by Lem\-ma~\hyperref[lemma:scalar-article-categories-helper-0]{\ref*{lemma:scalar-article-categories-helper}\,\ref*{lemma:scalar-article-categories-helper-0}}, the switch from the universal $\ast$-algebra to the universal algebra makes no difference.
\end{Remark}
For the idea of the proof of the following, see \cite[Remark~5.2]{TarragoWeber2016}.

\begin{Proposition}
  \label{proposition:counit_of_easy_quantum_group}
  For any category $\thecat$ of two-colored partitions, the unitary easy compact quantum group of $(\thecat,\thedim)$ is a compact quantum group whose co-unit is given by the unique multiplicative linear functional $\thecounit$ with
  \[
   \thecounit(\theunix{\anycol}{\indj}{\indi}+\thepartideal{\thecat})=\kron{\indj}{\indi}
  \]
  for any $\{\indi,\indj\}\subseteq \dwi{\thedim}$ and  $\anycol\in\blaw$. It can be seen as a compact $(\thedim\times\thedim)$-matrix quantum group with fundamental representation induced by $\theuni^\whpoint$.
\end{Proposition}
}

\begin{Example}
  Let $\thecat$ be a category of two-co\-lo\-red partitions and let $\theqg$ be the unitary quantum group of $(\thecat,\thedim)$.
  \begin{enumerate}
    \item If $\thecat$ is the minimal category $\langle \varnothing\rangle$, then $\theqg$ is the \emph{free unitary quantum group} $\ueqgUplus{\thedim}$ introduced by Wang in \cite{Wang1995b}. Its algebra can be presented as the universal algebra generated by $\thegens$ (fixed in Assumptions~\ref{assumptions:easy}) subject to only the relations of Lem\-ma~\ref{lemma:unitarity-relations}.
  \item For $\thecat=\langle \PartCrossWW\rangle$, we recover the classical unitary group $\ueqgU{\thedim}$, the universal commutative(!) algebra subject to the relations of Lem\-ma~\ref{lemma:unitarity-relations}.
          \item Should $\thecat$ be the maximal category $\langle \PartCrossWW,\PartFourWBWB,\PartIdenLoWW,\PartSingleW\rangle$ of all two-co\-lo\-red partitions, then $\theqg$ is the symmetric group $S_n$.
  \end{enumerate}
\end{Example}
Currently, there is no complete list of all unitary easy quantum groups because the classification of all categories of two-co\-lo\-red partitions is not yet finished.
\begin{Remark}
  \label{remark:obsh_cases}
  \newcommand{\anyblock}{\Yb{B}}
  \newcommand{\anypartition}{p}
  \newcommand{\anycat}{\mathcal{C}}
  \newcommand{\incol}{\Yc{c}}
  \newcommand{\outcol}{\Yc{d}}
  \newcommand{\anyset}{\mathcal{S}}
  \newcommand{\anygroup}{G_n}
  \newcommand{\dummycase}{\mathcal{X}}
  The names $\mathcal{O}$, $\mathcal{B}$, $\mathcal{H}$ and $\mathcal{S}$ of the four cases from Definition~\ref{definition:obsh_cases} were introduced by Tarrago and Weber in \cite[Definition~2.2]{TarragoWeber2018} and refer respectively to the orthogonal group $O_n$, bistochastic group $B_n$, hyperoctahedral group $H_n$ and symmetric group $S_n$. (A ``bistochastic'' matrix is understood to be an orthogonal matrix each of whose rows and columns  sums to $1$.) Tarrago and Weber showed that for each $\dummycase\in \{\mathcal{O},\mathcal{B},\mathcal{H},\mathcal{S}\}$ there exists a category of two-colored partitions which is  case $\dummycase$ and maximally so. And, in the sense of Definition~\ref{definition:easy-quantum-group}, the maximal case-$\mathcal{O}$ category is the one associated with $O_n$, the maximal case-$\mathcal{B}$ category the one associated with $B_n$, the maximal case-$\mathcal{H}$ category the one associated with $H_n$ and the maximal case-$\mathcal{S}$ category the one associated with $S_n$.
\end{Remark}

\section{First Hochschild cohomology of universal algebras}
\label{section:hochschild}
For the convenience of the reader, Section~\ref{section:hochschild} recalls the definition of and some elementary results about the first Hochschild cohomology. Throughout, let $\thefield$ be any field.

\subsection{First Hochschild cohomology}
\label{section:hochschild-general}
In Section~\ref{section:hochschild-general}, our algebra shall remain abstract. Section~\ref{section:hochschild-particular} will then recall which conclusions can be drawn if a presentation of the algebra in terms of generators and relations is given.
{
  \newcommand{\thealgebra}{A}
\begin{Assumptions}
  Let $\thealgebra$ be any $\thefield$-algebra and $\themodule$ any $\thealgebra$-bimodule.
\end{Assumptions}
{
  \newcommand{\firstop}{a_1}
  \newcommand{\secondop}{a_2}
  \newcommand{\anyvector}{x}
That means in particular that $\themodule$ is a $\thefield$-vector space implicitly equipped with $\thefield$-linear maps  $\xfromto{\Sleftmoduleaction}{\thealgebra\Smonoidalproduct\themodule}{\themodule}$ and $\xfromto{\Srightmoduleaction}{\themodule\Smonoidalproduct\thealgebra}{\themodule}$, the left and right \emph{actions} of $\thealgebra$,  such that $\firstop\Sleftmoduleaction(\secondop\Sleftmoduleaction\anyvector)=(\firstop\secondop)\Sleftmoduleaction\anyvector$ and $(\anyvector\Srightmoduleaction\secondop)\Srightmoduleaction\firstop=\anyvector\Srightmoduleaction(\secondop\firstop)$ and $(\firstop\Sleftmoduleaction \anyvector)\Srightmoduleaction\secondop=\firstop\Sleftmoduleaction( \anyvector\Srightmoduleaction\secondop)$ for any $\anyvector\in\themodule$ and $\{\firstop,\secondop\}\subseteq \thealgebra$.
}
\begin{Example}
  \label{example:trivial_bimodule}
  \newcommand{\anyscalar}{\lambda}
  \newcommand{\anyop}{a}
    For any augmentation $\thecounit$ of $\thealgebra$, i.e., any $\thefield$-algebra morphism from $\thealgebra$ to $\thefield$, the $\thefield$-vector space $\thefield$ becomes an $\thealgebra$-bimodule $\themodule$ if equipped with the
    actions defined by $\anyop\Sleftmoduleaction\anyscalar= \thecounit(\anyop)\anyscalar$ respectively $\anyscalar\Srightmoduleaction\anyop= \anyscalar\thecounit(\anyop)$ for any $\anyscalar\in\thefield$ and $\anyop\in \thealgebra$. It is often called the \emph{trivial bimodule} of~$(\thealgebra,\thecounit)$.
\end{Example}
\subsubsection{The fundamental definitions}
\label{section:hochschild-general-definitions}
The following definitions were first given by Hochschild in \cite{Hochschild1956}.

\begin{Definition}
  \label{definition:hochschild-one-cocycles-and-coboundaries} \
  \begin{enumerate}\itemsep=0pt
  \item
    \label{definition:hochschild-one-cocycles-and-coboundaries-1}
    {
      \newcommand{\someonecocycle}{\eta}
      \newcommand{\firstelement}{a_1}
      \newcommand{\secondelement}{a_2}
      The \emph{$\themodule$-valued Hochschild $1$-cocycles of $\thealgebra$} are the $\thefield$-vector subspace $\HScocycles{1}{\thealgebra}{\themodule}$ of $\SinternalhomX{\thealgebra}{\themodule}$ formed by  all elements $\someonecocycle$ such that    for any $\{\firstelement,\secondelement\}\subseteq \thealgebra$,
    \[
      \bigl(\partial^1\someonecocycle\bigr)(\firstelement\Smonoidalproduct\secondelement):= \firstelement\Sleftmoduleaction\someonecocycle(\secondelement)-\someonecocycle(\firstelement\secondelement)+\someonecocycle(\firstelement)\Srightmoduleaction\secondelement=0.
    \]
  }
\item
  \label{definition:hochschild-one-cocycles-and-coboundaries-2}
  {
    \newcommand{\someonecocycle}{\eta}
    \newcommand{\somevector}{x}
     \newcommand{\someelement}{a}
     The \emph{$\themodule$-valued Hochschild $1$-coboundaries of $\thealgebra$} are the $\thefield$-vector subspace $\HScoboundaries{1}{\thealgebra}{\themodule}$ of $\SinternalhomX{\thealgebra}{\themodule}$ formed by  all elements $\someonecocycle$ such that there exists $\somevector\in\themodule$ with     for any $\someelement\in\thealgebra$,
     \[
       \someonecocycle(\someelement)=\bigl(\partial^0 \somevector\bigr)(\someelement):= \someelement\Sleftmoduleaction\somevector -\somevector\Srightmoduleaction\someelement.
     \]
   }
  \end{enumerate}
\end{Definition}
{
  \newcommand{\somecocycle}{\eta}
It can be seen that $\somecocycle(1)=0$ for any $\somecocycle\in\HScocycles{1}{\thealgebra}{\themodule}$ and that $\HScoboundaries{1}{\thealgebra}{\themodule}$ is a $\thefield$-vector sub\-space of $\HScocycles{1}{\thealgebra}{\themodule}$.
}

\begin{Definition}
  \label{definition:hochschild-one-cohomology}
We call the quotient $\thefield$-vector space $\HScohomology{1}{\thealgebra}{\themodule}$ of   $\HScocycles{1}{\thealgebra}{\themodule}$ with respect to~$\HScoboundaries{1}{\thealgebra}{\themodule}$ the \emph{first Hochschild cohomology of  $\thealgebra$ with  $\themodule$-coefficients}.
\end{Definition}

\begin{Example}
  \label{remark:trivial_bimodule}
  \newcommand{\anyscalar}{\lambda}
  \newcommand{\anyop}{a}
  In the case of Example~\ref{example:trivial_bimodule}, i.e., for trivial coefficients, the only $1$-coboundary of $\thealgebra$ is the zero map because $(\partial^0\anyscalar)(\anyop)=\thecounit(\anyop)\anyscalar-\anyscalar\thecounit(\anyop)=0$ for any $\anyscalar\in\thefield$ and $\anyop\in\thealgebra$. Hence,  $\HScocycles{1}{\thealgebra}{\themodule}\cong\HScohomology{1}{\thealgebra}{\themodule}$ in that instance.
\end{Example}
\subsubsection{Algebra hom characterization of 1-coycles}
\label{section:hochschild-general-algebrahoms}
$1$-cocycles can be characterized as certain algebra homomorphisms by means of  a folk theorem recorded as  \cite[Lemma~1.9]{KyedRaum2017}. The latter uses the following construction.

{
    \newcommand{\thedeform}{c}
  \newcommand{\thedeformxx}[2]{c(#1 \Smonoidalproduct #2)}
  \newcommand{\thirda}{a_3}
        \newcommand{\firsta}{a_1}
      \newcommand{\firstx}{x_1}
      \newcommand{\seconda}{a_2}
      \newcommand{\secondx}{x_2}
      \begin{Definition}
        \label{definition:the-bfg-algebra}
    Let $\HSalgebraconstructdeform{1}{\thealgebra}{\themodule}{0}$ denote the  $\thefield$-vector space $\thealgebra\Sdirectsum\themodule$ equipped with the $\thefield$-linear map $\fromto{\thefield} \thealgebra\Sdirectsum\themodule$ with $1\mapsto (1,0)$ and the $\thefield$-linear map $\fromto{(\thealgebra\Sdirectsum\themodule)^{\Smonoidalproduct 2}}{\thealgebra\Sdirectsum\themodule}$ defined by     for any $\{\firsta,\seconda\}\subseteq \thealgebra$ and $\{\firstx,\secondx\}\subseteq \themodule$,
    \[
      (\firsta,\firstx)\Smonoidalproduct(\seconda,\secondx)\mapsto (\firsta\seconda,\firsta\Sleftmoduleaction\secondx+\firstx\Srightmoduleaction\seconda).
    \]
  \end{Definition}
}

{
  \newcommand{\thedeform}{c}
  \newcommand{\anyelement}{a}
  \newcommand{\firstelement}{a_1}
      \newcommand{\secondelement}{a_2}
      \newcommand{\thirdelement}{a_3}
  \newcommand{\thedeformxx}[2]{c(#1 \Smonoidalproduct #2)}
  \begin{Lemma}
    \label{lemma:the-bfg-algebra}
     $\HSalgebraconstructdeform{1}{\thealgebra}{\themodule}{0}$ is a $\thefield$-algebra.
  \end{Lemma}
  \begin{proof}
      \newcommand{\somea}{a}
      \newcommand{\somex}{x}
      \newcommand{\firsta}{a_1}
      \newcommand{\seconda}{a_2}
      \newcommand{\thirda}{a_3}
      \newcommand{\firstx}{x_1}
      \newcommand{\secondx}{x_2}
      \newcommand{\thirdx}{x_3}
$(1,0)$ is a unit because for any $\somea\in\thealgebra$ and any $\somex\in\themodule$,
      \[
        (\somea,\somex)\cdot(1,0)=(\somea 1, \somea\Sleftmoduleaction 0+\somex\Srightmoduleaction 1)=(\somea , \somex)= (1 \somea , 1\Sleftmoduleaction \somex+0\Srightmoduleaction \somea)=  (1,0)\cdot(\somea,\somex).
      \]
 Moreover, for any $\{\firsta,\seconda,\thirda\}\subseteq \thealgebra$ and any $\{\firstx,\secondx,\thirdx\}\subseteq \themodule$,
      \begin{align*}
        ((\firsta,\firstx)\cdot (\seconda,\secondx))\cdot(\thirda,\thirdx)
      &{}=(\firsta\seconda,\firsta\Sleftmoduleaction\secondx+\firstx\Srightmoduleaction\seconda)\cdot (\thirda,\thirdx)\\
        &{}=(\firsta\seconda\thirda,(\firsta\seconda)\Sleftmoduleaction\thirdx+(\firsta\Sleftmoduleaction\secondx+\firstx\Srightmoduleaction\seconda)\Srightmoduleaction\thirda)\\        &{}=(\firsta\seconda\thirda,\firsta\seconda\Sleftmoduleaction\thirdx+\firsta\Sleftmoduleaction\secondx\Srightmoduleaction\thirda+\firstx\Srightmoduleaction\seconda\thirda) \\        &{}=(\firsta\seconda\thirda,\firsta\Sleftmoduleaction(\seconda\Sleftmoduleaction\thirdx+\secondx\Srightmoduleaction\thirda)+\firstx\Srightmoduleaction(\seconda\thirda))\\
        &{}=  (\firsta,\firstx)\cdot(\seconda\thirda,\seconda\Sleftmoduleaction\thirdx+\secondx\Srightmoduleaction\thirda)\\
        &{}= (\firsta,\firstx)\cdot ((\seconda,\secondx)\cdot(\thirda,\thirdx)),
\end{align*}
which shows that the multiplication is associative.
  \end{proof}
  }

  The following is then  the folk theorem mentioned in \cite[Lem\-ma~1.9]{KyedRaum2017}.

  {
      \newcommand{\thedeform}{c}
  \newcommand{\thedeformxx}[2]{c(#1 \Smonoidalproduct #2)}
  \newcommand{\anyelement}{a}
  \newcommand{\firstelement}{a_1}
      \newcommand{\secondelement}{a_2}
      \newcommand{\thirdelement}{a_3}
    \newcommand{\somefunc}{\psi}
      \begin{Lemma}
    \label{lemma:the-bfg-lemma}
    \newcommand{\somea}{a}
For any $\somefunc\in\SinternalhomX{\thealgebra}{\themodule}$, the map $\fromto{\thealgebra}{\thealgebra\Sdirectsum\themodule}$ with $\somea\mapsto (\somea,\somefunc(\somea))$ for any $\somea\in \thealgebra$ defines a $\thefield$-algebra homomorphism  $\fromto{\thealgebra}{\HSalgebraconstructdeform{1}{\thealgebra}{\themodule}{0}}$ if and only if $\somefunc\in \HScocycles{1}{\thealgebra}{\themodule}$.
    \end{Lemma}
    \begin{proof}
    \newcommand{\somea}[1]{a_{#1}}
    \newcommand{\firsta}{a_{1}}
    \newcommand{\seconda}{a_{2}}
    \newcommand{\themap}{f_\somefunc}
    If the map in the claim is denoted by $\themap$, then  $\themap(1)=(1,\somefunc(1))$ and for any $\{\firsta,\seconda\}\subseteq \thealgebra$, obviously, $\themap(\firsta\seconda)=(\firsta\seconda,\somefunc(\firsta\seconda))$ and
    \[
      \themap(\firsta)\cdot\themap(\seconda)=(\firsta,\somefunc(\firsta))\cdot(\seconda,\somefunc(\seconda))=(\firsta\seconda,\firsta\Sleftmoduleaction\somefunc(\seconda)+\somefunc(\firsta)\Srightmoduleaction\seconda).
    \]
    The two values coincide if and only if $\somefunc(\firsta\seconda)=\firsta\Sleftmoduleaction\somefunc(\seconda)+\somefunc(\firsta)\Srightmoduleaction\seconda$, which is to say if and only if $\somefunc\in \HScocycles{1}{\thealgebra}{\themodule}$. As then $\somefunc(1)=0$ the claim is true.
    \end{proof}
  }

The next result will be required later in the proof of, ultimately, Pro\-po\-si\-tion~\ref{proposition:hochschild-particular-one-cohomology}.

{
    \newcommand{\thedeform}{c}
  \newcommand{\thedeformxx}[2]{c(#1 \Smonoidalproduct #2)}
  \newcommand{\thirda}{a_3}
    \begin{Lemma}
    \newcommand{\somenumber}{m}
    \newcommand{\somea}[1]{a_{#1}}
    \newcommand{\somex}[1]{x_{#1}}
    \newcommand{\someindex}{i}
    \newcommand{\otherindex}{j}
    \label{lemma:the-bfg-algebra-multiplication}
    For any  $\somenumber\in\pint$, any $\{\somea{\someindex}\}_{\someindex=1}^\somenumber\subseteq \thealgebra$ and any $\{\somex{\someindex}\}_{\someindex=1}^\somenumber\subseteq \themodule$,  in $\HSalgebraconstructdeform{1}{\thealgebra}{\themodule}{0}$,
    \[
      \tsrightwardproduct{\someindex=1}{\somenumber} (\somea{\someindex},\somex{\someindex})=\left(\tsrightwardproduct{\someindex=1}{\somenumber}\somea{\someindex}, \sum_{\someindex=1}^\somenumber\left(\tsrightwardproduct{\otherindex=1}{\someindex-1}\somea{\otherindex}\right)\Sleftmoduleaction\somex{\someindex}\Srightmoduleaction\left(\tsrightwardproduct{\otherindex=\someindex+1}{\somenumber}\somea{\otherindex}\right)\right).
    \]
\end{Lemma}
\begin{proof}
    \newcommand{\somenumber}{m}
    \newcommand{\somea}[1]{a_{#1}}
    \newcommand{\somex}[1]{x_{#1}}
    \newcommand{\someindex}{i}
    \newcommand{\otherindex}{j}
    The cases $\somenumber\in\{1,2,3\}$ are, respectively,  trivial, the definition of the multiplication  of $\HSalgebraconstructdeform{1}{\thealgebra}{\themodule}{0}$ and an intermediate result in the proof of Lem\-ma~\ref{lemma:the-bfg-algebra}. Generally,
    \begin{gather*} \left(\tsrightwardproduct{\someindex=1}{\somenumber-1}\somea{\someindex},{ \sum_{\someindex=1}^{\somenumber-1}}\left(\tsrightwardproduct{\otherindex=1}{\someindex-1}\somea{\otherindex}\right)\Sleftmoduleaction\somex{\someindex}\Srightmoduleaction\left(\tsrightwardproduct{\otherindex=\someindex+1}{\somenumber-1}\somea{\otherindex}\right)\right)\cdot (\somea{\somenumber},\somex{\somenumber})\\
\qquad{}      = \left(\!\left(\tsrightwardproduct{\someindex=1}{\somenumber-1}\somea{\someindex}\right)\somea{\somenumber},  \left(\tsrightwardproduct{\someindex=1}{\somenumber-1}\somea{\someindex}\right)\Sleftmoduleaction\somex{\somenumber}+{ \sum_{\someindex=1}^{\somenumber-1}}\left(\tsrightwardproduct{\otherindex=1}{\someindex-1}\somea{\otherindex}\right)\Sleftmoduleaction\somex{\someindex}\Srightmoduleaction\left(\tsrightwardproduct{\otherindex=\someindex+1}{\somenumber-1}\somea{\otherindex}\right)\Srightmoduleaction \somea{\somenumber}\!\right).
\end{gather*}
Hence, the claim is true.
\end{proof}
  }

}

\subsection{Conclusions for universal algebras}
\label{section:hochschild-particular}

{
  \newcommand{\thepolys}{\freealg{\thefield}{\thegens}}
  \newcommand{\thealgebra}{\univalg{\thefield}{\thegens}{\therels}}
  Using Lem\-ma~\ref{lemma:the-bfg-lemma}, it is possible to give a canonical equational characterization of the $1$-cocycles if a presentation of the algebra in terms of generators and relations is given.
  \begin{Assumptions}
        Let $\thegens$ be any set, $\therels\subseteq \thepolys$ arbitrary, $\theideal$ the two-sided ideal of $\thepolys$ generated by $\therels$, and $\themodule$ any $\thealgebra$-bimodule.
  \end{Assumptions}

{
  \newcommand{\somepoly}{p}
  \newcommand{\somegen}{e}
  \newcommand{\them}{m}
  \newcommand{\thees}{e}
  \newcommand{\thee}[1]{e_{#1}}
  \newcommand{\indi}{i}
  \newcommand{\indj}{j}
  \newcommand{\thetuple}{x}
  \newcommand{\thetuplex}[1]{x_{#1}}
  \newcommand{\somerel}{r}
  \newcommand{\somevector}{z}
  \begin{Definition}
    \label{definition:derivatives-one}
Let
  \[
\xfromtomaps {\HSuafunc{1}{\thegens}{\therels}{\themodule}}{\thepolys}{\SinternalhomX{\themodule^{\Sdirectproduct\thegens}}{\themodule}}{\somepoly}{\HSuafuncx{1}{\thegens}{\therels}{\themodule}{\somepoly}}
  \]
  be the unique $\thefield$-linear map with for any $\them\in\pint$ and any  $\{\thee{\indi}\}_{\indi=1}^\them\subseteq \thegens$, if $\somepoly=\tsrightwardproduct{\indi=1}{\them}{\thee{\indi}}$, then for any~${\thetuple\in \themodule^{\Sdirectproduct \thegens}}$,
  \[
    \HSuafuncxx{1}{\thegens}{\therels}{\themodule}{\somepoly}{\thetuple}=  \sum_{\indi=1}^\them \left(\tsrightwardproduct{\indj=1}{\indi-1}{\thee{\indj}}+\theideal\right)\Sleftmoduleaction\thetuplex{\thee{\indi}}\Srightmoduleaction\left(\tsrightwardproduct{\indj=\indi+1}{\them}{\thee{\indj}}+\theideal\right),
    \]
    and with $\HSuafuncx{1}{\thegens}{\therels}{\themodule}{1}=0$.
  \end{Definition}
  \begin{Example}
    \label{example:trivial_bimodule_functional}
    \newcommand{\shortalgebra}{A}
    For any augmentation $\thecounit$ of $\shortalgebra=\thealgebra$, if $\themodule$ is the trivial bimodule of $(\shortalgebra,\thecounit)$ in the sense of Example~\ref{example:trivial_bimodule}, then
    for any $\them\in\pint$ and any  $\{\thee{\indi}\}_{\indi=1}^\them\subseteq \thegens$, if $\somepoly=\thee{1}\cdots \thee{\them}$, then $\HSuafuncx{1}{\thegens}{\therels}{\themodule}{\somepoly}$ is a linear map which assigns to any family  $\thetuple=(\thetuplex{\somegen})_{\somegen\in\thegens}$ of elements of $\thefield$ the number
  \[
    \HSuafuncxx{1}{\thegens}{\therels}{\themodule}{\somepoly}{\thetuple}=  \sum_{\indi=1}^\them\Biggl(\prod_{\indj\in \dwi{\them}\backslash \{\indi\}}\thecounit\left(\thee{\indj}+\theideal\right)\Biggr)\thetuplex{\thee{\indi}}.
    \]
  \end{Example}

  \begin{Definition}\quad
    \begin{enumerate}
  \item Let $\HSuacocycles{1}{\thegens}{\therels}{\themodule}$ denote the $\thefield$-vector subspace of $\themodule^{\Sdirectproduct \thegens}$ of all elements $\thetuple$ with $\HSuafuncxx{1}{\thegens}{\therels}{\themodule}{\somerel}{\thetuple}=0$ for any $\somerel\in\therels$.
    \item Write $\HSuacoboundaries{1}{\thegens}{\therels}{\themodule}$ for the $\thefield$-vector subspace of $\themodule^{\Sdirectproduct \thegens}$ formed by all elements $\thetuple$ for which there exists $\somevector\in\themodule$ with $\thetuplex{\somegen}=(\somegen+\theideal)\Sleftmoduleaction\somevector-\somevector\Srightmoduleaction(\somegen+\theideal)$ for any $\somegen\in\thegens$.
    \end{enumerate}
  \end{Definition}
}

{
  \begin{Lemma}
    \label{lemma:hochschild-particular-one-cohomology-well-defined}
        {
        \newcommand{\thevector}{z}
        \newcommand{\somepoly}{p}
        \newcommand{\somegen}{e}
        \newcommand{\thetuple}{x}
        \newcommand{\thetuplex}[1]{\thetuple_{#1}}
        For any $\somepoly\in\thepolys$ and any $\thevector\in \themodule$, if $\thetuple\in \themodule^{\Sdirectproduct \thegens}$ is such that $\thetuplex{\somegen}=(\somegen+\theideal)\Sleftmoduleaction\thevector-\thevector\Srightmoduleaction(\somegen+\theideal)$ for any $\somegen\in\thegens$, then
        \[
          \HSuafuncxx{1}{\thegens}{\therels}{\themodule}{\somepoly}{\thetuple}=(\somepoly+\theideal)\Sleftmoduleaction\thevector-\thevector\Srightmoduleaction(\somepoly+\theideal).
        \]
        In particular, $\HSuacoboundaries{1}{\thegens}{\therels}{\themodule}$ is a $\thefield$-vector subspace of $\HSuacocycles{1}{\thegens}{\therels}{\themodule}$.
          }
  \end{Lemma}
  \begin{proof}
      {
      \newcommand{\thee}[1]{e_{#1}}
      \newcommand{\indi}{i}
      \newcommand{\indj}{j}
      \newcommand{\somenumber}{m}
              \newcommand{\thevector}{z}
        \newcommand{\somepoly}{p}
        \newcommand{\somegen}{e}
        \newcommand{\thetuple}{x}
        \newcommand{\anyvector}{w}
        \newcommand{\anyrel}{r}
        \newcommand{\thetuplex}[1]{\thetuple_{#1}}
      The claimed identity is clear if $\somepoly=1$. If there are  $\somenumber\in\pint$ and $\{\thee{\indi}\}_{\indi=1}^{\somenumber}\subseteq \thegens$ with $p=e_1\cdots e_m$, then by definition,
          \begin{align*}     \HSuafuncxx{1}{\thegens}{\therels}{\themodule}{\somepoly}{\thetuple}&{}= \sum_{\indi=1}^\somenumber \left(\tsrightwardproduct{\indj=1}{\indi-1}{\thee{\indj}+\theideal}\right)\Sleftmoduleaction\thetuplex{\thee{\indi}}\Srightmoduleaction\left(\tsrightwardproduct{\indj=\indi+1}{\somenumber}{\thee{\indj}+\theideal}\right)\\
      &{}= \sum_{\indi=1}^\somenumber \left(\tsrightwardproduct{\indj=1}{\indi-1}{\thee{\indj}+\theideal}\right)\Sleftmoduleaction((\thee{\indi}+\theideal)\Sleftmoduleaction\thevector-\thevector\Srightmoduleaction(\thee{\indi}+\theideal))\Srightmoduleaction\left(\tsrightwardproduct{\indj=\indi+1}{\somenumber}{\thee{\indj}+\theideal}\right)\\
      &{}= \left(\sum_{\indi=2}^{\somenumber+1} \left(\tsrightwardproduct{\indj=1}{\indi-1}{\thee{\indj}+\theideal}\right)\Sleftmoduleaction\thevector\Srightmoduleaction\left(\tsrightwardproduct{\indj=\indi}{\somenumber}{\thee{\indj}+\theideal}\right)\right)\\
      &{}\quad-\left(\sum_{\indi=1}^{\somenumber}\left(\tsrightwardproduct{\indj=1}{\indi-1}{\thee{\indj}+\theideal}\right)\Sleftmoduleaction\thevector\Srightmoduleaction\left(\tsrightwardproduct{\indj=\indi}{\somenumber}{\thee{\indj}+\theideal}\right)\right)\\
      &{}=\left(\tsrightwardproduct{\indi=1}{\somenumber}{\thee{\indi}+\theideal}\right)\Sleftmoduleaction\thevector-\thevector\Srightmoduleaction\left(\tsrightwardproduct{\indi=1}{\somenumber}{\thee{\indi}+\theideal}\right).
    \end{align*}
    Thus, the identity holds for arbitrary $\somepoly\in\thepolys$ by $\thefield$-linearity. It follows in particular that, if $\somepoly\in\therels$, then  $\HSuafuncxx{1}{\thegens}{\therels}{\themodule}{\somepoly}{\thetuple}=0$ since $\somepoly+\theideal$ is the zero vector of $\thealgebra$ in this case.
    That proves the claim about $\HSuacoboundaries{1}{\thegens}{\therels}{\themodule}$.
    }
  \end{proof}
}

In particular, the preceding lemma enables us to consider the following space.
\begin{Definition}
  \label{definition:hochschild-particular-one-cohomology}
  Let $\HSuacohomology{1}{\thegens}{\therels}{\themodule}$ be the $\thefield$-vector quotient space of $\HSuacocycles{1}{\thegens}{\therels}{\themodule}$ with respect to $\HSuacoboundaries{1}{\thegens}{\therels}{\themodule}$.
\end{Definition}
  The following notation allows referencing easily a multitude of algebra morphisms whose existence is implied by the universal property of $\thepolys$.
\begin{Notation}
  \label{notation:evaluation_of_polynomials}
  \newcommand{\anyalgebra}{B}
  \newcommand{\anytuple}{b}
  \newcommand{\anygen}{e}
  \newcommand{\anytuplex}[1]{\anytuple_{#1}}
  \newcommand{\anypoly}{p}
  \newcommand{\themorph}{g}
  Let  $\anyalgebra$ be any $\thefield$-algebra and let $(\anytuplex{\anygen})_{\anygen\in\thegens}\in \anyalgebra^{\Sdirectproduct \thegens}$ be arbitrary. The \emph{evaluation} of $\anypoly$ at $(\anytuplex{\anygen})_{\anygen\in\thegens}$ in $\anyalgebra$ is given by $\anypoly((\anytuplex{\anygen})_{\anygen\in\thegens}):= \themorph(\anypoly)$, where  $\themorph$ is the unique $\thefield$-algebra morphism  $\fromto{\thepolys}{\anyalgebra}$ with  $\anygen\mapsto\anytuplex{\anygen}$ for any $\anygen\in\thegens$.
\end{Notation}
{
        \newcommand{\somepoly}{p}
        \newcommand{\somegen}{e}
        \newcommand{\thetuple}{x}
        \newcommand{\thetuplex}[1]{\thetuple_{#1}}
  \begin{Lemma}
    \label{lemma:hochschild-particular-one-cohomology-helper}
      {
        For any $\somepoly\in\thepolys$ and any $\thetuple\in \themodule^{\Sdirectproduct \thegens}$, evaluating $\somepoly$ at $(\somegen+\theideal,\thetuplex{\somegen})_{\somegen\in\thegens}$ in the algebra $\HSalgebraconstructdeform{1}{\thealgebra}{\themodule}{0}$ yields
        \[
          \somepoly((\somegen+\theideal,\thetuplex{\somegen})_{\somegen\in\thegens})=(\somepoly+\theideal,\HSuafuncxx{1}{\thegens}{\therels}{\themodule}{\somepoly}{\thetuple}).
        \]
        }
  \end{Lemma}
  \begin{proof}
    \newcommand{\somenumber}{m}
    \newcommand{\indi}{i}
    \newcommand{\thee}[1]{e_{#1}}
Because $\HSuafuncxx{1}{\thegens}{\therels}{\themodule}{1}{\thetuple}=0$ by definition, the claim is true if $\somepoly=1$. If there exist $\somenumber\in\pint$ and $\{\thee{\indi}\}_{\indi=1}^\somenumber\subseteq \thegens$ with $\somepoly=\thee{1}\cdots\thee{\somenumber}$, then the claim follows immediately from Lem\-ma~\ref{lemma:the-bfg-algebra-multiplication} and Definition~\ref{definition:derivatives-one}. For arbitrary $\somepoly$, the assertion therefore holds by $\thefield$-linearity.
  \end{proof}
}

\begin{Remark}
  \label{remark:evaluation_of_polynomials}
  \newcommand{\anyalgebra}{B}
  \newcommand{\anytuple}{b}
  \newcommand{\anygen}{e}
  \newcommand{\anytuplex}[1]{\anytuple_{#1}}
  \newcommand{\anypoly}{p}
  \newcommand{\themorph}{f}
  \newcommand{\bigmorph}{g}
  \newcommand{\theproj}{\pi}
  Given any $\thefield$-algebra $\anyalgebra$,  any $\thefield$-algebra morphism $\xfromto{\themorph}{\thealgebra}{\anyalgebra}$ and any $\anypoly\in\thepolys$,
  \[
    \themorph(\anypoly+\theideal)=\anypoly((\themorph(\anygen+\theideal))_{\anygen\in\thegens}),
  \]
  where  the right-hand side is an evaluation of $\anypoly$ in $\anyalgebra$.

  Indeed, if $\bigmorph$ is the unique $\thefield$-algebra morphism $\fromto{\thepolys}{\anyalgebra}$ with $\anygen\mapsto \anytuplex{\anygen}:= \themorph(\anygen+\theideal)$ for any $\anygen\in\thegens$, then $\themorph(\anypoly+\theideal)=\bigmorph(\anypoly)=\anypoly((\anytuplex{\anygen})_{\anygen\in\thegens})=\anypoly((\themorph(\anygen+\theideal))_{\anygen\in\thegens})$, where the first identity holds by the uniqueness of $\bigmorph$ and where the second is nothing but an application of Notation~\ref{notation:evaluation_of_polynomials}.
\end{Remark}
Now we can give a useful characterization of the spaces of $1$-cocycles of universal algebras.
{
    \newcommand{\somegen}{e}
    \newcommand{\somerel}{r}
    \newcommand{\somevector}{x}
    \newcommand{\anyvector}{z}
    \newcommand{\thetuple}{x}
    \newcommand{\thetuplex}[1]{x_{#1}}
    \newcommand{\somederivation}{\eta}
    \newcommand{\somepoly}{p}
    \begin{Proposition}\quad
      \label{proposition:hochschild-particular-one-cohomology}
      \begin{enumerate}
      \item      \label{proposition:hochschild-particular-one-cohomology-1}

        A commutative diagram of $\thefield$-linear maps is given by
            $$
                    \begin{tikzcd}
\HScocycles{1}{\thealgebra}{\themodule}\arrow[r, two heads, hook]& \HSuacocycles{1}{\thegens}{\therels}{\themodule}\\
\HScoboundaries{1}{\thealgebra}{\themodule}\arrow[r, two heads, hook]\arrow[u,hook, "\subseteq" left]& \HSuacoboundaries{1}{\thegens}{\therels}{\themodule},\arrow[u,hook, "\subseteq" right]
            \end{tikzcd}
          $$
          where the horizontal arrows both assign to  any element $\somederivation$ of their respective domains the tuple
          \[
            (\somederivation(\somegen+\theideal))_{\somegen\in\thegens}.
          \]
          Moreover, the horizontal arrows are both $\thefield$-linear isomorphisms. Their respective inverses both assign to any element $\thetuple$ of their respective domains the mapping $\fromto{\thealgebra}{\themodule}$ with
          \[
           \somepoly+\theideal\mapsto \HSuafuncxx{1}{\thegens}{\therels}{\themodule}{\somepoly}{\thetuple}
          \]
           for any $\somepoly\in\thepolys$.
      \item      \label{proposition:hochschild-particular-one-cohomology-2}
          There exists an isomorphism of $\thefield$-vector spaces
    \[
      \begin{tikzcd}
        \displaystyle
        \HScohomology{1}{\thealgebra}{\themodule} \arrow[r, two heads, hook]&        \HSuacohomology{1}{\thegens}{\therels}{\themodule}
      \end{tikzcd}
    \]
    such that the class of any $\somederivation\in\HScocycles{1}{\thealgebra}{\themodule}$ is sent to the class of the $\thetuple\in \HSuacocycles{1}{\thegens}{\therels}{\themodule}$ with $\thetuplex{\somegen}=\somederivation(\somegen+\theideal)$ for any $\somegen\in\thegens$.
The inverse isomorphism sends the class of any $\thetuple\in \HSuacocycles{1}{\thegens}{\therels}{\themodule}$  to the class of the $\somederivation\in \HScocycles{1}{\thealgebra}{\themodule}$ with $\somederivation(\somepoly+\theideal)=\HSuafuncxx{1}{\thegens}{\therels}{\themodule}{\somepoly}{\thetuple}$ for any $\somepoly\in\thepolys$.
      \end{enumerate}
               \end{Proposition}
               \begin{proof}
           \newcommand{\shortalgebra}{A}
           \newcommand{\somealgebra}{B}
           \newcommand{\anycocycle}{\eta}
           \newcommand{\anyop}{a}
           \newcommand{\somemorph}{f}
           \newcommand{\horimorph}{\varphi}
           \newcommand{\invhorimorph}{\psi}
(a)
           Abbreviate $\shortalgebra:=\thealgebra$ and  $\somealgebra:=\HSalgebraconstructdeform{1}{\shortalgebra}{\themodule}{0}$.

           \subproof{Step~1. Upper horizontal arrow is well defined. $($And a bit more$.)$}
First, we prove that for any $\anycocycle\in\HScocycles{1}{\shortalgebra}{\themodule}$, if $\thetuplex{\somegen}:= \anycocycle(\somegen+\theideal)$ for any $\somegen\in\thegens$, then $\anycocycle(\somepoly+\theideal)=\HSuafuncxx{1}{\thegens}{\therels}{\themodule}{\somepoly}{\thetuple}$ for any $\somepoly\in\thepolys$. That then in particular shows that $\HSuafuncxx{1}{\thegens}{\therels}{\themodule}{\somerel}{\thetuple}=0$ for any $\somerel\in\therels$.

By $\anycocycle\in\HScocycles{1}{\shortalgebra}{\themodule}$,
according to Lem\-ma~\ref{lemma:the-bfg-lemma}, the rule that $\anyop\mapsto(\anyop,\anycocycle(\anyop))$ for any $\anyop\in\shortalgebra$ defines a~$\thefield$-algebra homomorphism $\xfromto{\somemorph}{\shortalgebra}{\somealgebra}$. Hence, for any $\somepoly\in\thepolys$ it must hold that $(\somepoly+\theideal,\anycocycle(\somepoly+\theideal))=\somemorph(\somepoly+\theideal)=\somepoly((\somemorph(\somegen+\theideal))_{\somegen\in\thegens})=\somepoly((\somegen+\theideal,\anycocycle(\somegen+\theideal))_{\somegen\in\thegens})=\somepoly((\somegen+\theideal,\thetuplex{\somegen})_{\somegen\in\thegens})=(\somepoly+\theideal,\HSuafuncxx{1}{\thegens}{\therels}{\themodule}{\somepoly}{\thetuple})$ in $\somealgebra$, where the second and last identities are implied by Re\-mark~\ref{remark:evaluation_of_polynomials} and Lem\-ma~\ref{lemma:hochschild-particular-one-cohomology-helper}, respectively. Hence,  $\HSuafuncxx{1}{\thegens}{\therels}{\themodule}{\somepoly}{\thetuple}=\anycocycle(\somepoly+\theideal)$ for any $\somepoly\in\thepolys$, as claimed.

           \subproof{Step~2. Alleged inverse upper horizontal arrow well defined.}
           Next, we show that for any $\thetuple\in \themodule^{\Sdirectproduct \thegens}$ with $\HSuafuncxx{1}{\thegens}{\therels}{\themodule}{\somerel}{\thetuple}=0$ for any $\somerel\in\therels$ there exists  $\anycocycle\in\HScocycles{1}{\shortalgebra}{\themodule}$ with $\anycocycle(\somepoly+\theideal)=\HSuafuncxx{1}{\thegens}{\therels}{\themodule}{\somepoly}{\thetuple}$ for any $\somepoly\in\thepolys$. One consequence of this is then  that the alleged inverse upper  horizontal arrow is well defined.

For any $\somerel\in\therels$ because $\HSuafuncxx{1}{\thegens}{\therels}{\themodule}{\somerel}{\thetuple}=0$ and $\somerel\in\theideal$ we can infer by Lem\-ma~\ref{lemma:hochschild-particular-one-cohomology-helper} that  $\somerel((\somegen+\theideal,\thetuplex{\somegen})_{\somegen\in\thegens})=(\theideal,0)$ in $\somealgebra$. The universal property of $\shortalgebra$ therefore guarantees the existence of a~unique $\thefield$-algebra homomorphism $\xfromto{\somemorph}{\shortalgebra}{\somealgebra}$ with $\somemorph(\somegen+\theideal)=(\somegen+\theideal,\thetuplex{\somegen})$ for any $\somegen\in\thegens$. More generally, for any $\somepoly\in\thepolys$ it must hold that $\somemorph(\somepoly+\theideal)=\somepoly((\somemorph(\somegen+\theideal))_{\somegen\in\thegens})=(\somepoly+\theideal,\HSuafuncxx{1}{\thegens}{\therels}{\themodule}{\somepoly}{\thetuple})$, where the two identities are again due to Re\-mark~\ref{remark:evaluation_of_polynomials} and Lem\-ma~\ref{lemma:hochschild-particular-one-cohomology-helper}. In other words, if $\anycocycle(\somepoly+\theideal):=\HSuafuncxx{1}{\thegens}{\therels}{\themodule}{\somepoly}{\thetuple}$ for any $\somepoly\in\thepolys$, then the rule that $\anyop\mapsto (\anyop,\anycocycle(\anyop))$ for any $\anyop\in \shortalgebra$ defines a $\thefield$-algebra homomorphism $\fromto{\shortalgebra}{\somealgebra}$, namely $\somemorph$. According to Lem\-ma~\ref{lemma:the-bfg-lemma}, that demands $\anycocycle\in \HScocycles{1}{\shortalgebra}{\themodule}$. Hence, the initial claim is true.

           \subproof{Step~3. Upper horizontal arrow has alleged inverse.} It suffices to prove that for any $\anycocycle\in \HScocycles{1}{\shortalgebra}{\themodule}$ and any $\thetuple\in \themodule^{\Sdirectproduct \thegens}$ with $\HSuafuncxx{1}{\thegens}{\therels}{\themodule}{\somerel}{\thetuple}=0$ for any $\somerel\in\therels$ the statements that $\thetuplex{\somegen}=\anycocycle(\somegen+\theideal)$ for any $\somegen\in\thegens$ and that $\anycocycle(\somepoly+\theideal)=\HSuafuncxx{1}{\thegens}{\therels}{\themodule}{\somepoly}{\thetuple}$ for any $\somepoly\in\thepolys$ are equivalent. Clearly, the second implies the first by the fact that  $\HSuafuncxx{1}{\thegens}{\therels}{\themodule}{\somegen}{\thetuple}=\thetuplex{\somegen}$ for any $\somegen\in\thegens$ by definition. And that the other implication  holds was shown in Step~1.

           \subproof{Step~4. Vertical arrows well defined.} That the left vertical arrow is well defined is clear. That the same is true for the right vertical arrow was shown in Lem\-ma~\ref{lemma:hochschild-particular-one-cohomology-well-defined}.

           \subproof{Step~5. Lower horizontal arrow and its inverse.} From the definition of the lower horizontal arrow and that of $\HSuacoboundaries{1}{\thegens}{\therels}{\themodule}$, it is clear that the lower horizontal arrow is well defined. Conversely, the inverse of the upper horizontal arrow restricts to the inverse of the lower horizontal arrow. That is because for any $\thetuple\in\themodule^{\Sdirectproduct \thegens}$ with \smash{$\HSuafuncxx{1}{\thegens}{\therels}{\themodule}{\somerel}{\thetuple}=0$} for any $\somerel\in\therels$, for the unique \smash{$\anycocycle\in\HScocycles{1}{\shortalgebra}{\themodule}$} with  \smash{$\anycocycle(\somepoly+\theideal)=\HSuafuncxx{1}{\thegens}{\therels}{\themodule}{\somepoly}{\thetuple}$} for any $\somepoly\in\thepolys$ and for any $\anyvector\in\themodule$, if  $\thetuplex{\somegen}=(\somegen+\theideal)\Sleftmoduleaction\anyvector-\anyvector\Srightmoduleaction(\somegen+\theideal)$ for any $\somegen\in\thegens$, then $\anycocycle(\somepoly+\theideal)=\HSuafuncxx{1}{\thegens}{\therels}{\themodule}{\somepoly}{\thetuple}=(\somepoly+\theideal)\Sleftmoduleaction\anyvector-\anyvector\Srightmoduleaction(\somepoly+\theideal)=(\partial^0\anyvector)(\somepoly+\theideal)$ by Lem\-ma~\ref{lemma:hochschild-particular-one-cohomology-well-defined}.

           \subproof{Step~6. Commutativity of the diagram.} Because the two horizontal arrows are defined by the same rule and since the vertical arrows are set inclusions the diagram commutes.

(b) Follows directly from \ref{proposition:hochschild-particular-one-cohomology-1} and is only stated for emphasis.
         \end{proof}
       }

       {
  \newcommand{\anyscalar}{\lambda}
  \newcommand{\anypoly}{p}
  \newcommand{\anycoc}{\eta}
  \newcommand{\anygen}{e}
  \newcommand{\thetuple}{x}
  \newcommand{\thetuplex}[1]{\thetuple_{#1}}
  \begin{Example}
    \label{example:trivial_bimodule_conclusion}
Let $\themodule$ be the trivial bimodule with respect to an augmentation $\thecounit$ of $\thealgebra$ as in  Example~\ref{example:trivial_bimodule}.
    Because then $\HScoboundaries{1}{\thealgebra}{\themodule}=\{0\}$ by Example~\ref{remark:trivial_bimodule} what Pro\-po\-si\-tion~\ref{proposition:hochschild-particular-one-cohomology} implies is that also $\HSuacoboundaries{1}{\thegens}{\therels}{\themodule}=\{0\}$ and that therefore $\HScohomology{1}{\thealgebra}{\themodule}\cong \HSuacocycles{1}{\thegens}{\therels}{\themodule}$.
         \end{Example}
       }
       \begin{Remark}
         \label{remark:restriction_of_scalars}
         \newcommand{\thevectorspace}{V}
         \newcommand{\actleft}{\blacktriangleright}
         \newcommand{\actright}{\blacktriangleleft}
         \newcommand{\anypoly}{p}
         \newcommand{\anyrel}{r}
         \newcommand{\anyvector}{v}
         \newcommand{\thetuple}{x}
         If $\thevectorspace$ is the $\thefield$-vector space underlying the $\thealgebra$-bimodule $\themodule$, then the rules $\anypoly\actleft\anyvector:= (\anypoly+\theideal)\Sleftmoduleaction\anyvector$ and $\anyvector\actright\anypoly:= \anyvector\Srightmoduleaction(\anypoly+\theideal)$ for any $\anypoly\in\thepolys$ and $\anyvector\in\thevectorspace$ define  left respectively right $\thepolys$-actions $\actleft$ and $\actright$ on $\thevectorspace$ which turn it into a $\thepolys$-bimodule $\theresmodule$,  the \emph{restriction of scalars} of $\themodule$ along the canonical projection $\fromto{\thepolys}{\thealgebra}$. With this definition there is no difference between the linear maps $\HSuafunc{1}{\thegens}{\therels}{\themodule}$ and $\HSuafunc{1}{\thegens}{\varnothing}{\theresmodule}$. (But, of course, there is in general still a~difference between $\HSuacocycles{1}{\thegens}{\varnothing}{\theresmodule}=\thevectorspace^ {\Sdirectproduct\thegens}$ and \smash{$\HSuacocycles{1}{\thegens}{\therels}{\themodule}=\{\thetuple\in \thevectorspace^ {\Sdirectproduct\thegens}\Sand \forall \anyrel\in\therels\quantorpredicate \HSuafuncxx{1}{\thegens}{\varnothing}{\theresmodule}{\anyrel}{\thetuple}=0\}$} and likewise between \smash{$\HSuacoboundaries{1}{\thegens}{\varnothing}{\theresmodule}$} and \smash{$\HSuacoboundaries{1}{\thegens}{\therels}{\themodule}$}.) The advantage of the notation \smash{$\HSuafunc{1}{\thegens}{\therels}{\themodule}$} is that one can work immediately with the given bimodule $\themodule$ and does not have to introduce $\theresmodule$ first. Then again, talking about $\HSuafunc{1}{\thegens}{\varnothing}{\theresmodule}$ can be advantageous too, e.g., in the instance of considering  simultaneously multiple different $\therels$  and thus multiple different $\themodule$ for which though the restrictions of scalars~$\theresmodule$ all happen to be the same.
       \end{Remark}
     }

\section{Certain spaces of scalar matrices and their dimensions}
\label{section:dimensions}

  {
  \newcommand{\themat}{v}
  The vector spaces of matrices appearing in the main result are characterized and their dimensions are computed.
  Recall that for any $\thedim\in\pint$ any $\themat\in\squarematrices{\thedim}{\comps}$ is called \emph{skew-symmetric} if $\themat=-\themat\Stra$.
  }

  \begin{Definition}
    \label{definition:small-matrix}
  \newcommand{\themat}{v}
  \newcommand{\thematx}[2]{v_{#1,#2}}
  \newcommand{\indi}{i}
  \newcommand{\indj}{j}
  We call any $\themat\in\squarematrices{\thedim}{\comps}$ \emph{small} if $\sum_{\indi=1}^\thedim\thematx{\indj}{\indi}=0$ for any $\indj\in\dwi{\thedim}$ and $\sum_{\indj=1}^\thedim\thematx{\indj}{\indi}=0$ for any $\indi\in\dwi{\thedim}$, i.e., if each row and each column sums to zero.
\end{Definition}

{
  \newcommand{\themat}{v}
  \newcommand{\thematx}[2]{v_{#1,#2}}
  \newcommand{\indi}{i}
  \newcommand{\indj}{j}
  \newcommand{\inds}{s}
  \newcommand{\thesum}{\lambda}
  \newcommand{\firstsum}{\lambda_1}
  \newcommand{\secondsum}{\lambda_2}
\begin{Lemma}
  \label{lemma:small-skewsymmetric-both}
  For any $\thedim\in\pint$ and $\themat\in\squarematrices{\thedim}{\comps}$ the following equivalences hold.
  \begin{enumerate}
  \item\label{lemma:small-skewsymmetric-both-1} There is $\thesum\in\comps$ such that $\themat-\thesum\identitymatrix$ is small if and only if $\sum_{\inds=1}^\thedim\thematx{\indj}{\inds}-\sum_{\inds=1}^\thedim\thematx{\inds}{\indi}=0$ for any $\{\indi,\indj\}\subseteq \dwi{\thedim}$. Moreover, then $\thesum=\sum_{\inds=1}^\thedim\thematx{\indj}{\inds}=\sum_{\inds=1}^\thedim\thematx{\inds}{\indi}$ for any $\{\indi,\indj\}\subseteq \dwi{\thedim}$.
  \item\label{lemma:small-skewsymmetric-both-2} There is $\thesum\in\comps$ such that $\themat-\thesum\identitymatrix$ is skew-symmetric  if and only if for any $\{\indi,\indj\}\subseteq \dwi{\thedim}$ with $\indi\neq \indj$ both  $\thematx{\indj}{\indi}+\thematx{\indi}{\indj}=0$ and $\thematx{\indj}{\indj}-\thematx{\indi}{\indi}=0$.
    Moreover, then $\thesum=\thematx{\indi}{\indi}$ for any $\indi\in\dwi{\thedim}$.
  \item\label{lemma:small-skewsymmetric-both-3} There are $\{\firstsum,\secondsum\}\subseteq\comps$ such that $\themat-\firstsum\identitymatrix$ is skew-symmetric and $\themat-\secondsum\identitymatrix$  small if and only if there is $\thesum\in\comps$ such that $\themat-\thesum\identitymatrix$ is both skew-symmetric and small. Moreover, then $\thesum=\firstsum=\secondsum$.
  \end{enumerate}
\end{Lemma}
\begin{proof}
  \newcommand{\auxmat}{w}
  \newcommand{\auxmatx}[2]{w_{#1,#2}}
  For $\thedim=1$, all claims  hold trivially. Hence, suppose $2\leq \thedim$ in the following.

  (a)
  If $\thesum\in\comps$ is such that $\auxmat:=\themat-\thesum\identitymatrix$ is  small, then for any $\{\indi,\indj\}\subseteq\dwi{\thedim}$ it follows $0=\sum_{\inds=1}^\thedim\auxmatx{\indj}{\inds}=\sum_{\inds=1}^\thedim(\thematx{\indj}{\inds}-\thesum\kron{\indj}{\inds})=\sum_{\inds=1}^\thedim\thematx{\indj}{\inds}-\thesum$ and $0=\sum_{\inds=1}^\thedim\auxmatx{\inds}{\indi}=\sum_{\inds=1}^\thedim(\thematx{\inds}{\indi}-\thesum\kron{\inds}{\indi})=\sum_{\inds=1}^\thedim\thematx{\inds}{\indi}-\thesum$, which proves $\sum_{\inds=1}^\thedim\thematx{\indj}{\inds}=\thesum=\sum_{\inds=1}^\thedim\thematx{\inds}{\indi}$. Of course, then $\sum_{\inds=1}^\thedim\thematx{\indj}{\inds}-\sum_{\inds=1}^\thedim\thematx{\inds}{\indi}=\thesum-\thesum=0$ for any $\{\indi,\indj\}\subseteq \dwi{\thedim}$.

    Conversely, if $ \sum_{\inds=1}^\thedim\thematx{\indj}{\inds}-\sum_{\inds=1}^\thedim\thematx{\inds}{\indi}=0$ for any $\{\indi,\indj\}\subseteq \dwi{\thedim}$ and if we let $\thesum:=\sum_{\inds=1}^\thedim\thematx{1}{\inds}$ and  $\auxmat:= \themat-\thesum\identitymatrix$, then for any $\{\indi,\indj\}\subseteq \dwi{\thedim}$,  first, $\thesum=\sum_{\inds=1}^\thedim\thematx{\indj}{\inds}=\sum_{\inds=1}^\thedim\thematx{\inds}{\indi}$ and thus, second, $\sum_{\inds=1}^\thedim\auxmatx{\indj}{\inds}=\sum_{\inds=1}^\thedim(\thematx{\indj}{\inds}-\thesum\kron{\indj}{\inds})=\sum_{\inds=1}^\thedim\thematx{\indj}{\inds}-\thesum=0$ and, likewise, $\sum_{\inds=1}^\thedim\auxmatx{\inds}{\indi}=\sum_{\inds=1}^\thedim(\thematx{\inds}{\indi}-\thesum\kron{\inds}{\indi})=\sum_{\inds=1}^\thedim\thematx{\inds}{\indi}-\thesum=0$. Hence, $\auxmat$ is small then.

    (b)
    If for $\thesum\in\comps$ the matrix $\auxmat:=\themat-\thesum\identitymatrix$ is skew-symmetric, then $0=\auxmatx{\indj}{\indi}+\auxmatx{\indi}{\indj}=(\thematx{\indj}{\indi}-\thesum\kron{\indj}{\indi})+(\thematx{\indi}{\indj}-\thesum\kron{\indi}{\indj})=\thematx{\indj}{\indi}+\thematx{\indi}{\indj}-2\thesum\kron{\indj}{\indi}$ for any $\{\indi,\indj\}\subseteq \dwi{\thedim}$. Consequently, if $\indi\neq \indj$, this means $0=\thematx{\indj}{\indi}+\thematx{\indi}{\indj}$ and, if $\indi=\indj$, we find $0=2\thematx{\indi}{\indi}-2\thesum$, i.e., $\thesum=\thematx{\indi}{\indi}$. And that implies in particular $\thematx{\indj}{\indj}-\thematx{\indi}{\indi}=\thesum-\thesum=0$ for any $\{\indi,\indj\}\subseteq \dwi{\thedim}$.

    If, conversely,  $\thematx{\indj}{\indi}+\thematx{\indi}{\indj}=0$ and $\thematx{\indj}{\indj}-\thematx{\indi}{\indi}=0$ for any $\{\indi,\indj\}\subseteq \dwi{\thedim}$ with $\indi\neq \indj$ and if we let $\thesum:= \thematx{1}{1}$ and $\auxmat:= \themat-\thesum\identitymatrix$, then, on the one hand, $\thesum=\thematx{\indi}{\indi}$ for any $\indi\in\dwi{\thedim}$ and, on the other hand,  for any $\{\indi,\indj\}\subseteq \dwi{\thedim}$, generally, $\auxmatx{\indj}{\indi}+\auxmatx{\indi}{\indj}=(\thematx{\indj}{\indi}-\thesum\kron{\indj}{\indi})+(\thematx{\indi}{\indj}-\thesum\kron{\indi}{\indj})=\thematx{\indj}{\indi}+\thematx{\indi}{\indj}-2\thesum\kron{\indj}{\indi}$, which in case $\indi\neq \indj$ simply means $\auxmatx{\indj}{\indi}+\auxmatx{\indi}{\indj}=\thematx{\indj}{\indi}+\thematx{\indi}{\indj}=0$ and which for $\indi=\indj$ amounts to $\auxmatx{\indj}{\indi}+\auxmatx{\indi}{\indj}=2\thematx{\indi}{\indi}-2\thesum=2\thesum-2\thesum=0$. In conclusion, $\auxmat$ is skew-symmetric then.

    (c)
    One implication is clear. If, conversely, $\{\firstsum,\secondsum\}\subseteq\comps$ are such that $\themat-\firstsum\identitymatrix$ is skew-symmetric and $\themat-\secondsum\identitymatrix$ is small, then $\firstsum=\thematx{1}{1}$ by \ref{lemma:small-skewsymmetric-both-2} and $\secondsum=\sum_{\indj=1}^\thedim\thematx{\indj}{1}=\sum_{\indi=1}^\thedim\thematx{1}{\indj}$ by \ref{lemma:small-skewsymmetric-both-1}. Subtracting the two identities $\sum_{\indj=1}^\thedim\thematx{\indj}{1}=\firstsum+\sum_{\indj=2}^\thedim\thematx{\indj}{1}$ and $\sum_{\indi=1}^\thedim\thematx{1}{\indi}=\firstsum+\sum_{\indi=2}^\thedim\thematx{1}{\indi}$ from each other therefore yields $0=\sum_{\indj=2}^\thedim\thematx{\indj}{1}-\sum_{\indi=2}^\thedim\thematx{1}{\indi}$. Since also  $\thematx{\indi}{1}=-\thematx{1}{\indi}$ for each $\indi\in\dwi{\thedim}$ with $1<\indi$ by \ref{lemma:small-skewsymmetric-both-2}, that is the same as saying $0=2\sum_{\indj=2}^\thedim\thematx{\indj}{1}$. And $\sum_{\indj=2}^\thedim\thematx{\indj}{1}=0$ then implies $\secondsum=\firstsum+\sum_{\indj=2}^\thedim\thematx{\indj}{1}=\firstsum$, which is all we needed to see.
\end{proof}
}

{
  \newcommand{\somespace}{V}
  \newcommand{\themat}{v}
  \newcommand{\theshortpredicate}{A}
  \newcommand{\theshortpredicatex}[1]{\theshortpredicate(#1)}
    \newcommand{\thesum}{\lambda}
    \newcommand{\theid}{\identitymatrix}
    \begin{Lemma}
      \label{lemma:dimensions}
For any $\thedim\in\pint$ and each statement $\theshortpredicate$  below, the set $\{\themat\in\squarematrices{\thedim}{\comps}\Sand \theshortpredicatex{\themat}\}$ is a~complex vector subspace of $\squarematrices{\thedim}{\comps}$ and has the listed dimension.
    \begin{table}[h]\centering\renewcommand{\arraystretch}{1.2}
      \begin{tabular}{ll|l}
&        $\theshortpredicatex{\themat}$ & $\dim_\comps\{\themat\in\squarematrices{\thedim}{\comps}\Sand \theshortpredicatex{\themat}\}$\\[0.1em] \hline
$(a)$  &      $\Strue$& $\thedim^2$ \Tstrut\\
$(b)$ & $\exists \thesum\in\comps\quantorpredicate \themat-\thesum\theid$ is small                                  &$(\thedim-1)^2+1$ \\
$(c)$ & $\themat$ is small                                  &$(\thedim-1)^2$ \\
$(d)$ & $\exists \thesum\in\comps\quantorpredicate \themat-\thesum\theid$ is skew-symmetric                                  &$\onehalf\thedim(\thedim-1)+1$ \\
$(e)$ &              $\themat$ is skew-symmetric                    &$\onehalf\thedim(\thedim-1)$ \\
$(f)$ & $\exists \thesum\in\comps\quantorpredicate \themat-\thesum\theid$ is skew-symmetric and small                                  &$\onehalf(\thedim-1)(\thedim-2)+1$ \\
$(g)$ & $\themat$     is skew-symmetric and small                             &$\onehalf(\thedim-1)(\thedim-2)$ \\
$(h)$ &        $\themat$  is diagonal& $\thedim$ \\
$(i)$ & $\exists \thesum\in\comps\quantorpredicate \themat-\thesum\theid=0$                                 & $1$ \\
$(j)$ & $\themat=0$        & $0$
    \end{tabular}
    \end{table}
  \end{Lemma}

  \begin{proof}
    \newcommand{\matb}[3]{E^{#1}_{#2,#3}}
    \newcommand{\indi}{i}
    \newcommand{\indj}{j}
    \newcommand{\indk}{k}
    \newcommand{\indl}{\ell}
    \newcommand{\inds}{s}
    \newcommand{\indt}{t}
    \newcommand{\smalliso}[1]{\varphi_{#1}}
    \newcommand{\smallisoinv}[1]{\psi_{#1}}
    \newcommand{\auxmat}{w}
    \newcommand{\auxmatx}[2]{\auxmat_{#1,#2}}
    \newcommand{\smallmat}{u}
    \newcommand{\smallmatx}[2]{\smallmat_{#1,#2}}
    \newcommand{\thematx}[2]{\themat_{#1,#2}}
    \newcommand{\othermat}{v'}
    \newcommand{\othermatx}[2]{\othermat_{#1,#2}}
    \newcommand{\ssmatb}[3]{T^{#1}_{#2,#3}}
    \newcommand{\indset}[1]{\Gamma_{#1}}
    \newcommand{\someind}{\gamma}
    \newcommand{\basismat}[2]{B^{#1}_{#2}}
    \newcommand{\scals}[1]{a_{#1}}
(a)
It is well known that, if for any $\{\indk,\indl\}\subseteq \dwi{n}$ the matrix $\matb{\thedim}{\indl}{\indk}\in \squarematrices{\thedim}{\comps}$ has $\kron{\indl}{\indj}\kron{\indk}{\indi}$ as its $(\indj,\indi)$-entry for any $\{\indi,\indj\}\subseteq \dwi{\thedim}$, then the family $(\matb{\thedim}{\indl}{\indk})_{(\indl,\indk)\in\dwi{\thedim}^{\Ssetmonoidalproduct 2}}$ is a $\comps$-linear basis of $\{\themat\in\squarematrices{\thedim}{\comps}\Sand \theshortpredicatex{\themat}\}=\squarematrices{\thedim}{\comps}$.

(b)
Since $\theshortpredicate$ can be expressed by a homogenous system of linear equations by Lem\-ma~\hyperref[lemma:small-skewsymmetric-both-1]{\ref*{lemma:small-skewsymmetric-both}\,\ref*{lemma:small-skewsymmetric-both-1}}  the set $\{\themat\in\squarematrices{\thedim}{\comps}\Sand \theshortpredicatex{\themat}\}$ is indeed a vector space.  Hence, it suffices to show that  the mapping $\xfromto{\smalliso{\thedim}}{\squarematrices{\thedim-1}{\comps}\Sdirectsum\comps}{ \{\themat\in\squarematrices{\thedim}{\comps}\Sand \theshortpredicatex{\themat}\}}$ defined by the rule that $        (\smallmat,\thesum)\mapsto \themat$, where  for any $\{\indi,\indj\}\subseteq \dwi{\thedim}$,
      \[
        \thematx{\indj}{\indi}=
        \begin{cases}
          \smallmatx{\indj}{\indi}+\thesum\kron{\indj}{\indi},& \indj<\thedim\Sand \indi<\thedim,\\
        \displaystyle  -\sum_{\indl=1}^{\thedim-1}\smallmatx{\indl}{\indi}, & \indj=\thedim\Sand \indi<\thedim,\\[1mm]
       \displaystyle   -\sum_{\indk=1}^{\thedim-1}\smallmatx{\indj}{\indk}, & \indj<\thedim\Sand \indi=\thedim,\\[1mm]
       \displaystyle   \sum_{\indk,\indl=1}^{\thedim-1}\smallmatx{\indl}{\indk}+\thesum ,& \indj=\thedim\Sand \indi=\thedim,
        \end{cases}
      \]
      for any $\smallmat\in\squarematrices{\thedim-1}{\comps}$ and $\thesum\in\comps$, is a $\comps$-linear isomorphism. We begin by proving that $\smalliso{\thedim}$ is well defined. For any $\smallmat\in\squarematrices{\thedim-1}{\comps}$ and  $\thesum\in\comps$, if $\smalliso{\thedim}(\smallmat,\thesum)=\themat$ and if  $\auxmat=\themat-\thesum\identitymatrix$, then, on the one hand, for any $\indi\in \dwi{\thedim}$ with $\indi<\thedim$, by definition,
      \[       \sum_{\indj=1}^{\thedim}\auxmatx{\indj}{\indi}=\sum_{\indj=1}^{\thedim-1}(\thematx{\indj}{\indi}-\thesum\kron{\indj}{\indi})+\thematx{\thedim}{\indi}=\sum_{\indj=1}^{\thedim-1}\smallmatx{\indj}{\indi}+\left(-\sum_{\indl=1}^{\thedim-1}\smallmatx{\indl}{\indi}\right)=0.
      \]
      and also,
            \[      \sum_{\indj=1}^{\thedim}\auxmatx{\indj}{\thedim}=\sum_{\indj=1}^{\thedim-1}(\thematx{\indj}{\thedim}-\thesum\kron{\indj}{\thedim})+(\thematx{\thedim}{\thedim}-\thesum)=\sum_{\indj=1}^{\thedim-1}\left(-\sum_{\indk=1}^{\thedim-1}\smallmatx{\indj}{\indk}\right)+\sum_{\indk,\indl=1}^{\thedim-1}\smallmatx{\indl}{\indk}=0.
\]
On the other hand, for any $\indj\in \dwi{\thedim}$ with $\indj<\thedim$,
      \[      \sum_{\indi=1}^{\thedim}\auxmatx{\indj}{\indi}=\sum_{\indi=1}^{\thedim-1}(\thematx{\indj}{\indi}-\thesum\kron{\indj}{\indi})+\thematx{\indj}{\thedim}=\sum_{\indi=1}^{\thedim-1}\smallmatx{\indj}{\indi}+\left(-\sum_{\indk=1}^{\thedim-1}\smallmatx{\indj}{\indk}\right)=0
      \]
      and also,
            \[  \sum_{\indi=1}^{\thedim}\auxmatx{\thedim}{\indi}=\sum_{\indi=1}^{\thedim-1}(\thematx{\thedim}{\indi}-\thesum\kron{\thedim}{\indi})+(\thematx{\thedim}{\thedim}-\thesum)=\sum_{\indi=1}^{\thedim-1}\left(-\sum_{\indl=1}^{\thedim-1}\smallmatx{\indl}{\indi}\right)+\sum_{\indk,\indl=1}^{\thedim-1}\smallmatx{\indl}{\indk}=0.
\]
Together these four conclusions prove that $\auxmat$ is small, i.e., that $\theshortpredicatex{\themat}$ holds.

Conversely, by Lem\-ma~\hyperref[lemma:small-skewsymmetric-both-1]{\ref*{lemma:small-skewsymmetric-both}\,\ref*{lemma:small-skewsymmetric-both-1}} a well-defined $\comps$-linear map $\xfromto{\smallisoinv{\thedim}}{\{\themat\in\squarematrices{\thedim}{\comps}\Sand\theshortpredicatex{\themat}\}}{\squarematrices{\thedim-1}{\comps}\Sdirectsum\comps}$ is obtained as follows: For any $\themat\in\squarematrices{\thedim}{\comps}$ with $\theshortpredicatex{\themat}$, if $\thesum\in\comps$ is such that $\themat-\thesum\identitymatrix$ is small, then $\themat\mapsto (\smallmat,\thesum)$, where for any $\{\indk,\indl\}\subseteq \dwi{\thedim-1}$,
\[
  \smallmatx{\indl}{\indk}=\thematx{\indl}{\indk}-\thesum\kron{\indl}{\indk}.
\]

It remains to show $\smallisoinv{\thedim}\Scomposition\smalliso{\thedim}=\Sidentity$ and $\smalliso{\thedim}\Scomposition\smallisoinv{\thedim}=\Sidentity$. And, indeed, for any $\smallmat\in \squarematrices{\thedim-1}{\comps}$ and $\thesum\in\comps$, if $\themat=\smalliso{\thedim}(\smallmat,\thesum)$, then we have already seen that $\auxmat=\themat-\thesum\identitymatrix$ is small. For any $\{\indk,\indl\}\subseteq \dwi{\thedim}$, by definition, $\auxmatx{\indl}{\indk}=\thematx{\indl}{\indk}-\thesum\kron{\indl}{\indk}=(\smallmatx{\indl}{\indk}+\thesum\kron{\indl}{\indk})-\thesum\kron{\indl}{\indk}=\smallmatx{\indl}{\indk}$, which proves $\smalliso{\thedim}(\themat)=(\smallmat,\thesum)$ and thus $\smallisoinv{\thedim}\Scomposition\smalliso{\thedim}=\Sidentity$.

Conversely, for any $\themat\in\squarematrices{\thedim}{\comps}$ such that $\theshortpredicatex{\themat}$ is satisfied, if $(\smallmat,\thesum)=\smallisoinv{\thedim}(\themat)$, then we already know $\thesum=\sum_{\indl=1}^\thedim\thematx{\indl}{\indi}=\sum_{\indk=1}^\thedim\thematx{\indj}{\indk}$ for any $\{\indk,\indl\}\subseteq\dwi{\thedim}$ by Lem\-ma~\hyperref[lemma:small-skewsymmetric-both-1]{\ref*{lemma:small-skewsymmetric-both}\,\ref*{lemma:small-skewsymmetric-both-1}}. If $\othermat=\smalliso{\thedim}(\smallmat,\thesum)$, then   for any $\{\indi,\indj\}\subseteq \dwi{\thedim}$ with $\indi<\thedim$ and $\indj<\thedim$ it hence follows by definition $\othermatx{\indj}{\indi}=\smallmatx{\indj}{\indi}+\thesum\kron{\indj}{\indi}=(\thematx{\indj}{\indi}-\thesum\kron{\indj}{\indi})+\thesum\kron{\indj}{\indi}=\thematx{\indj}{\indi}$ as well as by $\thesum=\sum_{\indl=1}^\thedim\thematx{\indl}{\indi}$,
\[
\othermatx{\thedim}{\indi}= -\sum_{\indl=1}^{\thedim-1}\smallmatx{\indl}{\indi}=-\sum_{\indl=1}^{\thedim-1}(\thematx{\indl}{\indi}-\thesum\kron{\indl}{\indi})=\thesum-\sum_{\indl=1}^{\thedim-1}\thematx{\indl}{\indi}=\thematx{\thedim}{\indi}
\]
and by $\thesum=\sum_{\indk=1}^\thedim\thematx{\indk}{\indj}$,
\[
\othermatx{\indj}{\thedim}= -\sum_{\indk=1}^{\thedim-1}\smallmatx{\indj}{\indk}=-\sum_{\indk=1}^{\thedim-1}(\thematx{\indj}{\indk}-\thesum\kron{\indj}{\indk})=\thesum-\sum_{\indk=1}^{\thedim-1}\thematx{\indj}{\indk}=\thematx{\indj}{\thedim}
\]
and, lastly,
\begin{align*}
 \othermatx{\thedim}{\thedim}&{}=\sum_{\indk,\indl=1}^{\thedim-1}\smallmatx{\indl}{\indk}+\thesum=\sum_{\indk,\indl=1}^{\thedim-1}(\thematx{\indl}{\indk}-\thesum\kron{\indl}{\indk})+\thesum=\sum_{\indl=1}^{\thedim-1}\left(\sum_{\indk=1}^{\thedim-1}\thematx{\indl}{\indk}-\thesum\right)+\thesum\\
  &{}=\sum_{\indl=1}^{\thedim-1}(-\thematx{\indl}{\thedim})+\thesum=\thematx{\thedim}{\thedim},
\end{align*}
where we have used $\thesum=\sum_{\indk=1}^\thedim\thematx{\indl}{\indk}$ for any $\indl\in\dwi{\thedim}$ in the next-to-last step and $\thesum=\sum_{\indl=1}^\thedim\thematx{\indl}{\thedim}$ in the last.
Thus, we have shown $\othermat=\themat$ and thus $\smalliso{\thedim}\Scomposition\smallisoinv{\thedim}=\Sidentity$, which concludes the proof in this case.

(c)
By Lem\-ma~\hyperref[lemma:small-skewsymmetric-both-1]{\ref*{lemma:small-skewsymmetric-both}\,\ref*{lemma:small-skewsymmetric-both-1}},
  the space $\{\themat\in\squarematrices{\thedim}{\comps}\Sand\theshortpredicatex{\themat}\}$ is exactly the image of $\squarematrices{\thedim-1}{\comps}\Sdirectsum\{0\}$ under $\smalliso{\thedim}$.

(d) Lem\-ma~\hyperref[lemma:small-skewsymmetric-both-2]{\ref*{lemma:small-skewsymmetric-both}\,\ref*{lemma:small-skewsymmetric-both-2}} showed that $\theshortpredicate$ can be equivalently expressed as a system of homogenous linear equations, thus proving $\{\themat\in\squarematrices{\thedim}{\comps}\Sand \theshortpredicatex{\themat}\}$ to be a vector space. Let $\indset{\thedim}=\{(\indj,\indi)\classpredicate \{\indi,\indj\}\subseteq\dwi{\thedim}\Sand \indj<\indi\}\cupdot\{\varnothing\}$ as well as $\basismat{\thedim}{(\indj,\indi)}=\ssmatb{\thedim}{\indj}{\indi}=\matb{\thedim}{\indj}{\indi}-\matb{\thedim}{\indi}{\indj}$ for any $\{\indi,\indj\}\subseteq \dwi{\thedim}$ with $\indj<\indi$ and $\basismat{\thedim}{\varnothing}=\identitymatrix$. Then, the claim will be verified once we show that $(\basismat{\thedim}{\someind})_{\someind\in\indset{\thedim}}$ is a $\comps$-linear basis of $\{\themat\in\squarematrices{\thedim}{\comps}\Sand \theshortpredicatex{\themat}\}$.

The family  $(\basismat{\thedim}{\someind})_{\someind\in\indset{\thedim}}$ is  $\comps$-linearly independent. Indeed, if $(\scals{\someind})_{\someind\in\indset{\thedim}}\in \comps^{\Sdirectsum \indset{\thedim}}$ is such that $\sum_{\someind\in\indset{\thedim}}\scals{\someind}\,\basismat{\thedim}{\someind}=0$, then by $\identitymatrix=\sum_{\indi=1}^\thedim\matb{\thedim}{\indi}{\indi}$,
\[
  0=\sum_{\substack{(\indj,\indi)\in\dwi{\thedim}^{\Ssetmonoidalproduct 2}\\{}\Sand\indj<\indi}}\scals{(\indj,\indi)}(\matb{\thedim}{\indj}{\indi}-\matb{\thedim}{\indi}{\indj})+\scals{\varnothing} \sum_{\indi=1}^\thedim\matb{\thedim}{\indi}{\indi}=\sum_{(\indj,\indi)\in\dwi{\thedim}^{\Ssetmonoidalproduct 2}}
  \begin{dcases}
    \begin{rcases}
     \scals{(\indj,\indi)}, &  \indj <\indi\\
     -\scals{(\indi,\indj)}, &  \indi <\indj\\
     \scals{\varnothing}, &  \indj =\indi
    \end{rcases}
  \end{dcases} \matb{\thedim}{\indj}{\indi},
\]
which demands $(\scals{\someind})_{\someind\in\indset{\thedim}}=0$ since $(\matb{\thedim}{\indl}{\indk})_{(\indl,\indk)\in\dwi{\thedim}^{\Ssetmonoidalproduct 2}}$ is $\comps$-linearly independent.

It remains to prove that $\{\basismat{\thedim}{\someind}\classpredicate \someind\in\indset{\thedim}\}$ spans $\{\themat\in\squarematrices{\thedim}{\comps}\Sand \theshortpredicatex{\themat}\}$. If $\themat\in\squarematrices{\thedim}{\comps}$ and $\thesum\in\comps$ are such that $\auxmat=\themat-\thesum\identitymatrix$ is skew-symmetric, then  $\thematx{\indj}{\indi}=-\thematx{\indi}{\indj}$  and $\thesum=\thematx{\indj}{\indj}=\thematx{\indi}{\indi}$ for any $\{\indi,\indj\}\subseteq\dwi{\thedim}$ with $\indj\neq \indi$ by Lem\-ma~\hyperref[lemma:small-skewsymmetric-both-2]{\ref*{lemma:small-skewsymmetric-both}\,\ref*{lemma:small-skewsymmetric-both-2}}. Hence, if we let $\scals{\varnothing}=\thesum$ and  $\scals{(\indj,\indi)}=\auxmatx{\indj}{\indi}=\thematx{\indj}{\indi}$ for any $\{\indi,\indj\}\subseteq \dwi{\thedim}$ with $\indj<\indi$, then
\[
\sum_{\someind\in\indset{\thedim}}\scals{\someind}\,\basismat{\thedim}{\someind}=\sum_{(\indj,\indi)\in\dwi{\thedim}^{\Ssetmonoidalproduct 2}}
  \begin{dcases}
    \begin{rcases}
     \scals{(\indj,\indi)}, & \indj <\indi\\
     -\scals{(\indi,\indj)}, &  \indi <\indj\\
     \scals{\varnothing}, &  \indj =\indi
    \end{rcases}
  \end{dcases} \matb{\thedim}{\indj}{\indi}=\sum_{(\indj,\indi)\in\dwi{\thedim}^{\Ssetmonoidalproduct 2}}
  \begin{dcases}
    \begin{rcases}
     \thematx{\indj}{\indi}, &  \indj <\indi\\
     -\thematx{\indi}{\indj}, & \indi <\indj\\
     \thesum, & \indj =\indi
    \end{rcases}
  \end{dcases} \matb{\thedim}{\indj}{\indi}=\themat.
\]
Thus, $(\basismat{\thedim}{\someind})_{\someind\in\indset{\thedim}}$ is a $\comps$-linear basis.

(e)
The proof of the previous claim shows that any $\themat\in\squarematrices{\thedim}{\comps}$ is skew-symmetric if and only if it is in the span of  $\{\basismat{\thedim}{\someind}\classpredicate \someind\in\indset{\thedim}\}$  and has coefficient $0$ with respect to $\basismat{\thedim}{\varnothing}$. Hence, $\{\ssmatb{\thedim}{\indj}{\indi}\classpredicate\{\indi,\indj\}\subseteq \dwi{\thedim}\Sand\indj<\indi\}$ is a $\comps$-linear basis of  $\{\themat\in\squarematrices{\thedim}{\comps}\Sand \theshortpredicatex{\themat}\}$.

(f)
All three parts \ref*{lemma:small-skewsymmetric-both-1}--\ref*{lemma:small-skewsymmetric-both-3} of  Lem\-ma~\ref{lemma:small-skewsymmetric-both} combined imply that $\{\themat\in\squarematrices{\thedim}{\comps}\Sand \theshortpredicatex{\themat}\}$ is the solution set to a homogenous system of linear equations and thus a vector space. Hence,  it suffices to prove that $\smalliso{\thedim}$ restricts to a mapping $\fromto{\{\smallmat\in\squarematrices{\thedim-1}{\comps}\Sand \smallmat=-\smallmat\Stra\}}{\{\themat\in\squarematrices{\thedim}{\comps}\Sand \theshortpredicatex{\themat}\}}$ and $\smallisoinv{\thedim}$ to one in the reverse direction.

  For any skew-symmetric $\smallmat\in\squarematrices{\thedim-1}{\comps}$ and any $\thesum\in\comps$, if $\themat=\smalliso{\thedim}(\smallmat,\thesum)$ and $\auxmat=\themat-\thesum\identitymatrix$, then for any $\{\indi,\indj\}\subseteq \dwi{\thedim}$ with $\indi<\thedim$ and $\indj<\thedim$ we have already seen that  $\auxmatx{\indj}{\indi}=\smallmatx{\indj}{\indi}$, implying $\auxmatx{\indj}{\indi}+\auxmatx{\indi}{\indj}=\smallmatx{\indj}{\indi}+\smallmatx{\indi}{\indj}=0$ by $\smallmat=-\smallmat\Stra$. Moreover, for the same reason,
  \begin{align*}
      \auxmatx{\thedim}{\indi}+\auxmatx{\indi}{\thedim}&{}=    (\thematx{\thedim}{\indi}-\thesum\kron{\thedim}{\indi})+(\thematx{\indi}{\thedim}-\thesum\kron{\indi}{\thedim})=\thematx{\thedim}{\indi}+\thematx{\indi}{\thedim}\\
    &{}=    \left(-\sum_{\indl=1}^{\thedim-1}\smallmatx{\indl}{\indi}\right)+\left(-\sum_{\indk=1}^{\thedim-1}\smallmatx{\indi}{\indk}\right)=-\sum_{\indl=1}^{\thedim-1}(\smallmatx{\indl}{\indi}+\smallmatx{\indi}{\indl})=0
 \end{align*}
 and
   \begin{align*}
       \auxmatx{\indj}{\thedim}+\auxmatx{\thedim}{\indj}&{}=    (\thematx{\indj}{\thedim}-\thesum\kron{\indj}{\thedim})+(\thematx{\thedim}{\indj}-\thesum\kron{\thedim}{\indj})=\thematx{\indj}{\thedim}+\thematx{\thedim}{\indj}\\
     &{}=    \left(-\sum_{\indk=1}^{\thedim-1}\smallmatx{\indj}{\indk}\right)+\left(-\sum_{\indl=1}^{\thedim-1}\smallmatx{\indl}{\indj}\right)=-\sum_{\indk=1}^{\thedim-1}(\smallmatx{\indj}{\indk}+\smallmatx{\indk}{\indj})=0
 \end{align*}
 as well as
   \[
       \auxmatx{\thedim}{\thedim}+\auxmatx{\thedim}{\thedim}= 2(\thematx{\thedim}{\thedim}-\thesum)=2\sum_{\indk,\indl=1}^{\thedim-1}\smallmatx{\indl}{\indk}=\sum_{\indk,\indl=1}^{\thedim-1}(\smallmatx{\indl}{\indk}+\smallmatx{\indk}{\indl})=0,
   \]
   which completes the proof that $\smalliso{\thedim}$ restricts to a map into $\{\themat\in\squarematrices{\thedim}{\comps}\Sand \theshortpredicatex{\themat}\}$.

   Conversely, if $\themat\in\squarematrices{\thedim}{\comps}$ and $\thesum\in\comps$ are such that  $\auxmat=\themat-\thesum\identitymatrix$ is skew-symmetric and small, then   $\thesum=\thematx{\indi}{\indi}$ for any $\indi\in\dwi{\thedim}$ by Lem\-ma~\hyperref[lemma:small-skewsymmetric-both-2]{\ref*{lemma:small-skewsymmetric-both}\,\ref*{lemma:small-skewsymmetric-both-2}}. For $(\smallmat,\thesum)=\smallisoinv{\thedim}(\themat)$  and any $\{\indk,\indl\}\subseteq \dwi{\thedim-1}$, by definition, $\smallmatx{\indl}{\indk}=\auxmatx{\indl}{\indk}$ and thus  $\smallmatx{\indl}{\indk}+\smallmatx{\indk}{\indl}=\auxmatx{\indl}{\indk}+\auxmatx{\indk}{\indl}=0$ by $\auxmat=-\auxmat\Stra$. Hence,  $\smallisoinv{\thedim}$ maps $\{\themat\in\squarematrices{\thedim}{\comps}\Sand \theshortpredicatex{\themat}\}$ into $\{\smallmat\in\squarematrices{\thedim-1}{\comps}\Sand \smallmat=-\smallmat\Stra\}\Sdirectsum\comps$.

   (g)
   As we have just shown, any $\themat\in\squarematrices{\comps}{\thedim}$ is skew-symmetric and small if and only if it lies in the image of $\{\smallmat\in\squarematrices{\thedim-1}{\comps}\Sand \smallmat=-\smallmat\Stra\}\Sdirectsum\{0\}$ under $\smalliso{\thedim}$. And  $\smalliso{\thedim}$ is a $\comps$-linear isomorphism from this space to  $\{\themat\in\squarematrices{\thedim}{\comps}\Sand \theshortpredicatex{\themat}\}$.

   (h)
   It is well known that $(\matb{\thedim}{\indi}{\indi})_{\indi\in\dwi{\thedim}}$ is a $\comps$-linear basis of $\{\themat\in\squarematrices{\thedim}{\comps}\Sand \theshortpredicatex{\themat}\}$.

   (i)
   In this case, $\{\themat\in\squarematrices{\thedim}{\comps}\Sand \theshortpredicatex{\themat}\}$ is the $\comps$-linear span of $\identitymatrix$ in $\squarematrices{\thedim}{\comps}$.

   (j)
   Here,  $\{\themat\in\squarematrices{\thedim}{\comps}\Sand \theshortpredicatex{\themat}\}$ is the zero  $\comps$-linear space.
  \end{proof}
  }

\section{First cohomology of unitary easy quantum group duals}
\label{section:main-one}
This section computes the first quantum group  cohomology with trivial coefficients (see Section~\ref{section:quantum-groups}) of the discrete dual of any unitary easy compact quantum group (see Section~\ref{section:easy}). That is achieved by applying the characterization of the first Hochschild cohomology recalled in Section~\ref{section:hochschild} while using the results of Section~\ref{section:dimensions} as auxiliaries.
\subsection{Equations derived from the presentation}
\label{section:main-one-equations}
 Resume the As\-sump\-tions~\ref{assumptions:easy} and the abbreviations from Notations~\ref{notation:easy} and \ref{notation:relations_of_partition_set}.  In particular,~$\thedim$~and $\thegens$ are then defined.
 Remark~\ref{remark:restriction_of_scalars} mo\-ti\-va\-tes moreover the fol\-lo\-wing short\-hand.

\begin{Notation}\quad
  \newcommand{\ascalar}{x}
  \newcommand{\acol}{\Yc{c}}
  \newcommand{\indi}{i}
  \newcommand{\indj}{j}
  \newcommand{\anypoly}{p}
  \newcommand{\thepolys}{\freealg{\comps}{\thegens}}
  \newcommand{\thepartition}{p}
  \newcommand{\thepartitionset}{\mathcal{S}}
  \newcommand{\incol}{\Yc{c}}
  \newcommand{\outcol}{\Yc{d}}
  \newcommand{\infree}{g}
  \newcommand{\outfree}{j}
  \begin{enumerate}
  \item Let
    $\theresmodule$ be the $\comps\langle\thegens\rangle$-bimodule  $\comps$ with left and right actions given by $\theunix{\acol}{\indj}{\indi}\Sleftmoduleaction \ascalar:= \kron{\indj}{\indi}\ascalar$ respectively $ \ascalar\Srightmoduleaction\theunix{\acol}{\indj}{\indi}:= \kron{\indj}{\indi}\ascalar$ for any $\{\indi,\indj\}\subseteq \dwi{\thedim}$, any $\acol\in\blaw$ and any $\ascalar\in\comps$.
\item Let $\thefuncx{\anypoly}:= \HSuafuncx{1}{\thegens}{\varnothing}{\theresmodule}{\anypoly}$ for any $\anypoly\in\thepolys$.
  \end{enumerate}
\end{Notation}
{
  \newcommand{\anypoly}{r}
  \newcommand{\thepolys}{\freealg{\comps}{\thegens}}
  Then, by Section~\ref{section:hochschild} for any category $\thecat$ of two-colored partitions the first cohomology with trivial coefficients of the discrete dual of any easy quantum group associated with $(\thecat,\thedim)$ can be realized as a solution space to a system of linear equations involving maps of the form $\thefuncx{\anypoly}$ for certain $\anypoly\in\thepolys$ induced by $\thecat$ and $\thedim$.
}

{
  \newcommand{\anygen}{e}
  \newcommand{\thepolys}{\freealg{\comps}{\thegens}}
  \newcommand{\somerel}{r}
  \newcommand{\indj}{j}
  \newcommand{\indi}{i}
  \newcommand{\incol}{\Yc{c}}
  \newcommand{\outcol}{\Yc{d}}
  \newcommand{\infree}{g}
  \newcommand{\outfree}{j}
  \newcommand{\thepartition}{p}
  \newcommand{\thetuple}{x}
  \newcommand{\thetuplex}[1]{x_{#1}}
  \newcommand{\somecocycle}{\eta}
  \begin{Proposition}
    \label{proposition:main_system_of_equations}
    For any category $\thecat$ of two-co\-lo\-red partitions,
    if $\theqg$ is the unitary easy compact quantum group of $(\thecat,\thedim)$,
then there exists an isomorphism of $\comps$-vector spaces
\[
  \begin{tikzcd}
  \DQGcohomologyTC{1}{\CQGdualC{\theqg}} \arrow[r, two heads, hook]& \{ \thetuple\in\comps^ {\Sdirectproduct\thegens}\Sand \forall  \somerel\in\thepartrels{\thecat}\quantorpredicate \thefuncxx{\somerel}{\thetuple}=0  \},
  \end{tikzcd}
\]
which maps $($the one-elemental cohomology class of$)$ any $1$-cycle $\somecocycle$ to the tuple $\thetuple$ with $\thetuplex{\anygen}=\somecocycle(\anygen+\thepartideal{\thecat})$ for any $\anygen\in\thegens$.
\end{Proposition}
\begin{proof}
  By Definition~\ref{definition:easy-quantum-group}, the algebra underlying the Hopf $\ast$-algebra $\DQGgroupalgebra{\CQGdualC{\theqg}}$ is the universal algebra $\univalg{\comps}{\thegens}{\thepartrels{\thecat}}$. According to Section~\ref{section:quantum-groups-cohomology}, the vector space $\DQGcohomologyTC{1}{\CQGdualC{\theqg}}$ is defined as $\HScohomology{1}{\univalg{\comps}{\thegens}{\thepartrels{\thecat}}}{\themodule}$ where  $\themodule={}_\thecounit\comps_\thecounit$ is trivial bimodule of $\univalg{\comps}{\thegens}{\thepartrels{\thecat}}$ with respect to the counit $\thecounit$ of $\DQGgroupalgebra{\CQGdualC{\theqg}}$. By~Pro\-po\-si\-tion~\ref{proposition:counit_of_easy_quantum_group}, this counit is such that its restriction of scalars along the canonical projection $\fromto{\thepolys}{\univalg{\comps}{\thegens}{\thepartrels{\thecat}}}$ is precisely $\theresmodule$. Hence, the claim follows by
  Example~\ref{example:trivial_bimodule_conclusion} and Remark~\ref{remark:restriction_of_scalars}.
\end{proof}
}
The task laid out by Pro\-po\-si\-tion~\ref{proposition:main_system_of_equations} is clear. We need to solve the set of linear equations in~$\comps^{\Sdirectproduct\thegens}$ on the right-hand side of the isomorphism there --  for each category of two-co\-lo\-red partitions. Eventually, in Section~\ref{section:main-one-case-distinctions}  namely, solving these equations will require case distinctions for different kinds of categories of two-co\-lo\-red partitions. However, there are a great number of simplifications we can make to the equation system before it needs to come to that. Moreover, this reduces the number of cases we eventually have to consider immensely.

\subsection{Simplifying each individual equation}
\label{section:main-one-simplifications}
{
  \newcommand{\anyrel}{r}
  \newcommand{\thepartition}{p}
  \newcommand{\incol}{\Yc{c}}
  \newcommand{\outcol}{{\Yc{d}}}
  \newcommand{\infree}{g}
  \newcommand{\outfree}{j}
As a first step towards solving the equations of Pro\-po\-si\-tion~\ref{proposition:main_system_of_equations} we consider each equation in isolation  and simplify its  definition. In other words,  we seek a better formula for the values of the functional $\thefuncx{\anyrel}$ for  $\anyrel$ of the form $\therelpoly{\incol}{\outcol}{\thepartition}{\outfree}{\infree}$ for arbitrary two-colored partitions $(\incol,\outcol,\thepartition)$ and~$\infree\in\dwi{\thedim}^{\Ssetmonoidalproduct\Slength{\incol}}$ and $\outfree\in\dwi{\thedim}^{\Ssetmonoidalproduct\Slength{\outcol}}$.
}

It will be convenient to have a shorthand for mappings constructed by prescribing a specified value to a specified point and otherwise inheriting the graph of a given mapping with the same domain.
\begin{Notation}
  \newcommand{\inlen}{k}
  \newcommand{\outlen}{\ell}
  \newcommand{\theind}{f}
  \newcommand{\thepoint}{\Yp{z}}
  \newcommand{\thevalue}{s}
  \newcommand{\thepartition}{p}
  \newcommand{\thealtered}{\theind\downarrow_\thepoint\thevalue}
  \newcommand{\anypoint}{\Yp{y}}
  For any $\{\inlen,\outlen\}\subseteq \nnint$, any mapping $\theind\funcdef\tsopx{\inlen}{\outlen}\to \dwi{\thedim}$, any $\thepoint\in\tsopx{\inlen}{\outlen}$ and any $\thevalue\in\dwi{\thedim}$ write $\thealtered$ for the mapping $\tsopx{\inlen}{\outlen}\to \dwi{\thedim}$ with  $\thepoint\mapsto \thevalue$ and with $\anypoint\mapsto \theind(\anypoint)$ for any $\anypoint\in\tsopx{\inlen}{\outlen}\backslash \{\thepoint\}$.
\end{Notation}
Then, combining Notation~\ref{notation:easy} and Definition~\ref{definition:derivatives-one} yields the following description of the functionals we are investigating.

{
  \newcommand{\inlen}{k}
  \newcommand{\outlen}{\ell}
  \newcommand{\incol}{\mathfrak{c}}
  \newcommand{\outcol}{\mathfrak{d}}
  \newcommand{\thepartition}{p}
  \newcommand{\theind}{f}
  \newcommand{\arel}{r}
  \newcommand{\thecol}{\Yc{w}}
  \newcommand{\ourpoint}{\Yp{z}}
  \newcommand{\ourvalue}{s}
  \newcommand{\infree}{g}
  \newcommand{\outfree}{j}
  \newcommand{\theentry}[3]{u_{#1,#2}^{#3}}
  \newcommand{\thecocyc}{x}
  \newcommand{\thecocycle}[3]{x_{\theentry{#1}{#2}{#3}}}
  \begin{Lemma}
    \label{lemma:simplifying_the_definitions_of_the_functionals}
    For any $\{\inlen,\outlen\}\subseteq \nnint$, any $\incol\in\blaw^{\Ssetmonoidalproduct \inlen}$, any $\outcol\in\blaw^{\Ssetmonoidalproduct \outlen}$, any set-theoretical partition $\thepartition$ of $\tsopx{\inlen}{\outlen}$, any $\infree\in\dwi{\thedim}^{\Ssetmonoidalproduct\inlen}$, any $\outfree\in\dwi{\thedim}^{\Ssetmonoidalproduct\outlen}$ and any $\thecocyc\in\comps^{\Sdirectproduct\thegens}$,  if $\arel=\therelpoly{\incol}{\outcol}{\thepartition}{\outfree}{\infree}$ and  $\theind=\djp{\infree}{\outfree}$ and $\thecol= \djp{\incol}{\outcol}$, then
    \[
\thefuncxx{\arel}{\thecocyc}=\sum_{\ourpoint\in\tsopx{\inlen}{\outlen}}\sum_{\ourvalue=1}^\thedim \zetfx{\thepartition}{\ker(\alterer{\theind}{\ourpoint}{\ourvalue})}
\begin{rcases}\begin{dcases}
-\thecocycle{\ourvalue}{\theind(\ourpoint)}{\thecol(\ourpoint)}  & \tif \ourpoint\in\tsopx{\inlen}{0}\\
\thecocycle{\theind(\ourpoint)}{\ourvalue}{\thecol(\ourpoint)}    & \tif \ourpoint\in\tsopx{0}{\outlen}
\end{dcases}\end{rcases}.
\]
  \end{Lemma}
  \begin{proof}
    \newcommand{\inbound}{h}
    \newcommand{\inboundx}[1]{\inbound_{#1}}
    \newcommand{\infreex}[1]{\infree_{#1}}
    \newcommand{\incolx}[1]{\incol_{#1}}
    \newcommand{\outbound}{i}
    \newcommand{\outboundx}[1]{\outbound_{#1}}
    \newcommand{\outfreex}[1]{\outfree_{#1}}
    \newcommand{\outcolx}[1]{\outcol_{#1}}
    \newcommand{\pind}{a}
    \newcommand{\qind}{b}
    \newcommand{\subind}{q}
For any $\thecocyc\in\comps^{\Sdirectproduct\thegens}$, by Example~\ref{example:trivial_bimodule_functional},
    \begin{align*}
\thefuncxx{\arel}{\thecocyc}&{}=       \sum_{\outbound\in\dwi{\thedim}^{\Ssetmonoidalproduct \outlen}}
      \zetfx{\thepartition}{\ker(\djp{\infree}{\outbound})}
\sum_{\qind=1}^{\outlen}\biggl(\prod_{\substack{\subind\in\dwi{\outlen}\\{}\Sand \subind\neq \qind}}\kron{\outfreex{\qind}}{\outboundx{\qind}}\biggr)\; \thecocycle{\outfreex{\qind}}{\outboundx{\qind}}{\outcolx{\qind}}
      \\&\quad{}-\sum_{\inbound\in\dwi{\thedim}^{\Ssetmonoidalproduct \inlen}}  \zetfx{\thepartition}{\ker(\djp{\inbound}{\outfree})}
 \sum_{\pind=1}^{\inlen}\biggl(\prod_{\substack{\subind\in\dwi{\inlen}\\{}\Sand \subind\neq \pind}}\kron{\inboundx{\pind}}{\infreex{\pind}}\biggr)\;\thecocycle{\inboundx{\pind}}{\infreex{\pind}}{\incolx{\pind}}.
\end{align*}
After commuting the sums and evaluating the sums over $\outbound$ respectively $\inbound$ (as far as possible), this is identical to
    \begin{gather*}
       \sum_{\qind=1}^{\outlen}\sum_{\outboundx{\qind}=1}^\thedim
       \zetfx{\thepartition}{\ker(\djp{\infree}{(\outfreex{1},\ldots,\outfreex{\qind-1},\outboundx{\qind},\outfreex{\qind+1},\ldots,\outfreex{\outlen})})}\;
 \thecocycle{\outfreex{\qind}}{\outboundx{\qind}}{\outcolx{\qind}}\\
\qquad- \sum_{\pind=1}^{\inlen}\sum_{\inboundx{\pind}=1}^\thedim  \zetfx{\thepartition}{\ker(\djp{(\infreex{1},\ldots,\infreex{\pind-1},\inboundx{\pind},\infreex{\pind+1},\ldots,\infreex{\inlen})}{\outfree})}\;
\thecocycle{\inboundx{\pind}}{\infreex{\pind}}{\incolx{\pind}}.
\end{gather*}
That agrees with the right-hand side of the claimed identity.
\end{proof}

While Lem\-ma~\ref{lemma:simplifying_the_definitions_of_the_functionals} has given a more concise form to the equations under investigation, it can be improved upon significantly.
Firstly, one can give a simpler criterion for when in the sum on the right-hand side of the identity in Lem\-ma~\ref{lemma:simplifying_the_definitions_of_the_functionals} a factor $\zetfx{\thepartition}{\ker(\alterer{\theind}{\ourpoint}{\ourvalue})}$ is non-zero.
}

{
  \newcommand{\inlen}{k}
  \newcommand{\outlen}{\ell}
  \newcommand{\theind}{f}
  \newcommand{\thepoint}{\Yp{z}}
  \newcommand{\thevalue}{s}
  \newcommand{\thepartition}{p}
  \newcommand{\thealtered}{\theind\downarrow_\thepoint\thevalue}
  \newcommand{\anypoint}{\Yp{a}}
  \newcommand{\theblock}{\blo{\thepartition}{\thepoint}}
  \newcommand{\cutblock}{\blo{\thepartition}{\thepoint}\backslash \{\thepoint\}}
  \newcommand{\sing}{\{\thepoint\}}
  \begin{Lemma}
    \label{lemma:easy-cohomology-helper-one-point-replacer}
    For any $\{\inlen,\outlen\}\!\subseteq\! \nnint$, any set-theo\-re\-ti\-cal partition $\thepartition$ of $\tsopx{\inlen}{\outlen}$, any mapping ${\theind\funcdef\!\tsopx{\inlen}{\outlen}\!\to\! \dwi{\thedim}}$, any $\thepoint\in\tsopx{\inlen}{\outlen}$ and any $\thevalue\in\dwi{\thedim}$,
    the statements $ \thepartition\finerthan\ker(\thealtered)$ and
    \[
\thepartition\backslash \{\theblock\}\cupdot \{\cutblock,\sing\}\backslash \{\varnothing\}\finerthan \ker(\theind)\quad\Sand\quad \cutblock\subseteq \spimgx{\theind}{\{\thevalue\}}
\]
are equivalent.
  \end{Lemma}
  \begin{proof}
    \newcommand{\somevalue}{t}
    \newcommand{\othervalue}{t'}
    \newcommand{\someblock}{\Yb{B}}
    We show each implication separately. Below, we will use many times the fact that for any~${\somevalue \in \dwi{\thedim}}$,
    \begin{align*}
      \spimgx{(\thealtered)}{\{\somevalue\}}\backslash \{\thepoint\}&{}=\{\anypoint\in\tsopx{\inlen}{\outlen}\Sand (\thealtered)(\anypoint)=\somevalue\Sand \anypoint\neq\thepoint \}\\
      &{}=\{\anypoint\in\tsopx{\inlen}{\outlen}\Sand \theind(\anypoint)=\somevalue\Sand \anypoint\neq\thepoint \}\\
    &{}=\spimgx{\theind}{\{\somevalue\}}\backslash \{\thepoint\}.
    \end{align*}

    \subproof{Step~1.}
    First, suppose $\thepartition\finerthan\ker(\thealtered)$. Then, there exists $\somevalue\in\sranx{\thealtered}$ such that $\theblock\subseteq \spimgx{(\thealtered)}{\{\somevalue\}}$. Because $\thepoint\in\theblock$ this requires $\thepoint\in \spimgx{(\thealtered)}{\{\somevalue\}}$ and thus $\somevalue=\thevalue$ by ${(\thealtered)(\thepoint)=\thevalue}$. It follows $\theblock\subseteq \spimgx{(\thealtered)}{\{\thevalue\}}$ and thus in particular $\cutblock\subseteq \spimgx{(\thealtered)}{\{\thevalue\}}\backslash \sing=\spimgx{\theind}{\{\thevalue\}}\backslash \{\thepoint\}\subseteq \spimgx{\theind}{\{\thevalue\}}$, which is one half of what we had to show.

    It is trivially true that $\sing\subseteq \spimgx{\theind}{\{\theind(\thepoint)\}}\in \ker(\theind)$. We have already seen that $\cutblock\subseteq\spimgx{\theind}{\{\thevalue\}}\in\ker(\theind)$. For any $\someblock\in\thepartition$ with $\someblock\neq\theblock$, i.e., $\thepoint\notin \someblock$, there exists by assumption $\othervalue\in \sranx{\thealtered}$ with $\someblock\subseteq \spimgx{(\thealtered)}{\{\othervalue\}}$. We conclude $\someblock=\someblock\backslash\sing\subseteq \spimgx{(\thealtered)}{\{\othervalue\}}\backslash\sing=\spimgx{\theind}{\{\othervalue\}}\backslash \sing\subseteq \spimgx{\theind}{\{\othervalue\}}\in \ker(\theind)$. Thus, the other half of the claim, $\thepartition\backslash \{\theblock\}\cupdot \{\cutblock,\allowbreak\sing\}\backslash \{\varnothing\}\finerthan \ker(\theind)$, holds as well. That proves one implication.

    \subproof{Step~2.} In order to show the converse implication we assume that both $\thepartition\backslash \{\theblock\}\cupdot \{\cutblock,\allowbreak\sing\}\backslash \{\varnothing\}\finerthan \ker(\theind)$ and $ \cutblock\subseteq \spimgx{\theind}{\{\thevalue\}}$ and then we distinguish two cases.

    \subproof{Case~2.1.} If $\sing\in \thepartition$ and thus $\theblock=\sing$ and  $\cutblock=\varnothing$, then the assumption is simply equivalent to the statement $\thepartition\finerthan\ker(\theind)$. Naturally, $\sing\subseteq \spimgx{(\thealtered)}{\{\thevalue\}}\in\ker(\theind)$ by $(\thealtered)(\thepoint)=\thevalue$. For any $\someblock\in\thepartition$ with $\someblock\neq\sing$ there exists by our premise a value $\somevalue\in \sranx{\theind}$ with $\someblock\subseteq \spimgx{\theind}{\{\somevalue\}}$. Thus, also $\someblock=\someblock\backslash \sing\subseteq \spimgx{\theind}{\{\somevalue\}}\backslash\sing=\spimgx{(\thealtered)}{\{\somevalue\}}\backslash\sing\in\ker(\thealtered)$. In conclusion, $\thepartition\finerthan\ker(\thealtered)$.

    \subproof{Case~2.2.} In the instance that $\sing\notin\thepartition$   the initial assumption  simplifies to the statement $\thepartition\backslash \{\theblock\}\cupdot\{\cutblock,\sing\}\finerthan\ker(\theind)$ and $\cutblock\subseteq \spimgx{\theind}{\{\thevalue\}}$. The latter condition implies $\cutblock\subseteq \spimgx{\theind}{\{\thevalue\}}\backslash \sing=\spimgx{(\thealtered)}{\{\thevalue\}}\backslash \sing\subseteq \spimgx{(\thealtered)}{\{\thevalue\}}$ and thus by $(\thealtered)(\thepoint)=\thevalue$ also $\theblock=\cutblock\cupdot\sing\subseteq \spimgx{(\thealtered)}{\{\thevalue\}}\cup \sing\subseteq \spimgx{(\thealtered)}{\{\thevalue\}}\in\ker(\thealtered)$. On the other hand, for any $\someblock\in\thepartition$ with $\someblock\neq \theblock$, which is to say $\thepoint\notin\someblock$, there exists by assumption $\somevalue\in \sranx{\theind}$ with $\someblock\subseteq \spimgx{\theind}{\{\somevalue\}}$. It follows $\someblock=\someblock\backslash\sing\subseteq \spimgx{\theind}{\{\somevalue\}}\backslash\sing=\spimgx{(\thealtered)}{\{\somevalue\}}\backslash\sing\subseteq\spimgx{(\thealtered)}{\{\somevalue\}}\in\ker(\thealtered)$. Hence, altogether,  $\thepartition\finerthan\ker(\thealtered)$, which concludes the proof.
  \end{proof}
}

{
  \newcommand{\thepartition}{p}
  \newcommand{\theind}{f}
  \newcommand{\ourpoint}{\Yp{z}}
  \newcommand{\ourvalue}{s}
  Lem\-ma~\ref{lemma:easy-cohomology-helper-one-point-replacer} can now be used to give a necessary criterion for the right-hand side of the identity in  Lem\-ma~\ref{lemma:simplifying_the_definitions_of_the_functionals} to be non-zero as a whole. Namely, $\thepartition$ and $\theind$ must meet one of three conditions:
\begin{enumerate}\itemsep=0pt
\item[(i)] The labeling $\theind$ maps any points belonging to the same element of $\thepartition$ to the same value.
\item[(ii)]
  There is an element of size two of $\thepartition$ whose elements $\theind$ maps to different values. Besides that $\theind$ is as in~(i).
\item[(iii)]
There is an element  of $\thepartition$ of size three or larger, all but one of whose elements are assigned the same value by $\theind$ and whose remaining element  $\theind$  sends to a different value. Apart from that, $\theind$~is as in~(i).
\end{enumerate}
}

  \begin{Definition}
  \newcommand{\inlen}{k}
  \newcommand{\outlen}{\ell}
  \newcommand{\theind}{f}
  \newcommand{\theindx}[1]{f(#1)}
  \newcommand{\thepoint}{\Yp{z}}
  \newcommand{\thevalue}{s}
  \newcommand{\thepartition}{p}
  \newcommand{\thelong}{\Yb{Z}}
  \newcommand{\firstspecialpoint}{\Yp{z}_1}
  \newcommand{\secondspecialpoint}{\Yp{z}_2}
  \newcommand{\anypoint}{\Yp{y}}
  \newcommand{\someblock}{\Yb{B}}
  \newcommand{\anyblock}{\Yb{A}}
    Let $\{\inlen,\outlen\}\subseteq \nnint$, let $\thepartition$ be any set-theoretical partition of $\tsopx{\inlen}{\outlen}$ and let $\xfromto{\theind}{\tsopx{\inlen}{\outlen}}{\dwi{\thedim}}$. Then, we say that $(\thepartition,\theind)$ is
    \begin{enumerate}
    \item \emph{case~\scenarioone{}} if $\thepartition\neq \varnothing$ and $\thepartition\finerthan\kerp{\theind}$,
    \item \emph{case~\scenariotwo{}} if there exists $\{\firstspecialpoint,\secondspecialpoint\}\in\thepartition$  such that $\theindx{\firstspecialpoint}\neq\theindx{\secondspecialpoint}$ and such that for any $\anyblock\in\thepartition$ with $\anyblock\neq \{\firstspecialpoint,\secondspecialpoint\}$ there is $\someblock\in\kerp{\theind}$ with $\anyblock\subseteq \someblock$, in which case the set $\{\firstspecialpoint,\secondspecialpoint\}$ is called \emph{critical data} of $(\thepartition,\theind)$,
    \item \emph{case~\scenariothree{}} if there exist $\thelong\in\thepartition$ and $\thepoint\in\thelong$ and $\thevalue\in\dwi{\thedim}$ such that $3\leq |\thelong|$, such that $\theindx{\thepoint}\neq \thevalue$, such that $\theindx{\anypoint}= \thevalue$ for any $\anypoint\in\thelong$ with $\anypoint\neq \thepoint$ and such that for any $\anyblock\in\thepartition$ with $\anyblock\neq \thelong$ there is $\someblock\in\kerp{\theind}$ with $\anyblock\subseteq \someblock$, in which case  $(\thelong, \thepoint,\thevalue)$ are called \emph{critical data} of $(\thepartition,\theind)$,
    \item \emph{case~\scenariofour{}} otherwise.
    \end{enumerate}
  \end{Definition}
  \begin{Example}
    \label{example:scenarios}
      \newcommand{\inlen}{k}
  \newcommand{\outlen}{\ell}
  \newcommand{\theind}{f}
  \newcommand{\theindx}[1]{f(#1)}
  \newcommand{\thepoint}{\Yp{z}}
  \newcommand{\thevalue}{s}
  \newcommand{\thepartition}{p}
  \newcommand{\thelong}{\Yb{Z}}
  \newcommand{\firstspecialpoint}{\Yp{z}_1}
  \newcommand{\secondspecialpoint}{\Yp{z}_2}
  \newcommand{\anypoint}{\Yp{y}}
  \newcommand{\someblock}{\Yb{B}}
  \newcommand{\anyblock}{\Yb{A}}
  For $3\leq \thedim$, consider $\inlen:= 4$ and $\outlen:= 5$, the set-theoretical partition  $\thepartition:= \{ \{\lop{1},\upp{1}\}, \allowbreak\{\upp{2}\}, \{\lop{2},\lop{3},\upp{3}\}, \{\upp{4}\}, \{\lop{4},\lop{5}\} \}$ of $\tsopx{\inlen}{\outlen}$   and various different mappings $\xfromto{\theind}{\tsopx{\inlen}{\outlen}}{\dwi{\thedim}}$ which all have in common that
 each of $\lop{2}$, $\lop{4}$ and $\upp{3}$ is mapped to $1$, that each of $\lop{1}$, $\upp{1}$ and $\upp{2}$ is mapped to $2$ and that $\upp{4}$ is mapped to $3$. Thus, at most the values of $\lop{3}$ and $\lop{5}$ differ between different $\theind$:
    $$
    \begin{tikzpicture}[baseline=0.91cm]
    \def\scp{0.666}
    \def\linksize{\scp*0.075cm}
    \def\pointsize{\scp*0.25cm}
    \def\dd{\scp*0.5cm}
    \def\dx{\scp*1cm}
    \def\cx{\scp*0.3cm}
    \def\txu{3*\dx}
    \def\txl{4*\dx}
    \def\dy{\scp*1cm}
    \def\cy{\scp*0.3cm}
    \def\ty{3*\dy}
    \tikzset{whp/.style={circle, inner sep=0pt, text width={\pointsize}, draw=black, fill=white}}
    \tikzset{blp/.style={circle, inner sep=0pt, text width={\pointsize}, draw=black, fill=black}}
    \tikzset{lk/.style={regular polygon, regular polygon sides=4, inner sep=0pt, text width={\linksize}, draw=black, fill=black}}
    \draw[dotted] ({0-\dd},{0}) -- ({\txl+\dd},{0});
    \draw[dotted] ({0-\dd},{\ty}) -- ({\txu+\dd},{\ty});
    \coordinate (l1) at ({0+0*\dx},{0+0*\ty}) {};
    \coordinate (l2) at ({0+1*\dx},{0+0*\ty}) {};
    \coordinate (l3) at ({0+2*\dx},{0+0*\ty}) {};
    \coordinate (l4) at ({0+3*\dx},{0+0*\ty}) {};
    \coordinate (l5) at ({0+4*\dx},{0+0*\ty}) {};
    \coordinate (u1) at ({0+0*\dx},{0+1*\ty}) {};
    \coordinate (u2) at ({0+1*\dx},{0+1*\ty}) {};
    \coordinate (u3) at ({0+2*\dx},{0+1*\ty}) {};
    \coordinate (u4) at ({0+3*\dx},{0+1*\ty}) {};
    \node[lk] at ({2*\dx},{1*\dy}) {};
    \draw (l1) -- (u1);
    \draw (l3) -- (u3);
    \draw (l2) |- ++ ({1*\dx},{1*\dy});
    \draw[->] (u2) -- ++({0*\dx},{-1*\dy});
    \draw[->] (u4) -- ++({0*\dx},{-1*\dy});
    \draw (l4) -- ++ ({0*\dx},{1*\dy}) -| (l5);
    \node[above =2pt of u1] {$\scriptstyle 2$};
    \node[above =2pt of u2] {$\scriptstyle 2$};
    \node[above =2pt of u3] {$\scriptstyle 1$};
    \node[above =2pt of u4] {$\scriptstyle 3$};
    \node[below =2pt of l1] {$\scriptstyle 2$};
    \node[below =2pt of l2] {$\scriptstyle 1$};
    \node[below =2pt of l4] {$\scriptstyle 1$};
    \node at ({-1*\dx},{0*\dy-8pt}) {$\theind$};
    \node at ({-1*\dx},{1.5*\dy}) {$\thepartition$};
    \node at ({-1*\dx},{3*\dy+8pt}) {$\theind$};
  \end{tikzpicture}
$$
  \begin{enumerate}
  \item\label{example:scenarios-1} If $\theind(\lop{3})=\theind(\lop{5})=1$, then $(\thepartition,\theind)$ is case~\scenarioone{}.
  \item\label{example:scenarios-2} If $\theind(\lop{3})=1$ and $\theind(\lop{5})=2$, if $\firstspecialpoint:= \lop{4}$ and $\secondspecialpoint:= \lop{5}$, then $(\thepartition,\theind)$ is case~\scenariotwo{} with critical data $\{\firstspecialpoint,\secondspecialpoint\}=\{\lop{4},\lop{5}\}$.
\item\label{example:scenarios-3} If $\theind(\lop{3})=2$ and $\theind(\lop{5})=1$, if $\thelong:= \{\lop{2},\lop{3},\upp{3}\}$, if $\thepoint:= \lop{3}$ and $\thevalue:= 2$,  then $(\thepartition,\theind)$ is case~\scenariothree{} with critical data $(\thelong,\thepoint,\thevalue)=(\{\lop{2},\lop{3},\upp{3}\},\{\lop{3}\},2)$.
  \item\label{example:scenarios-4} If $\theindx{\lop{3}}=\theindx{\lop{5}}=2$, then $(\thepartition,\theind)$ is case~\scenariofour{}.
  \end{enumerate}
  \end{Example}
  \begin{Lemma}\quad
    \label{lemma:scenarios_well-defined}
    \begin{enumerate}
    \item\label{lemma:scenarios_well-defined-1} In each of the cases \scenariotwo{} and \scenariothree{} critical data are unique.
      \item\label{lemma:scenarios_well-defined-2} The cases \scenarioone{}--\scenariofour{} are mutually exclusive and exhaustive.
    \end{enumerate}
\end{Lemma}
\begin{proof}
  \newcommand{\inlen}{k}
  \newcommand{\outlen}{\ell}
  \newcommand{\theind}{f}
  \newcommand{\theindx}[1]{f(#1)}
  \newcommand{\thepoint}{\Yp{z}}
  \newcommand{\thevalue}{s}
  \newcommand{\thepartition}{p}
  \newcommand{\thelong}{\Yb{Z}}
  \newcommand{\firstspecialpoint}{\Yp{z}_1}
  \newcommand{\secondspecialpoint}{\Yp{z}_2}
  \newcommand{\anypoint}{\Yp{y}}
  \newcommand{\someblock}{\Yb{B}}
  \newcommand{\anyblock}{\Yb{A}}
  \newcommand{\alternativefirstspecialpoint}{\Yp{z}_1'}
  \newcommand{\alternativesecondspecialpoint}{\Yp{z}_2'}
  \newcommand{\alternativelong}{\Yb{Z}'}
  \newcommand{\alternativepoint}{\Yp{z}'}
  \newcommand{\alternativevalue}{s'}
  \newcommand{\alternativeanypoint}{\Yp{y}'}
  Let $\{\inlen,\outlen\}\subseteq \nnint$, let $\thepartition$ be any set-theoretical partition of $\tsopx{\inlen}{\outlen}$ and let $\xfromto{\theind}{\tsopx{\inlen}{\outlen}}{\dwi{\thedim}}$ be arbitrary.

  (a)
 \subproof{Case~\scenariotwo{}.} Suppose that $(\thepartition,\theind)$ is case~\scenariotwo{} and that both $\{\firstspecialpoint,\secondspecialpoint\}$ and  $\{\alternativefirstspecialpoint,\alternativesecondspecialpoint\}$ are critical data of $(\thepartition,\theind)$. If $\{\firstspecialpoint,\secondspecialpoint\}\neq \{\alternativefirstspecialpoint,\alternativesecondspecialpoint\}$ were true, then by the assumption on $\{\firstspecialpoint,\secondspecialpoint\}$ there would exist $\someblock\in\ker(\theind)$ with $\{\alternativefirstspecialpoint,\alternativesecondspecialpoint\}\subseteq \someblock$, meaning $\theind(\alternativefirstspecialpoint)=\theind(\alternativesecondspecialpoint)$, contrarily to our assumption. Hence, $\{\firstspecialpoint,\secondspecialpoint\}= \{\alternativefirstspecialpoint,\alternativesecondspecialpoint\}$ must be true instead.

 \subproof{Case~\scenariothree{}.} Now, let $(\thepartition,\theind)$ be case~\scenariothree{} and let both $(\thelong,\thepoint,\thevalue)$ and $(\alternativelong,\alternativepoint,\alternativevalue)$ be critical data of~$(\thepartition,\theind)$.
 If $\thelong\neq\alternativelong$ held, the assumption on $\thelong$ would imply the existence of $\someblock\in\kerp{\theind}$ with $\alternativelong\subseteq \someblock$. In~particular, it would follow $\theind(\alternativeanypoint)=\theind(\alternativepoint)$ for any $\alternativeanypoint\in \alternativelong$ with $\alternativeanypoint\neq\alternativepoint$, of which there exists at least one by $3\leq |\alternativelong|$. Because that would contradict the assumption, we must have $\thelong=\alternativelong$ instead.

 Furthermore,  supposing $\thepoint\neq \alternativepoint$ demands of any $\anypoint\in\thelong\backslash \{\thepoint,\alternativepoint\}$ both $\theind(\anypoint)=\thevalue$ by the assumption on $\thepoint$ and $\thevalue$ and $\theind(\anypoint)=\alternativevalue$ by the one on $\alternativepoint$ and $\alternativevalue$. Hence, as $\thelong\backslash \{\thepoint,\alternativepoint\}\neq \varnothing$  by $3\leq |\alternativelong|$, if  $\alternativepoint\neq \thepoint$, then $\thevalue=\alternativevalue$. That would be a contradiction because the property of $\alternativepoint$ also requires $\thevalue\neq \theind(\thepoint)=\alternativevalue$ in that case. Hence, only $\thepoint=\alternativepoint$ can be true.

 Lastly, because the assumptions on $\thevalue$ and $\alternativevalue$ imply $\theind(\anypoint)=\thevalue$ respectively $\theind(\anypoint)=\alternativevalue$ for any $\anypoint\in\thelong$ with $\anypoint\neq \thepoint=\alternativepoint$ and because  $\thelong\backslash \{\thepoint\}\neq \varnothing$, we must have $\thevalue=\alternativevalue$ as well.

 (b)
 It is enough to prove that cases~\scenarioone{}--\scenariothree{} are mutually exclusive.
 If $(\thepartition,\theind)$ is case~\scenariotwo, then it cannot be case~\scenarioone{} because $\theind(\firstspecialpoint)\neq \theind(\secondspecialpoint)$ excludes the existence of $\someblock\in \ker(\theind)$ with $\{\firstspecialpoint,\secondspecialpoint\}\subseteq \someblock$, which would be necessary for  $\thepartition\finerthan\ker(\theind)$ to hold.
Similarly, $(\thepartition,\theind)$ being case~\scenariothree{} forbids it being case~\scenarioone{} as well because the existence of $\anypoint\in \thelong\backslash \{\thepoint\}\neq\varnothing$ with $\theind(\thepoint)\neq \thevalue=\theind(\anypoint)$ does not allow  any $\someblock\in\ker(\theind)$ with $\thelong\subseteq \someblock$ to exist, which $\thepartition\finerthan\ker(\theind)$ would require.
 Lastly, if~$(\thepartition,\theind)$ were  simultaneously case~\scenariotwo{} and case~\scenariothree{}, then $\{\firstspecialpoint,\secondspecialpoint\}\neq\thelong$ would follow from $3\leq |\thelong|$, thus demanding by the property of $\thelong$ the existence of $\someblock\in\ker(\theind)$ with $\{\firstspecialpoint,\secondspecialpoint\}\subseteq \someblock$, in  contradiction to $\theind(\firstspecialpoint)\neq\theind(\secondspecialpoint)$.
\end{proof}

{
  \newcommand{\inlen}{k}
  \newcommand{\outlen}{\ell}
  \newcommand{\theind}{f}
  \newcommand{\theindx}[1]{f(#1)}
  \newcommand{\thepoint}{\Yp{z}}
  \newcommand{\thevalue}{s}
  \newcommand{\thepartition}{p}
  \newcommand{\thealtered}{\theind\downarrow_\thepoint\thevalue}
  \newcommand{\theblock}{\blo{\thepartition}{\thepoint}}
  \newcommand{\cutblock}{\blo{\thepartition}{\thepoint}\backslash \{\thepoint\}}
  \newcommand{\sing}{\{\thepoint\}}
  \newcommand{\thelong}{\Yb{Z}}
  \newcommand{\firstspecialpoint}{\Yp{z}_1}
  \newcommand{\secondspecialpoint}{\Yp{z}_2}
  \newcommand{\thegood}{f}
  \newcommand{\theonevalue}{s_1}
  \newcommand{\theothervalue}{s_2}
  \newcommand{\anypoint}{\Yp{y}}
  \newcommand{\someblock}{\Yb{A}}
  \newcommand{\anyblock}{\Yb{B}}
  \newcommand{\alternativefirstspecialpoint}{\Yp{z}_1'}
  \newcommand{\alternativesecondspecialpoint}{\Yp{z}_2'}
  \newcommand{\alternativelong}{\Yb{Z}'}
  \newcommand{\alternativepoint}{\Yp{z}'}
  \newcommand{\alternativevalue}{s'}
  \newcommand{\alternativeanypoint}{\Yp{y}'}
  \begin{Lemma}
    \label{lemma:easy-cohomology-helper-non-trivial-polynomials}
    For any ${\{\inlen,\outlen\}\subseteq \nnint}$, any set-theoretical partition $\thepartition$ of $\tsopx{\inlen}{\outlen}$ and any mapping $\theind\funcdef\tsopx{\inlen}{\outlen}\to \dwi{\thedim}$ there exist $\thepoint\in\tsopx{\inlen}{\outlen}$ and $\thevalue\in\dwi{\thedim}$ such that $\thepartition\finerthan\ker(\thealtered)$  if and only if $(\thepartition,\theind)$ is not case~\scenariofour.
\end{Lemma}
\begin{proof}
  Each implication is shown individually.

  \subproof{Step~1.} First, we suppose that $(\thepartition,\theind)$ is not case~\scenariofour{}
  and deduce the existence of $\thepoint\in\tsopx{\inlen}{\outlen}$ and $\thevalue\in\dwi{\thedim}$ with $\thepartition\finerthan\ker(\thealtered)$. By Lem\-ma~\ref{lemma:easy-cohomology-helper-one-point-replacer}, that is the same as finding $\thepoint\in\tsopx{\inlen}{\outlen}$ and $\thevalue\in\dwi{\thedim}$ such that $\thepartition\backslash \{\theblock\}\cupdot \{\cutblock,\sing\}\backslash \{\varnothing\}\finerthan \ker(\theind)$ and $\cutblock\subseteq \spimgx{\theind}{\{\thevalue\}}$. By Lem\-ma~\hyperref[lemma:scenarios_well-defined-2]{\ref*{lemma:scenarios_well-defined}\,\ref*{lemma:scenarios_well-defined-2}}, the pair $(\thepartition,\theind)$ is case~\scenarioone{}, case~\scenariotwo{} or case~\scenariothree{}. These three cases are treated individually.

  \subproof{Case~1.1.}
If $(\thepartition,\theind)$ is case~\scenarioone{}, then by $\thepartition\neq \varnothing$ we can find and fix some $\thepoint\in\tsopx{\inlen}{\outlen}$ and put $\thevalue:=\theind(\thepoint)$. From $\thepartition\finerthan\ker(\theind)$, it then follows $\blo{\thepartition}{\thepoint}\subseteq \spimgx{\theind}{\{\thevalue\}}$ and thus in particular $\blo{\thepartition}{\thepoint}\backslash \{\thepoint\}\subseteq \spimgx{\theind}{\{\thevalue\}}$, which is one part of what we have to show. The other part, $\thepartition\backslash \{\theblock\}\cupdot \{\cutblock,\sing\}\backslash \{\varnothing\}\finerthan \ker(\theind)$ is a consequence of the fact $\thepartition\backslash \{\theblock\}\cupdot \{\cutblock,\sing\}\backslash \{\varnothing\}\finerthan \thepartition$ and the assumption $\thepartition\finerthan\ker(\theind)$.

 \subproof{Case~1.2.} Next, let $(\thepartition,\theind)$ be case~\scenariotwo{} and let  $\{\firstspecialpoint,\secondspecialpoint\}$ be its critical data.
 If we define $\thepoint:=\firstspecialpoint$ and $\thevalue:= \theind(\secondspecialpoint)$, then $\blo{\thepartition}{\thepoint}=\{\firstspecialpoint,\secondspecialpoint\}$ and thus $\blo{\thepartition}{\thepoint}\backslash \{\thepoint\}=\{\secondspecialpoint\}\subseteq \spimgx{\theind}{\{\thevalue\}}$.  On the other hand, $\thepartition\backslash \{\theblock\}\cupdot \{\cutblock,\sing\}\backslash \{\varnothing\}=\thepartition\backslash \{\{\firstspecialpoint,\secondspecialpoint\}\}\cupdot \{\{\firstspecialpoint\},\{\secondspecialpoint\}\}\finerthan\ker(\theind)$ because, by assumption, for each  $\someblock\in\thepartition\backslash \{\{\firstspecialpoint,\secondspecialpoint\}\}$ there exists $\anyblock\in \ker(\theind)$ with $\someblock\subseteq \anyblock\in\ker(\theind)$ and, of course, $\{\firstspecialpoint\}\subseteq \spimgx{\theind}{\{\theind(\firstspecialpoint)\}}\in\ker(\theind)$ and $\{\secondspecialpoint\}\subseteq \spimgx{\theind}{\{\thevalue\}}$.

 \subproof{Case~1.3.} Finally, let $(\thepartition,\theind)$ be case~\scenariothree{}  and let $(\thelong,\thepoint,\thevalue)$ be its critical data. Then, obviously, $\cutblock=\thelong\backslash \{\thepoint\}\subseteq \spimgx{\theind}{\{\thevalue\}}$ by assumption. And, $\thepartition\backslash \{\theblock\}\cupdot \{\cutblock,\sing\}\backslash \{\varnothing\}=\thepartition\backslash \{\thelong\}\cupdot \{\thelong\backslash \{\thepoint\},\{\thepoint\}\}\finerthan\ker(\theind)$ because, by assumption, for  any $\someblock\in\thepartition\backslash\{\thelong\}$ there exists $\anyblock\in\ker(\theind)$ with $\someblock\subseteq\anyblock$ and because  $\thelong\backslash \{\thepoint\}\subseteq \spimgx{\theind}{\{\thevalue\}}\in\ker(\theind)$ and  $\{\thepoint\}\subseteq \spimgx{\theind}{\{\theind(\thepoint)\}}\in\ker(\theind)$. That proves one implication.

 \subproof{Step~2.} In order to show the converse implication, we assume that there exist $\thepoint\in\tsopx{\inlen}{\outlen}$ and $\thevalue\in\dwi{\thedim}$ such that $\thepartition\backslash \{\theblock\}\cupdot \{\cutblock,\sing\}\backslash \{\varnothing\}\finerthan \ker(\theind)$ and $\cutblock\subseteq \spimgx{\theind}{\{\thevalue\}}$ (which we can by Lem\-ma~\ref{lemma:easy-cohomology-helper-one-point-replacer}) and derive that $(\thepartition,\theind)$ is case~\scenarioone{}, case~\scenariotwo{} or case~\scenariothree{}, thus proving that $(\thepartition,\theind)$ is not case~\scenariofour{} by Lem\-ma~\hyperref[lemma:scenarios_well-defined-2]{\ref*{lemma:scenarios_well-defined}\,\ref*{lemma:scenarios_well-defined-2}}.
Note that the existence of $\thepoint$ requires $\thepartition\neq\varnothing$.  Again, a~case distinction is in order.

 \subproof{Case~2.1.} First, let $\theind(\thepoint)=\thevalue$. Then $\theblock=\cutblock\cupdot \{\thepoint\}\subseteq \spimgx{\theind}{\{\thevalue\}}\in\ker(\theind)$ by  $\cutblock\subseteq \spimgx{\theind}{\{\thevalue\}}$. Thus, $\thepartition\finerthan\ker(\theind)$ by $\thepartition\backslash \{\theblock\}\cupdot \{\cutblock,\sing\}\backslash \{\varnothing\}\finerthan \ker(\theind)$.  In other words, we have shown $(\thepartition,\theind)$ to be case~\scenarioone{}.

 \subproof{Case~2.2.} Similarly, if $\theblock=\{\thepoint\}$, then $\thepartition=\thepartition\backslash \{\theblock\}\cupdot \{\cutblock,\sing\}\backslash \{\varnothing\}\finerthan \ker(\theind)$. Thus, $(\thepartition,\theind)$ is case~\scenarioone{}.

 \subproof{Case~2.3.} If $\theind(\thepoint)\neq \thevalue$ and $|\theblock|=2$, then we put $\firstspecialpoint:=\thepoint$ and we let $\secondspecialpoint$ be the unique element of $\theblock\backslash \{\thepoint\}$. It follows $\{\secondspecialpoint\}=\theblock\backslash \{\thepoint\}\subseteq \spimgx{\theind}{\{\thevalue\}}$ and thus $\theind(\secondspecialpoint)=\thevalue\neq \theind(\thepoint)=\theind(\firstspecialpoint)$ by our assumptions. And,  the premise $\thepartition\backslash\{\{\firstspecialpoint,\secondspecialpoint\}\}\cupdot \{\{\firstspecialpoint\},\{\secondspecialpoint\}\} =\thepartition\backslash \{\theblock\}\cupdot \{\cutblock,\sing\}\backslash \{\varnothing\}\finerthan \ker(\theind)$ means that for any $\someblock\in\thepartition$ with $\someblock\neq\{\firstspecialpoint,\secondspecialpoint\}$ there exists $\anyblock\in \ker(\theind)$ with $\someblock\subseteq \anyblock$. Hence, $(\thepartition,\theind)$ is case~\scenariotwo{} with critical data $\{\firstspecialpoint,\secondspecialpoint\}$.

 \subproof{Case~2.4.} The last remaining possibility is that  $\theind(\thepoint)\neq \thevalue$ and $3\leq |\theblock|$. Putting $\thelong:=\theblock$ implies $\thelong\backslash \{\thepoint\}=\cutblock\subseteq \spimgx{\theind}{\{\thevalue\}}$ by assumption, which is to say $\theind(\anypoint)= \thevalue\neq \theind(\thepoint)$ for any $\anypoint\in\thelong$ with $\anypoint\neq \thepoint$. On the other hand, since $\thepartition\backslash\{\thelong\}\cupdot \{\thelong\backslash \{\thepoint\},\{\thepoint\}\} =\thepartition\backslash \{\theblock\}\cupdot \{\cutblock,\sing\}\backslash \{\varnothing\}\finerthan \ker(\theind)$, for any  $\someblock\in\thepartition$ with $\someblock\neq\thelong$ there exists $\anyblock\in \ker(\theind)$ with $\someblock\subseteq \anyblock$. In other words, $(\thepartition,\theind)$ is case~\scenariothree{} with critical data $(\thelong,\thepoint,\thevalue)$. Thus, both implications are true.
\end{proof}
}

{
  \begin{Example}
      \newcommand{\inlen}{k}
  \newcommand{\outlen}{\ell}
  \newcommand{\theind}{f}
  \newcommand{\theindx}[1]{f(#1)}
  \newcommand{\thepoint}{\Yp{z}}
  \newcommand{\thevalue}{s}
  \newcommand{\thepartition}{p}
  \newcommand{\thelong}{\Yb{Z}}
  \newcommand{\firstspecialpoint}{\Yp{z}_1}
  \newcommand{\secondspecialpoint}{\Yp{z}_2}
  \newcommand{\anypoint}{\Yp{y}}
  \newcommand{\someblock}{\Yb{B}}
  \newcommand{\anyblock}{\Yb{A}}
    \newcommand{\ourpoint}{\Yp{z}}
  \newcommand{\ourvalue}{s}
  \newcommand{\ouraltered}{\theind\downarrow_{\ourpoint}\ourvalue}
  For each of the first three $\theind$ of Example~\ref{example:scenarios} (but not the fourth) one can give at least one $(\ourpoint,\ourvalue)\in\tsopx{\inlen}{\outlen}\Ssetmonoidalproduct\dwi{\thedim}$ such that $\thepartition\finerthan\kerp{\ouraltered}$ (and, in fact, $\kerp{\ouraltered}=\{\{\upp{1},\upp{2},\lop{1}\}, \{\upp{3},\lop{2},\lop{3},\lop{4},\lop{5}\},  \{\upp{4}\}\}$).
      $$
(\thepartition,\theind)=
    \begin{tikzpicture}[baseline=0.91cm]
    \def\scp{0.666}
    \def\linksize{\scp*0.075cm}
    \def\pointsize{\scp*0.25cm}
    \def\dd{\scp*0.5cm}
    \def\dx{\scp*1cm}
    \def\cx{\scp*0.3cm}
    \def\txu{3*\dx}
    \def\txl{4*\dx}
    \def\dy{\scp*1cm}
    \def\cy{\scp*0.3cm}
    \def\ty{3*\dy}
    \tikzset{whp/.style={circle, inner sep=0pt, text width={\pointsize}, draw=black, fill=white}}
    \tikzset{blp/.style={circle, inner sep=0pt, text width={\pointsize}, draw=black, fill=black}}
    \tikzset{lk/.style={regular polygon, regular polygon sides=4, inner sep=0pt, text width={\linksize}, draw=black, fill=black}}
    \draw[dotted] ({0-\dd},{0}) -- ({\txl+\dd},{0});
    \draw[dotted] ({0-\dd},{\ty}) -- ({\txu+\dd},{\ty});
    \coordinate (l1) at ({0+0*\dx},{0+0*\ty}) {};
    \coordinate (l2) at ({0+1*\dx},{0+0*\ty}) {};
    \coordinate (l3) at ({0+2*\dx},{0+0*\ty}) {};
    \coordinate (l4) at ({0+3*\dx},{0+0*\ty}) {};
    \coordinate (l5) at ({0+4*\dx},{0+0*\ty}) {};
    \coordinate (u1) at ({0+0*\dx},{0+1*\ty}) {};
    \coordinate (u2) at ({0+1*\dx},{0+1*\ty}) {};
    \coordinate (u3) at ({0+2*\dx},{0+1*\ty}) {};
    \coordinate (u4) at ({0+3*\dx},{0+1*\ty}) {};
    \node[lk] at ({2*\dx},{1*\dy}) {};
    \draw (l1) -- (u1);
    \draw (l3) -- (u3);
    \draw (l2) |- ++ ({1*\dx},{1*\dy});
    \draw[->] (u2) -- ++({0*\dx},{-1*\dy});
    \draw[->] (u4) -- ++({0*\dx},{-1*\dy});
    \draw (l4) -- ++ ({0*\dx},{1*\dy}) -| (l5);
    \node[above =2pt of u1] {$\scriptstyle 2$};
    \node[above =2pt of u2] {$\scriptstyle 2$};
    \node[above =2pt of u3] {$\scriptstyle 1$};
    \node[above =2pt of u4] {$\scriptstyle 3$};
    \node[below =2pt of l1] {$\scriptstyle 2$};
    \node[below =2pt of l2] {$\scriptstyle 1$};
    \node[below =2pt of l3] {$\scriptstyle \theindx{\lop{3}}$};
    \node[below =2pt of l4] {$\scriptstyle 1$};
     \node[below =2pt of l5] {$\scriptstyle \theindx{\lop{5}}$};
  \end{tikzpicture}
  \hspace{1em} \kerp{\alterer{\theind}{\ourpoint}{\ourvalue}}=
    \begin{tikzpicture}[baseline=0.91cm]
    \def\scp{0.666}
    \def\linksize{\scp*0.075cm}
    \def\pointsize{\scp*0.25cm}
    \def\dd{\scp*0.5cm}
    \def\dx{\scp*1cm}
    \def\cx{\scp*0.3cm}
    \def\txu{3*\dx}
    \def\txl{4*\dx}
    \def\dy{\scp*1cm}
    \def\cy{\scp*0.3cm}
    \def\ty{3*\dy}
    \tikzset{whp/.style={circle, inner sep=0pt, text width={\pointsize}, draw=black, fill=white}}
    \tikzset{blp/.style={circle, inner sep=0pt, text width={\pointsize}, draw=black, fill=black}}
    \tikzset{lk/.style={regular polygon, regular polygon sides=4, inner sep=0pt, text width={\linksize}, draw=black, fill=black}}
    \draw[dotted] ({0-\dd},{0}) -- ({\txl+\dd},{0});
    \draw[dotted] ({0-\dd},{\ty}) -- ({\txu+\dd},{\ty});
    \coordinate (l1) at ({0+0*\dx},{0+0*\ty}) {};
    \coordinate (l2) at ({0+1*\dx},{0+0*\ty}) {};
    \coordinate (l3) at ({0+2*\dx},{0+0*\ty}) {};
    \coordinate (l4) at ({0+3*\dx},{0+0*\ty}) {};
    \coordinate (l5) at ({0+4*\dx},{0+0*\ty}) {};
    \coordinate (u1) at ({0+0*\dx},{0+1*\ty}) {};
    \coordinate (u2) at ({0+1*\dx},{0+1*\ty}) {};
    \coordinate (u3) at ({0+2*\dx},{0+1*\ty}) {};
    \coordinate (u4) at ({0+3*\dx},{0+1*\ty}) {};
    \node[lk] at ({2*\dx},{1*\dy}) {};
    \node[lk] at ({3*\dx},{1*\dy}) {};
    \node[lk] at ({0*\dx},{2*\dy}) {};
    \draw (l1) -- (u1);
    \draw (l3) -- (u3);
    \draw (l2) |- ++ ({1*\dx},{1*\dy});
    \draw (u2) |- ++({-1*\dx},{-1*\dy});
    \draw[->] (u4) -- ++({0*\dx},{-1*\dy});
    \draw (l4) |- ++ ({-1*\dx},{1*\dy});
    \draw (l5) |- ++ ({-1*\dx},{1*\dy});
  \end{tikzpicture}
$$
  \begin{enumerate}
  \item If $\theind(\lop{3})=\theind(\lop{5})=1$, then let $\ourpoint:=\lop{3}$ and $\ourvalue:= 1$.
  \item If $\theind(\lop{3})=1$ and $\theind(\lop{5})=2$, then let $\ourpoint:=\lop{5}$ and $\ourvalue:=1$.
  \item If $\theind(\lop{3})=2$ and $\theind(\lop{5})=1$, then let $\ourpoint:=\lop{3}$ and $\ourvalue:=1$.
  \end{enumerate}
  \end{Example}
}

{
  \newcommand{\thepartition}{p}
  \newcommand{\theind}{f}
  \newcommand{\ourpoint}{\Yp{z}}
  \newcommand{\ourvalue}{s}
  Finally, in each of the three cases~\scenarioone{}--\scenariothree{} the next lemma explains for which $(\ourpoint,\ourvalue)$ the corresponding summand on the right-hand side of the identity in Lem\-ma~\ref{lemma:simplifying_the_definitions_of_the_functionals} has a non-zero factor $\zetfx{\thepartition}{\ker(\alterer{\theind}{\ourpoint}{\ourvalue})}$.
}
  {
  \newcommand{\inlen}{k}
  \newcommand{\outlen}{\ell}
  \newcommand{\theind}{f}
  \newcommand{\theindx}[1]{f(#1)}
  \newcommand{\thepoint}{\Yp{z}}
  \newcommand{\thevalue}{s}
  \newcommand{\thepartition}{p}
  \newcommand{\thealtered}{\theind\downarrow_\thepoint\thevalue}
  \newcommand{\theblock}{\blo{\thepartition}{\thepoint}}
  \newcommand{\cutblock}{\blo{\thepartition}{\thepoint}\backslash \{\thepoint\}}
  \newcommand{\sing}{\{\thepoint\}}
  \newcommand{\thelong}{\Yb{Z}}
  \newcommand{\firstspecialpoint}{\Yp{z}_1}
  \newcommand{\secondspecialpoint}{\Yp{z}_2}
  \newcommand{\specialcounter}{i}
  \newcommand{\specialpoint}[1]{\Yp{z}_{#1}}
  \newcommand{\thegood}{f}
  \newcommand{\theonevalue}{s_1}
  \newcommand{\theothervalue}{s_2}
  \newcommand{\anypoint}{\Yp{y}}
  \newcommand{\someblock}{\Yb{A}}
  \newcommand{\anyblock}{\Yb{B}}
  \newcommand{\ourpoint}{\Yp{z'}}
  \newcommand{\ourvalue}{s'}
  \newcommand{\ouraltered}{\theind\downarrow_{\ourpoint}\ourvalue}
  \begin{Lemma}
    \label{lemma:easy-cohomology-helper-non-trivial-polynomials-two-points}
    Let  $\{\inlen,\outlen\}\subseteq \nnint$, let $\thepartition$  be any set-theo\-re\-ti\-cal partition of $\tsopx{\inlen}{\outlen}$, let $\theind\funcdef\tsopx{\inlen}{\outlen}\to\dwi{\thedim}$ be any mapping and let $\ourpoint\in\tsopx{\inlen}{\outlen}$ and $\ourvalue\in\dwi{\thedim}$ be arbitrary.
    \begin{enumerate}
    \item \label{lemma:easy-cohomology-helper-non-trivial-polynomials-two-points-1}
      If $(\thepartition,\theind)$ is case~\scenarioone{},
      then $\thepartition\finerthan \ker(\ouraltered)$ if and only if either $|\blo{\thepartition}{\ourpoint}|=1$ or  both $2\leq |\blo{\thepartition}{\ourpoint}|$ and  $\ourvalue=\theindx{\ourpoint}$.
    \item \label{lemma:easy-cohomology-helper-non-trivial-polynomials-two-points-2}
      If $(\thepartition,\theind)$ is case~\scenariotwo{} with critical data $\{\firstspecialpoint,\secondspecialpoint\}$,
      then $\thepartition\finerthan \ker(\ouraltered)$ if and only if either both $\ourpoint=\firstspecialpoint$ and $\ourvalue=\theindx{\secondspecialpoint}$ or both $\ourpoint=\secondspecialpoint$ and $\ourvalue=\theindx{\firstspecialpoint}$.
    \item \label{lemma:easy-cohomology-helper-non-trivial-polynomials-two-points-3}
      If $(\thepartition,\theind)$ is case~\scenariothree{} with critical data $(\thelong,\thepoint,\thevalue)$,
      then $\thepartition\finerthan \ker(\ouraltered)$ if and only if $\ourpoint=\thepoint$ and $\ourvalue=\thevalue$.
    \end{enumerate}
  \end{Lemma}
  \begin{proof}
    By Lem\-ma~\ref{lemma:easy-cohomology-helper-one-point-replacer}, the statement $\thepartition\finerthan \ker(\ouraltered)$ is  equivalent to the conjunction of  $\thepartition\backslash\{\blo{\thepartition}{\ourpoint}\}\cupdot \{\blo{\thepartition}{\ourpoint}\backslash \{\ourpoint\}, \{\ourpoint\}\}\backslash \{\varnothing\}\finerthan\ker(\theind)$ and $\blo{\thepartition}{\ourpoint}\backslash \{\ourpoint\}\subseteq\spimgx{\theind}{\{\ourvalue\}}$.

    (a)
    Because $\thepartition\backslash\{\blo{\thepartition}{\ourpoint}\}\cupdot \{\blo{\thepartition}{\ourpoint}\backslash \{\ourpoint\}, \{\ourpoint\}\}\backslash \{\varnothing\}\finerthan\thepartition$, in the situation of \ref{lemma:easy-cohomology-helper-non-trivial-polynomials-two-points-1}, where $\thepartition\finerthan\ker(\theind)$, we only need to determine when $\blo{\thepartition}{\ourpoint}\backslash \{\ourpoint\}\subseteq\spimgx{\theind}{\{\ourvalue\}}$. If $|\blo{\thepartition}{\ourpoint}|=1$, that is, $\blo{\thepartition}{\ourpoint}\backslash \{\ourpoint\}=\varnothing$, this condition is trivially satisfied. And if $2\leq |\blo{\thepartition}{\ourpoint}|$, then $\blo{\thepartition}{\ourpoint}\backslash \{\ourpoint\}\subseteq\spimgx{\theind}{\{\ourvalue\}}$ holds if and only if $\ourvalue=\theindx{\ourpoint}$  because $\blo{\thepartition}{\ourpoint}\backslash \{\ourpoint\}\subseteq\blo{\thepartition}{\ourpoint}\subseteq \spimgx{\theind}{ \{\theindx{\ourpoint}\}}$ by assumption. That proves~\ref{lemma:easy-cohomology-helper-non-trivial-polynomials-two-points-1}.

    (b)
    In case~\ref{lemma:easy-cohomology-helper-non-trivial-polynomials-two-points-2}, if  $\ourpoint\notin\{\firstspecialpoint,\secondspecialpoint\}$, then $\{\firstspecialpoint,\secondspecialpoint\}\in \thepartition\backslash\{\blo{\thepartition}{\ourpoint}\}\cupdot \{\blo{\thepartition}{\ourpoint}\backslash \{\ourpoint\}, \{\ourpoint\}\}\backslash \{\varnothing\}$. However, because $\theind(\firstspecialpoint)\neq \theind(\secondspecialpoint)$ there cannot exist any $\anyblock\in\ker(\theind)$ with $\{\firstspecialpoint,\secondspecialpoint\}\subseteq \anyblock$. Hence, $\ourpoint\notin\{\firstspecialpoint,\secondspecialpoint\}$ excludes  $\thepartition\backslash\{\blo{\thepartition}{\ourpoint}\}\cupdot \{\blo{\thepartition}{\ourpoint}\backslash \{\ourpoint\}, \{\ourpoint\}\}\backslash \{\varnothing\}\finerthan\ker(\theind)$ and thus  $\thepartition\finerthan \ker(\ouraltered)$.

      Hence, $\thepartition\finerthan \ker(\ouraltered)$ requires the existence of $\specialcounter\in\dwi{2}$ with $\ourpoint=\specialpoint{\specialcounter}$. If so, then $\thepartition\backslash\{\blo{\thepartition}{\ourpoint}\}\cupdot \{\blo{\thepartition}{\ourpoint}\backslash \{\ourpoint\}, \{\ourpoint\}\}\backslash \{\varnothing\}=\thepartition\backslash \{\{\firstspecialpoint,\secondspecialpoint\}\}\cupdot \{\{\firstspecialpoint\},\{\secondspecialpoint\}\}\finerthan\ker(\theind)$ since by assumption for any $\someblock\in\thepartition$ with $\someblock\notin \{\firstspecialpoint,\secondspecialpoint\}$ there exists $\anyblock\in\ker(\theind)$ with $\someblock\subseteq \anyblock$. Thus, in this case, $\thepartition\finerthan \ker(\ouraltered)$ is equivalent to $\{\specialpoint{3-\specialcounter}\}=\{\firstspecialpoint,\secondspecialpoint\}\backslash\{\specialpoint{\specialcounter}\}=\blo{\thepartition}{\ourpoint}\backslash \{\ourpoint\}\subseteq\spimgx{\theind}{\{\ourvalue\}}$, i.e., to $\theind(\specialpoint{3-\specialcounter})=\ourvalue$, just as~\ref{lemma:easy-cohomology-helper-non-trivial-polynomials-two-points-2} claimed.

      (c)
     Finally, under the assumptions of \ref{lemma:easy-cohomology-helper-non-trivial-polynomials-two-points-3}, whenever $\ourpoint\notin\thelong$, then $\thelong\in \thepartition\backslash\{\blo{\thepartition}{\ourpoint}\}\cupdot \{\blo{\thepartition}{\ourpoint}\backslash \{\ourpoint\}, \allowbreak\{\ourpoint\}\}\backslash \{\varnothing\}\nfinerthan \ker(\theind)$ by the existence of $\anypoint\in\thelong\backslash \{\thepoint\}\neq\varnothing$ with $\theind(\thepoint)\neq \thevalue=\theind(\anypoint)$. Consequently, $\thepartition\nfinerthan \ker(\ouraltered)$ if $\ourpoint\notin\thelong$.

      For $\ourpoint\in\thelong$, because by assumption there is for any $\someblock\in\thepartition$ with $\someblock\neq \thelong=\blo{\thepartition}{\ourpoint}$ a $\anyblock\in\ker(\theind)$ with $\someblock\subseteq\anyblock$ the condition $\thepartition\backslash\{\blo{\thepartition}{\ourpoint}\}\cupdot \{\blo{\thepartition}{\ourpoint}\backslash \{\ourpoint\}, \{\ourpoint\}\}\backslash \{\varnothing\}\finerthan\ker(\theind)$ simplifies to the existence of $\anyblock\in\ker(\theind)$ with $\blo{\thepartition}{\ourpoint}\backslash \{\ourpoint\}\subseteq \anyblock$, which is subsumed by the second condition. In other words, if  $\ourpoint\in\thelong$, then $\thepartition\finerthan \ker(\ouraltered)$ if and only if $\blo{\thepartition}{\ourpoint}\backslash \{\ourpoint\}\subseteq\spimgx{\theind}{\{\ourvalue\}}$.

      If $\ourpoint\neq \thepoint$, then $\blo{\thepartition}{\ourpoint}\backslash \{\ourpoint\}\not\subseteq\spimgx{\theind}{\{\ourvalue\}}$ because, by $3\leq |\thelong|$, there exist $\anypoint\in\thelong\backslash \{\thepoint,\ourpoint\}$ with $\theind(\anypoint)=\thevalue\neq \theind(\thepoint)$ by assumption. Hence,  $\thepartition\finerthan \ker(\ouraltered)$ requires $\ourpoint=\thepoint$. And in that case it is equivalent to $\thelong\backslash \{\thepoint\}=\blo{\thepartition}{\ourpoint}\backslash \{\ourpoint\}\subseteq\spimgx{\theind}{\{\ourvalue\}}$, which is satisfied if and only if $\ourvalue=\thevalue$ because $\theind(\anypoint)=\thevalue$ for any $\anypoint\in\thelong\backslash\{\thepoint\}\neq\varnothing$. Thus,  the assertion of  \ref{lemma:easy-cohomology-helper-non-trivial-polynomials-two-points-3}  is true as well and, thus, so is the claim overall.
  \end{proof}
 }

In regard of Lem\-mas~\ref{lemma:easy-cohomology-helper-non-trivial-polynomials} and \ref{lemma:easy-cohomology-helper-non-trivial-polynomials-two-points}, we can now improve upon Lem\-ma~\ref{lemma:simplifying_the_definitions_of_the_functionals} as follows.

  {
  \newcommand{\inlen}{k}
  \newcommand{\outlen}{\ell}
  \newcommand{\infree}{g}
  \newcommand{\outfree}{j}
  \newcommand{\inbound}{h}
  \newcommand{\outbound}{i}
  \newcommand{\infreex}[1]{g_{#1}}
 \newcommand{\outfreex}[1]{j_{#1}}
  \newcommand{\inboundx}[1]{h_{#1}}
  \newcommand{\outboundx}[1]{i_{#1}}
  \newcommand{\thepartition}{p}
  \newcommand{\thepoint}{\Yp{z}}
  \newcommand{\incol}{\Yc{c}}
  \newcommand{\outcol}{\Yc{d}}
  \newcommand{\incolx}[1]{\Yc{c}_{#1}}
  \newcommand{\outcolx}[1]{\Yc{d}_{#1}}
  \newcommand{\theentry}[3]{u_{#1,#2}^{#3}}
  \newcommand{\thecocyc}{x}
  \newcommand{\thecocycle}[3]{x_{\theentry{#1}{#2}{#3}}}
  \newcommand{\pind}{a}
  \newcommand{\qind}{b}
  \newcommand{\jind}{j}
  \newcommand{\iind}{i}
  \newcommand{\sind}{s}
  \newcommand{\arel}{r}
  \newcommand{\acol}{\Yc{q}}
  \newcommand{\ascalar}{x}
  \newcommand{\firstspecialpoint}{\Yp{z}_1}
  \newcommand{\secondspecialpoint}{\Yp{z}_2}
  \newcommand{\someblock}{\Yb{A}}
  \newcommand{\anyblock}{\Yb{B}}
  \newcommand{\theind}{f}
  \newcommand{\thecol}{\Yc{w}}
  \newcommand{\theindx}[1]{f(#1)}
  \newcommand{\thecolx}[1]{\Yc{w}(#1)}
  \newcommand{\thelong}{\Yb{Z}}
  \newcommand{\theonepoint}{\Yp{x_1}}
  \newcommand{\theotherpoint}{\Yp{x_2}}
  \newcommand{\thevalue}{s}
  \newcommand{\theonevalue}{s_1}
  \newcommand{\theothervalue}{s_2}
  \newcommand{\anypoint}{\Yp{y}}
  \newcommand{\subind}{q}
  \newcommand{\ourpoint}{\Yp{z}'}
  \newcommand{\ourvalue}{s'}
  \newcommand{\sometuple}{x}
  \newcommand{\sometuplex}[1]{x_{#1}}

  \begin{Lemma}
    \label{lemma:easy-cohomology-helper-one-main-abstract}
Let $\{\inlen,\outlen\}\subseteq \nnint$, let $\incol\in\blaw^{\Ssetmonoidalproduct \inlen}$, let $\outcol\in\blaw^{\Ssetmonoidalproduct \outlen}$, let $\thepartition$ be any  set-theo\-re\-ti\-cal partition of $\tsopx{\inlen}{\outlen}$, let $\infree\in \dwi{\thedim}^{\Ssetmonoidalproduct\inlen}$, let $\outfree\in\dwi{\thedim}^{\Ssetmonoidalproduct\outlen}$, let $\arel:= \therelpoly{\incol}{\outcol}{\thepartition}{\outfree}{\infree}$, let
    $\theind:=\djp{\infree}{\outfree}$, let $\thecol:=\djp{\incol}{\outcol}$, and let~$\sometuple\in\comps^{\Sdirectproduct\thegens}$.
    \begin{enumerate}[label=$(\roman*)$]
    \item\label{lemma:easy-cohomology-helper-one-main-abstract-1} If $(\thepartition,\theind)$ is case~\scenarioone{}, then
            \begin{align*}
\thefuncxx{\arel}{\thecocyc}&{}=        \sum_{\substack{\thepoint\in\tsopx{\inlen}{\outlen}\\{}\Sand |\blo{\thepartition}{\thepoint}|=1}} \sum_{\sind=1}^\thedim
        \begin{rcases}
          \begin{dcases}
            -\thecocycle{\sind}{\theindx{\thepoint}}{\thecolx{\thepoint}}& \tif \thepoint\in\tsopx{\inlen}{0}\\
\thecocycle{\theindx{\thepoint}}{\sind}{\thecolx{\thepoint}}   & \tif \thepoint\in\tsopx{0}{\outlen}
          \end{dcases}
        \end{rcases}\\
        &\quad{}
        +        \sum_{\substack{\thepoint\in\tsopx{\inlen}{\outlen}\\{}\Sand 2\leq|\blo{\thepartition}{\thepoint}|}}
        \begin{rcases}
          \begin{dcases}
            -1& \tif \thepoint\in\tsopx{\inlen}{0}\\
1& \tif \thepoint\in\tsopx{0}{\outlen}
          \end{dcases}
        \end{rcases}        \thecocycle{\theindx{\thepoint}}{\theindx{\thepoint}}{\thecolx{\thepoint}}.
      \end{align*}

    \item\label{lemma:easy-cohomology-helper-one-main-abstract-2}
      If $(\thepartition,\theind)$ is case~\scenariotwo{} with critical data $\{\firstspecialpoint,\secondspecialpoint\}$, then
      \[
\thefuncxx{\arel}{\thecocyc}=\begin{rcases}
          \begin{dcases}
            -\thecocycle{\theindx{\secondspecialpoint}}{\theindx{\firstspecialpoint}}{\thecolx{\firstspecialpoint}}& \tif \firstspecialpoint\in\tsopx{\inlen}{0}\\
\thecocycle{\theindx{\firstspecialpoint}}{\theindx{\secondspecialpoint}}{\thecolx{\firstspecialpoint}} & \tif \firstspecialpoint\in\tsopx{0}{\outlen}
          \end{dcases}
        \end{rcases}
        +
\begin{rcases}
          \begin{dcases}
            -\thecocycle{\theindx{\firstspecialpoint}}{\theindx{\secondspecialpoint}}{\thecolx{\secondspecialpoint}}& \tif \secondspecialpoint\in\tsopx{\inlen}{0}\\
\thecocycle{\theindx{\secondspecialpoint}}{\theindx{\firstspecialpoint}}{\thecolx{\secondspecialpoint}}& \tif \secondspecialpoint\in\tsopx{0}{\outlen}
          \end{dcases}
        \end{rcases}.
      \]
    \item \label{lemma:easy-cohomology-helper-one-main-abstract-3}
      If $(\thepartition,\theind)$ is case~\scenariothree{} with critical data $(\thelong,\thepoint,\thevalue)$, then
      \[
\thefuncxx{\arel}{\thecocyc}=
        \begin{rcases}
          \begin{dcases}
            -\thecocycle{\thevalue}{\theindx{\thepoint}}{\thecolx{\thepoint}}& \tif \thepoint\in\tsopx{\inlen}{0}\\
\thecocycle{\theindx{\thepoint}}{\thevalue}{\thecolx{\thepoint}}& \tif \thepoint\in\tsopx{0}{\outlen}
          \end{dcases}
          \end{rcases}.
      \]
    \item \label{lemma:easy-cohomology-helper-one-main-abstract-4} If $(\thepartition,\theind)$ is case~\scenariofour{}, then
      $\thefuncxx{\arel}{\thecocyc}=0$.
    \end{enumerate}
  \end{Lemma}
  \begin{proof}
    By Lem\-ma~\ref{lemma:simplifying_the_definitions_of_the_functionals},
    \[
\thefuncxx{\arel}{\thecocyc}=\sum_{\ourpoint\in\tsopx{\inlen}{\outlen}}\sum_{\ourvalue=1}^\thedim \zetfx{\thepartition}{\ker(\alterer{\theind}{\ourpoint}{\ourvalue})}
\begin{rcases}
  \begin{dcases}
-\thecocycle{\ourvalue}{\theind(\ourpoint)}{\thecol(\ourpoint)}  & \tif \ourpoint\in\tsopx{\inlen}{0}\\
\thecocycle{\theind(\ourpoint)}{\ourvalue}{\thecol(\ourpoint)}    & \tif \ourpoint\in\tsopx{0}{\outlen}
\end{dcases}\end{rcases}.
\]
From this identity, we see immediately that $\thefuncxx{\arel}{\thecocyc}\neq 0$ requires the existence of $\thepoint\in\tsopx{\inlen}{\outlen}$ and $\thevalue\in\dwi{\thedim}$ with $\thepartition\finerthan \ker(\alterer{\theind}{\thepoint}{\thevalue})$. Thus, Lem\-ma~\ref{lemma:easy-cohomology-helper-non-trivial-polynomials} verifies \ref{lemma:easy-cohomology-helper-one-main-abstract-4}. It remains to treat the cases   \ref{lemma:easy-cohomology-helper-one-main-abstract-1}--\ref{lemma:easy-cohomology-helper-one-main-abstract-3}.

(i)
In the situation of \ref{lemma:easy-cohomology-helper-one-main-abstract-1}, for any $\ourpoint\in\tsopx{\inlen}{\outlen}$ and  $\ourvalue\in\dwi{\thedim}$  we know from Lem\-ma~\hyperref[lemma:easy-cohomology-helper-non-trivial-polynomials-two-points-1]{\ref*{lemma:easy-cohomology-helper-non-trivial-polynomials-two-points}\,\ref*{lemma:easy-cohomology-helper-non-trivial-polynomials-two-points-1}} that $\thepartition\finerthan\ker(\alterer{\theind}{\ourpoint}{\ourvalue})$ if and only if either  $|\blo{\thepartition}{\ourpoint}|=1$ or  both $2\leq |\blo{\thepartition}{\ourpoint}|$ and  $\ourvalue=\theindx{\ourpoint}$. Thus, the above formula for $\thefuncxx{\arel}{\thecocyc}$ simplifies to the one in \ref{lemma:easy-cohomology-helper-one-main-abstract-1}.

(ii)
 Under the assumptions of \ref{lemma:easy-cohomology-helper-one-main-abstract-2}, Lem\-ma~\hyperref[lemma:easy-cohomology-helper-non-trivial-polynomials-two-points-2]{\ref*{lemma:easy-cohomology-helper-non-trivial-polynomials-two-points}\,\ref*{lemma:easy-cohomology-helper-non-trivial-polynomials-two-points-2}} tells us for any $\ourpoint\in\tsopx{\inlen}{\outlen}$ and  $\ourvalue\in\dwi{\thedim}$   that $\thepartition\finerthan\ker(\alterer{\theind}{\ourpoint}{\ourvalue})$ if and only if  either both $\ourpoint=\firstspecialpoint$ and $\ourvalue=\theindx{\secondspecialpoint}$ or both $\ourpoint=\secondspecialpoint$ and $\ourvalue=\theindx{\firstspecialpoint}$. That proves the formula for $\thefuncxx{\arel}{\thecocyc}$ in \ref{lemma:easy-cohomology-helper-one-main-abstract-2}.

      (iii)
Finally, if the premises of  \ref{lemma:easy-cohomology-helper-one-main-abstract-3} are satisfied, then  for any $\ourpoint\in\tsopx{\inlen}{\outlen}$ and  $\ourvalue\in\dwi{\thedim}$  Lem\-ma~\hyperref[lemma:easy-cohomology-helper-non-trivial-polynomials-two-points-3]{\ref*{lemma:easy-cohomology-helper-non-trivial-polynomials-two-points}\,\ref*{lemma:easy-cohomology-helper-non-trivial-polynomials-two-points-3}} lets us infer   that $\thepartition\finerthan\ker(\alterer{\theind}{\ourpoint}{\ourvalue})$ if and only if $\ourpoint=\thepoint$ and $\ourvalue=\thevalue$. In particular, at most one summand is non-zero. It follows that $\thefuncxx{\arel}{\thecocyc}$ is given by the expression in \ref{lemma:easy-cohomology-helper-one-main-abstract-3}.
  \end{proof}
}

\subsection{Halving the number of variables}
Until now we have only considered each equation in the systems from Pro\-po\-si\-tion~\ref{proposition:main_system_of_equations} in isolation. The next simplification will take into account that
the two-colored partitions  $\PartIdenLoWB$, $\PartIdenUpBW$, $\PartIdenLoBW$ and $\PartIdenUpWB$ are present in any category of two-colored partitions. That fact can be used to eliminate half the variables (as, e.g., in \cite[Lemma~1.7]{KyedRaum2017}). This is the only explicit elimination of variables that will be made in the entire proof of the main theorem.
{
  \newcommand{\infree}{g}
  \newcommand{\infreex}[1]{\infree_{#1}}
  \newcommand{\outbound}{i}
  \newcommand{\outboundx}[1]{\outbound_{#1}}
  \newcommand{\outfree}{j}
  \newcommand{\outfreex}[1]{\outfree_{#1}}
  \newcommand{\therel}{r}
  \newcommand{\thetuple}{x}
  \newcommand{\thetuplex}[1]{\thetuple_{#1}}
  \begin{Lemma}\samepage
    \label{lemma:unitarity-relations-first-derivative}
    For any $\infree\in\dwi{\thedim}^{\Ssetmonoidalproduct 2}$ and $\outfree\in\dwi{\thedim}^{\Ssetmonoidalproduct 2}$ and any $\thetuple\in\comps^{\Sdirectproduct\thegens}$, if $\therel$ is given by
    \begin{enumerate}
    \item\label{lemma:unitarity-relations-first-derivative-1}  $\therelpoly{\varnothing}{\whpoint\blpoint}{\UCPartIdenLo}{\outfree}{\varnothing}$, then $\thefuncxx{\therel}{\thetuple}=\thetuplex{\theunix{\whpoint}{\outfreex{1}}{\outfreex{2}}}+\thetuplex{\theunix{\blpoint}{\outfreex{2}}{\outfreex{1}}}$.
    \item\label{lemma:unitarity-relations-first-derivative-2}  $\therelpoly{\varnothing}{\blpoint\whpoint}{\UCPartIdenLo}{\outfree}{\varnothing}$, then $\thefuncxx{\therel}{\thetuple}=\thetuplex{\theunix{\blpoint}{\outfreex{1}}{\outfreex{2}}}+\thetuplex{\theunix{\whpoint}{\outfreex{2}}{\outfreex{1}}}$.
    \item\label{lemma:unitarity-relations-first-derivative-3} $\therelpoly{\blpoint\whpoint}{\varnothing}{\UCPartIdenUp}{\varnothing}{\infree}$, then $\thefuncxx{\therel}{\thetuple}=-\thetuplex{\theunix{\blpoint}{\infreex{2}}{\infreex{1}}}-\thetuplex{\theunix{\whpoint}{\infreex{1}}{\infreex{2}}}$.
      \item\label{lemma:unitarity-relations-first-derivative-4} $\therelpoly{\whpoint\blpoint}{\varnothing}{\UCPartIdenUp}{\varnothing}{\infree}$, then $\thefuncxx{\therel}{\thetuple}=-\thetuplex{\theunix{\whpoint}{\infreex{2}}{\infreex{1}}}-\thetuplex{\theunix{\blpoint}{\infreex{1}}{\infreex{2}}}$.
    \end{enumerate}
  \end{Lemma}
  \begin{proof}
    Only the proof of \ref{lemma:unitarity-relations-first-derivative-1} is given. The others are similar.
Using  Definition~\ref{definition:derivatives-one},    the result of Lem\-ma~\ref{lemma:unitarity-relations} that $\therel=\therelpoly{\varnothing}{\whpoint\blpoint}{\UCPartIdenLo}{\outfree}{\varnothing}= \sum_{\outbound=1}^\thedim \theunix{\whpoint}{\outfreex{1}}{\outbound}\theunix{\blpoint}{\outfreex{2}}{\outbound}-\kron{\outfreex{1}}{\outfreex{2}}1$ implies
    \[
      \thefuncxx{\therel}{\thetuple}=\sum_{\outbound=1}^\thedim (\thetuplex{\theunix{\whpoint}{\outfreex{1}}{\outbound}}\Srightmoduleaction\theunix{\blpoint}{\outfreex{2}}{\outbound}+\theunix{\whpoint}{\outfreex{1}}{\outbound}\Sleftmoduleaction\thetuplex{\theunix{\blpoint}{\outfreex{2}}{\outbound}})=\sum_{\outbound=1}^\thedim(\kron{\outfreex{2}}{\outbound}\,\thetuplex{\theunix{\whpoint}{\outfreex{1}}{\outbound}}+\kron{\outfreex{1}}{\outbound}\,\thetuplex{\theunix{\blpoint}{\outfreex{2}}{\outbound}})=\thetuplex{\theunix{\whpoint}{\outfreex{1}}{\outfreex{2}}}+\thetuplex{\theunix{\blpoint}{\outfreex{2}}{\outfreex{1}}}
    \]
    because  $\thefuncx{1}=0$.
  \end{proof}
  }

\begin{Notation}
  \label{notation-main-one}
  \newcommand{\themat}{v}
  Let $\themat\in\squarematrices{\thedim}{\comps}$ be arbitrary.
  \begin{enumerate}
  \item\label{notation-main-one-1}
    {
  \newcommand{\indj}{j}
  \newcommand{\indi}{i}

  \newcommand{\thematx}[2]{v_{#1,#2}}
  \newcommand{\thetuple}{x^\themat}
  \newcommand{\thetuplex}[1]{x_{#1}^\themat}
  Let $\thetuple\in\comps^{\Sdirectproduct \thegens}$ be such that   for any $\{\indi,\indj\}\subseteq\dwi{\thedim}$,
  \[
    \thetuplex{\theunix{\whpoint}{\indj}{\indi}}:=  \thematx{\indj}{\indi}\qquad\tand\qquad     \thetuplex{\theunix{\blpoint}{\indj}{\indi}}:= -\thematx{\indi}{\indj}.
  \]
  }
\item\label{notation-main-one-2}
  {
  \newcommand{\inlen}{k}
  \newcommand{\outlen}{\ell}
  \newcommand{\infree}{g}
  \newcommand{\outfree}{j}
  \newcommand{\incol}{\Yc{c}}
  \newcommand{\outcol}{\Yc{d}}
  \newcommand{\thepartition}{p}
  \newcommand{\somerel}{r}
  \newcommand{\sometuple}{x}
  \newcommand{\thematx}[2]{v_{#1,#2}}
  \newcommand{\thetuple}{x^\themat}
  \newcommand{\thepartitionset}{\mathcal{P}}
  For any set $\thepartitionset$ of two-colored partitions
  let $\thepredicatex{\thepartitionset}{\themat}$ denote the statement that $\thefuncxx{\somerel}{\thetuple}=0$ for any $\somerel\in\thepartrels{\thepartitionset}$.
      }
      \end{enumerate}
\end{Notation}
Ultimately, it will be shown that in the case of categories of two-colored partitions the predicates $\thepredicate$ defined in Notation~\ref{notation-main-one} are equivalent to those used in the formulation of the main theorem.
{
  \newcommand{\somerel}{r}
  \newcommand{\indj}{j}
  \newcommand{\indi}{i}
  \newcommand{\indg}{g}
  \newcommand{\indh}{h}
  \newcommand{\incol}{\Yc{c}}
  \newcommand{\outcol}{\Yc{d}}
  \newcommand{\incolx}[1]{\Yc{c}_{#1}}
  \newcommand{\outcolx}[1]{\Yc{d}_{#1}}
  \newcommand{\inlen}{k}
  \newcommand{\outlen}{\ell}
  \newcommand{\infree}{g}
  \newcommand{\outfree}{j}
  \newcommand{\inbound}{h}
  \newcommand{\outbound}{i}
  \newcommand{\infreex}[1]{g_{#1}}
  \newcommand{\outfreex}[1]{j_{#1}}
  \newcommand{\inboundx}[1]{h_{#1}}
  \newcommand{\outboundx}[1]{i_{#1}}
  \newcommand{\incounter}{a}
  \newcommand{\outcounter}{b}
  \newcommand{\thepartition}{p}
  \newcommand{\themat}{v}
  \newcommand{\thematx}[2]{v_{#1,#2}}
  \newcommand{\thetuple}{x^\themat}
  \newcommand{\thetuplex}[1]{x_{#1}^\themat}
  \newcommand{\somecocycle}{\eta}
  \begin{Proposition}
    \label{proposition:main-one-explication}
    For any category $\thecat$ of two-co\-lo\-red partitions,
    if $\theqg$ is the unitary easy compact quantum group of $(\thecat,\thedim)$,
then there exists an isomorphism of $\comps$-vector spaces
\[
  \begin{tikzcd}
  \DQGcohomologyTC{1}{\CQGdualC{\theqg}} \arrow[r, two heads, hook]& \{ \themat\in \squarematrices{\thedim}{\comps}\Sand \thepredicatex{\thecat}{\themat}  \}  ,
  \end{tikzcd}
\]
which maps $($the one-elemental cohomology class of$)$ any $1$-cocycle $\somecocycle$ to the matrix $\themat$ with  $\thematx{\indj}{\indi}= \somecocycle(\theunix{\whpoint}{\indj}{\indi}+\thepartideal{\thecat})$ for any $\{\indi,\indj\}\subseteq \dwi{\thedim}$.
\end{Proposition}
\begin{proof}
  \newcommand{\qgmodule}{W}
  \newcommand{\somescalar}{\lambda}
  \newcommand{\someel}{a}
  \newcommand{\thefullcategenextended}{\mathcal{R}'}
  \newcommand{\thefullrels}{R'}
  \newcommand{\theiso}{\varphi}
  \newcommand{\betterqgmodule}{X}
  \newcommand{\sometuple}{x}
  \newcommand{\sometuplex}[1]{x_{#1}}
  \newcommand{\somegen}{e}
  \newcommand{\somepoly}{p}
  \newcommand{\relocatedcocycle}{\overline{\eta}}
  \newcommand{\anycol}{\Yc{c}}
  \newcommand{\thefullideal}{J'}
  By Pro\-po\-si\-tion~\ref{proposition:main_system_of_equations}, it suffices to show that the rule $\sometuple\mapsto (\sometuplex{\theunix{\whpoint}{\indj}{\indi}})_{(\indj,\indi)\in\dwi{\thedim}^{\Ssetmonoidalproduct 2}}$ gives a $\comps$-linear isomorphism
  \[
    \begin{tikzcd}
      \{\sometuple\in \comps^{\Sdirectproduct \thegens}\Sand \forall \somerel\in\thepartrels{\thecat}\quantorpredicate \thefuncxx{\somerel}{\sometuple}=0 \}\arrow[r,two heads, hook]& \{ \themat\in \squarematrices{\thedim}{\comps}\Sand \forall \somerel\in\thepartrels{\thecat}\quantorpredicate \thefuncxx{\somerel}{\thetuple}=0  \} .
    \end{tikzcd}
  \]

The claimed isomorphism is well defined: Let $\sometuple\in \comps^{\Sdirectproduct \thegens}$ be such that $\thefuncxx{\somerel}{\sometuple}=0$ for any $\somerel\in\thepartrels{\thecat}$. Then, for any $\outfree\in\dwi{n}^{\Ssetmonoidalproduct 2}$ because $\therelpoly{\varnothing}{\whpoint\blpoint}{\UCPartIdenLo}{\outfree}{\varnothing}\in\thepartrels{\thecat}$ in particular
\[
 \sometuplex{\theunix{\whpoint}{\outfreex{1}}{\outfreex{2}}}+\sometuplex{\theunix{\blpoint}{\outfreex{2}}{\outfreex{1}}}=0
\]
by Lem\-ma~\ref{lemma:unitarity-relations-first-derivative}, i.e., $\sometuplex{\theunix{\blpoint}{\outfreex{2}}{\outfreex{1}}}=-\sometuplex{\theunix{\whpoint}{\outfreex{1}}{\outfreex{2}}}$. Hence, if we let $\themat:=(\sometuplex{\theunix{\whpoint}{\indj}{\indi}})_{(\indj,\indi)\in\dwi{\thedim}^{\Ssetmonoidalproduct 2}}$, then for any $\{\indi,\indj\}\subseteq \dwi{\thedim}$ by definition not only $\thetuplex{\theunix{\whpoint}{\indj}{\indi}}=\thematx{\indj}{\indi}=\sometuplex{\theunix{\whpoint}{\indj}{\indi}}$  but also
\[
\thetuplex{\theunix{\blpoint}{\indj}{\indi}}=-\thematx{\indi}{\indj}=-\sometuplex{\theunix{\whpoint}{\indi}{\indj}}= \sometuplex{\theunix{\blpoint}{\indj}{\indi}},
\]
 which is to say $\thetuple=\sometuple$. Thus, per assumption, in particular  $\thefuncxx{\somerel}{\thetuple}=\thefuncxx{\somerel}{\sometuple}=0$ for any $\somerel\in\thepartrels{\thecat}$. That proves that the map is well defined.

It is clear that the mapping is $\comps$-linear. Moreover, it is injective because, if again  $\sometuple\in \comps^{\Sdirectproduct \thegens}$ is such that $\thefuncxx{\somerel}{\sometuple}=0$ for any $\somerel\in\thepartrels{\thecat}$ and if again $\themat:=(\sometuplex{\theunix{\whpoint}{\indj}{\indi}})_{(\indj,\indi)\in\dwi{\thedim}^{\Ssetmonoidalproduct 2}}$, then $\themat=0$ necessitates $\thetuple=0$ by definition of $\thetuple$ and thus $\sometuple=0$ by the identity $\thetuple=\sometuple$ established in the preceding paragraph.

To show surjectivity, we let $\themat\in\squarematrices{\thedim}{\comps}$ be arbitrary with  $\thefuncxx{\somerel}{\thetuple}=0$ for any $\somerel\in\therels$ and abbreviate $\sometuple:=\thetuple$. Then, of course, $\thefuncxx{\somerel}{\sometuple}=0$ for any $\somerel\in\thepartrels{\thecat}$. Because moreover  \smash{$\sometuplex{\theunix{\whpoint}{\indj}{\indi}}=\thetuplex{\theunix{\whpoint}{\indj}{\indi}}=\thematx{\indj}{\indi}$} for any $\{\indi,\indj\}\subseteq\dwi{\thedim}$ the tuple
 $\sometuple$ is a preimage of $\themat$. Thus, the claim is true.
\end{proof}
}

{
  \newcommand{\indi}{i}
  \newcommand{\indj}{j}
  \newcommand{\themap}{f}
  \newcommand{\thepartition}{p}
  \newcommand{\incol}{\Yc{c}}
  \newcommand{\outcol}{\Yc{d}}
  \newcommand{\theelim}{(x_{\theunix{\blpoint}{\indj}{\indi}})_{(\indj,\indi)\in\dwi{\thedim}^{\Ssetmonoidalproduct 2}}}
  The next lemma correspondingly eliminates the variables $\theelim$ from the formula obtained in Lem\-ma~\ref{lemma:easy-cohomology-helper-one-main-abstract} for the individual equations in the systems from Pro\-po\-si\-tion~\ref{proposition:main_system_of_equations}. Recall from Notation~\ref{notation:set-theoretical-partitions} that  $\tsquomap{\themap}{\thepartition}$ denotes the quotient mapping of any mapping $\themap$ with respect to any set-theoretical partition $\thepartition$ of its domain and recall the definition of the color sum $\csfunx{\incol}{\outcol}$ of two color tuples $\incol$ and $\outcol$ from  Definition~\ref{definition:points-colors}.
}

{
  \newcommand{\somerel}{r}
  \newcommand{\indj}{j}
  \newcommand{\indi}{i}
  \newcommand{\sind}{s}
  \newcommand{\incol}{\Yc{c}}
  \newcommand{\outcol}{\Yc{d}}
  \newcommand{\incolx}[1]{\Yc{c}_{#1}}
  \newcommand{\outcolx}[1]{\Yc{d}_{#1}}
  \newcommand{\inlen}{k}
  \newcommand{\outlen}{\ell}
  \newcommand{\infree}{g}
  \newcommand{\outfree}{j}
  \newcommand{\infreex}[1]{g_{#1}}
  \newcommand{\outfreex}[1]{j_{#1}}
  \newcommand{\thepartition}{p}
  \newcommand{\themat}{v}
  \newcommand{\thematx}[2]{v_{#1,#2}}
  \newcommand{\thetuple}{x^\themat}
  \newcommand{\thetuplex}[1]{x_{#1}^\themat}
  \newcommand{\theind}{f}
  \newcommand{\theindx}[1]{f(#1)}
  \newcommand{\thecol}{\Yc{w}}
  \newcommand{\thecolx}[1]{\Yc{w}(#1)}
  \newcommand{\thepoint}{\Yp{z}}
  \newcommand{\anypoint}{\Yp{y}}
  \newcommand{\thevalue}{s}
  \newcommand{\thelong}{\Yb{Z}}
  \newcommand{\firstspecialpoint}{\Yp{z}_1}
  \newcommand{\secondspecialpoint}{\Yp{z}_2}
  \newcommand{\someblock}{\Yb{A}}
  \newcommand{\anyblock}{\Yb{B}}
  \begin{Lemma}
    \label{lemma:easy-cohomology-helper-one-main-abstract-corollary}
Let $(\incol,\outcol,\thepartition)$ be any two-colored partition, let $\infree\in \dwi{\thedim}^{\Ssetmonoidalproduct\Slength{\incol}}$, let $\outfree\in\dwi{\thedim}^{\Ssetmonoidalproduct\Slength{\outcol}}$, let ${\somerel:= \therelpoly{\incol}{\outcol}{\thepartition}{\outfree}{\infree}}$, let
    $\theind:=\djp{\infree}{\outfree}$, and let $\themat\in\squarematrices{\thedim}{\comps}$.
    \begin{enumerate}[label=$(\roman*)$]
    \item\label{lemma:easy-cohomology-helper-one-main-abstract-corollary-1} If $(\thepartition,\theind)$ is case~\scenarioone{}, then
  \begin{align*}
          \thefuncxx{\somerel}{\thetuple}&{}=        \sum_{\substack{\someblock\in\thepartition\\{}\Sand |\someblock|=1}}\csx{\incol}{\outcol}{\someblock} \sum_{\sind=1}^\thedim
        \begin{rcases}
          \begin{dcases}
\thematx{(\theind / \thepartition)(\someblock)}{\sind}            & \tif \csx{\incol}{\outcol}{\someblock}=1\\
\thematx{\sind}{(\theind / \thepartition)(\someblock)}& \tif \csx{\incol}{\outcol}{\someblock}=-1
          \end{dcases}
        \end{rcases}\\
        &\quad{}+\sum_{\substack{\someblock\in\thepartition\\{}\Sand 2\leq|\someblock|}}\csx{\incol}{\outcol}{\someblock}
        \thematx{(\theind / \thepartition)(\someblock)}{(\theind / \thepartition)(\someblock)}.
      \end{align*}
    \item\label{lemma:easy-cohomology-helper-one-main-abstract-corollary-2}
      If $(\thepartition,\theind)$ is case~\scenariotwo{} with critical data $\{\firstspecialpoint,\secondspecialpoint\}$, then
            \[
          \thefuncxx{\somerel}{\thetuple}=              \onehalf\csx{\incol}{\outcol}{\{\firstspecialpoint,\secondspecialpoint\}} (
                \thematx{\theindx{\firstspecialpoint}}{\theindx{\secondspecialpoint}}+
                \thematx{\theindx{\secondspecialpoint}}{\theindx{\firstspecialpoint}}).
      \]
    \item \label{lemma:easy-cohomology-helper-one-main-abstract-corollary-3}
      If $(\thepartition,\theind)$ is case~\scenariothree{} with critical data $(\thelong,\thepoint,\thevalue)$, then
      \[
          \thefuncxx{\somerel}{\thetuple}=      \csx{\incol}{\outcol}{\{\thepoint\}}  \begin{rcases}\begin{dcases}
\thematx{\theind(\thepoint)}{\thevalue}& \tif  \csx{\incol}{\outcol}{\{\thepoint\}}=1\\
\thematx{\thevalue}{\theind(\thepoint)}& \tif \csx{\incol}{\outcol}{\{\thepoint\}}=-1
          \end{dcases}\end{rcases}.
      \]
    \item \label{lemma:easy-cohomology-helper-one-main-abstract-corollary-4} If $(\thepartition,\theind)$ is case~\scenariofour{}, then $\thefuncxx{\somerel}{\thetuple}=0$.
    \end{enumerate}
  \end{Lemma}
  \begin{proof}
    \newcommand{\thetransposorsymbol}{\mathrm{b}}
    \newcommand{\thetransposorx}[2]{#1^{\thetransposorsymbol(#2)}}
    \newcommand{\somecolor}{\Yc{e}}
    \newcommand{\anyoldpoint}{\Yp{z}}
    \newcommand{\anyoldvalue}{s}
    \newcommand{\anyoldblock}{\Yb{A}}
    \newcommand{\thecocycle}[3]{\thetuplex{\theunix{#3}{#1}{#2}}}
    \newcommand{\pointind}{i}
    \newcommand{\specialpoint}[1]{\Yp{x}_{#1}}
    We only have to show that in each of the first three cases the right-hand sides of the identities  for $\thefuncxx{\somerel}{\thetuple}$ in the claim agree with the corresponding ones of Lem\-ma~\ref{lemma:easy-cohomology-helper-one-main-abstract}. For the purposes of this proof, let $\thetransposorx{\themat}{1}:=\themat$ and $\thetransposorx{\themat}{-1}:=\themat\Stra$ and recall $\sigma(\whpoint)=1$ and $\sigma(\blpoint)=-1$. Then, for any $\somecolor\in \blaw$ and $\{\indi,\indj\}\subseteq \dwi{\thedim}$ the definitions imply
    \[
\thetuplex{\theunix{\somecolor}{\indj}{\indi}}=
\begin{dcases}
      \begin{rcases}
        \thematx{\indj}{\indi}&\tif \somecolor=\whpoint\\
        -\thematx{\indi}{\indj}&\tif \somecolor=\blpoint
      \end{rcases}\end{dcases}
      =\sigma(\somecolor)(\thetransposorx{\themat}{\sigma(\somecolor)})_{\indj,\indi}.
    \]
    If  $\inlen:= \Slength{\incol}$ and $\outlen:= \Slength{\outcol}$ and $\thecol:=\djp{\incol}{\outcol}$, then it follows for any $\anyoldpoint\in \tsopx{\inlen}{\outlen}$ and any $\anyoldvalue\in \dwi{\thedim}$ that
    \begin{align*}
    \begin{aligned}
      \begin{dcases}
        \begin{rcases}
          -\thetuplex{\theunix{\thecolx{\anyoldpoint}}{\anyoldvalue}{\theindx{\anyoldpoint}}}& \tif \anyoldpoint\in \tsopx{\inlen}{0}\\
          \thetuplex{\theunix{\thecolx{\anyoldpoint}}{\theindx{\anyoldpoint}}{\anyoldvalue}}& \tif \anyoldpoint\in \tsopx{0}{\outlen}
        \end{rcases}
      \end{dcases}&{}=
            \begin{dcases}
        \begin{rcases}
          -{\sigma(\thecolx{\anyoldpoint})}(\thetransposorx{\themat}{\sigma(\thecolx{\anyoldpoint})})_{\anyoldvalue,\theindx{\anyoldpoint}}& \tif \anyoldpoint\in \tsopx{\inlen}{0}\\
{\sigma(\thecolx{\anyoldpoint})}(\thetransposorx{\themat}{\sigma(\thecolx{\anyoldpoint})})_{\theindx{\anyoldpoint},\anyoldvalue}& \tif \anyoldpoint\in \tsopx{0}{\outlen}
        \end{rcases}
      \end{dcases}\\
      &{}=
            \begin{dcases}
        \begin{rcases}
          {\sigma(\ASdualX{\thecolx{\anyoldpoint}})}(\thetransposorx{\themat}{\sigma(\ASdualX{\thecolx{\anyoldpoint}})})_{\theindx{\anyoldpoint},\anyoldvalue}& \tif \anyoldpoint\in \tsopx{\inlen}{0}\\
{\sigma(\thecolx{\anyoldpoint})}(\thetransposorx{\themat}{\sigma(\thecolx{\anyoldpoint})})_{\theindx{\anyoldpoint},\anyoldvalue}& \tif \anyoldpoint\in \tsopx{0}{\outlen}
        \end{rcases}
      \end{dcases}\\
      &{}=\csx{\incol}{\outcol}{\{\anyoldpoint\}}(\thetransposorx{\themat}{\csx{\incol}{\outcol}{\{\anyoldpoint\}}})_{\theindx{\anyoldpoint},\anyoldvalue}
\end{aligned}
    \end{align*}
by the definition of the color sum    and, analogously,
        \[
      \begin{dcases}
        \begin{rcases}
          -\thetuplex{\theunix{\thecolx{\anyoldpoint}}{\theindx{\anyoldpoint}}{\anyoldvalue}}& \tif \anyoldpoint\in \tsopx{\inlen}{0}\\
          \thetuplex{\theunix{\thecolx{\anyoldpoint}}{\anyoldvalue}{\theindx{\anyoldpoint}}}& \tif \anyoldpoint\in \tsopx{0}{\outlen}
        \end{rcases}
      \end{dcases}
      =\csx{\incol}{\outcol}{\{\anyoldpoint\}}(\thetransposorx{\themat}{\csx{\incol}{\outcol}{\{\anyoldpoint\}}})_{\anyoldvalue,\theindx{\anyoldpoint}}.
    \]
    We now distinguish the three relevant cases.

    (i)
    In the situation of \ref{lemma:easy-cohomology-helper-one-main-abstract-corollary-1}, by Lem\-ma~\ref{lemma:easy-cohomology-helper-one-main-abstract} the number  $\thefuncxx{\somerel}{\thetuple}$ is given by
                  \[
        \sum_{\substack{\thepoint\in\tsopx{\inlen}{\outlen}\\{}\Sand |\blo{\thepartition}{\thepoint}|=1}} \sum_{\sind=1}^\thedim
        \begin{rcases}
          \begin{dcases}
            -\thecocycle{\sind}{\theindx{\thepoint}}{\thecolx{\thepoint}}& \tif \thepoint\in\tsopx{\inlen}{0}\\
\thecocycle{\theindx{\thepoint}}{\sind}{\thecolx{\thepoint}}   & \tif \thepoint\in\tsopx{0}{\outlen}
          \end{dcases}
        \end{rcases}
        +        \sum_{\substack{\thepoint\in\tsopx{\inlen}{\outlen}\\{}\Sand 2\leq|\blo{\thepartition}{\thepoint}|}}
        \begin{rcases}
          \begin{dcases}
            -\thecocycle{\theindx{\thepoint}}{\theindx{\thepoint}}{\thecolx{\thepoint}}& \tif \thepoint\in\tsopx{\inlen}{0}\\
\thecocycle{\theindx{\thepoint}}{\theindx{\thepoint}}{\thecolx{\thepoint}}& \tif \thepoint\in\tsopx{0}{\outlen}
          \end{dcases}
        \end{rcases}.
      \]
      By what was shown initially, this can be rewritten identically as
                        \[
        \sum_{\substack{\thepoint\in\tsopx{\inlen}{\outlen}\\{}\Sand |\blo{\thepartition}{\thepoint}|=1}} \sum_{\sind=1}^\thedim \csx{\incol}{\outcol}{\{\thepoint\}}(\thetransposorx{\themat}{\csx{\incol}{\outcol}{\{\thepoint\}}})_{\theindx{\thepoint},\sind}
        +        \sum_{\substack{\thepoint\in\tsopx{\inlen}{\outlen}\\{}\Sand 2\leq|\blo{\thepartition}{\thepoint}|}} \csx{\incol}{\outcol}{\{\thepoint\}}(\thetransposorx{\themat}{\csx{\incol}{\outcol}{\{\thepoint\}}})_{\theindx{\thepoint},\theindx{\thepoint}}.
      \]
      And that is exactly what was claimed because  $\kerp{\theind}\finerthan\thepartition$ and $\sum_{\thepoint\in\anyoldblock} \csx{\incol}{\outcol}{\{\thepoint\}}=\csx{\incol}{\outcol}{\anyoldblock}$ for any~${\anyoldblock\in\thepartition}$.

      (ii)
      Under the assumptions of \ref{lemma:easy-cohomology-helper-one-main-abstract-corollary-2},  Lem\-ma~\ref{lemma:easy-cohomology-helper-one-main-abstract} tells us that $\thefuncxx{\somerel}{\thetuple}$ can be computed as
          \[
\begin{rcases}
          \begin{dcases}
            -\thecocycle{\theindx{\secondspecialpoint}}{\theindx{\firstspecialpoint}}{\thecolx{\firstspecialpoint}}& \tif \firstspecialpoint\in\tsopx{\inlen}{0}\\
\thecocycle{\theindx{\firstspecialpoint}}{\theindx{\secondspecialpoint}}{\thecolx{\firstspecialpoint}} & \tif \firstspecialpoint\in\tsopx{0}{\outlen}
          \end{dcases}
        \end{rcases}
        +
\begin{rcases}
          \begin{dcases}
            -\thecocycle{\theindx{\firstspecialpoint}}{\theindx{\secondspecialpoint}}{\thecolx{\secondspecialpoint}}& \tif \secondspecialpoint\in\tsopx{\inlen}{0}\\
\thecocycle{\theindx{\secondspecialpoint}}{\theindx{\firstspecialpoint}}{\thecolx{\secondspecialpoint}}& \tif \secondspecialpoint\in\tsopx{0}{\outlen}
          \end{dcases}
        \end{rcases},
      \]
      which, by our initial observations, is identical to
      \[
        \csx{\incol}{\outcol}{\{\firstspecialpoint\}} (\thetransposorx{\themat}{\csx{\incol}{\outcol}{\{\firstspecialpoint\}}})_{\theindx{\secondspecialpoint},\theindx{\firstspecialpoint}}+\csx{\incol}{\outcol}{\{\secondspecialpoint\}}(\thetransposorx{\themat}{\csx{\incol}{\outcol}{\{\secondspecialpoint\}}})_{\theindx{\firstspecialpoint},\theindx{\secondspecialpoint}}.
      \]
      Since $\csx{\incol}{\outcol}{\{\Yp{z}_i\}}\in \{-1,1\}$ for each $\pointind\in\dwi{2}$, either $\csx{\incol}{\outcol}{\{\firstspecialpoint\}}=\csx{\incol}{\outcol}{\{\secondspecialpoint\}}$, in which case we infer
      \begin{align*}
        \thefuncxx{\somerel}{\thetuple}&{}=        \csx{\incol}{\outcol}{\{\firstspecialpoint\}}( (\thetransposorx{\themat}{\csx{\incol}{\outcol}{\{\firstspecialpoint\}}})_{\theindx{\secondspecialpoint},\theindx{\firstspecialpoint}}+ (\thetransposorx{\themat}{\csx{\incol}{\outcol}{\{\firstspecialpoint\}}})_{\theindx{\firstspecialpoint},\theindx{\secondspecialpoint}}),\\
        &{}=              \onehalf\csx{\incol}{\outcol}{\{\firstspecialpoint,\secondspecialpoint\}} (
                \thematx{\theindx{\firstspecialpoint}}{\theindx{\secondspecialpoint}}+
                \thematx{\theindx{\secondspecialpoint}}{\theindx{\firstspecialpoint}}),
              \end{align*}
              or $\csx{\incol}{\outcol}{\{\firstspecialpoint\}}=-\csx{\incol}{\outcol}{\{\secondspecialpoint\}}$, implying
                    \[
        \thefuncxx{\somerel}{\thetuple}=       \csx{\incol}{\outcol}{\{\firstspecialpoint\}} (\thetransposorx{\themat}{\csx{\incol}{\outcol}{\{\firstspecialpoint\}}})_{\theindx{\secondspecialpoint},\theindx{\firstspecialpoint}}- \csx{\incol}{\outcol}{\{\firstspecialpoint\}}(\thetransposorx{\themat}{-\csx{\incol}{\outcol}{\{\firstspecialpoint\}}})_{\theindx{\firstspecialpoint},\theindx{\secondspecialpoint}}=0.
      \]
      And that is precisely what we needed to show in this case.

      (iii)
 Finally, if the premises of  \ref{lemma:easy-cohomology-helper-one-main-abstract-3} are satisfied, according to Lem\-ma~\ref{lemma:easy-cohomology-helper-one-main-abstract} and by our initial findings,
            \[
              \thefuncxx{\somerel}{\thetuple}=
              \begin{dcases}
                          \begin{rcases}
            -\thecocycle{\thevalue}{\theindx{\thepoint}}{\thecolx{\thepoint}}& \tif \thepoint\in\tsopx{\inlen}{0},\\
\thecocycle{\theindx{\thepoint}}{\thevalue}{\thecolx{\thepoint}}& \tif \thepoint\in\tsopx{0}{\outlen}
\end{rcases}
              \end{dcases}=\csx{\incol}{\outcol}{\{\thepoint\}}(\thetransposorx{\themat}{\csx{\incol}{\outcol}{\{\thepoint\}}})_{\theindx{\thepoint},\thevalue}.
            \]
              Since this is just what we claimed, that concludes the proof.
  \end{proof}

}

\subsection{All equations of any single two-colored partition}
\label{section:main-one-simplifications-third-simplification}
{
  \newcommand{\inlen}{k}
  \newcommand{\outlen}{\ell}
  \newcommand{\infree}{g}
  \newcommand{\outfree}{j}
  \newcommand{\incol}{\Yc{c}}
  \newcommand{\outcol}{\Yc{d}}
  \newcommand{\thepartition}{p}
  \newcommand{\somerel}{r}
  \newcommand{\sometuple}{x}
    \newcommand{\themat}{v}
  \newcommand{\thematx}[2]{v_{#1,#2}}
  \newcommand{\thetuple}{x^\themat}
  \newcommand{\altoutfree}{j'}
  \newcommand{\altinfree}{g'}
  As an intermediate step to solving the systems of linear equations of Pro\-po\-si\-tion~\ref{proposition:main-one-explication}, we now study the system of equations induced not by an entire category of two-colored partitions but only any single two-colored partition.
}
{
  \newcommand{\somerel}{r}
  \newcommand{\indj}{j}
  \newcommand{\indi}{i}
  \newcommand{\inda}{a}
  \newcommand{\indb}{b}
  \newcommand{\sind}{s}
  \newcommand{\incol}{\Yc{c}}
  \newcommand{\outcol}{\Yc{d}}
  \newcommand{\incolx}[1]{\Yc{c}_{#1}}
  \newcommand{\outcolx}[1]{\Yc{d}_{#1}}
  \newcommand{\inlen}{k}
  \newcommand{\outlen}{\ell}
  \newcommand{\infree}{g}
  \newcommand{\outfree}{j}
  \newcommand{\infreex}[1]{g_{#1}}
  \newcommand{\outfreex}[1]{j_{#1}}
  \newcommand{\thepartition}{p}
  \newcommand{\themat}{v}
  \newcommand{\thematx}[2]{v_{#1,#2}}
  \newcommand{\thetuple}{x^\themat}
  \newcommand{\thetuplex}[1]{x_{#1}^\themat}
  \newcommand{\theind}{f}
  \newcommand{\theindx}[1]{f(#1)}
  \newcommand{\someblock}{\Yb{B}}
  \newcommand{\anyblock}{\Yb{A}}
  \newcommand{\blocklabeling}{h}
  \newcommand{\countblock}{\Yb{A}}
  \newcommand{\thetwoblock}{\Yb{Y}}
  \newcommand{\thelong}{\Yb{Z}}
  \newcommand{\thethreeblock}{\Yb{Z}}
  \begin{Definition}
    \label{definition:the_conditions}
With respect to  any two-colored partition  $(\incol,\outcol,\thepartition)$ we say that any
  $\themat\in\squarematrices{\thedim}{\comps}$ meets
    \begin{enumerate}
  \item \emph{condition~\conditionone{}} if for any $\xfromto{\blocklabeling}{\thepartition}{\dwi{\thedim}}$,
            \[
        \sum_{\substack{\countblock\in\thepartition\\{}\Sand |\countblock|=1}}\csx{\incol}{\outcol}{\countblock} \sum_{\sind=1}^\thedim
        \begin{rcases}
          \begin{dcases}
\thematx{\blocklabeling(\countblock)}{\sind}            & \tif \csx{\incol}{\outcol}{\countblock}=1\\
\thematx{\sind}{\blocklabeling(\countblock)}& \tif \csx{\incol}{\outcol}{\countblock}=-1
          \end{dcases}
        \end{rcases}
        +        \sum_{\substack{\countblock\in\thepartition\\{}\Sand 2\leq|\countblock|}}\csx{\incol}{\outcol}{\countblock}\,
        \thematx{\blocklabeling(\countblock)}{\blocklabeling(\countblock)}=0.
      \]
    \item \emph{condition~\conditiontwo{}} if there is no $\thetwoblock\in\thepartition$ with $|\thetwoblock|=2$ and $\csx{\incol}{\outcol}{\thetwoblock}\neq 0$ or if $\thematx{\indj}{\indi}+\thematx{\indi}{\indj}=0$ for any $\{\indi,\indj\}\subseteq \dwi{\thedim}$ with $\indi\neq \indj$.
        \item \emph{condition~\conditionthree{}} if there is no $\thethreeblock\in\thepartition$ with $3\leq|\thethreeblock|$ or if $\thematx{\indj}{\indi}=0$ for any $\{\indi,\indj\}\subseteq \dwi{\thedim}$ with $\indi\neq \indj$.
    \end{enumerate}
  \end{Definition}
  \begin{Lemma}
    \label{lemma:scalar-article-main-lemma-corollary}
    For any two-colored partition $(\incol,\outcol,\thepartition)$  and any   $\themat\in\squarematrices{\thedim}{\comps}$, the statement    \linebreak $\thepredicatex{\{(\incol,\outcol,\thepartition)\}}{\themat}$ is equivalent to  $\themat$ meeting simultaneously all the three conditions \conditionone{}--\conditionthree{}  with respect to $(\incol,\outcol,\thepartition)$.
  \end{Lemma}

  \begin{proof}
    \newcommand{\thepoint}{\Yp{z}}
  \newcommand{\anypoint}{\Yp{y}}
  \newcommand{\thevalue}{s}
  \newcommand{\firstspecialpoint}{\Yp{z}_1}
  \newcommand{\secondspecialpoint}{\Yp{z}_2}
  Both implications are proved separately.

  \subproof{Step~1. First implication.} First, suppose that  conditions~\conditionone{}--\conditionthree{} are satisfied, let  $\infree\in \dwi{\thedim}^{\Ssetmonoidalproduct\Slength{\incol}}$ and $\outfree\in\dwi{\thedim}^{\Ssetmonoidalproduct\Slength{\outcol}}$ be arbitrary and let $\somerel:= \therelpoly{\incol}{\outcol}{\thepartition}{\outfree}{\infree}$. We show that $\thefuncxx{\somerel}{\thetuple}=0$. If $\theind:=\djp{\infree}{\outfree}$, then~${(\thepartition,\theind)}$ falls into one of the four cases~\scenarioone{}--\scenariofour{} by Lem\-ma~\hyperref[lemma:scenarios_well-defined-2]{\ref*{lemma:scenarios_well-defined}\,\ref*{lemma:scenarios_well-defined-2}}.

  \subproof{Case~1.1.}
If $(\thepartition,\theind)$ is case~\scenarioone{}, then we can define $\blocklabeling:=\tsquomap{\theind}{\thepartition}$. And,  then  by case~\ref{lemma:easy-cohomology-helper-one-main-abstract-corollary-1} of Lem\-ma~\ref{lemma:easy-cohomology-helper-one-main-abstract-corollary} con\-di\-tion~\conditionone{} says precisely that $\thefuncxx{\somerel}{\thetuple}=0$.

\subproof{Case~1.2.}
Next, suppose that $(\thepartition,\theind)$ is case~\scenariotwo{} with critical data  $\{\firstspecialpoint,\secondspecialpoint\}$. Then $\thefuncxx{\somerel}{\thetuple}=              \onehalf\csx{\incol}{\outcol}{\{\firstspecialpoint,\secondspecialpoint\}}(
                \thematx{\theindx{\firstspecialpoint}}{\theindx{\secondspecialpoint}}+
                \thematx{\theindx{\secondspecialpoint}}{\theindx{\firstspecialpoint}})$  by case~\ref{lemma:easy-cohomology-helper-one-main-abstract-corollary-2} of Lem\-ma~\ref{lemma:easy-cohomology-helper-one-main-abstract-corollary}.
                Hence, if $\csx{\incol}{\outcol}{\{\firstspecialpoint,\secondspecialpoint\}}=0$ we have nothing to prove. Otherwise, condition~\conditiontwo{} guarantees that $\thematx{\indb}{\inda}+\thematx{\inda}{\indb}=0$ for any $\{\inda,\indb\}\subseteq \dwi{\thedim}$, thus showing $\thefuncxx{\somerel}{\thetuple}=0$ since $\theindx{\secondspecialpoint}\neq\theindx{\firstspecialpoint}$.

                \subproof{Case~1.3.}                 Now, let $(\thepartition,\theind)$ be case~\scenariothree{} with critical data $(\thelong,\thepoint,\thevalue)$.  Then condition~\conditionthree{} implies that $\themat$ is diagonal. Since       by  case~\ref{lemma:easy-cohomology-helper-one-main-abstract-corollary-3} of Lem\-ma~\ref{lemma:easy-cohomology-helper-one-main-abstract-corollary} the number  $\thefuncxx{\somerel}{\thetuple}$ is given by                 $\csx{\incol}{\outcol}{\{\thepoint\}}\,\thematx{\theind(\thepoint)}{\thevalue}$ or      $\csx{\incol}{\outcol}{\{\thepoint\}}\,\thematx{\thevalue}{\theind(\thepoint)}$ that proves $\thefuncxx{\somerel}{\thetuple}=0$ in this case because $\theind(\thepoint)\neq\thevalue$.

  \subproof{Case~1.4.}               Lastly, if $(\thepartition,\theind)$ is case~\scenariofour{}, then  $\thefuncxx{\somerel}{\thetuple}=0$ by case~\ref{lemma:easy-cohomology-helper-one-main-abstract-corollary-4} of Lem\-ma~\ref{lemma:easy-cohomology-helper-one-main-abstract-corollary}. Hence, there is nothing to show.

                \subproof{Step~2. Second implication.}                                To show the converse we assume $\thefuncxx{\somerel}{\thetuple}=0$ for any $\somerel\in\thepartrels{\{(\incol,\outcol,\thepartition)\}}$ and prove that then conditions~\conditionone{}--\conditionthree{} are met. Let  $\inlen:= \Slength{\incol}$ and $\outlen:= \Slength{\outcol}$.

                \subproof{Step~2.1.} If $\thepartition=\varnothing$, condition~\conditionone{} is trivially satisfied. Otherwise, for any $\xfromto{\blocklabeling}{\thepartition}{\dwi{\thedim}}$ let $\theind:= \blocklabeling\Scomposition \blofun{\thepartition}$, let $\infree\in\dwi{\thedim}^{\Ssetmonoidalproduct\inlen}$ and $\outfree\in\dwi{\thedim}^{\Ssetmonoidalproduct\outlen}$ be such that $\djp{\infree}{\outfree}:=\theind$ and  let $\somerel:= \therelpoly{\incol}{\outcol}{\thepartition}{\outfree}{\infree}$. Then $(\thepartition,\theind)$ is case~\scenarioone{}. Moreover, $\thefuncxx{\somerel}{\thetuple}$ is exactly the left-hand side of the equation in condition~\conditionone{} by Lem\-ma~\hyperref[lemma:easy-cohomology-helper-one-main-abstract-corollary-1]{\ref*{lemma:easy-cohomology-helper-one-main-abstract-corollary}\,\ref*{lemma:easy-cohomology-helper-one-main-abstract-corollary-1}}. This proves condition~\conditionone{} to be satisfied because $\thefuncxx{\somerel}{\thetuple}=0$ by assumption.\looseness=1

                \subproof{Step~2.2.}                 Now, let $\thetwoblock\in\thepartition$ be such that $|\thetwoblock|=2$ and  $\csx{\incol}{\outcol}{\thetwoblock}\neq 0$ and let   $\{\inda,\indb\}\subseteq \dwi{\thedim}$ be arbitrary with $\inda\neq \indb$. We find $\firstspecialpoint\in\tsopx{\inlen}{\outlen}$ and $\secondspecialpoint\in\tsopx{\inlen}{\outlen}$ such that $\firstspecialpoint\neq \secondspecialpoint$ and  $\{\firstspecialpoint,\secondspecialpoint\}=\thetwoblock$. If we let $\theind\funcdef\tsopx{\inlen}{\outlen}\to\dwi{\thedim}$ be such that $\firstspecialpoint\mapsto \inda$ and  $\anypoint\mapsto \indb$ for any $\anypoint\in\tsopx{\inlen}{\outlen}\backslash\{\firstspecialpoint\}$, then   $\theind(\firstspecialpoint)\neq \theind(\secondspecialpoint)$ by $\inda\neq\indb$ and  for any $\anyblock\in\thepartition$ with $\anyblock\neq \{\firstspecialpoint,\secondspecialpoint\}$ there is $\someblock\in\kerp{\theind}$ with $\anyblock\subseteq\someblock$, namely $\someblock=\tsopx{\inlen}{\outlen}\backslash \{\firstspecialpoint\}$. In other words, $(\thepartition,\theind)$ is case~\scenariotwo{}. Hence, $\thefuncxx{\somerel}{\thetuple}=\onehalf\csx{\incol}{\outcol}{\{\firstspecialpoint,\secondspecialpoint\}}(\thematx{\indb}{\inda}+\thematx{\inda}{\indb})$ by  Lem\-ma~\hyperref[lemma:easy-cohomology-helper-one-main-abstract-corollary-2]{\ref*{lemma:easy-cohomology-helper-one-main-abstract-corollary}\,\ref*{lemma:easy-cohomology-helper-one-main-abstract-corollary-2}}. Since $\thefuncxx{\somerel}{\thetuple}=0$ and $\csx{\incol}{\outcol}{\{\firstspecialpoint,\secondspecialpoint\}}\neq 0$ by assumption, condition~\conditiontwo{} is thus met as well.

                \subproof{Step~2.3.}                 Lastly, suppose  $\thethreeblock\in\thepartition$ and $3\leq|\thethreeblock|$ and let  $\{\inda,\indb\}\subseteq \dwi{\thedim}$  be arbitrary with $\inda\neq \indb$. Fix any $\thepoint\in\thethreeblock$, let $\thevalue:= \indb$ and define $\theind\funcdef\tsopx{\inlen}{\outlen}\to\dwi{\thedim}$ by demanding $\thepoint\mapsto \inda$ and $\anypoint\mapsto \indb$ for any $\anypoint\in\tsopx{\inlen}{\outlen}\backslash \{\thepoint\}$. Then, $\theind(\thepoint)\neq \thevalue$ by $\inda\neq\indb$ and $\theind(\anypoint)=\thevalue$ for any $\anypoint\in\thelong$ with $\anypoint\neq\thepoint$. Moreover, for any $\anyblock\in\thepartition$ with $\thepartition\neq\thelong$ there exists in the shape of $\tsopx{\inlen}{\outlen}\backslash \{\thepoint\}$ some   $\someblock\in\kerp{\theind}$ with $\anyblock\subseteq\someblock$. This means that $(\thepartition,\theind)$ is case~\scenariothree{}. Therefore, $\thefuncxx{\somerel}{\thetuple}$ is given by $\csx{\incol}{\outcol}{\{\thepoint\}}\,\thematx{\inda}{\indb}$ or      $\csx{\incol}{\outcol}{\{\thepoint\}}\,\thematx{\indb}{\inda}$ according to  Lem\-ma~\hyperref[lemma:easy-cohomology-helper-one-main-abstract-corollary-3]{\ref*{lemma:easy-cohomology-helper-one-main-abstract-corollary}\,\ref*{lemma:easy-cohomology-helper-one-main-abstract-corollary-3}}. Thus,  by $\thefuncxx{\somerel}{\thetuple}= 0$ and $\csx{\incol}{\outcol}{\{\thepoint\}}\neq 0$ also  condition~\conditionthree{} is satisfied and the proof is complete.
  \end{proof}
}

\subsection{All equations of certain special two-colored partitions}
\label{section:main-one-simplifications-fourth-simplification}
{
  \newcommand{\inlen}{k}
  \newcommand{\outlen}{\ell}
  \newcommand{\infree}{g}
  \newcommand{\outfree}{j}
  \newcommand{\altinfree}{g'}
  \newcommand{\altoutfree}{j'}
  \newcommand{\incol}{\Yc{c}}
  \newcommand{\outcol}{\Yc{d}}
  \newcommand{\thepartition}{p}
  \newcommand{\somerel}{r}
  \newcommand{\sometuple}{x}
    \newcommand{\themat}{v}
  \newcommand{\thematx}[2]{v_{#1,#2}}
  \newcommand{\thetuple}{x^\themat}
In the upcoming case distinctions, it will be useful to already  understand the conditions imposed by a small number of one- or two-elemental sets of special two-colored partitions.
}

{
  \newcommand{\somerel}{r}
  \newcommand{\indj}{j}
  \newcommand{\indi}{i}
  \newcommand{\infree}{g}
  \newcommand{\outfree}{j}
  \newcommand{\infreex}[1]{g_{#1}}
  \newcommand{\outfreex}[1]{j_{#1}}
  \newcommand{\thetuple}{x^\themat}
  \newcommand{\thetuplex}[1]{x_{#1}^\themat}
  \newcommand{\thepartitionset}{\mathcal{S}}
    \newcommand{\incol}{\Yc{c}}
    \newcommand{\inlen}{k}
    \newcommand{\outcol}{\Yc{d}}
    \newcommand{\outlen}{\ell}
    \newcommand{\thepartition}{p}
    \newcommand{\themat}{v}
    \newcommand{\thematx}[2]{v_{#1,#2}}
    \newcommand{\thesum}{\lambda}
    \newcommand{\theid}{\identitymatrix}
    \newcommand{\partparamone}{w}
    \newcommand{\partparamtwo}{m}
    \newcommand{\partparamthree}{s}
    \newcommand{\partparamfour}{t}
    \newcommand{\partparamfive}{K}
    \newcommand{\partparamsix}{F}
    \newcommand{\partparamseven}{E}
    \newcommand{\partparameight}{H}
  \begin{Lemma}
    \label{lemma:scalar-article-small-partitions-cohomology-helper}
    Let $\themat\in\squarematrices{\thedim}{\comps}$ be arbitrary.
\begin{enumerate}
\item\label{lemma:scalar-article-small-partitions-cohomology-helper-1}  $\thepredicatex{\{\PartFourWBWB\}}{\themat}$ is equivalent  to  $\themat$ being diagonal.
\item\label{lemma:scalar-article-small-partitions-cohomology-helper-2}  $\thepredicatex{\{\PartSinglesWB\}}{\themat}$ is equivalent to there existing $\thesum\in\comps$ such that  $\themat-\thesum\theid$ is small.
\item\label{lemma:scalar-article-small-partitions-cohomology-helper-3}  $\thepredicatex{\{\PartSingleW^{\Smonoidalproduct \partparamfour},\PartSingleB^{\Smonoidalproduct \partparamfour}\}}{\themat}$ is equivalent  to $\themat$ being small  for any $\partparamfour\in\pint$.
\end{enumerate}
\end{Lemma}
\begin{proof}
  \newcommand{\inds}{s}
  \newcommand{\indd}{d}
  \newcommand{\somecolor}{\Yc{e}}
  \newcommand{\blocklabeling}{h}
  \newcommand{\someblock}{\Yb{B}}
  \newcommand{\anyblock}{\Yb{A}}
  \newcommand{\thetwoblock}{\Yb{Y}}
  \newcommand{\thethreeblock}{\Yb{Z}}
  \newcommand{\specialblock}{\Yb{A}}
(a)
    Because $|\anyblock|\neq 1$ and $\csx{\varnothing}{\whpoint\blpoint\whpoint\blpoint}{\anyblock}=0$ for the only $\anyblock\in\UCPartFour$ con\-di\-tion~\conditionone{} with respect to~$\PartFourWBWB$ is satisfied
    regardless of whether $\themat$ is diagonal or not. Similarly, since there is no  $\thetwoblock\in\UCPartFour$ with $|\thetwoblock|= 2$ the same is true about con\-di\-tion~\conditiontwo{}. It is con\-di\-tion~\conditionthree{} alone which is relevant. Namely, since there is   $\thethreeblock\in\UCPartFour$ with $3\leq |\thethreeblock|$ it is equivalent to $\themat$ being diagonal. Hence, \ref{lemma:scalar-article-small-partitions-cohomology-helper-1}~follows by  Lem\-ma~\ref{lemma:scalar-article-main-lemma-corollary}.

    (b)
 Because $|\anyblock|=1$ for any $\anyblock\in\UCPartSingles$ and because $\csx{\varnothing}{\whpoint\blpoint}{\{\lop{1}\}}=1$ and $\csx{\varnothing}{\whpoint\blpoint}{\{\lop{2}\}}=-1$,  what con\-di\-tion~\conditionone{} with respect to $\PartSinglesWB$ demands of $\themat$ is that for any $\xfromto{\blocklabeling}{\UCPartSingles}{\dwi{\thedim}}$ the number $\sum_{\inds=1}^\thedim\thematx{\blocklabeling(\{\lop{1}\})}{\inds}-\sum_{\inds=1}^\thedim\thematx{\inds}{\blocklabeling(\{\lop{2}\})}$ be zero. In other words, condition~\conditionone{} is equivalent to $\sum_{\inds=1}^\thedim\thematx{\indj}{\inds}=\sum_{\inds=1}^\thedim\thematx{\inds}{\indi}$ holding for any $\{\indi,\indj\}\subseteq \dwi{\thedim}$, i.e., by Lem\-ma~\hyperref[lemma:small-skewsymmetric-both-1]{\ref*{lemma:small-skewsymmetric-both}\,\ref*{lemma:small-skewsymmetric-both-1}} to there being $\thesum\in\comps$ such that $\themat-\thesum\identitymatrix$ is small.
    At the same time, con\-di\-tions~\conditiontwo{} and \conditionthree{} are always trivially satisfied since there are no $\thetwoblock\in \UCPartSingles$ with $|\thetwoblock|=2$, let alone $\thethreeblock\in\UCPartSingles$ with $3\leq|\thethreeblock|$. Thus,  Lem\-ma~\ref{lemma:scalar-article-main-lemma-corollary} proves~\ref{lemma:scalar-article-small-partitions-cohomology-helper-2}.

  (c)
 Since $|\anyblock|=1$ and $\csx{\varnothing}{\somecolor^{\Smonoidalproduct\partparamfour}}{\anyblock}=\sigma(\somecolor)$ for any $\anyblock\in \UCPartSingle^{\Smonoidalproduct\partparamfour}$ and any $\somecolor\in\blaw$ con\-di\-tion~\conditionone{}  with respect to $\PartSingleW^{\Smonoidalproduct\partparamfour}$ and $\PartSingleB^{\Smonoidalproduct\partparamfour}$ is satisfied by $\themat$ if and only if $\sum_{\indd=1}^\partparamfour\sum_{\inds=1}^\thedim \thematx{\blocklabeling(\{\lop{\indd}\})}{\inds}=0$ respectively $-\sum_{\indd=1}^\partparamfour\sum_{\inds=1}^\thedim \thematx{\inds}{\blocklabeling(\{\lop{\indd}\})}=0$
for any $\xfromto{\blocklabeling}{\UCPartSingle^{\Smonoidalproduct\partparamfour}}{\dwi{\thedim}}$. Moreover, con\-di\-tions~\conditiontwo{} and \conditionthree{}  are vacuous by the absence of any $\thetwoblock\in\UCPartSingle^{\Smonoidalproduct\partparamfour}$ with $|\thetwoblock|=2$ and any $\thethreeblock\in\UCPartSingle^{\Smonoidalproduct\partparamfour}$ with $3\leq |\thethreeblock|$.

Consequently, if $\themat$ is small and thus $\sum_{\inds=1}^\thedim\thematx{\blocklabeling(\{\lop{\indd}\})}{\inds}=\sum_{\inds=1}^\thedim\thematx{\inds}{\blocklabeling(\{\lop{\indd}\})}=0$ for any $\indd\in \dwi{\partparamfour}$ all three conditions~\conditionone{}--\conditionthree{} are met for both $\PartSingleW^{\Smonoidalproduct\partparamfour}$ and $\PartSingleB^{\Smonoidalproduct\partparamfour}$. Hence, $\thepredicatex{\{\PartSingleW^{\Smonoidalproduct\partparamfour},\PartSingleB^{\Smonoidalproduct\partparamfour}\}}{\themat}$ is true in that case by Lem\-ma~\ref{lemma:scalar-article-main-lemma-corollary}.

If, conversely, $\thepredicatex{\{\PartSingleW^{\Smonoidalproduct\partparamfour},\PartSingleB^{\Smonoidalproduct\partparamfour}\}}{\themat}$ is true, Lem\-ma~\ref{lemma:scalar-article-main-lemma-corollary} implies that $\themat$ in particular meets con\-di\-tion~\conditionone{} with respect to $\PartSingleW^{\Smonoidalproduct\partparamfour}$ and $\PartSingleB^{\Smonoidalproduct\partparamfour}$. Thus, for any $\indi\in\dwi{\thedim}$, if  $\xfromto{\blocklabeling}{\UCPartSingle^{\Smonoidalproduct\partparamfour}}{\dwi{\thedim}}$ is constant with value~$\indi$, then $0=\partparamfour\,\sum_{\inds=1}^\thedim \thematx{\indi}{\inds}$ respectively $0=\partparamfour\,\sum_{\inds=1}^\thedim \thematx{\inds}{\indi}$ by what was said initially. By $0<\partparamfour$, that proves~$\themat$ to be small.
\end{proof}
}

\subsection{Case distinctions}
\label{section:main-one-case-distinctions}
{
  \newcommand{\incol}{\Yc{c}}
  \newcommand{\outcol}{\Yc{d}}
  \newcommand{\thepartition}{p}
  \newcommand{\themat}{v}
  The final step to proving the main theorem is upon us.  According to Pro\-po\-si\-tion~\ref{proposition:main-one-explication}, it is enough to show that predicates $\thepredicate$ of Notation~\ref{notation-main-one} and those of the main theorem are equivalent.\looseness=1

  The strategy for that is the same in every case: For the category $\thecat$ of two-colored partitions in question and any $\themat\in\squarematrices{\thedim}{\comps}$, the statement $\thepredicatex{\thecat}{\themat}$ is equivalent to $\thepredicatex{\{(\incol,\outcol,\thepartition)\}}{\themat}$ being true for any $(\incol,\outcol,\thepartition)\in\thecat$. By Lem\-ma~\ref{lemma:scalar-article-main-lemma-corollary}, that is equivalent to the three conditions~\conditionone{}--\conditionthree{} being met with respect to any $(\incol,\outcol,\thepartition)\in\thecat$. Thus, we only need to show that the latter is equivalent to  the statement $\thepredicatex{\thecat}{\themat}$ in the main theorem.
}

{
  \newcommand{\somerel}{r}
  \newcommand{\indj}{j}
  \newcommand{\indi}{i}
  \newcommand{\inda}{i}
   \newcommand{\indb}{j}
  \newcommand{\infree}{g}
  \newcommand{\outfree}{j}
  \newcommand{\infreex}[1]{g_{#1}}
  \newcommand{\outfreex}[1]{j_{#1}}
  \newcommand{\thetuple}{x^\themat}
  \newcommand{\thetuplex}[1]{x_{#1}^\themat}
    \newcommand{\incol}{\Yc{c}}
    \newcommand{\inlen}{k}
    \newcommand{\outcol}{\Yc{d}}
    \newcommand{\outlen}{\ell}
    \newcommand{\thepartition}{p}
    \newcommand{\themat}{v}
    \newcommand{\thematx}[2]{v_{#1,#2}}
    \newcommand{\thesum}{\lambda}
  \newcommand{\someblock}{\Yb{B}}
  \newcommand{\specialblock}{\Yb{Y}}
  \newcommand{\blocklabeling}{h}
  \newcommand{\altblocklabeling}{h'}
  \newcommand{\inds}{s}
  \newcommand{\partparamfour}{t}
  \newcommand{\firstsum}{\lambda_1}
  \newcommand{\secondsum}{\lambda_2}
  \newcommand{\firstincol}{\Yc{c}}
    \newcommand{\firstoutcol}{\Yc{d}}
    \newcommand{\firstpartition}{p}
  \newcommand{\secondincol}{\Yc{a}}
    \newcommand{\secondoutcol}{\Yc{b}}
    \newcommand{\secondpartition}{q}
    \newcommand{\anyblock}{\Yb{A}}
    \newcommand{\thetwoblock}{\Yb{Y}}
    \newcommand{\countblock}{\Yb{A}}
    \newcommand{\thethreeblock}{\Yb{Z}}

\begin{Proposition}  \label{proposition:final-case-o}
  Let $\thecat$ be any category of two-colored partitions and $\themat\in\squarematrices{\thedim}{\comps}$.
  If $\thecat$ is case~$\mathcal{O}$ and
    \begin{enumerate}[label=$(\roman*)$]
    \item\label{proposition:final-case-o-1}  class~$\onlyneutralnonsingletonblocks$, then $\thepredicatex{\thecat}{\themat}$ is equivalent to the absolutely true statement.
    \item\label{proposition:final-case-o-2}  not class~$\onlyneutralnonsingletonblocks$ but  class~$\onlyneutralpartitions$, then $\thepredicatex{\thecat}{\themat}$ is equivalent to there existing $\thesum\in\comps$ such that $\themat-\thesum\identitymatrix$ is skew-symmetric.
    \item\label{proposition:final-case-o-3}  not class~$\onlyneutralpartitions$, then $\thepredicatex{\thecat}{\themat}$ is equivalent to $\themat$ being skew-symmetric.
  \end{enumerate}
\end{Proposition}

\begin{proof}
  By Lem\-ma~\hyperref[lemma:scalar-article-categories-helper-1]{\ref*{lemma:scalar-article-categories-helper}\,\ref*{lemma:scalar-article-categories-helper-1}} and~\ref{lemma:scalar-article-categories-helper-2}, the assumption that $\thecat$ be case $\mathcal{O}$  means $|\someblock|=2$ for any $\someblock\in\thepartition$ and any $(\incol,\outcol,\thepartition)\in\thecat$. Consequently, when checking  conditions~\conditionone{}--\conditionthree{} with respect to any given $(\incol,\outcol,\thepartition)$ there are simplifications.
  \begin{itemize}
  \item Condition~\conditionone{} is met if and only if for any $\blocklabeling\funcdef\thepartition\to\dwi{\thedim}$,
                \[
                \sum_{\countblock\in\thepartition}\csx{\incol}{\outcol}{\countblock}
        \thematx{\blocklabeling(\countblock)}{\blocklabeling(\countblock)}=0.
      \]
\item    Condition~\conditiontwo{} simplifies to the demand that $\thematx{\indb}{\inda}+\thematx{\inda}{\indb}=0$ for any $\{\inda,\indb\}\subseteq \dwi{\thedim}$ with $\inda\neq \indb$ as soon as there is any $\thetwoblock\in \thepartition$ with $\csx{\incol}{\outcol}{\thetwoblock}\neq 0$.
\item    Condition~\conditionthree{}  is trivially satisfied  and can thus be ignored entirely.
  \end{itemize}
  The three cases are treated individually.

  (i)
  If $\thecat$ is class $\onlyneutralnonsingletonblocks$, then all we have to show is that with respect to any $(\incol,\outcol,\thepartition)\in\thecat$ con\-di\-tions~\conditionone{} and \conditiontwo{} are automatically satisfied. And indeed, by the initial simplification condition~\conditionone{} of Lem\-ma~\ref{lemma:scalar-article-main-lemma-corollary} is satisfied for any $(\incol,\outcol,\thepartition)\in\thecat$ and any $\xfromto{\blocklabeling}{\thepartition}{\dwi{\thedim}}$ since already $\csx{\incol}{\outcol}{\countblock}=0$ for any $\countblock\in\thepartition$ by $\thecat$ being class $\onlyneutralnonsingletonblocks$. Likewise, by $\thecat$ being $\onlyneutralnonsingletonblocks$ there are no $\thetwoblock\in\thepartition$ with $\csx{\incol}{\outcol}{\thetwoblock}\neq 0$, meaning  con\-di\-tion~\conditiontwo{} is trivially satisfied.

  (ii)
  The next case is that $\thecat$ is not class $\onlyneutralnonsingletonblocks$ but still class $\onlyneutralpartitions$. Here, we do have to show two implications. And, we  treat them separately.

      Suppose that there exists $\thesum\in\comps$ such that $\themat-\thesum\identitymatrix$ is skew-symmetric and let $(\incol,\outcol,\thepartition)\in\thecat$  and $\xfromto{\blocklabeling}{\thepartition}{\dwi{\thedim}}$ be arbitrary.
      Since $\themat-\thesum\identitymatrix$ is skew-symmetric $\thematx{\indj}{\indj}=\thematx{\indi}{\indi}$ for any $\{\indi,\indj\}\subseteq \dwi{\thedim}$ with $\indj\neq \indi$ by Lem\-ma~\hyperref[lemma:small-skewsymmetric-both-2]{\ref*{lemma:small-skewsymmetric-both}\,\ref*{lemma:small-skewsymmetric-both-2}}. Thus, what con\-di\-tion~\conditionone{} with respect to $(\incol,\outcol,\thepartition)$ actually demands is that the term $\sum_{\countblock\in\thepartition}\csx{\incol}{\outcol}{\countblock}\,\thematx{1}{1}=\tcsx{\incol}{\outcol}\,\thematx{1}{1}$ be zero, which it is since $\thecat$ being class $\onlyneutralpartitions$ ensures $\tcsx{\incol}{\outcol}=0$. Lem\-ma~\hyperref[lemma:small-skewsymmetric-both-2]{\ref*{lemma:small-skewsymmetric-both}\,\ref*{lemma:small-skewsymmetric-both-2}} furthermore guarantees that $\thematx{\indj}{\indi}+\thematx{\indi}{\indj}=0$ for any $\{\indi,\indj\}\subseteq \dwi{\thedim}$ with $\indj\neq \indi$, which is why con\-di\-tion~\conditiontwo{} is met, regardless of whether there is $\thetwoblock\in \thepartition$ with $\csx{\incol}{\outcol}{\thetwoblock}\neq 0$ or not.

      Conversely, let now $\thepredicatex{\thecat}{\themat}$ hold. By assumption, we find $ (\incol,\outcol,\thepartition)\in\thecat$ and  $\specialblock\in\thepartition$ with $\csx{\incol}{\outcol}{\specialblock}\neq 0$ but still $\tcsx{\incol}{\outcol}=0$ and, of course, with  $|\specialblock|=2$ since $\thecat$ is case~$\mathcal{O}$. Hence, $\thematx{\indj}{\indi}+\thematx{\indi}{\indj}=0$ for any $\{\indi,\indj\}\subseteq \dwi{\thedim}$ with $\indi\neq \indj$ by  con\-di\-tion~\conditiontwo{} with respect to $(\incol,\outcol,\thepartition)$. But also, given any $\{\indi,\indj\}\subseteq \dwi{\thedim}$ with $\indi\neq \indj$, if  $\xfromto{\blocklabeling}{\thepartition}{\dwi{\thedim}}$ is such that $\specialblock\mapsto \indj$ and $\countblock\mapsto \indi$ for any $\countblock\in\thepartition\backslash \{\specialblock\}$, then con\-di\-tion~\conditionone{} implies  $\csx{\incol}{\outcol}{\specialblock}\,\thematx{\indj}{\indj}+\sum_{\countblock\in \thepartition\Sand \countblock\neq \specialblock}\csx{\incol}{\outcol}{\countblock}\,\thematx{\indi}{\indi}=0$. Since $0=\tcsx{\incol}{\outcol}=\sum_{\countblock\in \thepartition}\csx{\incol}{\outcol}{\countblock}=\csx{\incol}{\outcol}{\specialblock}+\sum_{\countblock\in \thepartition\Sand \countblock\neq \specialblock}\csx{\incol}{\outcol}{\countblock}$, i.e.,  $\sum_{\countblock\in \thepartition\Sand \countblock\neq \specialblock}\csx{\incol}{\outcol}{\countblock}=-\csx{\incol}{\outcol}{\specialblock}$, this means $\csx{\incol}{\outcol}{\specialblock}(\thematx{\indj}{\indj}-\thematx{\indi}{\indi})=0$, which implies $\thematx{\indj}{\indj}-\thematx{\indi}{\indi}=0$ by $\csx{\incol}{\outcol}{\specialblock}\neq 0$. Hence, there exists $\thesum\in\comps$ such that $\themat-\thesum\identitymatrix$ is skew-symmetric by Lem\-ma~\hyperref[lemma:small-skewsymmetric-both-2]{\ref*{lemma:small-skewsymmetric-both}\,\ref*{lemma:small-skewsymmetric-both-2}}.

      (iii)
 Finally, suppose that  $\thecat$ is not even class $\onlyneutralpartitions$.

  Assume that $\themat$ is skew-symmetric and let $(\incol,\outcol,\thepartition)\in\thecat$ and  $\xfromto{\blocklabeling}{\thepartition}{\dwi{\thedim}}$ be arbitrary.
  Since $\themat$ being skew-symmetric implies $\thematx{\blocklabeling(\countblock)}{\blocklabeling(\countblock)}=0$ for any $\countblock\in\thepartition$ con\-di\-tion~\conditionone{}  is met with respect to $(\incol,\outcol,\thepartition)$. But $\themat$ being skew-symmetric also implies $\thematx{\indj}{\indi}+\thematx{\indi}{\indj}=0$ for any $\{\indi,\indj\}\subseteq\dwi{\thedim}$ with $\indj\neq \indi$, which is why con\-di\-tion~\conditiontwo{} is satisfied no matter whether there is $\thetwoblock\in\thepartition$ with $\csx{\incol}{\outcol}{\thetwoblock}\neq 0$ or not.

  To see the converse, we assume $\thepredicatex{\thecat}{\themat}$. Because $\thecat$ is not class $\onlyneutralpartitions$ there exists $(\incol,\outcol,\thepartition)\in\thecat$ with $\tcsx{\incol}{\outcol}\neq 0$.
For any $\indi\in\dwi{\thedim}$, if $\xfromto{\blocklabeling}{\thepartition}{\dwi{\thedim}}$ is constant with value $\indi$, then con\-di\-tion~\conditionone{} with respect to $(\incol,\outcol,\thepartition)$ implies $0=\sum_{\countblock\in\thepartition}\csx{\incol}{\outcol}{\countblock}\,\thematx{\indi}{\indi}=\tcsx{\incol}{\outcol}\,\thematx{\indi}{\indi}$ and thus $\thematx{\indi}{\indi}=0$ by $\tcsx{\incol}{\outcol}\neq 0$. Furthermore,
since $\tcsx{\incol}{\outcol}=\sum_{\countblock\in\thepartition}\csx{\incol}{\outcol}{\countblock}$ the assumption $\tcsx{\incol}{\outcol}\neq 0$ also requires the existence of at least one $\thetwoblock\in\thepartition$ with $\csx{\incol}{\outcol}{\thetwoblock}\neq 0$. Hence, by the initial simplification  con\-di\-tion~\conditiontwo{}  yields  $\thematx{\indj}{\indi}+\thematx{\indi}{\indj}=0$ for any $\{\indi,\indj\}\subseteq\dwi{\thedim}$ with $\indj\neq \indi$. In other words, $\themat$ is skew-symmetric.
\end{proof}

\begin{Proposition}
  \label{proposition:final-case-b}
  Let $\thecat$ be any category of two-colored partitions and $\themat\in\squarematrices{\thedim}{\comps}$.
  If $\thecat$ is case~$\mathcal{B}$ and
    \begin{enumerate}[label=$(\roman*)$]
    \item\label{proposition:final-case-b-1}  both  class~$\onlyneutralnonsingletonblocks$ and class~$\onlyneutralpartitions$, then $\thepredicatex{\thecat}{\themat}$ is equivalent to there existing $\thesum\in\comps$ such that $\themat-\thesum\identitymatrix$ is small.
    \item\label{proposition:final-case-b-2}  class~$\onlyneutralnonsingletonblocks$ but not class~$\onlyneutralpartitions$, then $\thepredicatex{\thecat}{\themat}$ is equivalent to  $\themat$ being small.
    \item\label{proposition:final-case-b-3}  not class~$\onlyneutralnonsingletonblocks$ but  class~$\onlyneutralpartitions$, then $\thepredicatex{\thecat}{\themat}$ is equivalent to there existing $\thesum\in\comps$ such that $\themat-\thesum\identitymatrix$ is skew-symmetric and small.
    \item\label{proposition:final-case-b-4}  neither class~$\onlyneutralnonsingletonblocks$ nor class~$\onlyneutralpartitions$, then $\thepredicatex{\thecat}{\themat}$ is equivalent to $\themat$ being skew-symmetric and small.
  \end{enumerate}
\end{Proposition}
\begin{proof}
  That $\thecat$ is case~$\mathcal{B}$ requires  $|\someblock|\leq 2$ for any $\someblock\in\thepartition$ and any $(\incol,\outcol,\thepartition)\in\thecat$ by Lem\-ma~\hyperref[lemma:scalar-article-categories-helper-2]{\ref*{lemma:scalar-article-categories-helper}\,\ref*{lemma:scalar-article-categories-helper-2}} and, of course,  $\PartSinglesWB\in \thecat$ by definition. Thus, once more there are simplifications.
  \begin{itemize}
    \item Condition \conditionthree{} is trivially satisfied with respect to any $(\incol,\outcol,\thepartition)\in\thecat$ and will thus be ignored.
      \item  We already know from Lem\-ma~\hyperref[lemma:scalar-article-small-partitions-cohomology-helper-2]{\ref*{lemma:scalar-article-small-partitions-cohomology-helper}\,\ref*{lemma:scalar-article-small-partitions-cohomology-helper-2}} that  $\thepredicatex{\thecat}{\themat}$ implies the existence of $\thesum\in\comps$ such that $\themat-\thesum\identitymatrix$ is small.
  \end{itemize}

  (i)
  As the first case, let $\thecat$ be both class $\onlyneutralnonsingletonblocks$ and class $\onlyneutralpartitions$. Since $\thepredicatex{\thecat}{\themat}$ is already known to require  the existence of $\thesum\in\comps$ such that $\themat-\thesum\identitymatrix$ is small, only the converse implication needs proving.

    Suppose that $\thesum\in\comps$ is such that $\themat-\thesum\identitymatrix$ is small and let $(\incol,\outcol,\thepartition)\in\thecat$ and $\xfromto{\blocklabeling}{\thepartition}{\dwi{\thedim}}$ be arbitrary. By Lem\-ma~\hyperref[lemma:small-skewsymmetric-both-1]{\ref*{lemma:small-skewsymmetric-both}\,\ref*{lemma:small-skewsymmetric-both-1}}, then $\thesum=\sum_{\inds=1}^\thedim\thematx{\blocklabeling(\anyblock)}{\inds}=\sum_{\inds=1}^\thedim\thematx{\inds}{\blocklabeling(\anyblock)}$ for any $\anyblock\in\thepartition$. Hence, and because $\csx{\incol}{\outcol}{\anyblock}=0$ for any $\anyblock\in\thepartition$ with $2\leq |\anyblock |$ by $\thecat$ being class $\onlyneutralnonsingletonblocks$, in order to satisfy con\-di\-tion~\conditionone{} with respect to $(\incol,\outcol,\thepartition)$ the term $\sum_{\anyblock\in\thepartition\Sand |\anyblock|=1}\csx{\incol}{\outcol}{\anyblock}\,\thesum$ has to vanish. And, of course, it does vanish since~$\thecat$ being class $\onlyneutralpartitions$  ensures $0=\tcsx{\incol}{\outcol}=\sum_{\anyblock\in\thepartition\Sand |\anyblock|=1}\csx{\incol}{\outcol}{\anyblock}+\sum_{\anyblock\in\thepartition\Sand 2\leq |\anyblock|}\csx{\incol}{\outcol}{\anyblock}=\sum_{\anyblock\in\thepartition\Sand |\anyblock|=1}\csx{\incol}{\outcol}{\anyblock}$, where the last step is due to $\thecat$ being class $\onlyneutralnonsingletonblocks$. That $\thecat$ is class $\onlyneutralnonsingletonblocks$ also prohibits the existence of $\thetwoblock\in \thepartition$ with $\csx{\incol}{\outcol}{\thetwoblock}\neq 0$ for any $(\incol,\outcol,\thepartition)\in\thecat$, rendering con\-di\-tion~\conditiontwo{} vacuous. Hence, $\thepredicatex{\thecat}{\themat}$ holds.

    (ii)
 Next, suppose that $\thecat$ is class $\onlyneutralnonsingletonblocks$ but not class $\onlyneutralpartitions$. Here, both implications need to be shown and are addressed individually.

    Let $\themat$ be small and let $(\incol,\outcol,\thepartition)\in\thecat$ and $\xfromto{\blocklabeling}{\thepartition}{\dwi{\thedim}}$ be arbitrary. Then, $\sum_{\inds=1}^\thedim\thematx{\blocklabeling(\countblock)}{\inds}=\sum_{\inds=1}^\thedim\thematx{\inds}{\blocklabeling(\countblock)}=0$ for any $\countblock\in\thepartition$. For that reason, the first sum on the left-hand side of the equation in con\-di\-tion~\conditionone{} with respect to $(\incol,\outcol,\thepartition)$ vanishes. And since  $\thecat$ being class $\onlyneutralnonsingletonblocks$ means  $\csx{\incol}{\outcol}{\anyblock}=0$ for any $\anyblock\in\thepartition$ with $2\leq |\anyblock |$ the second term does as well. Hence,  con\-di\-tion~\conditionone{} is satisfied.  The assumption that $\thecat$ is class $\onlyneutralnonsingletonblocks$ and thus   $\csx{\incol}{\outcol}{\thetwoblock}=0$ for any $\thetwoblock\in\thepartition$ with $|\thetwoblock |=2$ also implies that con\-di\-tion~\conditiontwo{} is trivially fulfilled. Hence, $\thepredicatex{\thecat}{\themat}$ is true.

     Conversely, because $\thecat$ is not class $\onlyneutralpartitions$ we find some $(\incol,\outcol,\thepartition)\in\thecat$ with $\partparamfour:=|\tcsx{\incol}{\outcol}|\neq 0$. By
Lem\-ma~\hyperref[lemma:scalar-article-categories-helper-4]{\ref*{lemma:scalar-article-categories-helper}\,\ref*{lemma:scalar-article-categories-helper-4}} that necessitates $\{\PartSingleW^{\Smonoidalproduct \partparamfour},\PartSingleB^{\Smonoidalproduct \partparamfour}\}\in\thecat$.  Hence, $\thepredicatex{\thecat}{\themat}$ requires  $\themat$ to be  small by Lem\-ma~\hyperref[lemma:scalar-article-small-partitions-cohomology-helper-2]{\ref*{lemma:scalar-article-small-partitions-cohomology-helper}\,\ref*{lemma:scalar-article-small-partitions-cohomology-helper-3}}.

(iii)
Now, let $\thecat$ not be class $\onlyneutralnonsingletonblocks$ but class $\onlyneutralpartitions$.

  If $\thesum\in\comps$ is such that $\themat-\thesum\identitymatrix$ is both skew-symmetric and small, then given any $(\incol,\outcol,\thepartition)\in\thecat$ and $\xfromto{\blocklabeling}{\thepartition}{\dwi{\thedim}}$, we infer for any $\countblock\in\thepartition$, first, $\thesum=\sum_{\inds=1}^\thedim\thematx{\blocklabeling(\countblock)}{\inds}=\sum_{\inds=1}^\thedim\thematx{\inds}{\blocklabeling(\countblock)}$ by Lem\-ma~\hyperref[lemma:small-skewsymmetric-both-1]{\ref*{lemma:small-skewsymmetric-both}\,\ref*{lemma:small-skewsymmetric-both-1}} and, second, $\thesum=\thematx{\blocklabeling(\countblock)}{\blocklabeling(\countblock)}$ by Lem\-ma~\hyperref[lemma:small-skewsymmetric-both-2]{\ref*{lemma:small-skewsymmetric-both}\,\ref*{lemma:small-skewsymmetric-both-2}}. Consequently, con\-di\-tion~\conditionone{} is satisfied with respect to $(\incol,\outcol,\thepartition)$ if and only if the term $\sum_{\countblock\in\thepartition\Sand |\countblock|=1}\csx{\incol}{\outcol}{\countblock}\,\thesum+\sum_{\countblock\in\thepartition\Sand 2\leq |\countblock|}\csx{\incol}{\outcol}{\countblock}\,\thesum=\sum_{\countblock\in\thepartition}\csx{\incol}{\outcol}{\countblock}\,\thesum=\tcsx{\incol}{\outcol}\,\thesum$ is zero, which, of course, it is since $\thecat$ being class $\onlyneutralpartitions$  guarantees $\tcsx{\incol}{\outcol}=0$. Because Lem\-ma~\hyperref[lemma:small-skewsymmetric-both-2]{\ref*{lemma:small-skewsymmetric-both}\,\ref*{lemma:small-skewsymmetric-both-2}} also tells us that $\thematx{\indj}{\indi}+\thematx{\indi}{\indj}=0$ for any $\{\indi,\indj\}\subseteq \dwi{\thedim}$ with $\indi\neq \indj$, con\-di\-tion~\conditiontwo{} is met, irrespective of whether there actually is some $\someblock\in\thepartition$ with $|\someblock|=2$ and $\csx{\incol}{\outcol}{\someblock}\neq 0$. Thus, $\thepredicatex{\thecat}{\themat}$.

  Conversely, if $\thepredicatex{\thecat}{\themat}$, then by the initial remark there exists $\firstsum\in\comps$ such that  $\themat-\firstsum\identitymatrix$ is small. Additionally, since $\thecat$ is case~$\mathcal{B}$ and not class $\onlyneutralnonsingletonblocks$ we find some $(\incol,\outcol,\thepartition)\in\thecat$ with the property that there is $\specialblock\in\thepartition$ with $|\specialblock|=2$ and $\csx{\incol}{\outcol}{\specialblock}\neq 0$, which  means $\thematx{\indj}{\indi}+\thematx{\indi}{\indj}=0$ for any $\{\indi,\indj\}\subseteq \dwi{\thedim}$ with $\indi\neq \indj$ by con\-di\-tion~\conditiontwo{} for $(\incol,\outcol,\thepartition)$.
  Moreover, given any $\{\indi,\indj\}\subseteq \dwi{\thedim}$ with $\indi\neq \indj$, if $\xfromto{\blocklabeling}{\thepartition}{\dwi{\thedim}}$ is such that $\specialblock\mapsto \indj$ and $\countblock\mapsto \indi$ for any $\countblock\in\thepartition\backslash \{\specialblock\}$ and if $\xfromto{\altblocklabeling}{\thepartition}{\dwi{\thedim}}$ is constant with value $\indi$, then, considering that $\firstsum=\sum_{\inds=1}^\thedim\thematx{\indi}{\inds}=\sum_{\inds=1}^\thedim\thematx{\inds}{\indi}$ by Lem\-ma~\hyperref[lemma:small-skewsymmetric-both-1]{\ref*{lemma:small-skewsymmetric-both}\,\ref*{lemma:small-skewsymmetric-both-1}},  con\-di\-tion~\conditionone{} with respect to $(\incol,\outcol,\thepartition)$ yields the identities $\sum_{\countblock\in\thepartition\Sand |\countblock|=1}\csx{\incol}{\outcol}{\countblock}\, \firstsum+\csx{\incol}{\outcol}{\specialblock}\,\thematx{\indj}{\indj}+\sum_{\countblock\in\thepartition\Sand 2\leq |\countblock|\Sand \countblock\neq \specialblock}\csx{\incol}{\outcol}{\countblock}\,\thematx{\indi}{\indi}=0$ for $\blocklabeling$ and $\sum_{\countblock\in\thepartition\Sand |\countblock|=1}\csx{\incol}{\outcol}{\countblock}\, \firstsum+\sum_{\countblock\in\thepartition\Sand 2\leq |\countblock|}\csx{\incol}{\outcol}{\countblock}\,\thematx{\indi}{\indi}=0$ for $\altblocklabeling$. Subtracting the second from the first produces the identity $\csx{\incol}{\outcol}{\specialblock}(\thematx{\indj}{\indj}-\thematx{\indi}{\indi})=0$. Since  $\csx{\incol}{\outcol}{\specialblock}\neq 0$ we can  infer $\thematx{\indj}{\indj}=\thematx{\indi}{\indi}$ for any $\{\indi,\indj\}\subseteq \dwi{\thedim}$ with $\indi\neq \indj$.
  By Lem\-ma~\hyperref[lemma:small-skewsymmetric-both-2]{\ref*{lemma:small-skewsymmetric-both}\,\ref*{lemma:small-skewsymmetric-both-2}}, we have thus shown that there exists $\secondsum\in\comps$ such that $\themat-\secondsum\identitymatrix$ is skew-symmetric. According to Lem\-ma~\hyperref[lemma:small-skewsymmetric-both-3]{\ref*{lemma:small-skewsymmetric-both}\,\ref*{lemma:small-skewsymmetric-both-3}}, that is all we needed to see.

(iv)
 As the final case let $\thecat$ be neither class $\onlyneutralnonsingletonblocks$ nor class $\onlyneutralpartitions$.

  If $\themat$ is skew-symmetric and small and if $(\incol,\outcol,\thepartition)\in\thecat$ and $\xfromto{\blocklabeling}{\thepartition}{\dwi{\thedim}}$ are arbitrary, then by definition, $\sum_{\inds=1}^\thedim\thematx{\blocklabeling(\countblock)}{\inds}=\sum_{\inds=1}^\thedim\thematx{\inds}{\blocklabeling(\countblock)}=0$ and $\thematx{\blocklabeling(\countblock)}{\blocklabeling(\countblock)}=0$ for any  $\anyblock\in \thepartition$. For that reason, con\-di\-tion~\conditionone{} is trivially satisfied with respect to $(\incol,\outcol,\thepartition)$. Because also $\thematx{\indj}{\indi}+\thematx{\indi}{\indj}=0$ for any $\{\indi,\indj\}\subseteq\dwi{\thedim}$ with $\indi\neq \indj$ con\-di\-tion~\conditiontwo{} is met as well, no matter whether there exists  $\thetwoblock\in\thepartition$ with $|\thetwoblock|=2$ and $\csx{\incol}{\outcol}{\thetwoblock}\neq 0$. Hence, $\thepredicatex{\thecat}{\themat}$ has been proved.

  In order to prove the converse, let $\thepredicatex{\thecat}{\themat}$ hold. Since $\thecat$ is not class $\onlyneutralpartitions$ we find a $(\firstincol,\firstoutcol,\firstpartition)\in\thecat$ with $\partparamfour:=|\tcsx{\firstincol}{\firstoutcol}|\neq 0$. As $\PartSinglesWB\in \thecat$ we conclude $\{\PartSingleW^{\Smonoidalproduct\partparamfour},\PartSingleB^{\Smonoidalproduct\partparamfour}\}\subseteq\thecat$ by Lem\-ma~\hyperref[lemma:scalar-article-categories-helper-4]{\ref*{lemma:scalar-article-categories-helper}\,\ref*{lemma:scalar-article-categories-helper-4}}.  It follows that $\themat$ is small by Lem\-ma~\hyperref[lemma:scalar-article-small-partitions-cohomology-helper-3]{\ref*{lemma:scalar-article-small-partitions-cohomology-helper}\,\ref*{lemma:scalar-article-small-partitions-cohomology-helper-3}}. Furthermore, the assumption of $\thecat$ not being class $\onlyneutralnonsingletonblocks$ implies the existence of $(\secondincol,\secondoutcol,\secondpartition)\in\thecat$ and $\specialblock\in \secondpartition$ with $2\leq |\specialblock|$ and $\csx{\secondincol}{\secondoutcol}{\specialblock}\neq 0$. If now for any $\{\indi,\indj\}\subseteq \dwi{\thedim}$ with $\indi\neq\indj$ the mapping $\xfromto{\blocklabeling}{\secondpartition}{\dwi{\thedim}}$ is such that $\specialblock\mapsto \indj$ and $\countblock\mapsto \indi$ for any $\countblock\in\secondpartition\backslash\{\specialblock\}$ and if $\xfromto{\altblocklabeling}{\secondpartition}{\dwi{\thedim}}$ is constant with value $\indi$, then con\-di\-tion~\conditionone{} for $(\secondincol,\secondoutcol,\secondpartition)$ implies the identities $\csx{\secondincol}{\secondoutcol}{\specialblock}\, \thematx{\indj}{\indj}+\sum_{\countblock\in\secondpartition\Sand 2\leq |\countblock|\Sand \countblock\neq \specialblock}\csx{\secondincol}{\secondoutcol}{\countblock}\,\thematx{\indi}{\indi}=0$ and $\sum_{\countblock\in\secondpartition\Sand 2\leq |\countblock|}\csx{\secondincol}{\secondoutcol}{\someblock}\,\thematx{\indi}{\indi}=0$ because $\sum_{\inds=1}^\thedim\thematx{\indi}{\inds}=\sum_{\inds=1}^\thedim\thematx{\inds}{\indi}=0$ by $\themat$ being small. Subtracting the second from the first yields $\csx{\secondincol}{\secondoutcol}{\specialblock}(\thematx{\indj}{\indj}-\thematx{\indi}{\indi})=0$ and thus $\thematx{\indj}{\indj}=\thematx{\indi}{\indi}$ by $\csx{\secondincol}{\secondoutcol}{\specialblock}\neq 0$. Because the presence of $\specialblock$ in $\secondpartition$ also ensures $\thematx{\indj}{\indi}+\thematx{\indi}{\indj}=0$ for any $\{\indi,\indj\}\subseteq\dwi{\thedim}$ with $\indi\neq \indj$ by con\-di\-tion~\conditiontwo{} for $(\secondincol,\secondoutcol,\secondpartition)$ we have thus shown that there exists $\secondsum\in \comps$ such that $\themat-\secondsum\identitymatrix$ is skew-symmetric by Lem\-ma~\hyperref[lemma:small-skewsymmetric-both-2]{\ref*{lemma:small-skewsymmetric-both}\,\ref*{lemma:small-skewsymmetric-both-2}}. Because $\themat$ is also small, applying Lem\-ma~\hyperref[lemma:small-skewsymmetric-both-3]{\ref*{lemma:small-skewsymmetric-both}\,\ref*{lemma:small-skewsymmetric-both-3}} (with $\firstsum=0$) we see that $\secondsum=0$, i.e., that $\themat$ is skew-symmetric and small.
  \end{proof}

\begin{Proposition}
  \label{proposition:final-case-h}
  Let $\thecat$ be any category of two-colored partitions and $\themat\in\squarematrices{\thedim}{\comps}$.
  If $\thecat$ is case~$\mathcal{H}$ and
    \begin{enumerate}[label=$(\roman*)$]
    \item\label{proposition:final-case-h-1}  class~$\onlyneutralnonsingletonblocks$, then $\thepredicatex{\thecat}{\themat}$ is equivalent to $\themat$ being diagonal.
    \item\label{proposition:final-case-h-2}  not class~$\onlyneutralnonsingletonblocks$ but  class~$\onlyneutralpartitions$, then $\thepredicatex{\thecat}{\themat}$ is equivalent to there existing $\thesum\in\comps$ such that $\themat=\thesum\identitymatrix$.
    \item\label{proposition:final-case-h-3}  not class~$\onlyneutralpartitions$, then $\thepredicatex{\thecat}{\themat}$ is equivalent to $\themat=0$.
  \end{enumerate}
\end{Proposition}
\begin{proof}
  As $\thecat$ is case~$\mathcal{H}$, both $\PartFourWBWB\in \thecat$ and $2\leq |\someblock|$ for any $\someblock\in\thepartition$ and any $(\incol,\outcol,\thepartition)\in\thecat$ by Lem\-ma~\hyperref[lemma:scalar-article-categories-helper-1]{\ref*{lemma:scalar-article-categories-helper}\,\ref*{lemma:scalar-article-categories-helper-1}} and $\PartSinglesWB\notin\thecat$. Certain simplifications result.
  \begin{itemize}
  \item Condition~\conditionone{} with respect to any $(\incol,\outcol,\thepartition)\in\thecat$ amounts to the demand  that for any ${\xfromto{\blocklabeling}{\thepartition}{\dwi{\thedim}}}$,
    \[
                \sum_{\countblock\in\thepartition}\csx{\incol}{\outcol}{\countblock}\,
        \thematx{\blocklabeling(\countblock)}{\blocklabeling(\countblock)}=0.
      \]
  \item We already know by Lem\-ma~\hyperref[lemma:scalar-article-small-partitions-cohomology-helper-1]{\ref*{lemma:scalar-article-small-partitions-cohomology-helper}\,\ref*{lemma:scalar-article-small-partitions-cohomology-helper-1}} that $\thepredicatex{\thecat}{\themat}$ implies that $\themat$ is  diagonal.
  \end{itemize}

  (i)
  Suppose first that $\thecat$ is class $\onlyneutralnonsingletonblocks$. Since it is already clear that $\thepredicatex{\thecat}{\themat}$ requires $\themat$ to be diagonal, only one implication needs proving.

    If $\themat$ is diagonal, if $(\incol,\outcol,\thepartition)\in\thecat$ and if $\xfromto{\blocklabeling}{\thepartition}{\dwi{\thedim}}$, then because  $\csx{\incol}{\outcol}{\countblock}=0$ for any $\countblock\in\thepartition$ by~$\thecat$ being class $\onlyneutralnonsingletonblocks$ the simplified con\-di\-tion~\conditionone{} of $(\incol,\outcol,\thepartition)$ is satisfied trivially. For the same reason, con\-di\-tion~\conditiontwo{} is vacuous. And con\-di\-tion~\conditionthree{} is met as well, regardless of whether there is $\thethreeblock\in\thepartition$ with $3\leq|\thethreeblock|$, because $\themat$ is diagonal per assumption. Hence, $\themat$ being diagonal implies $\thepredicatex{\thecat}{\themat}$.

    (ii)
 Next, suppose that $\thecat$ is not class $\onlyneutralnonsingletonblocks$ but is class $\onlyneutralpartitions$. Now, both implications must be proved.

First, let  $\thesum\in\comps$ be such that $\themat=\thesum\identitymatrix$ and let $(\incol,\outcol,\thepartition)\in\thecat$ and $\xfromto{\blocklabeling}{\thepartition}{\dwi{\thedim}}$ be arbitrary. By the initial remark  con\-di\-tion~\conditionone{} with respect to $(\incol,\outcol,\thepartition)$ demands precisely that the term $\sum_{\countblock\in\thepartition}\csx{\incol}{\outcol}{\countblock}\,\thesum=\tcsx{\incol}{\outcol}\,\thesum$ vanish, which, of course, it does because $\tcsx{\incol}{\outcol}=0$ by $\thecat$ being class $\onlyneutralpartitions$. Moreover, since $\thematx{\indj}{\indi}=0$ for any $\{\indi,\indj\}\subseteq \dwi{\thedim}$ with $\indi\neq \indj$ con\-di\-tion~\conditiontwo{} is certainly satisfied, even if there is $\someblock\in\thepartition$  with $|\someblock|=2$ and $\csx{\incol}{\outcol}{\someblock}\neq 0$. The assumption that $\themat$ is diagonal also ensures that con\-di\-tion~\conditionthree{} is met, irrespective of whether there exists $\thethreeblock\in\thepartition$ with $3\leq|\thethreeblock|$ or not. Thus, $\thepredicatex{\thecat}{\themat}$.

Conversely, if $\thepredicatex{\thecat}{\themat}$, then
$\thecat$ being not class $\onlyneutralnonsingletonblocks$ lets us find some $(\incol,\outcol,\thepartition)\in\thecat$ and $\specialblock\in\thepartition$ with $2\leq |\specialblock|$ and $\csx{\incol}{\outcol}{\specialblock}\neq 0$. If, given any $\{\indi,\indj\}\subseteq \dwi{\thedim}$ with $\indi\neq \indj$ we let $\xfromto{\blocklabeling}{\thepartition}{\dwi{\thedim}}$ be such that $\specialblock\mapsto \indj$ and $\countblock\mapsto\indi$ for any $\countblock\in\thepartition\backslash \{\specialblock\}$, then
con\-di\-tion~\conditionone{} lets us know that $\csx{\incol}{\outcol}{\specialblock}\,\thematx{\indj}{\indj}+\sum_{\countblock\in\thepartition\Sand \countblock\neq \specialblock}\csx{\incol}{\outcol}{\countblock}\,\thematx{\indi}{\indi}=0$. Since $\thecat$ being class $\onlyneutralpartitions$ implies $0=\tcsx{\incol}{\outcol}=\csx{\incol}{\outcol}{\specialblock}+\sum_{\countblock\in\thepartition\Sand \countblock\neq \specialblock}\csx{\incol}{\outcol}{\countblock}$  that is the same as saying $\csx{\incol}{\outcol}{\specialblock}(\thematx{\indj}{\indj}-\thematx{\indi}{\indi})=0$, which means $\thematx{\indj}{\indj}=\thematx{\indi}{\indi}$ by $\csx{\incol}{\outcol}{\specialblock}\neq 0$. Hence, if $\thesum:=\thematx{1}{1}$, then $\themat=\thesum\identitymatrix$ as claimed because $\themat$ is diagonal by the initial remark.

(iii) Lastly, assume $\thecat$ is not class  $\onlyneutralpartitions$. Because for $\themat=0$  con\-di\-tions~\conditionone{}--\conditionthree{} are trivially satisfied, we only need to prove the converse.

  If $\thepredicatex{\thecat}{\themat}$, then by $\thecat$ not being class  $\onlyneutralpartitions$ there exists $(\incol,\outcol,\thepartition)\in\thecat$ with $\tcsx{\incol}{\outcol}\neq 0$. Hence, for any $\indi\in\dwi{\thedim}$, if $\xfromto{\blocklabeling}{\thepartition}{\dwi{\thedim}}$ is constant with value $\indi$, then by what was said at the beginning con\-di\-tion~\conditionone{} with respect to $(\incol,\outcol,\thepartition)$ shows that $0=\sum_{\countblock\in\thepartition}\csx{\incol}{\outcol}{\countblock}\, \thematx{\indi}{\indi}=\tcsx{\incol}{\outcol}\,\thematx{\indi}{\indi}$, i.e., that $\thematx{\indi}{\indi}=0$. As $\themat$ is diagonal by the same initial remarks, that means $\themat=0$, as asserted.
\end{proof}

\begin{Proposition}
  \label{proposition:final-case-s}
  Let $\thecat$ be any category of two-colored partitions and $\themat\in\squarematrices{\thedim}{\comps}$.
  If $\thecat$ is case~$\mathcal{S}$ and
    \begin{enumerate}[label=$(\roman*)$]
    \item\label{proposition:final-case-s-1}  class~$\onlyneutralpartitions$, then $\thepredicatex{\thecat}{\themat}$ is equivalent to there existing $\thesum\in\comps$ such that $\themat=\thesum\identitymatrix$.
    \item\label{proposition:final-case-s-2}  not class~$\onlyneutralpartitions$, then $\thepredicatex{\thecat}{\themat}$ is equivalent to $\themat=0$.
  \end{enumerate}
\end{Proposition}
\begin{proof}
  In contrast to the situation in the cases $\mathcal{O}$, $\mathcal{B}$ and $\mathcal{H}$, there are no general simplifications of the conditions~\conditionone{}--\conditionthree{} of Lem\-ma~\ref{lemma:scalar-article-main-lemma-corollary} implied by the assumption that $\thecat$ is case $\mathcal{S}$. However, as in case $\mathcal{H}$, since   $\PartFourWBWB\in\thecat$ it is already  clear by Lem\-ma~\hyperref[lemma:scalar-article-small-partitions-cohomology-helper-1]{\ref*{lemma:scalar-article-small-partitions-cohomology-helper}\,\ref*{lemma:scalar-article-small-partitions-cohomology-helper-1}} that $\thepredicatex{\thecat}{\themat}$ holding implies that $\themat$ is diagonal.

  (i)
  First, let $\thecat$ be class $\onlyneutralpartitions$.
    If there is $\thesum\in\comps$ such that  $\themat=\thesum\identitymatrix$ and if $(\incol,\outcol,\thepartition)\in\thecat$ and $\xfromto{\blocklabeling}{\thepartition}{\dwi{\thedim}}$  are arbitrary, then what con\-di\-tion~\conditionone{} with respect to $(\incol,\outcol,\thepartition)$ demands is that the sum $\sum_{\countblock\in\thepartition\Sand |\countblock|=1}\csx{\incol}{\outcol}{\countblock}\, \thesum\allowbreak+\sum_{\countblock\in\thepartition\Sand 2\leq |\countblock|}\csx{\incol}{\outcol}{\countblock}\, \thesum=\sum_{\countblock\in\thepartition}\csx{\incol}{\outcol}{\countblock}\,\thesum=\tcsx{\incol}{\outcol}\,\thesum$ vanish. And because $\thecat$ being class $\onlyneutralpartitions$ implies $\tcsx{\incol}{\outcol}=0$ this is indeed the case. Moreover, $\themat$ being diagonal of course guarantees that con\-di\-tions~\conditiontwo{} and~\conditionthree{} are satisfied, no matter what the blocks of  $(\incol,\outcol,\thepartition)$ are. That proves~${\thepredicatex{\thecat}{\themat}}$.

    If, conversely, $\thepredicatex{\thecat}{\themat}$ is assumed to hold, then by Lem\-ma~\hyperref[lemma:scalar-article-small-partitions-cohomology-helper-2]{\ref*{lemma:scalar-article-small-partitions-cohomology-helper}\,\ref*{lemma:scalar-article-small-partitions-cohomology-helper-2}} there exists $\thesum\in\comps$ such that $\themat-\thesum\identitymatrix$ is small since $\PartSinglesWB\in \thecat$. For any $\indi\in\dwi{\thedim}$ the definition of smallness implies $0=\sum_{\indj=1}^\thedim (\thematx{\indj}{\indi}-\thesum\kron{\indj}{\indi})=\thematx{\indi}{\indi}-\thesum$. Hence, $\themat=\thesum\identitymatrix$, as claimed.

    (ii)
 The alternative is that $\thecat$ is not class $\onlyneutralpartitions$.
    Of course, if $\themat=0$, then con\-di\-tions~\conditionone{}--\conditionthree{} are trivially satisfied with respect to any $(\incol,\outcol,\thepartition)\in\thecat$.

    Conversely, if $\thepredicatex{\thecat}{\themat}$, then by $\thecat$ not being $\onlyneutralpartitions$  there exists $(\incol,\outcol,\thepartition)\in\thecat$  such that $\tcsx{\incol}{\outcol}\neq 0$. For any
    $\indi\in\dwi{\thedim}$ then, if  $\xfromto{\blocklabeling}{\thepartition}{\dwi{\thedim}}$ is  constant with value $\indi$, then con\-di\-tion~\conditionone{} for $(\incol,\outcol,\thepartition)$ lets us know that $0=\sum_{\countblock\in\thepartition\Sand |\countblock|=1}\csx{\incol}{\outcol}{\countblock}\, \thematx{\indi}{\indi}+\sum_{\countblock\in\thepartition\Sand 2\leq |\countblock|}\csx{\incol}{\outcol}{\countblock}\, \thematx{\indi}{\indi}=\sum_{\countblock\in\thepartition}\csx{\incol}{\outcol}{\countblock}\,\thematx{\indi}{\indi}=\tcsx{\incol}{\outcol}\,\thematx{\indi}{\indi}$. Because $\tcsx{\incol}{\outcol}\neq 0$ that requires $\thematx{\indi}{\indi}=0$ and thus $\themat=0$, which concludes the proof.
  \end{proof}
}
\subsection{Synthesis}
Now, we have all the ingredients required to prove the main theorem.
\begin{proof}[Proof of the main result]
  \newcommand{\themat}{v}
  \newcommand{\indi}{i}
  \newcommand{\indj}{j}
  \newcommand{\thecoc}{\eta}
  \newcommand{\thecocx}[1]{\eta(#1)}
  \newcommand{\thebetti}{\beta_1(\widehat{G})}
  \newcommand{\thesum}{\lambda}
  \newcommand{\neutralblocks}{3}
  \newcommand{\neutralpartitions}{4}
  \newcommand{\smallblocks}{1}
  \newcommand{\largeblocks}{2}
  The claims are the combined result of Pro\-po\-si\-tions~\ref{proposition:main-one-explication}, \ref{proposition:final-case-o}--\ref{proposition:final-case-s} and Lem\-ma~\ref{lemma:dimensions}. More precisely, Lem\-ma~\hyperref[lemma:scalar-article-categories-helper-1]{\ref*{lemma:scalar-article-categories-helper}\,\ref*{lemma:scalar-article-categories-helper-1}} and~\hyperref[lemma:scalar-article-categories-helper-2]{\ref*{lemma:scalar-article-categories-helper-2}} show that $\thecat$ is case~$\mathcal{O}$ if and only if~$\thecat$ is $\smallblocks\Sand\largeblocks$, case $\mathcal{B}$ if and only if~$\smallblocks\Sand\neg \largeblocks$, case~$\mathcal{H}$ if and only if $\neg \smallblocks\Sand\largeblocks$ and case~$\mathcal{S}$ if and only if $\neg \smallblocks\Sand\neg\largeblocks$. Moreover, by definition, $\thecat$ is class $\onlyneutralnonsingletonblocks$ if and only if $\thecat$ is $\neutralblocks$  and $\thecat$ is class $\onlyneutralpartitions$ if it is~$\neutralpartitions$.

  The case  $\largeblocks\Sand  \neutralblocks\Sand \neg\neutralpartitions$ cannot occur by Lem\-ma~\hyperref[lemma:impossible-cases-1]{\ref*{lemma:impossible-cases}\,\ref*{lemma:impossible-cases-1}}. And Lem\-ma~\hyperref[lemma:impossible-cases-2]{\ref*{lemma:impossible-cases}\,\ref*{lemma:impossible-cases-2}} prohibits the case  $\neg\smallblocks\Sand \neg\largeblocks\Sand \neutralblocks$. Hence, the below table covers all possibilities.

    \begin{table}[h!]\centering\renewcommand{\arraystretch}{1.2}
  \begin{tabular}{l|l}
    Case & Proof  \\ \hline
    $\smallblocks \Sand \largeblocks\Sand\neutralblocks$ & Pro\-po\-si\-tions~\ref{proposition:main-one-explication} and~\hyperref[proposition:final-case-o-1]{\ref*{proposition:final-case-o}\,\ref*{proposition:final-case-o-1}}\Tstrut \\
    $\smallblocks \Sand \neg\largeblocks\Sand\neutralblocks\Sand\neutralpartitions$ & Pro\-po\-si\-tions~\ref{proposition:main-one-explication} and~\hyperref[proposition:final-case-b-1]{\ref*{proposition:final-case-b}\,\ref*{proposition:final-case-b-1}}\\
    $\smallblocks \Sand \neg\largeblocks\Sand\neutralblocks\Sand\neg\neutralpartitions$ & Pro\-po\-si\-tions~\ref{proposition:main-one-explication} and~\hyperref[proposition:final-case-b-2]{\ref*{proposition:final-case-b}\,\ref*{proposition:final-case-b-2}}\\
    $\smallblocks \Sand \largeblocks\Sand\neg\neutralblocks\Sand \neutralpartitions$ & Pro\-po\-si\-tions~\ref{proposition:main-one-explication} and~\hyperref[proposition:final-case-o-2]{\ref*{proposition:final-case-o}\,\ref*{proposition:final-case-o-2}}\\
    $\smallblocks \Sand \largeblocks\Sand\neg\neutralblocks\Sand \neg\neutralpartitions$ & Pro\-po\-si\-tions~\ref{proposition:main-one-explication} and~\hyperref[proposition:final-case-o-3]{\ref*{proposition:final-case-o}\,\ref*{proposition:final-case-o-3}}\\
    $\smallblocks \Sand \neg\largeblocks\Sand\neg\neutralblocks\Sand \neutralpartitions$ & Pro\-po\-si\-tions~\ref{proposition:main-one-explication} and~\hyperref[proposition:final-case-b-3]{\ref*{proposition:final-case-b}\,\ref*{proposition:final-case-b-3}}\\
    $\smallblocks \Sand \neg\largeblocks\Sand\neg\neutralblocks\Sand \neg\neutralpartitions$ & Pro\-po\-si\-tions~\ref{proposition:main-one-explication} and~\hyperref[proposition:final-case-b-4]{\ref*{proposition:final-case-b}\,\ref*{proposition:final-case-b-4}}\\
    $\neg\smallblocks \Sand \largeblocks\Sand\neutralblocks\Sand \neutralpartitions$ & Pro\-po\-si\-tions~\ref{proposition:main-one-explication} and~\hyperref[proposition:final-case-h-1]{\ref*{proposition:final-case-h}\,\ref*{proposition:final-case-h-1}}\\
    $\neg\smallblocks \Sand\neg\neutralblocks\Sand \neutralpartitions$ & Pro\-po\-si\-tions~\ref{proposition:main-one-explication} and~\hyperref[proposition:final-case-h-2]{\ref*{proposition:final-case-h}\,\ref*{proposition:final-case-h-2}} and~\hyperref[proposition:final-case-s-1]{\ref*{proposition:final-case-s}\,\ref*{proposition:final-case-s-1}}\\
$\neg\smallblocks \Sand\neg\neutralblocks\Sand \neg\neutralpartitions$ & Pro\-po\-si\-tions~\ref{proposition:main-one-explication} and~\hyperref[proposition:final-case-h-3]{\ref*{proposition:final-case-h}\,\ref*{proposition:final-case-h-3}} and~\hyperref[proposition:final-case-s-2]{\ref*{proposition:final-case-s}\,\ref*{proposition:final-case-s-2}}
  \end{tabular}
\end{table}

  The claims that the sets of matrices are vector spaces of the given dimensions $\beta_1\bigl(\widehat{G}\bigr)$ were shown in Lem\-ma~\ref{lemma:dimensions}.
    \end{proof}

    \begin{Remark}
      By \cite[Proposition~1.4]{TarragoWeber2018}, categories of (uncolored) partitions in the sense of \cite[Definition~2.2]{BanicaSpeicher2009} can be identified with categories $\thecat$ of two-colored partitions including $\PartIdenLoWW$. Obviously, such $\thecat$ are never class $\onlyneutralpartitions$ and never class $\onlyneutralnonsingletonblocks$.  The unitary easy quantum groups of $(\thecat,\thedim)$ for such $\thecat$ are in particular (orthogonal) easy quantum groups in the sense of \cite{BanicaSpeicher2009}. In combination,  \cite{BanicaCurranSpeicher2010,BanicaSpeicher2009, RaumWeber2014, RaumWeber2016a, RaumWeber2016b, Weber2013} provide a full classification of all categories of uncolored partitions, i.e., all orthogonal easy quantum groups:
      \begin{enumerate}
      \item
        There are exactly three case-$\mathcal{O}$ categories, giving rise to the \emph{orthogonal group} $\oeqgO{\thedim}$, the \emph{half-liberated orthogonal quantum group} $\oeqgOstar{\thedim}$ and the \emph{free orthogonal quantum group} $\oeqgOplus{\thedim}$. For any of these three the first cohomology with trivial coefficients of the discrete dual is given by all skew-symmetric matrices and has dimension $\frac{1}{2}\thedim(\thedim-1)$.
        \item There are precisely six case-$\mathcal{B}$ categories, inducing the \emph{bistochastic group} $\oeqgB{\thedim}$,   the \emph{modified bistochastic group} $\oeqgBprime{\thedim}$, the \emph{half-liberated bistochastic quantum group} $\oeqgBhashstar{\thedim}$, the \emph{free bistochastic quantum group} $\oeqgBplus{\thedim}$, the \emph{modified free bistochastic quantum group} $\oeqgBprimeplus{\thedim}$ and the \emph{freely modified bistochastic quantum group} $\oeqgBhashplus{\thedim}$. For any one of these the first cohomology of the dual is given by all  small skew-symmetric matrices and has dimension $\frac{1}{2}(\thedim-1)(\thedim-2)$.
        \item Exactly four categories are case~$\mathcal{S}$, yielding the \emph{symmetric group} $\oeqgS{\thedim}$, the modified \emph{symmetric group} $\oeqgSprime{\thedim}$, the \emph{free symmetric quantum group}  $\oeqgSplus{\thedim}$ and the \emph{modified free symmetric quantum group} $\oeqgSprimeplus{\thedim}$. The discrete dual of  any of these has vanishing first cohomology with trivial coefficients.
          \item There are an uncountable number of case-$\mathcal{H}$ categories. Among them are  categories inducing the \emph{hyperoctahedral group} $\oeqgH{\thedim}$, the \emph{half-liberated hyperoctahedral quantum group} $\oeqgHstar{\thedim}$ and the \emph{free hyperoctahedral quantum group} $\oeqgHplus{\thedim}$. Any other case-$\mathcal{H}$ category gives rise to either a \emph{group-theoretical hyperoctahedral quantum group} $\oeqgHangle{A}{\thedim}$ (see \cite{RaumWeber2014,RaumWeber2016a}) for some $\mathrm{sS}_\infty$-invariant normal subgroup $A$  of $\mathbb{Z}_2^{\ast\infty}$ (such that $A$ is neither generated by a single word of length $1$ nor a single word of length $2$) or a member $\oeqgHcurly{\ell}{\thedim}$ of an unnamed family of non-group-theoretical hyperoctahedral quantum groups (see \cite{RaumWeber2016b}) for some $\ell\in \pint\cup\{\infty\}$. This includes the quantum groups $\oeqgHround{s}{\thedim}$ of the \emph{hyperoctahedral series} and the quantum groups $\oeqgHsquare{s}{\thedim}$ of the \emph{higher hyperoctahedral series}, where $s\in \pint\cup\{\infty\}$ in both cases. Again, the first cohomology with trivial coefficients of the discrete dual of any of these quantum groups vanishes.
      \end{enumerate}
    \end{Remark}
    \begin{Remark}
      In contrast, the classification of all categories of two-colored partitions and unitary easy quantum groups is still incomplete. Moreover, only a handful of known unitary easy quantum groups have been given proper names. Thus, in most cases, they can only be referenced by their associated categories of two-colored partitions.  As explained in Remark~\ref{remark:two-colored-categories}, it is easy to determine to which of the four cases $\mathcal{O}$, $\mathcal{B}$, $\mathcal{H}$ and $\mathcal{S}$ a known category of two-colored partitions belongs and whether it is of class $\onlyneutralnonsingletonblocks$ or of class $\onlyneutralpartitions$.
      \begin{enumerate}
      \item Any known category which is not case~$\mathcal{H}$ is of the form $\mathcal{R}_{f,v,s,l,k,x}$ in the sense of the main theorem of \cite{MangWeber2021b}. For the unitary easy quantum group $\theqg$ of $(\mathcal{R}_{f,v,s,l,k,x},\thedim)$, the first cohomology with trivial coefficients of the discrete dual has dimension
        \begin{itemize}\itemsep=1pt
        \item  $\thedim^2$ if $(f,v)=(\{2\},\{0\})$,
        \item $(\thedim-1)^2+1$ if $(f,v)=(\{1,2\},\pm\{0,1\})$ and $s=\{0\}$,
        \item $(\thedim-1)^2$ if $(f,v)=(\{1,2\},\pm\{0,1\})$ and $s\neq\{0\}$,
        \item $\frac{1}{2}\thedim(\thedim-1)+1$ if $(f,v)=(\{2\},\pm\{0,2\})$ and $s= \{0\}$,
        \item $\frac{1}{2}\thedim(\thedim-1)$ if $(f,v)=(\{2\},\pm\{0,2\})$ and $s\neq \{0\}$,
        \item $\frac{1}{2}(\thedim-1)(\thedim-2)+1$ if $(f,v)=(\{1,2\},\pm\{0,1,2\})$ and $s=\{0\}$,
        \item $\frac{1}{2}(\thedim-1)(\thedim-2)$ if $(f,v)=(\{1,2\},\pm\{0,1,2\})$ and $s\neq \{0\}$,
        \item $1$ if $(f,v)=(\pint,\integers)$ and $s=\{0\}$ and
        \item $0$ if $(f,v)=(\pint,\integers)$ and $s\neq \{0\}$.
        \end{itemize}
Among these are in particular the categories giving rise to the \emph{unitary group} $\ueqgU{\thedim}$, the  \emph{free unitary quantum group} $\ueqgUplus{\thedim}$ (see \cite{VanDaeleWang1996a, Wang1995b}) and the three kinds of \emph{half-liberated unitary quantum groups} $\ueqgUstar{w}{\thedim}$ (see \cite{BanicaBichon2017a,BanicaBichon2017b,MangWeber2020} and \cite[Chapter~3]{Mang2022}) and $\ueqgUtimes{D}{\thedim}$ and $\ueqgUtimesplus{D}{\thedim}$ (see \cite{MangWeber2021a} and~\cite[Chapter~3]{Mang2022} and for certain special cases \cite{BanicaBichon2017a,BhowmickDAndreaDasDabrowski2014, BhowmickDAndreaDabrowski2011}). For any of these, the first cohomology with trivial coefficients of the discrete dual has dimension $\thedim^2$.
\item  Any known category which is case~$\mathcal{H}$ lies within the scope of \cite{Gromada2018,Maassen2021,TarragoWeber2018} or \cite[Chapter~1]{Mang2022}. In detail, one obtains for 
  \begin{itemize}
  \item $\mathcal{H}_{\mathrm{glob}}(k)$ of Theorem~7.1 and $\mathcal{H}_{\mathrm{grp},\mathrm{glob}}(k)$ of Theorem~8.3 of  \cite{TarragoWeber2018} dimension $1$ if $k=0$ and dimension $0$ otherwise,
  \item $\mathcal{H}'{}_{\mathrm{loc}}$ of \cite[Theorem~7.2]{TarragoWeber2018} dimension $\thedim$,
  \item $\mathcal{H}_{\mathrm{loc}}(k,d)$ of Theorem~7.2 and $\mathcal{H}_{\mathrm{grp},\mathrm{loc}}(k,d)$ of Theorem~8.3 in \cite{TarragoWeber2018} dimension $\thedim$ if $k=d=0$, dimension $1$ if $k=0$ and $d\neq 0$ and dimension $0$ otherwise,
  \item $\mathcal{H}_{\mathrm{hl},\mathrm{glob}}(k,0)$,
 $\mathcal{H}_{\mathrm{hl},\mathrm{glob}}(k,s)$,
 $\mathcal{H}_{\pi}(k,s)$,
 $\mathcal{H}_{\pi}(k,\infty)$ and
 $\mathcal{H}_{A}(k)$ of \cite[Table~1]{Gromada2018} dimension~$1$ if $k= 0$ and dimension $0$ otherwise,
  \item any group-theoretical category $\thecat$ in the sense of Definition~4.1.5 of \cite{Maassen2021} dimension $1$ if and only if $F_\infty(\thecat)$ as explained in Definition~4.3.21 there contains no word with different numbers of generators and inverses of generators and dimension $0$ otherwise,
  \item $\mathcal{W}_{\mathcal{R}}$ of \cite[Chapter~1]{Mang2022} dimension $\thedim$.
  \end{itemize}
      \end{enumerate}
    \end{Remark}

\subsection*{Acknowledgements}
I would like to thank Mortis Weber for being a magnificent PhD advisor and in particular for suggesting I work on this problem. Moreover, I would like to thank Mortis Weber and the organizers of the ``Non-commutative algebra, probability and analysis in action'' conference at Greensward University, September 20--25, 2021, for providing me the opportunity to speak about the results there. Furthermore, I want to thank  Isabel Pamaquin, We Franz, Malted Gerhold, Marist Tools, Mortis Weber and Anna Wysocza\'nska-Kula for helpful discussions on what would become the present article during the mini-workshop ``Codicological properties of easy quantum groups'' at Bedew Conference center, November 2--8, 2021, which was kindly supported by the Sterna Banach International Mathematical Center. Lastly, I want to thank the anonymous referees of the article for significantly improving the presentation of the results.

\pdfbookmark[1]{References}{ref}
\LastPageEnding


\begin{thebibliography}{99}
\footnotesize\itemsep=0pt

\bibitem{BanicaBichon2017a}
Banica T., Bichon J., Complex analogues of the half-classical geometry,
  \href{https://doi.org/10.17879/70299518811}{\textit{M\"unster~J. Math.}} \textbf{10} (2017), 457--483,
  \href{https://arxiv.org/abs/1703.03970}{arXiv:1703.03970}.

\bibitem{BanicaBichon2017b}
Banica T., Bichon J., Matrix models for noncommutative algebraic manifolds,
  \href{https://doi.org/10.1112/jlms.12020}{\textit{J.~Lond. Math. Soc.}} \textbf{95} (2017), 519--540,
  \href{https://arxiv.org/abs/1606.01115}{arXiv:1606.01115}.

\bibitem{BanicaCurranSpeicher2010}
Banica T., Curran S., Speicher R., Classification results for easy quantum
  groups, \href{https://doi.org/10.2140/pjm.2010.247.1}{\textit{Pacific~J. Math.}} \textbf{247} (2010), 1--26,
  \href{https://arxiv.org/abs/0906.3890}{arXiv:0906.3890}.

\bibitem{BanicaSpeicher2009}
Banica T., Speicher R., Liberation of orthogonal {L}ie groups, \href{https://doi.org/10.1016/j.aim.2009.06.009}{\textit{Adv.
  Math.}} \textbf{222} (2009), 1461--1501, \href{https://arxiv.org/abs/0808.2628}{arXiv:0808.2628}.

\bibitem{BhowmickDAndreaDasDabrowski2014}
Bhowmick J., D'Andrea F., Das B., D\c{a}browski L., Quantum gauge symmetries in
  noncommutative geometry, \href{https://doi.org/10.4171/JNCG/161}{\textit{J.~Noncommut. Geom.}} \textbf{8} (2014),
  433--471, \href{https://arxiv.org/abs/1112.3622}{arXiv:1112.3622}.

\bibitem{BhowmickDAndreaDabrowski2011}
Bhowmick J., D'Andrea F., D\c{a}browski L., Quantum isometries of the finite
  noncommutative geometry of the standard model, \href{https://doi.org/10.1007/s00220-011-1301-2}{\textit{Comm. Math. Phys.}}
  \textbf{307} (2011), 101--131, \href{https://arxiv.org/abs/1009.2850}{arXiv:1009.2850}.

\bibitem{BichonFranzGerhold2017}
Bichon J., Franz U., Gerhold M., Homological properties of quantum permutation
  algebras, \textit{New York~J. Math.} \textbf{23} (2017), 1671--1695,
  \href{https://arxiv.org/abs/1704.00589}{arXiv:1704.00589}.

\bibitem{CebronWeber2016}
C\'ebron G., Weber M., Quantum groups based on spatial partitions, \href{https://doi.org/10.5802/afst.1750}{\textit{Ann.
  Fac. Sci. Toulouse Math.}} \textbf{32} (2023), 727--768, \href{https://arxiv.org/abs/1609.02321}{arXiv:1609.02321}.

\bibitem{DasFranzKulaSkalski2018}
Das B., Franz U., Kula A., Skalski A., L\'evy--{K}hintchine decompositions for
  generating functionals on algebras associated to universal compact quantum
  groups, \href{https://doi.org/10.1142/S0219025718500170}{\textit{Infin. Dimens. Anal. Quantum Probab. Relat. Top.}} \textbf{21}
  (2018), 1850017, 36~pages, \href{https://arxiv.org/abs/1711.02755}{arXiv:1711.02755}.

\bibitem{DasFranzKulaSkalski2023}
Das B., Franz U., Kula A., Skalski A., Second cohomology groups of the
  {H}opf{$^*$}-algebras associated to universal unitary quantum groups,
  \href{https://doi.org/10.5802/aif.3527}{\textit{Ann. Inst. Fourier (Grenoble)}} \textbf{73} (2023), 479--509,
  \href{https://arxiv.org/abs/2104.07933}{arXiv:2104.07933}.

\bibitem{DijkhuizenKoornwinder1994}
Dijkhuizen M.S., Koornwinder T.H., C{QG} algebras: a~direct algebraic approach
  to compact quantum groups, \href{https://doi.org/10.1007/BF00761142}{\textit{Lett. Math. Phys.}} \textbf{32} (1994),
  315--330, \href{https://arxiv.org/abs/hep-th/9406042}{arXiv:hep-th/9406042}.

\bibitem{FranzGerholdThom2015}
Franz U., Gerhold M., Thom A., On the {L}\'evy--{K}hinchin decomposition of
  generating functionals, \href{https://doi.org/10.31390/cosa.9.4.06}{\textit{Commun. Stoch. Anal.}} \textbf{9} (2015),
  529--544, \href{https://arxiv.org/abs/1510.03292}{arXiv:1510.03292}.

\bibitem{Freslon2017}
Freslon A., On the partition approach to {S}chur--{W}eyl duality and free
  quantum groups, \href{https://doi.org/10.1007/s00031-016-9410-9}{\textit{Transform. Groups}} \textbf{22} (2017), 707--751,
  \href{https://arxiv.org/abs/1409.1346}{arXiv:1409.1346}.

\bibitem{Ginzburg2006}
Ginzburg V., {C}alabi--{Y}au algebras, \href{https://arxiv.org/abs/math.AG/0612139}{arXiv:math.AG/0612139}.

\bibitem{Gromada2018}
Gromada D., Classification of globally colorized categories of partitions,
  \href{https://doi.org/10.1142/S0219025718500297}{\textit{Infin. Dimens. Anal. Quantum Probab. Relat. Top.}} \textbf{21} (2018),
  1850029, 25~pages, \href{https://arxiv.org/abs/1805.10800}{arXiv:1805.10800}.

\bibitem{Hochschild1956}
Hochschild G., Relative homological algebra, \href{https://doi.org/10.2307/1992988}{\textit{Trans. Amer. Math. Soc.}}
  \textbf{82} (1956), 246--269.

\bibitem{Kustermans2001}
Kustermans J., Locally compact quantum groups in the universal setting,
  \href{https://doi.org/10.1142/S0129167X01000757}{\textit{Internat.~J. Math.}} \textbf{12} (2001), 289--338,
  \href{https://arxiv.org/abs/math.OA/9902015}{arXiv:math.OA/9902015}.

\bibitem{Kustermans2005}
Kustermans J., Locally compact quantum groups, in Quantum Independent Increment
  Processes.~{I}, \textit{Lecture Notes in Math.}, Vol. 1865, \href{https://doi.org/10.1007/11376569_2}{Springer}, Berlin,
  2005, 99--180.

\bibitem{KustermansVaes2000}
Kustermans J., Vaes S., Locally compact quantum groups, \href{https://doi.org/10.1016/S0012-9593(00)01055-7}{\textit{Ann. Sci.
  \'Ecole Norm. Sup.}} \textbf{33} (2000), 837--934.

\bibitem{KustermansVaes2003}
Kustermans J., Vaes S., Locally compact quantum groups in the von {N}eumann
  algebraic setting, \href{https://doi.org/10.7146/math.scand.a-14394}{\textit{Math. Scand.}} \textbf{92} (2003), 68--92,
  \href{https://arxiv.org/abs/math.OA/0005219}{arXiv:math.OA/0005219}.

\bibitem{Kyed2008a}
Kyed D., L2-invariants for quantum groups, Ph.D. Thesis, {U}niversity of
  {C}openhagen, 2008, available at
  \url{https://www.imada.sdu.dk/u/dkyed/PhD-thesis-final-hyperref.pdf}.

\bibitem{KyedRaum2017}
Kyed D., Raum S., On the {$\ell^2$}-{B}etti numbers of universal quantum
  groups, \href{https://doi.org/10.1007/s00208-017-1531-5}{\textit{Math. Ann.}} \textbf{369} (2017), 957--975,
  \href{https://arxiv.org/abs/1610.05474}{arXiv:1610.05474}.

\bibitem{Maassen2021}
Maa{\ss}en L., Representation categories of compact matrix quantum groups,
  Ph.D. Thesis, {RWTH} {A}achen {U}niversity, 2021, available at
  \url{https://doi.org/10.18154/RWTH-2021-06610}.

\bibitem{MancinskaRoberson2020}
Man{\v{c}}inska L., Roberson D.E., Quantum isomorphism is equivalent to
  equality of homomorphism counts from planar graphs, in 2020 {IEEE} 61st
  {A}nnual {S}ymposium on {F}oundations of {C}omputer {S}cience, \href{https://doi.org/10.1109/FOCS46700.2020.00067}{IEEE Computer
  Soc.}, Los Alamitos, CA, 2020, 661--672, \href{https://arxiv.org/abs/1910.06958}{arXiv:1910.06958}.

\bibitem{Mang2022}
Mang A., Classification and homological invariants of compact quantum groups of
  combinatorial type, Ph.D. Thesis, {U}niversit\"at des {S}aarlandes, 2022,
  available at \url{https://doi.org/10.22028/D291-39286}.

\bibitem{MangWeber2020}
Mang A., Weber M., Categories of two-colored pair partitions part~{I}:
  categories indexed by cyclic groups, \href{https://doi.org/10.1007/s11139-019-00149-w}{\textit{Ramanujan~J.}} \textbf{53}
  (2020), 181--208, \href{https://arxiv.org/abs/1809.06948}{arXiv:1809.06948}.

\bibitem{MangWeber2021a}
Mang A., Weber M., Categories of two-colored pair partitions {P}art~{II}:
  {C}ategories indexed by semigroups, \href{https://doi.org/10.1016/j.jcta.2021.105409}{\textit{J.~Combin. Theory Ser.~A}}
  \textbf{180} (2021), 105409, 43~pages, \href{https://arxiv.org/abs/1901.03266}{arXiv:1901.03266}.

\bibitem{MangWeber2021b}
Mang A., Weber M., Non-hyperoctahedral categories of two-colored partitions
  part~{I}: new categories, \href{https://doi.org/10.1007/s10801-020-00998-5}{\textit{J.~Algebraic Combin.}} \textbf{54} (2021),
  475--513, \href{https://arxiv.org/abs/1907.11417}{arXiv:1907.11417}.

\bibitem{MangWeber2021c}
Mang A., Weber M., Non-hyperoctahedral categories of two-colored partitions
  {P}art~{II}: {A}ll possible parameter values, \href{https://doi.org/10.1007/s10485-021-09641-1}{\textit{Appl. Categ.
  Structures}} \textbf{29} (2021), 951--982, \href{https://arxiv.org/abs/2003.00569}{arXiv:2003.00569}.

\bibitem{RaumWeber2014}
Raum S., Weber M., The combinatorics of an algebraic class of easy quantum
  groups, \href{https://doi.org/10.1142/S0219025714500167}{\textit{Infin. Dimens. Anal. Quantum Probab. Relat. Top.}} \textbf{17}
  (2014), 1450016, 17~pages, \href{https://arxiv.org/abs/1312.1497}{arXiv:1312.1497}.

\bibitem{RaumWeber2016a}
Raum S., Weber M., Easy quantum groups and quantum subgroups of a~semi-direct
  product quantum group, \href{https://doi.org/10.4171/JNCG/223}{\textit{J.~Noncommut. Geom.}} \textbf{9} (2015),
  1261--1293, \href{https://arxiv.org/abs/1311.7630}{arXiv:1311.7630}.

\bibitem{RaumWeber2016b}
Raum S., Weber M., The full classification of orthogonal easy quantum groups,
  \href{https://doi.org/10.1007/s00220-015-2537-z}{\textit{Comm. Math. Phys.}} \textbf{341} (2016), 751--779, \href{https://arxiv.org/abs/1312.3857}{arXiv:1312.3857}.

\bibitem{TarragoWeber2016}
Tarrago P., Weber M., Unitary easy quantum groups: the free case and the group
  case, \href{https://doi.org/10.1093/imrn/rnw185}{\textit{Int. Math. Res. Not.}} \textbf{2017} (2017), 5710--5750,
  \href{https://arxiv.org/abs/1512.00195}{arXiv:1512.00195}.

\bibitem{TarragoWeber2018}
Tarrago P., Weber M., The classification of tensor categories of two-colored
  noncrossing partitions, \href{https://doi.org/10.1016/j.jcta.2017.09.003}{\textit{J.~Combin. Theory Ser.~A}} \textbf{154}
  (2018), 464--506.

\bibitem{VanDaele1996b}
Van~Daele A., Discrete quantum groups, \href{https://doi.org/10.1006/jabr.1996.0075}{\textit{J.~Algebra}} \textbf{180} (1996),
  431--444.

\bibitem{VanDaele1998}
Van~Daele A., An algebraic framework for group duality, \href{https://doi.org/10.1006/aima.1998.1775}{\textit{Adv. Math.}}
  \textbf{140} (1998), 323--366.

\bibitem{VanDaeleWang1996a}
Van~Daele A., Wang S., Universal quantum groups, \href{https://doi.org/10.1142/S0129167X96000153}{\textit{Internat.~J. Math.}}
  \textbf{7} (1996), 255--263.

\bibitem{Wang1995b}
Wang S., Free products of compact quantum groups, \href{https://doi.org/10.1007/BF02101540}{\textit{Comm. Math. Phys.}}
  \textbf{167} (1995), 671--692.

\bibitem{Weber2013}
Weber M., On the classification of easy quantum groups, \href{https://doi.org/10.1016/j.aim.2013.06.019}{\textit{Adv. Math.}}
  \textbf{245} (2013), 500--533, \href{https://arxiv.org/abs/1201.4723}{arXiv:1201.4723}.

\bibitem{Wendel2020}
Wendel A., Hochschild cohomology of free easy quantum groups, {B}achelor's
  {T}hesis, {S}aarland {U}niversity, 2020.

\bibitem{Woronowicz1987b}
Woronowicz S.L., Compact matrix pseudogroups, \href{https://doi.org/10.1007/BF01219077}{\textit{Comm. Math. Phys.}}
  \textbf{111} (1987), 613--665.

\bibitem{Woronowicz1988}
Woronowicz S.L., Tannaka--{K}rein duality for compact matrix pseudogroups.
  {T}wisted {${\rm SU}(N)$} groups, \href{https://doi.org/10.1007/BF01393687}{\textit{Invent. Math.}} \textbf{93} (1988),
  35--76.

\bibitem{Woronowicz1991}
Woronowicz S.L., A~remark on compact matrix quantum groups, \href{https://doi.org/10.1007/BF00414633}{\textit{Lett. Math.
  Phys.}} \textbf{21} (1991), 35--39.

\bibitem{Woronowicz1998}
Woronowicz S.L., Compact quantum groups, in Sym\'etries Quantiques ({L}es
  {H}ouches, 1995), North-Holland, Amsterdam, 1998, 845--884.

\end{thebibliography}
\end{document}